\documentclass[
11pt,                          % standard font size
english                        % standard language
]{book}

\usepackage[english]{babel}    % with explicit language
\usepackage{amsmath}           % ams mathematical stuff
\usepackage[utf8]{inputenc}    % smart input of funny chars
\usepackage[T1]{fontenc}       % also for the font encoding
\usepackage{longtable}         % tables longer than one page
\usepackage{exscale}           % large summation signs in 11pt
\usepackage[final]{graphicx}   % to include pdf pictures
\usepackage[sort]{cite}        % nicer citations
\usepackage{array}             % nice tables
\usepackage{wasysym}           % smiley symbols
\usepackage[a4paper]{geometry} % geometry of page layout
\usepackage{xspace}            % better spacing after macros
\usepackage{tikz}              % for commutative diagrams and stuff
\usepackage{chairx}            % the Chair X style file
\usepackage[expansion=false    % no font expansion
           ]{microtype}        % only protrusion
\usepackage[nottoc]{tocbibind} % refs and index in the toc
\usepackage[backref=page,      % backrefs in the bibliography
           final=true,         % always treat as final
           pdfpagelabels       % use pdf page labels
           ]{hyperref}         % hyperrefs are cool!
\usepackage[toc,page]{appendix}

\graphicspath{{../tikz/}}

\usetikzlibrary{matrix}
\usetikzlibrary{arrows}
\usetikzlibrary{patterns}
\usetikzlibrary{decorations.pathreplacing}
\usetikzlibrary{decorations.text}

\newcommand{\Eta}{\operatorname{H}}

\newcommand{\coproduct}{\Delta}
\newcommand{\ocoproduct}{\Delta}
\newcommand{\ostar}{\star}
\newcommand{\bch}[2]{\mathrm{BCH}\left(#1, #2\right)}
\newcommand{\bchpart}[3]{\mathrm{BCH}_{#1}\left(#2, #3\right)}
\newcommand{\bchparts}[4]{\mathrm{BCH}_{#1, #2}\left(#3, #4\right)}
\newcommand{\bchtilde}[4]{\widetilde{\mathrm{BCH}}_{#1, #2}\left(#3; #4\right)}
\newcommand\ot[2]{\mathrel{\overset{\makebox[0pt]
	{\mbox{\normalfont\footnotesize\sffamily #1}}}{#2}}}
\renewcommand{\vector}[2]{\begin{pmatrix} #1 \\ #2 \end{pmatrix}}

\begin{document}

% ============================================================================
% 			Header part
% ============================================================================

% title page
%
\begin{titlepage}
	\begin{center}
	\vspace*{2.5cm}
		\Huge{
			Convergence of the Gutt star product
		}
		\\
		\vspace{1cm}
		\Large{
			\textbf{Paul Stapor,} \\
			October 27, 2015 \\
		}
		
		\vspace{7.5cm}
				
		\Large{
				Master Thesis in the Study Program
				Mathematical Physics, M. Sc.
				\\
				with
				\\
				Prof. Dr. Stefan Waldmann \\
		}
		\vspace{1cm}
		\Large{
			Chair X (Mathematical Physics), \\
			Department of Mathematics, \\
			Julius-Maximilians-University, 
			W\"urzburg
		}
	\end{center}
\end{titlepage}
\thispagestyle{empty}

% table of contents
%
\tableofcontents
\thispagestyle{empty}

% ============================================================================
% ////////////////////////////////////////////////////////////////////////////

% ============================================================================
% 			Main Part
% ============================================================================

% Introduction
%

%
% the Introduction of my master thesis
%

\chapter{Introduction}
\label{sec:Intro}

\section*{Mathematical Physics}
Throughout history, the fields of mathematics and physics have been 
closely linked to each other. The great physicists of the past have always been 
great mathematicians and vice versa: Carl Friedrich Gauss, one of the
most brilliant mathematicians, did not only find plenty of mathematical 
relations, prove a lot of theorems and develop many new ideas, which 
should become rich and fruitful new fields in later mathematics, but he also has a 
large number of credits in physics: the recovery of the dwarf planet Ceres in 
astronomy, new results in electromagnetism (like the Gauss's law or a 
representation for the unit of magnetism, which was named after him) and the 
Gaussian lens formula in geometric optics are just some of his best known 
merits. Isaac Newton on the other hand may have been rather a physicist, but 
it was his so-called second law of motion that used for the first time 
something like differential calculus and therefore opened up the door for new 
branches of mathematics. Even if one want claims that differential calculus 
was actually invented by Gottfried Wilhelm Leibniz, this does not change 
the basic observation, since Leibniz wrote a lot of essays on physics and can 
be considered as the inventor of the concept of kinetic energy (or, as he 
called it, the vis viva) and its conservation in certain mechanical 
systems. Of course one has to name Joseph-Louis Lagrange, who was an ingenious 
mathematician with rich contributions to number theory and algebra, but also to 
the fields of analytical mechanics and astronomy. Still today, every physics 
student has to learn his Lagrangian formalism in the lecture on theoretical 
mechanics and how one can derive the laws of motion for various mechanical 
systems from it. A last name we want to mention here is Paul 
Dirac, who is certainly one of the founding fathers of quantum mechanics. He 
provided the ideas for a lot of structures and relations in differential 
geometry, functional analysis and distribution theory. Many of the concept he 
introduced using his physicist's intuition have been later proven to be correct or 
used as starting points for new theories by mathematicians.

A lot of developments in mathematics can be seen as triggered by physics: they 
were necessary to describe the physical behaviour of our world and therefore 
pushed forward by scientists. We already mentioned differential calculus, 
without whom modern analysis, the theory of ordinary or partial differential 
equations or differential geometry would not be possible. Besides the also 
named field of functional analysis, also Lie theory and many parts of geometry 
provide examples for mathematics which was inspired by physics. Of course, this 
correspondence is not a one way street since the understanding of nature made 
great progress due to a better knowledge of the mathematical laws that describe 
it. A good example fir this is Lebesgue's theory of integration and its 
application to quantum mechanics: the space $L^2(\mathbb{R}^{3n})$ is the state 
space of standard $n$-particle system in quantum mechanics.

There are good reasons to say that this extremely tight binding of mathematics 
and physics persisted until the 20th century. Without any doubt, these two 
areas are still closely linked, but one could say that at a certain point in 
history they started walking away from each other. Of course, there have always 
been mathematicians who did not take their motivation from physics and 
physicists who did not use elaborate mathematics or even find new theorems to 
describe aspects of the world around them, but for a long time, the vast 
majority of both groups showed at least an interest for the other domain. This 
definitely changed during the 20th century. The main reason for this can 
surely be found in the extremely fast development which both of the domains 
experienced in this time. It is already impossible for one person to overview 
the whole field of mathematics or the one of physics, since there are too many 
new things coming up every day and one has to specialize for being able to do 
research. Another reason is surely the fact that modern mathematics is 
strongly influenced by the desire to formulate things as clean as possible, 
without using handwaving ``physical arguments''. This is a principle which 
surely allowed many new and fruitful evolutions in the last decades and which 
is mostly due to the Bourbaki movement in the middle of the last century, but 
it also forgets about the fact that physical intuition was often a powerful 
tool for new ideas or also for heuristics which led to proofs of important 
theorems. Another reason, which is more situated in the domain of physics, is 
certainly the incredibly fast development of the knowledge about 
semiconductors. This became possible due to quantum mechanics which forms the 
foundation of this theory, but for very most of modern applications a basic 
understanding of the quantum theory behind is enough or one can even get new 
results with so called semi-classical approaches. Here, a lot of new results 
can be established without going deep into mathematics and hence without giving 
a new stimulus to it. In this sense, it is enough for many modern physicists to 
acquire a certain amount of mathematical knowledge and then they never have to 
care about mathematical theories again.

Certainly, this situation is due to a natural development does not present a 
problem, although it is a bit of a pity. However, it would be too much to say 
that those fields are falling apart: there are still a lot of intersections of 
the two sciences and these contact areas provide rich domains of research. The 
range of topics, where either mathematics takes its motivation from physics, 
or where theoretical physics needs very elaborate mathematical methods, is 
usually grouped under the name \emph{mathematical physics}.

\section*{A Mathematical Theory of Quantization}

One of the younger branches in this field is the theory of quantization.
It belongs to the area of pure mathematics, but takes its inspiration from physics 
and is therefore a part of mathematical physics. The idea is, roughly spoken, to 
find a correspondence between the quantum and the classical world in physics. The 
mathematical description of their laws are different but yet there are a lot of 
similarities. It is more or less clear, how the classical world emerges out 
of a huge number of quantum objects and the mathematics of classical mechanics 
can be understood as a limit case the behaviour of $n$ quantum objects where 
$n \longrightarrow \infty$. The other way round, it is not clear how one can 
create the mathematical description of a quantum system out of the one of a 
classical system. This reversed process is usually called \emph{quantization} and 
its understanding is a mathematical task, not a physical one. 

One has to say that the question of \emph{how to quantize a given classical 
system} is still far from being well understood and there are various approaches 
to it. We will give an overview of them in Chapter 2.2.3. However, all these 
attempts to create a mathematical theory of quantization have led to many new 
developments in mathematics, also in other fields: they had a strong drawback on 
differential and algebraic geometry, Lie theory and functional analysis. In 
particular the theory of Poisson Lie groups is completely due to them.

\subsection*{Deformation Quantization and how this Master Thesis fits in}

Deformation quantization is not the only one, but maybe the most developed theory 
in this area. Roughly speaking, it tries to \emph{deform} the (idealized) 
commutative algebra of smooth functions on a manifold (of classical physical 
observables) by making it noncommutative and wants to get an (idealized) algebra 
(of quantum mechanical observables) this way. This is done by replacing the 
pointwise product of functions with a noncommutative product, which takes into 
account certain derivatives of the functions and plugs in a formal parameter which 
is called $\hbar$ and which should correspond to Planck's constant in physics. 
This new product is then called a star product. There are basically two different 
ways of constructing such a star product: either by integral formulas, or by 
formal power series in $\hbar$. These two approaches are linked by the fact that 
one gets a formal series out of the integral formulas by an asymptotic expansion. 
This mostly algebraic theory of formal power series is called formal deformation 
quantization, and it is quite well understood: although there are still a lot of 
open questions remaining, many existence and classification results for star 
products were found, there is also a light physical interpretation: The zeroth 
order of the formal parameter represents classical mechanics and the first order 
quantum mechanics. Different mechanical systems allow different star products, and 
of course it is interesting to ask whether some star products are more ``natural'' 
for a given system than others.

Besides the algebraic aspects of this theory, one can also speak about convergence 
of these formal series or study the integral formulas: so one can ask the 
question, if the deformation is continuous or smooth in some sense. This leads to 
the field of \emph{strict deformation quantization}. Here, no closed theory exists 
yet and there are just very few things known about \emph{continuous} star 
products. The approach via integral formulas is better developed at the moment, 
but it is only possible for finite-dimensional systems. For a formal power series, 
the first question is of course in what topology one wants it to converge. Once 
one has an answer to it, one has to do some more or less complicated analysis of 
the formal power series, in order to control it somehow.

For this second approach via formal power series, there are just very few examples 
of continuous star products known. They are all related to classical systems, 
which can be described as symplectic vector spaces. This master thesis presents a 
new example of a formal star product, which can be turned into a convergent series 
and which does not come from a symplectic vector spaces, but from a class of 
vector spaces with a more general (Poisson) structures. So it enlarges this theory 
of strict deformation quantization by a very substantial and big group of examples 
and provides new scientific results. The core part of it is hence made public as a 
preprint \cite{esposito.stapor.waldmann:2015a:pre} by now. This preprint will not 
be cited again in this thesis, but it is clear that the latter contains the whole 
content of the former and even a bit more.

Besides the fact, that this work presents a contribution to the theory of 
quantization, its main parts are also of independent interest. The topologized, 
deformed algebra, which is treated here is nothing else but the universal 
enveloping algebra $\algebra{U}(\lie{g})$ of a Lie algebra $\lie{g}$. So it also 
has an impact on Lie theory, since there is no canonical topology on $\algebra{U}
(\lie{g})$: we give a candidate for it by an explicit and functorial 
construction, which also applies to some infinite-dimensional Lie algebras. To the 
best of our knowledge, such a topology for infinite-dimensional cases has not been 
known before.

\subsection*{An Outlook}

Although (or maybe because) we found many new results, we have to face a lot of 
open questions. The first one is of course, how far we can generalize the 
construction of the topology on the universal enveloping algebra and how large its 
completion can be. It is almost answered by this thesis, but some limiting cases 
are left to consider, as pointed out in the Chapters 6.2 and 6.4. Another question 
considers the representation theory of $\algebra{U}(\lie{g})$: now we have a 
topology and can ask for continuous representations, for example. Some more 
questions rather touch the field of deformation quantization again. First: we can 
think of continuous star products not only on vector spaces, but also on some 
really geometric manifolds, namely on the coadjoint orbits of a Lie group $G$.
The results of this work are a very big step into this direction and may allow new 
examples of strict deformations. There is of course still work to be done, but 
this would be the very first class of such examples. Second: we can think of more 
complicated classes of quantized classical mechanical systems on vector spaces and 
try to make those more general star products convergent, too. The final aim would 
be to find results for the so-called Kontsevich construction, which is somehow the 
biggest result in formal deformation quantization. The question, how this can be 
transferred into a strict setting is completely open. Third: a big class of star 
products can be generated using twists. That is a technique, which uses the Hopf 
algebra structure of $\algebra{U}(\lie{g})$. Since we now have a topology, we can 
think of convergent twists and ask, if this leads to a convergent (non-formal) 
power series, hence to a continuous star product. This question is closely linked 
to symmetries of the quantized system and could lead to a better understanding of 
such things as ``quantized symmetries'' which have yet not been understood at all.

\section*{Summary and Organization}

\subsection*{Summary}

This work focusses on a particular star product, the so-called Gutt star 
product. It can be established on a certain class of classical mechanical systems, 
which turn out to be vector spaces with a particular type of biderivation (a 
Poisson bracket) on their algebra of smooth functions. The aim of this thesis
is to find a large subalgebra of the smooth functions and a (locally 
convex) topology on them, such that the Gutt star product is continuous and that 
the commutative classical algebra can be deformed smoothly into the 
noncommutative quantum algebra. Of course, we must give a proper definition 
how this smoothness is meant. We also establish many helpful
properties of this construction, like a certain functoriality. Moreover, this work 
relates to other fields of mathematics, such as Lie theory. For example, the 
quantized algebra of smooth functions is closely linked to a universal enveloping 
algebra of a Lie algebra, which naturally appears for the considered classical 
systems. Therefore, it carries the structure of a Hopf algebra, which is also 
continuous with respect to the constructed topology.

\subsection*{Organization}

This master thesis is organized as follows: In Chapter 2, we will
introduce the most important concepts of classical and quantum mechanics, 
explain their relations and give an overview over the field of deformation 
quantization, its history and its current state. We will also classify it this 
thesis into the theory of deformation quantization. We will also give an outlook 
on the next steps, which can be done using the results of this work.
In the third chapter, we will explain in more detail the kind of Poisson 
systems this thesis deals with and see that they are in fact Lie algebras. We 
will construct the Gutt star product, which is characteristic for those 
systems, in different ways and show that these constructions are equal. We will 
explain the link to Lie theory, the Poincar\'e-Birkhoff-Witt and the 
Baker-Campbell-Hausdorff theorem and will introduce those results on the 
Baker-Campbell-Hausdorff series, which we will need for our later work.
Chapter four is devoted to finding explicit formulas for the Gutt star 
product and explaining them using two examples, as well as to some easy 
conclusions one can draw from those formulas.
Chapter five is the core of this thesis. First, we will introduce briefly the 
concept of locally convex topologies and explain why they are the convenient 
setting for our task. We will show in detail how the locally convex topology for 
the Gutt star product is constructed, what their properties are and what kind of 
topology our Poisson system must have had such that the construction of our 
topology is possible. At this point, we will introduce the concept of 
\emph{asymptotic estimate algebras}, 
which can be seen as a concept between locally multiplicatively convex algebras 
a general locally convex ones. Then we will show that the Gutt star product is 
indeed continuous with respect to our topology, that the deformation is 
analytic (even entire, if the underlying field is $\mathbb{C}$), that the 
construction is functorial and we will analyse the completion of this algebra. 
We will also show that this topology is optimal in a certain way.
The sixth chapter is devoted to particular systems, namely to nilpotent Lie 
algebras. We will show how the results we found previously can be improved, 
but we will also see the limits of our construction. We will establish 
the link to a work of Stefan Waldmann's and we will show that we come to the 
same conclusion by taking a different way using the Gutt star product. Finally, 
we will see that those stronger results are not limited to the very case of 
strictly nilpotent Lie algebras, but that there are (in infinite dimensions) 
weaker notions of nilpotency which lead to the same result, when the 
construction of the topology is adapted a bit.
In the end, Chapter seven treats the Hopf algebraic part which is very short 
due to the algebraic properties of the deformation. We will see that 
the coproduct and the antipode remain undeformed and continuous with respect 
to our topology, too.

\subsection*{Summary and Organization (German)}

Diese Masterarbeit befasst sich mit dem Gutt-Sternprodukt. Es kann f\"ur bestimmte 
Poissonmannigfaltigkeiten definiert werden, n\"amlich f\"ur Vektorr\"aume mit 
linearem Poissontensor. Die Poissonmannigfaltigkeit stellt dann den Dualraum einer 
Lie-Algebra dar. Auf der Polynomalgebra \"uber diesem Dualraum, bzw. der 
symmetrischen Tensoralgebra \"uber der Lie-Algebra, l\"asst sich dann das 
Gutt-Sternprodukt definieren, welches als formale Deformation der symmetrischen 
Tensoralgebra aufgefasst werden kann. Im Rahmen dieser Arbeit wird mit Hilfe einer 
von Waldmann entwickelten Konstruktion die symmetrische Tensoralgebra 
lokal-konvex topologisiert und dann der Nachweis erbracht, dass das 
Gutt-Sternprodukt in dieser Topologie tats\"achlich stetig ist. Das besondere bei 
diesem Vorgehen ist, dass es sich insbesondere auch auf eine große Klasse von 
unendlich-dimensionalen Lie-Algebren anwenden l\"asst. Damit tr\"agt diese Arbeit 
neue wissenschaftliche Ergebnisse zum Themengebiet der strikten 
Deformationsquantisierung bei.

Im nächsten Abschnitt, Kapitel 2, wird der Begriff ''Quantisierung'' und 
insbesondere der Themenbereich der Deformationsquantisierung n\"aher erl\"autert 
und die Fragestellung motiviert. Im dritten Kapitel werden die notwendigen 
Vorkenntnisse \"uber Kirillov-Kostant-Souriau-Strukturen, das Gutt-Sternprodukt 
und die beiden Theoreme nach Poincar\'e, Birkhoff und Witt sowie Baker, Campbell 
und Hausdorff vermittelt. Kapitel 4 besch\"aftigt sich mit der Herleitung 
verschiedener Formeln f\"ur das Gutt-Sternprodukt. Im f\"unften Kapitel wird dann 
die lokal-konvexe Topologie auf der Tensoralgebra eingef\"uhrt sowie die 
notwendige Unterkategorie lokal-konvexer Lie-Algebren erkl\"art, mit der wir uns 
befassen. Danach wird die Stetigkeit des Sternprodukts bewiesen, ebenso wie die 
analytische Abh\"angigkeit vom Deformationsparameter, die Funktorialit\"at der 
Konstruktion und die Optimalit\"at des Ergebnisses. Kapitel 6 besch\"aftigt 
sich mit nilpotenten Lie-Algebren und damit, wie weit sich die bisherigen 
Ergebnisse in diesem Fall ausbauen lassen. Dar\"uber hinaus wird der Zusammenhang 
mit einer Arbeit von Waldmann bez\"uglich der Weyl-Algebra aufgezeigt. Im siebten 
Kapitel wird gezeigt, dass neben der assoziativen Algebra- auch eine Hopf-Struktur 
vorliegt, die ebenfalls lokal-konvex topologiert werden kann und stetig ist.

\section*{Thanks to \ldots}

I want to thank my advisor, Stefan Waldmann, for the time he invested in this 
work and for his intense supervision. I am also very grateful to Chiara 
Esposito, who also spent a lot of time in discussion with me and was so kind to 
offer me a place in her office, so that I had a place at the chair where I 
could work. It is mainly due to both of them, that this work progressed so 
quickly and finally became a preprint 
which will hopefully become a publication in the next time. I am also grateful 
to Matthias Sch\"otz and Thorsten Reichert for many fruitful discussions 
and their patience with me bothering them with questions. Moreover, I am 
grateful for the support of the Studienstiftung and the Max-Weber Programm, who 
financed my studies and without whom my life would have looked very 
differently during the past six years. Finally, I want to 
say thank you to my parents who supported me all the time of school and 
studies, who encouraged me and who awakened my interest for the understanding 
of nature by their way of educating me. A last thank you is dedicated to my 
girlfriend for encouraging and supporting me when I was unmotivated and for not 
being angry when I spent a lot of time working on this thesis. :-)

% Chapter 2
%

%
% Chapter 2 of my master thesis:
% Still quiete a lot introductory
%

\chapter{Deformation quantization}

The starting point for every theory of quantization is without any doubt the 
theory of mechanics. Here, mechanics means both, the classical and the quantum 
theory. Since we want to explain how one can link those two together, we will give 
a short overview of both theories and show their similarities and their 
differences in the first section of this chapter. Afterwards we will explain what 
a quantization should actually be and collect the existing approaches to it. In 
Section three, we will focus on Deformation quantization and give an overview of 
this rather young theory.

\section{Mechanics: The Classical and the Quantum World}
\label{sec:chap2_Mechanics}

\subsection{Classical Mechanical Systems}
\label{subsec:chap2_Classical}
We want to briefly recall the notions of classical mechanics. There are many 
good books on this subject and also the notation for the basic concepts is more or 
less uniform everywhere. Nevertheless, we want to refer to the books of Marsden 
and Ratiu \cite{marsden.ratiu:1999a} and Arnold \cite{arnold:1989a}, which give 
very good introductions and overviews of the theory.
Imagine the simplest model for a mechanical system, which is not trivial: a 
single particle with mass $m$ moving in $\mathbb{R}^3$ in a scalar potential $V$. 
We will denote its position by $q = (q^1, q^2, q^3) \in Q = \mathbb{R}^3$ and call 
the set $Q$ of all possible positions the \emph{configuration space}. Since we 
also have a time coordinate $t \in \mathbb{R}$, we can describe the path on which 
the particle moves by a parametrized curve $q(t)$. The state of the particle is 
completely described by its position and its velocity $(q(t), \dot q(t))$ and the 
velocity should be understood as a tangent vector $\dot q(t) \in T_{q(t)}Q = 
\mathbb{R}^3$. Therefore, the tangent bundle $TQ$, which is in this case 
$\mathbb{R}^6$, is sometimes called the state space. In classical mechanics, we 
can describe the movement of the particle using the Euler-Lagrange equations, 
which are derived from the so-called Lagrange function by a variational principle. 
The Lagrange function of this system reads
\begin{equation*}
	\mathcal{L}(q, \dot q)
	=
	T(q, \dot q) - V(q, \dot q)
\end{equation*}
and $T(q, \dot q) = \frac{m}{2}  \sum_{i=1}^3 \left( \dot{q}^i \right)^2$ is the 
kinetic energy. The action along a path is defined as the integral
\begin{equation*}
	\mathcal{S}(q, \dot q)
	=
	\int\limits_{t_0}^{t_1}
	\mathcal{L}(q, \dot q)
	dt.
\end{equation*}
It is a physical observation, similar to a mathematical axiom, that this action 
functional is stationary along the trajectories of the particle. So by fixing a 
starting point $q_0$ and an end point $q_1$ of a trajectory, we find
\begin{equation*}
	\delta \mathcal{S}
	\at{q(t_0) = q_0, q(t_1) = q_1}
	=
	0.
\end{equation*}
From this, one finds the Euler-Lagrange equations for $i = 1, 2, 3$
\begin{equation*}
	\frac{d}{dt} 
	\frac{\partial \mathcal{L}}{\partial \dot q^i}
	-
	\frac{\partial \mathcal{L}}{\partial q^i}
	=
	0.
\end{equation*}
In our case this means
\begin{equation*}
	\frac{d}{dt} 
	\frac{\partial T}{\partial \dot q^i}
	-
	\frac{\partial V}{\partial q^i}
	=
	m \ddot q
	-
	\frac{\partial V}{\partial q^i}
	=
	0
\end{equation*}
and we can integrate the equations to get the trajectory. So far, we described 
the system in the state space. This description is called 
the Lagrange formalism. It is not the only picture, which describes to 
behaviour of the system: one can go to the so-called Hamilton formalism, which 
is based on the conjugate momenta $p = (p_1, p_2, p_3) \in T_{q(t)}^*Q$
\begin{equation*}
	p_i(t)
	=
	\frac{\partial \mathcal{L}}{\partial \dot q^i}
\end{equation*}
and linked to the Lagrange formalism by a Legendre transformation. The set $T^*Q$ 
also describes all possible states of the system and is called the phase space. 
Now we can define the Hamilton function
\begin{equation*}
	H(q, p)
	=
	p_i q^i
	-
	\mathcal{L}(q, \dot q),
\end{equation*} 
which is the crucial quantity in this setting. It represents the energy of the 
system. One finds the equations
\begin{equation*}
	\frac{\partial H}{\partial p_i}
	=
	\frac{d q^i}{dt}
	\quad \text{ and } \quad
	- \frac{\partial H}{\partial q^i}
	=
	\frac{d p_i}{dt}.
\end{equation*}
They strongly remind of a symplectic structure and indeed this is the case: 
using the standard sympectic matrix
\begin{equation*}
	\omega
	=
	\begin{pmatrix}
		0 & \Unit
		\\
		- \Unit & 0
	\end{pmatrix}
\end{equation*}
one finds
\begin{equation*}
	\frac{d}{dt} \vector{q^i}{p_i}
	=
	\omega
	\vector{\frac{\partial H}{\partial q^i}}
	{\frac{\partial H}{\partial p_i}}.
\end{equation*}
More generally, one can define the Poisson bracket $\{f, g\}$ for two functions 
$\Cinfty(T^*Q)$ by
\begin{equation*}
	\{f, g\}
	=
	\sum\limits_{i=1}^3
	\left(
		\frac{\partial f}{\partial q^i}
		\frac{\partial g}{\partial p_i}
		-
		\frac{\partial f}{\partial p_i}
		\frac{\partial g}{\partial q^i}
	\right)
\end{equation*}
and finds the time evolution of a function $f$ given by
\begin{equation*}
	\frac{d}{dt} f
	=
	\{ f, H \}.
\end{equation*}
All of those objects have, of course, geometrical interpretations which 
allow generalizations from this very easy example. For a more 
complex system than one particle in a scalar potential, the configuration space 
$Q$ is a smooth manifold, the state space $TQ$ is described by the tangent 
bundle and the phase space $T^*Q$ by the cotangent bundle. Points in the phase 
space, which can also be understood as Dirac measures in $T^*Q$, describe the 
possible states of the system. This interpretation allows us to speak of positive 
Borel measures on $T^*Q$ as mixed states, which describe a probabilistic 
distribution of the state of the system. The transition from 
the Lagrange to the Hamilton formalism is done by a fiber derivative, the 
kinetic energy has a mathematical interpretation as a Riemannian metric and the 
Poisson bracket is the one which is due to the canonical symplectic form on the 
cotangent bundle. So in some sense, classical mechanics can be described by 
symplectic geometry.

The last conclusion, however, was a bit too fast. There are mechanical systems 
like the rigid body, which can not be described in the symplectic formalism, 
or which have certain symmetries and therefore allow a \emph{reduction}. If a 
given system has a symmetry (which is described by a certain type of Lie group 
action), there will be mathematical tools, which allow to divide out a part of 
phase space. One gets a reduced phase space, which is a quotient of $T^*Q$ 
and which, in general, is not symplectic any more. In these cases, we end up 
with a more general structure than a symplectic manifold. The Poisson bracket 
however is still be there. So the objects we actually want to use in order to 
describe classical mechanics are \emph{Poisson manifolds}. 
\begin{definition}[Poisson Manifold]
	A Poisson manifold is a pair $(M, \{ \cdot, \cdot \})$ of a smooth 
	manifold $M$ and a Poisson bracket $\{ \cdot, \cdot\}$. The bracket is a 
	biderivation
	\begin{equation*}
		\{ \cdot, \cdot \}
		\colon
		\Cinfty(M)
		\times
		\Cinfty(M)
		\longrightarrow
		\Cinfty(M),
	\end{equation*}
	which is anti-symmetric and fulfils the Jacobi identity
	\begin{equation*}
		\{ \{f, g\}, h \}
		+
		\{ \{g, h\}, f \}
		+
		\{ \{h, f\}, g \}
		=
		0.
	\end{equation*}
\end{definition}
It is worth mentioning that sometimes, it will be more convenient not to talk 
about the Poisson bracket $\{\cdot, \cdot\}$, but about the \emph{Poisson tensor} 
$\pi \in \Secinfty(\Anti^2(TV))$, which is a equivalent description. Like the 
bracket, the tensor has to fulfil the Jacobi identity. We have for 
all $f,g \in \Cinfty(M)$ and $x \in M$
\begin{equation*}
	\{f, g\}(x)
	=
	\pi\at{x}(\D f, \D g)(x).
\end{equation*}
The ``physical observables'' in this setting are the smooth functions on the 
Poisson manifold. Together with the Poisson bracket, they form the prototype of a 
\emph{Poisson algebra}: an algebra with a bracket, which is an antisymmetric 
biderivation and which fulfils the Jacobi identity. The range, sometimes also 
called the spectrum, of a given smooth function has the physical interpretation as 
the measurable values of this observable. So in full generality, we could say that 
a classical mechanical system is a Poisson manifold together with a Hamilton 
function $H$, and the time evolution of every observable $f \in \Cinfty(M)$ is 
given by the relation
\begin{equation*}
	\frac{d}{dt} f(t)
	=
	\{f(t), H\}.
\end{equation*}
We will soon see that this formalism is as close as we can get to the one of 
quantum mechanics.

\subsection{Quantum Mechanics}
\label{subsec:chap2_Quantum}

In quantum theory, the model of a physical system looks completely different on 
the first sight. There is indeed no obvious way to describe how this formalism was 
derived out of the one describing classical mechanics. It was a hard 
and non-straightforward way, taken in small steps of which some were 
more or less obvious and some were ingenious guesses. Built up on first ideas 
by Max Planck, who described the radiation of a black body, Albert Einstein, 
who explained the photo-electrical effect and Niels Bohr, who solved the 
problem of atomic spectra, a hand full of physicists including Max Born, 
Werner Heisenberg, Erwin Schr\"odinger, Paul Dirac, Wolfgang Pauli and Pascual 
Jordan developed a complex, counter-intuitive but extremely well working theory, 
which was able to give precise answers to very most of the open questions at that 
time. However, this theory, mainly risen between 1925 and 1928, did not have a 
satisfying mathematical fundament. It took around ten more years until 
mathematicians, mainly John von Neumann, worked out a mathematical theory for 
quantum mechanics and gave also a physical interpretation to their mathematical 
formulation \cite{vonneumann:1996a}. Today, quantum mechanics is usually 
introduced as a mathematical theory based on axioms and most of the textbooks 
(at least most of the more mathematical ones) use this axiomatical approach to 
explain the theory. Two very nice mathematical introductions are given by 
Bongaarts \cite{bongaarts:2015a} and Hall \cite{hall:2013a}.

The state of a physical system is described by a vector $\psi$ on a 
Hilbert space $\hilbert H$ and the observables are self-adjoint (and 
usually unbounded) operators on it. The spectra of these operators describe 
their measurable values. The theory is probabilistic, but has a deterministic 
time evolution: the pair $(A, \psi)$ gives a probabilistic distribution, 
from which we can calculate the probability to measure $A$ with the value $a 
\in \spec (A)$ in the state $\psi$ using the spectral resolution of $A$. We 
want to look again at the example of a particle in $\mathbb{R}^3$ 
moving in a potential. Let $x = (x_1, x_2, x_3) \in \mathbb{R}^3$ be a 
coordinate vector, then we have the so-called Schr\"odinger representation 
given by the Hilbert space $\hilbert{H} = \Lzwei(\mathbb{R}^3)$ of square 
integrable functions with the Lebesgue measure. The state $\psi$ of a particle 
is a time-dependent element of this Hilbert space which fulfils the 
Schr\"odinger equation
\begin{equation}
	\label{QM:Schroedinger}
	i \hbar \frac{d}{dt} \psi(x, t)
	=
	\hat{H}
	\psi(x, t)
\end{equation}
where $\hat{H}$ denotes the Hamilton operator which is given by
\begin{equation}
	\hat{H}
	=
	- \frac{\hbar^2}{2 m}
	\Laplace + V(x),
\end{equation}
where $\Laplace$ denotes the Laplace operators and $V$ the potential.
Similar to the Hamilton function in classical mechanics, the Hamilton operator 
(or more precisely: its spectrum) describes the energy of this system, $m$ is the 
mass of the particle and $V$ the potential. The operators which correspond to 
the coordinate $x_i$ and the momentum $p_j$ of the particle are given by
\begin{align*}
	\left( \hat{q}_i
	\psi \right) (x, t)
	&=
	x_i
	\psi(x, t)
	\\
	\left( \hat{p}_j
	\psi \right) (x, t)
	&=
	- i \hbar 
	\frac{\partial}{\partial x_j}
	\psi(x, t).
\end{align*}
This turns Equation \eqref{QM:Schroedinger} into
\begin{equation*}
	i \hbar \frac{d}{dt} \psi(x, t)
	=
	\left( \frac{\hat{p}^2}{2m} + V(x) \right)
	\psi(x, t)
\end{equation*}
and looks therefore similar to the classical time evolution.
Since the  state is a function of the coordinate $x$, one calls this the 
representation in \emph{position space}. Sometimes it is more convenient to 
describe a state in terms of its momentum and one changes to \emph{momentum 
space} via Fourier transformation. Both representations are equivalent, very 
similar to the variables $q$ and $p$ in the Hamilton form of classical 
mechanics, from which this picture is very much inspired. Unlike the classical 
position and momentum observables being functions, the quantum mechanical 
observables are \emph{operators} and do not commute any longer. One can check 
these by calculating commutators and finds the so-called Heisenberg or 
\emph{classical commutations relations}
\begin{align}
	[p_i, p_j]
	& =
	[q_i, q_j]
	=
	0,
	\\
	[q_i, p_j]
	& =
	i \hbar \delta_{ij}.
\end{align}
These relations are the starting point for all quantization theories. The main 
goal is always to find an algebra which fulfils these (or maybe equivalent) 
relations and this is what makes a quantum algebra so different from a 
classical one. Physically spoken, this noncommutativity means that a 
measurement of the position of a particle influences its momentum and vice 
versa, such that it becomes impossible to measure both, position and momentum, 
at the same time with arbitrarily high precision. This still very vague 
statement can be made precise using the Hilbert norm and the Cauchy-Schwarz 
inequality and yields the famous Heisenberg uncertainty relation
\begin{equation*}
	\Laplace A_\psi
	\Laplace B_\psi
	\geq
	\frac{1}{2}
	\left|
		\SP{\psi, [A, B] \psi}
	\right|.
\end{equation*}
Here $A, B$ are observables and $\Laplace A_\psi, \Laplace B_\psi$ denote their 
standard deviations in the state $\psi$. For the position and momentum of a 
particle, this gives
\begin{equation*}
	\Laplace \hat{p}
	\Laplace \hat{q}
	\geq
	\frac{\hbar}{2}.
\end{equation*}
However, the Schr\"odinger picture is not the only possibility to describe the 
behaviour of a quantum system. Heisenberg proposed a description in which not 
the particles, but the observables depend on the time. Since quantum mechanics,
time evolution is described by a one-parameter group of unitary operators, we have
\begin{equation*}
	\psi(x,t)
	=
	U(t - t_0) \psi(x, t_0),
\end{equation*}
at least in systems where the interactions are not explicitly time dependent. 
This makes it possible to apply unitary transformations to the states 
as well as to the operators 
\begin{equation*}
	\psi(x)
	=
	\psi(x, t_0)
	=
	U(t_0 - t) \psi(x, t)
	\quad \text{ and } \quad
	A(t)
	=
	U(t_0 - t)
	A
	U(t - t_0).
\end{equation*}
One gets the so called \emph{Heisenberg picture}, in which the 
time-dependence is not expressed in the states but in the operators. Both 
pictures are equivalent and it is just a question of convenience which one is 
used. In our example, the unitary operators are given by
\begin{equation*}
	U(t - t_0)
	=
	\E^{- \frac{i}{\hbar} \hat{H} t},
\end{equation*}
since they describe the time evolution of solutions of the Schr\"odinger 
equation. If one differentiates an operator with respect to $t$ in the 
Heisenberg picture, one finds the \emph{Heisenberg equation}
\begin{equation*}
	\frac{d}{dt}
	A(t)
	=
	\frac{1}{i \hbar}
	[A(t), H].
\end{equation*}
This reminds us of course of the time evolution in terms of Poisson brackets in 
classical mechanics. The idea that both objects, the quantum mechanical commutator 
and the Poisson bracket, should somehow be seen as counterparts to each other, is 
called the \emph{correspondence principle}. It is the second staring point for 
nearly all theories of quantization.

\section{Making Things Noncommutative: Quantization}
\label{sec:chap2_Quantization}

\subsection{The Task}
\label{subsec:chap2_Task}

We have had a brief overview over the mathematical formulation of classical and 
quantum mechanics. On one hand, if we compare for example the way in which states 
of a physical system are described (a point in the cotangent bundle of a manifold 
and a vector in a Hilbert space) we will see that they are very different. On the 
other hand, we have seen that both theories have certain things in common: the 
time evolution can be described in a similar way. We want to give a short list of 
the main concepts of both theories and compare them (a much more detailed list and 
discussion can be found in Chapter 5 of \cite{waldmann:2007a}):
\bgroup
\renewcommand{\arraystretch}{1.6}
\begin{center}
	\begin{tabular}
	{lll}
		~ 
		&
		\textbf{Classical} 
		&
		\textbf{Quantum}
		\\
		\parbox{4cm}
		{
			Observables
		}
		&
		\parbox{5cm}{
			Poisson(-*)-algebra $\algebra{A}_{Cl}$
			of \\
			smooth functions $\Cinfty(M)$ \\
			on a Poisson manifold
		}
		&
		\parbox{5cm}{
			*-algebra $\algebra{A}_{QM}$
			of self-adjoint operators
			on a Hilbert space $\hilbert{H}$
		}
		\\
		\parbox{4cm}
		{
			Measurable 
			Values
		}
		&
		$\spec (f) \subseteq \mathbb{R}$
		&
		$\spec (A) \subseteq \mathbb{R}$
		\\
		States
		&
		\parbox{5cm}{
			Points in the
			phase space
		}
		&
		\parbox{5cm}{
			Vectors in a
			Hilbert space
		}
		\\
		Time evolution
		&
		Hamilton function $H$
		&
		Hamilton operator $\hat{H}$
		\\
		\parbox{4cm}{
			Infinitesimal\\
			time evolution
		}
		&
		$
		\frac{d}{dt} f(t)
		=
		\{f(t), H\}
		$
		&
		$\frac{d}{dt} A(t)
		=
		\frac{1}{i \hbar}
		[A(t), H]
		$
	\end{tabular} 
\end{center}
\egroup
When we look at this table, we see that it is probably a good guideline to 
understand the observable algebras as the main concepts of mechanics, rather than 
the states. We know that the classical theory emerges as a limit of many quantum 
systems, such that the physical constant $\hbar$ becomes small in enough (in 
comparison to the characteristic action scale of the system) to be neglected. How 
to construct a correspondence the other way round is the question which the theory 
of quantization addresses. Physicists used to solve this by ``making a hat on the 
variables $p$ and $q$ and saying that they the do not commute any more''. 
Surprisingly enough, this approach, which is called \emph{canonical quantization} 
in physics, works for most of the simple examples. In particular, this is the 
case when we do not have high powers of $\hat p$ and $\hat q$ and no 
mixed terms. That this idea is however not canonical in a \emph{mathematical 
sense} is almost needless to say and of course it causes severe problems when 
higher and mixed terms in position and momentum appear. Briefly stated, the 
problem is that in the classical theory, the expressions $q^2p, qpq, pq^2$ 
describe the same polynomial, but in the quantum theory, they do not. So we must 
think about the question, to which operators we want to map mixed polynomials, and 
canonical quantization does not provide an answer. To create a mathematical theory 
of quantization, one needs to impose an ordering on the $\hat p$ and $\hat q$. 
Already Hermann Weyl knew about this problem and looked for ways how to solve it 
by proposing a totally symmetrized expression \cite{weyl:1931a} for such terms. 
This ordering is now known as the Weyl ordering and his formula became the 
starting point for the theory of deformation quantization many years later.

\subsection{A Definition}
First, we want to take a step back and try to formalize what we actually mean by a 
quantization. See also Section 3.2 in \cite{esposito:2015a} or Section 5.1.2 in
\cite{waldmann:2007a} for a good introduction to this topic. For us, a 
quantization should be a correspondence
\begin{equation*}
	\mathcal{Q}
	\colon
	\algebra{A}_{\mathrm{Classical}}
	\longrightarrow
	\algebra{A}_{\mathrm{Quantum}}
\end{equation*}
between a commutative algebra $\algebra{A}_{\mathrm{Classical}}$ and a 
noncommutative algebra $\algebra{A}_{\mathrm{Quantum}}$, which is a 
``bijection'' in some sense: all classical observables appear as classical 
limits from quantum observables, hence we should expect $\mathcal{Q}$ to be 
''injective''. On the other hand, if the quantum algebra was bigger than the 
classical algebra, the whole concept of quantization would be pointless, since 
we could never recover all observables, thus we want $\mathcal{Q}$ to be 
``surjective'', too. The correspondence should also keep somehow track of the 
physical meaning of our observables, or at least tell us, what the observable 
$\mathcal{Q}(f)$ should be for some $f \in \algebra{A}_{\mathrm{Classical}}$. 
Moreover, $\mathcal{Q}^{-1}$ should behave like the classical limit, since this 
is the rather well-understood part of both. Finally, we want the correspondence to
keep associative structures: the classical algebra is associative and we want 
the product in our noncommutative algebra to correspond to the concatenation of 
operators on a Hilbert space which is associative, too. One possibility of 
formulating this would be the following: in the very prototype of a classical 
observable algebra, a symplectic vector space of dimension $2n$, 
we want that the following three axioms are fulfilled:
\begin{enumerate}
	\item[(Q1)]
	We want a ``map'', which allows a representation in the standard 
	picture of quantum mechanics:
	$\mathcal{Q}(1) = \Unit$, 
	$\mathcal{Q}(q^i) = \hat{q}^i$, 
	$\mathcal{Q}(p_j) = - i \hbar \frac{\partial}{\partial q^j}$.
	
	\item[(Q2)]
	The correspondence principle should be fulfilled:
	$[\mathcal{Q}(f), \mathcal{Q}(g)] = i \hbar \mathcal{Q}(\{f, g\})$, 
	for all $f, g \in \algebra{A}_{\mathrm{Classical}}$.
	
	\item[(Q3)]
	For mathematical simplicity, $\mathcal{Q}$ should be linear.
\end{enumerate}
We finally split up the whole process: in a first step, we \emph{quantize} the 
system and in a second step, we look for \emph{representations} on a Hilbert 
space.
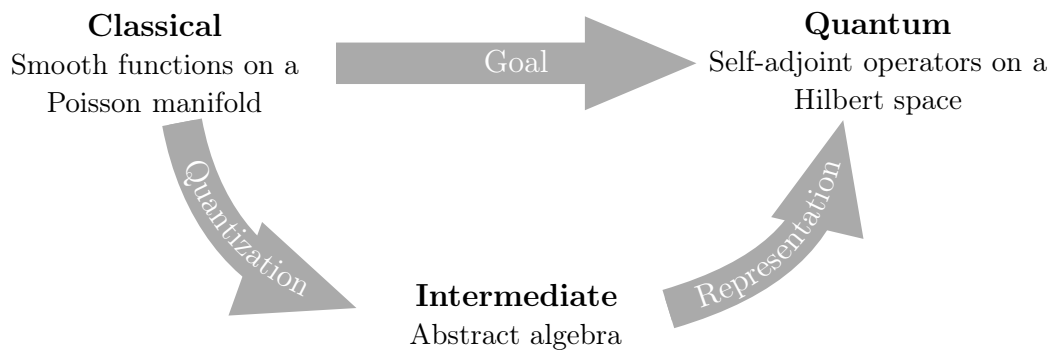
\begin{center}
    \begin{tikzpicture}
    	\tikzstyle{myarrowq} = 
    	[line width = 3.4mm,
    	draw = {rgb:black,1;white,2},
    	-triangle 45,
    	shorten >= -3mm,
    	postaction = {
    		draw, 
    		line width = 5.3mm, 
    		shorten >=4.5mm,
    		-},
    	postaction = {
    		decorate,
    		decoration = {
    			text along path, 
    			text = {|\large \color{white}| Quantization}, 
    			text align = {align = center},
    			raise = -0.7ex
    			}
    		},
    	]
    	\tikzstyle{myarrowg} = 
    	[line width = 3.2mm,
    	draw = {rgb:black,1;white,2},
    	-triangle 45,
    	postaction = {
    		draw, 
    		line width = 5.5mm, 
    		shorten >=6.5mm,
    		-},
    	postaction = {
    		decorate,
    		decoration = {
    			text along path, 
    			text = {|\large \color{white}| Goal }, 
    			text align = {align = center},
    			raise = -0.7ex
    			}
    		}
    	]
    	\tikzstyle{myarrowr} = 
    	[line width = 3.2mm,
    	draw = {rgb:black,1;white,2},
   		shorten >=-5mm,
   		shorten <=-4mm,
    	-triangle 45,
    	postaction = {
    		draw, 
    		line width = 5.3mm, 
    		shorten >=1.5mm,
    		-},
    	postaction = {
    		decorate,
    		decoration = {
    			text along path, 
    			text = {|\large \color{white}|Representations}, 
    			text align = {align = left},
    			raise = -0.6ex
    			}
    		}
    	]
    		    	
        \matrix (m)[
        matrix of nodes,
        row sep=3em,
        column sep=0pt,
		ampersand replacement=\&
        ]
        {
        	\parbox{4.5cm}
        	{
        	  \begin{center}
        		{\large
        			\textbf{Classical}
        		}
        		\\
        		{\normalsize
        			Smooth functions on
        			a Poisson manifold
        		}
        	  \end{center}
        	}
        	\&
        	\&
        	\parbox{4.5cm}
        	{
        	  \begin{center}
        		{\large
        			\textbf{Quantum}
        		}
        		\\
        		{\normalsize
        			Self-adjoint operators
        			on a Hilbert space
        		}
        	  \end{center}
        	}
        	\\
        	\&
        	\parbox{4.5cm}
        	{
        	  \begin{center}
        		{\large
        			\textbf{Intermediate}
        		}
        		\\
        		{\normalsize
        			Abstract algebra
        		}
        	  \end{center}
        	}
        	\&
        	\\
        };
        \draw
        [-stealth]
        (m-1-1) edge[myarrowg]					(m-1-3) 
        		edge[myarrowq, bend right=30]	(m-2-2)
        (m-2-2) edge[myarrowr, bend right=30]	(m-1-3)
		;
    \end{tikzpicture}
\end{center}

\subsection{Different approaches}
There exist various ideas about how to build up such a scheme. For example, in 
axiomatic quantum field theory one usually wants the quantum algebra to be a 
$C^*$-algebra (see for example \cite{haag:1993a} or 
\cite{baer.ginoux.pfaeffle:2007a}). The reason is that the bounded operators 
$\Bounded(\hilbert{H})$ naturally form such an algebra (and actually even more 
than that). Unfortunately, most of the operators in quantum mechanics are 
unbounded. This problem is cured by looking at the exponentiated operators
\begin{equation*}
	A
	\longmapsto
	\E^{i t A},
	\quad t \in \mathbb{R}.
\end{equation*}
This way, one gets a one parameter group for $t \in \mathbb{R}$ and unitary 
operators are clearly bounded. Moreover, there is a very nice correspondence 
between $C^*$-algebras and locally compact Hausdorff spaces, which is known as 
the Gelfand-Naimark theorem. Roughly stated, this means that there is an 
equivalence of categories between compact Hausdorff spaces and unital commutative 
$C^*$-algebras: every such space gives rise to a commutative algebra of 
complex-valued, continuous functions, which form a $C^*$-algebra. On the other 
hand, from every unital commutative $C^*$-algebra $\algebra{A}$ one can construct 
a compact Hausdorff space $X$ such that $\algebra{A} \cong \Stetig(X)$. This 
correspondence can be extended to open and closed subsets of the space, to 
homeomorphism, locally compact spaces, their compactification, vector bundles and 
even more. There is a strong link between commutative algebra and topology or 
geometry. One possibility to think of a quantized system is to think of continuous 
or smooth functions on a noncommutative space, which then should correspond to a 
noncommutative $C^*$-algebra. This idea leads to noncommutative geometry, which 
is mostly due to Alain Connes. A very detailed, but not necessarily easy to read 
book is \cite{connes:1994a} by Connes, another and rather brief introduction is 
for example \cite{varilly:2006a}.

A different approach, called geometric quantization, tries to fulfil all 
of the three axioms (Q1) - (Q3) from the previous part. Unfortunately, this 
causes problems: already for a symplectic vector space, it is impossible to have a 
one-to-one correspondence of the Poisson bracket and the quantum mechanical 
commutator. This is known as the Groenewold-van Hove theorem, which was found 
around 1950 \cite{vanhove:1951a, groenewold:1946a}. More precisely spoken: no 
representation of the Lie algebra, which is generated by the $q^i$ and $p_j$ and 
which is defined by the classical commutation relations, can be extended 
irreducibly and faithfully to the commutator Lie algebra which comes from the 
associative unital algebra which is generated by the $q^i$ and $p_j$. Thus 
geometric quantization restricts to smaller observable algebras, which are not 
problematic. At its present state, this approach is however still far from being a 
general theory, but it provide a procedure for number of exemplary cases. Its 
ideas are mostly due to Souriau \cite{souriau:1970a}, Kostant and Segal.

There are also different approaches. Berezin proposed a quantization scheme for 
particular K\"ahler manifolds, \cite{berezin:1975a, berezin:1975b, berezin:1975c}. 
Still, new ideas keep coming up, but in the following, we want to concentrate on a 
particular type of quantization, which is called \emph{deformation quantization}.

\section{Deformation Quantization}
\label{sec:chap2_DQ}

\subsection{The Concept}
\label{subsec:chap2_Concept}

Deformation quantization tries to realize the three points (Q1) and (Q2) from the 
previous section, but weakens the third. If it is not possible to have such a 
correspondence exactly, we will at least want to have it 
\emph{asymptotically}. The motivating example is the Weyl quantization which we 
already talked about. There are actually two such formulas, that can be given: the 
first maps a function $f$ in the variables $q, p \in \mathbb{R}^n$ to a 
differential operator on $\mathbb{R}^n$, which is actually a formal power series 
in the parameter $\hbar$. For polynomial functions, this series is a sum (more 
precisely: a polynomial again) and well defined, for general functions this will 
really be a formal power series, hence a truly infinite sequence, which a priori 
has no analytical, but just an algebraic meaning. 
The second is given by an integral formula and holds for another class of 
functions (Schwartz functions), but one gets the first formula out of the second 
as an asymptotic expansion for $\hbar \longrightarrow 0$. With these 
quantizations, one can also define a product of two functions $f, g$, which will 
necessarily take those two functions to a formal power series in $\hbar$. Moyal 
showed that the commutator of this product can be understood as a series which 
approximates the quantum mechanical commutator \cite{moyal:1949a}. The reason why 
seemingly all of a sudden power series appear is the following: if one wants the 
correspondence principle to be asymptotically fulfilled, i.e.
\begin{equation*}
	[\mathcal{Q}(f), \mathcal{Q}(g)]
	=
	i \hbar \mathcal{Q}(\{ f, g \})
	+ \mathcal{O}(\hbar^2)
\end{equation*}
and the multiplication of these quantized functions to be associative 
(as needed for representations on the Hilbert space), one will necessarily get 
higher and higher orders in $\hbar$. This iteration can never be stopped without 
loosing associativity. Motivated by this observation, a group of mathematicians, 
the so called Dijon-school, started working out this idea of products as formal 
power series. They understood quantization as a \emph{deformation} of the 
commutative product by a formal parameter (mostly called $\hbar$, $\lambda$ or 
$\nu$, in this work we will call it $z$ from now on), which controls the 
noncommutativity of the theory. These deformed products should moreover fulfil 
some compatibility conditions with the classical theory. This was the 
hour of birth of deformation quantization. The main characters of this group 
were Flato, Lichnerowicz, Bayen, Fr{\o}nsdal and Sternheimer, who published 
their ideas in the late 70's \cite{bayen.et.al:1977a, bayen.et.al:1978a} and 
gave a first definition of a star product. These two articles became the starting 
point for what has now become a rich and fruitful theory. The deformed products 
are tits corner stone and one defines them in the following way.
\begin{definition}[Star Product]
	\label{Def:StarProduct}
	Let $(M, \{\cdot, \cdot\})$ be a Poisson manifold over a field 
	$\mathbb{K}$ ($\mathbb{K} = \mathbb{R}$ or $\mathbb{C}$).
	A star product on $M$ is a bilinear map
	\begin{equation*}
	    \star_z \colon 
    	\Cinfty(M) 
    	\times 
	    \Cinfty(M) 
	    \longrightarrow
	    \Cinfty(M) \llbracket z \rrbracket
	    , \
	    (f,g) 
	    \longmapsto 
	    f \star_z g 
	    =
	    \sum\limits_{n = 0}^\infty 
	    z^n C_n(f,g)
	\end{equation*}
	such that its $\mathbb{K}\llbracket z \rrbracket$-linear extension to 
	$\Cinfty(M) \llbracket z \rrbracket$ fulfils the following properties:
	\begin{definitionlist}
		\item
		$\star_z$ is associative.
		
		\item
		$C_0(f, g) = f \cdot g,\ \forall_{f,g \in \Cinfty(M)}$ (Classical limit).
		
		\item
		$C_1(f, g) - C_1(g, f)= z \{f, g\},\ \forall_{f,g \in \Cinfty(M)}$
		(Semi-classical limit).
		
		\item
		$1 \star_z f = f \star_z 1 = f,\ \forall_{f \in \Cinfty(M)}$.
	\end{definitionlist}
	If the $C_n$ are bidifferential operators, the star product is said to be 
	differential and if the order of differentiation of the $C_n$ does not 
	exceed $n$ in both arguments, a differential star product is said to be 
	natural. Moreover, we will say that a star product $\star_z$ is of Weyl-type,
	if $\cc{f \star_z g} = \cc{g} \star_z \cc{f}$ for all $f,g \in \Cinfty(M)$
	where $\cc{\phantom{a}}$ denotes the complex conjugation.
\end{definition}
They also defined a notion of equivalence of star products. The idea behind is 
that two equivalent star products should give rise to the same physics.
\begin{definition}[Equivalence of Star Products]
	Two star products $\star_z$ and $\widehat{\star}_z$ for a Poisson manifold 
	$(M, \{\cdot, \cdot\})$ are said to be equivalent, if there is a formal 
	power series
	\begin{equation*}
		T
		=
		\id +
		\sum\limits_{n=0}^{\infty}
		z^n T_n
	\end{equation*}
	of linear maps $T_n \colon \Cinfty(M) \longrightarrow \Cinfty(M)$, which 
	extends $\mathbb{K}\llbracket z \rrbracket$-linearly to 
	$\Cinfty(M) \llbracket z \rrbracket$, such that the following statements hold:
	\begin{equation*}
		f \star_z g
		=
		T^{-1}
		\left(
			T(f) \widehat{\star}_z T(g)
		\right)
		, \ 
		\forall_{f,g \in \Cinfty(M) \llbracket z \rrbracket}
		\quad \text{ and } \quad
		T(1) 
		= 
		1.
	\end{equation*}
	For differential or natural star products, we accordingly speak of 
	differential or 	natural equivalences.
\end{definition}
Note that these definitions are purely algebraic, since we do not ask for the
convergence of those power series. Hence the theory which was developed from 
this in the following years is also a mostly algebraic theory. Like for the Weyl 
product, there were also integral formulas around for other types of star products 
for which one can also make sense of convergence. But as already pointed out, one 
has to strongly restrict the algebra of functions, for example to the Schwartz 
space.

\subsection{A Mathematical Theory}
\label{subsec:chap2_MathTheory}
Deformation quantization is a good example for a mathematical theory, which is 
motivated by a physical idea. It is not really talking about a physical problem, 
since the world does not need to be quantized -- it already is. It is talking 
about a mathematical problem: how to recover the quantized (and very 
counter-intuitive) mathematics, which describe the world on very small scales out 
of the classical (and much more intuitive) mathematics, which describe the world 
on our scale? Deformation quantization tries to give an answer to that, but is 
unfortunately (at least at its present state) still far from doing so completely. 
Like many such theories, it started of from a more or less precise physical 
background and developed into something very different: a mostly algebraic, purely 
mathematical theory. That is also due to the fact that the questions, which had to 
be answered in the beginning, were very hard and of mathematical nature. It took a 
lot of time to find answers and meanwhile, the mathematicians working on them
were interested other aspects of the theory. We want to give a short overview of
those questions and their answers very briefly and summarize a bit the history of 
deformation quantization next. A more detailed summary can be found in section 6.1 
of \cite{waldmann:2007a}.

The prototype is, as already mentioned, a symplectic vector space, for which 
Weyl proposed a star product with a certain (symmetric) ordering (although the 
definition of a star product did not exist at his time). However, other 
orderings are possible: one can have a standard or an anti-standard ordering, 
where the $\hat p$'s are all ordered to the right or to the left, respectively, or 
something in  between. One of the first questions was, if one could also construct 
star products on symplectic manifolds and if these products will be standard or 
Weyl ordered, if they will be differential or natural and so on. Locally, the 
answer was yes, but it took some time and many small steps, until DeWilde and 
Lecomte could show that every cotangent bundle of a smooth manifold has star 
products \cite{dewilde.lecomte:1983a} and then extended this result to arbitrary 
smooth manifolds \cite{dewilde.lecomte:1983b}. Another proof was given 
independently from that by Omori, Maeda and Yoshioka 
\cite{omori.maeda.yoshioka:1991a} and then by Fedosov, who presented a simple 
and very geometric construction \cite{fedosov:1994a}, which always give rise to 
natural star products, as shown in \cite{bordemann.waldmann:1997a} 
or more generally in \cite{gutt.rawnsley:2003a}. Moreover, every star product on a 
symplectic manifold is equivalent to a Fedosov star product 
\cite{bertelson.cahen.gutt:1997a}. A lot of 
results were found for K\"ahler manifolds and also the already mentioned 
procedure, which is due to Berezin, gives rise to star products. There are 
moreover standard, anti-standard ordered and many other types of star products 
on every cotangent bundle. The next question was the one concerning the 
equivalence classes of star products in the symplectic case. One can show that 
locally, two star products on a symplectic manifold are always equivalent. Hence 
a classification result should depend on global phenomena. Indeed, this is the 
case and it can be shown that star products on a symplectic manifold are 
classified by its second deRham cohomology $\HdR^2(M)$. This result is due to 
Deligne \cite{deligne:1995a}, a different proof was given by Nest and Tsygan 
\cite{nest.tsygan:1995a}, another one by Bertelson, Cahen and Gutt 
\cite{bertelson.cahen.gutt:1997a}. More precisely: the choice of an equivalence 
class of closed, nondenerate 2-forms $\omega \in \Formen^2(M)$ determines a 
Fedosov star product and from every Fedosov 
star product one can calculate such an equivalence class. This already 
determines all star products on symplectic manifolds, since every star product 
on such a manifold is equivalent to a Fedosov star product. The case of Poisson 
manifolds took longer and was much harder to solve, since associativity turned 
out to be a complicated condition to fulfil. There were some examples of star 
products known for particular Poisson structure, like the Gutt star product 
\cite{gutt:1983a}, which was also found by Drinfel'd \cite{drinfeld:1983a} 
independently, but the general existence (and also the classification) result was 
proven by Kontsevich \cite{kontsevich:1997:pre, kontsevich:2003a} many years 
later. His classification result is known as the formality theorem and needs the 
notion of $L_{\infty}$-algebras, which are fairly involved objects. He gave an 
explicit construction, how star products can be built out of Poisson brackets on 
$\mathbb{R}^d$. This construction was extended by Cattaneo, Felder and Tomassini 
to Poisson manifolds \cite{cattaneo.felder.tomassini:2002b} and indepedently from 
that by Dolgushev \cite{dolgushev:2005a}. Another and easier formulation of the 
Kontsevich construction on $\mathbb{R}^d$ in terms of operads was later given by 
Tamarkin \cite{tamarkin:2003a}.

\subsection{From Formal to Strict}
\label{subsec:chap2_Formal2Strict}

So far, one could say that the big cornerstones of the theory are already there 
and that it is somehow ``finished''. For two reasons, this is not the case. 
First, a mathematical theory is never actually ``finished'', since there are 
always a lot of new things which can be found. There are still many different 
types of star products to classify, like invariant or equivariant star products 
in the case that one has Lie group or Lie algebra actions. A very recent result 
concerning the classification of equivariant star products on symplectic 
manifolds is, for example, due to Reichert and Waldmann 
\cite{reichert.waldmann:2015a:pre}. Second, the theory of deformation 
quantization still has a different aspect: all we talked about so far was purely 
algebraic and there is no notions of convergence of these formal power series. 
If some day, this theory shall have a real drawback on physics, it will be 
necessary to talk about the convergence properties of these star products, since 
in physics $\hbar$ is \emph{not} a formal parameter but a constant with a 
dimension and a fixed value and therefore the question of convergence matters. 
When we dace those problems and speak about continuous star products, we leave the 
field of \emph{formal} deformation quantization and come to \emph{strict} 
deformation quantization.

Although it is closer to physics, strict deformation quantization is still a 
mathematical theory. There are two different approaches to it: we 
already mentioned integral formulas, which allow to speak of continuous star 
products. The second approach uses the formal power series instead and wants to 
construct a topology on the polynomial algebra, such that the star product becomes 
continuous. Then one completes the tensor algebra over the vector space to a 
subalgebra of the smooth functions, on which the star product will still be 
continuous.

The first approach is mostly due to Rieffel, who developed these ideas in some of 
his papers \cite{rieffel:1989a, rieffel:1990c, rieffel:1993a}. He wants to 
realize a strict deformation quantization by actions of an abelian Lie group on a 
$C^*$-algebra of classical observables. Later Rieffel formulated a list of open 
questions, which strict deformation quantization should try to answer 
\cite{rieffel:1998a} in the next years. His approach was carried on by Bieliavsky 
and Gayral \cite{bieliavsky:2002a, bieliavsky.gayral:2015a}, who extended these 
concepts to much more general Lie groups and different manifolds. To get 
reasonably big observable algebras, they used oscillatory integrals and pushed 
this theory forward. A similar idea was realized by Natsume, Nest and 
Peter \cite{natsume.nest.peter:2003a}, who could show that under certain 
topological conditions, symplectic manifolds always admit strict deformations.

The second approach is due to Beiser and Waldmann 
\cite{beiser:2011a, beiser.waldmann:2014a, waldmann:2014a}. They restrict to the 
local situation, that means to Poisson structures on vector spaces. Then, they 
look at the polynomial functions on this vector spaces and try to find 
continuity estimates for them by constructing an explicit locally convex 
topology on the symmetric tensor algebra (which is isomorphic to the polynomial 
algebra). The aim is to make the topology as coarse as possible, to get then a 
large completion and hence a big quantized algebra of 
observables. There are two big advantages of this approach: the first one is 
that it can be applied to infinite dimensional vector spaces, what is necessary 
for quantum field theory, which deals with infinitely many degrees of freedom. 
The second is that we can really speak about all observables, also those 
which will correspond to unbounded operators, without having to exponentiate. In 
this sense, this idea is somehow more fundamental. The disadvantage is, however, 
that it is just a local theory at the moment. The idea is worked out just for 
one type of star products by now, which are star products of exponential type like 
the Weyl product. This means until now one can only control star products on 
symplectic vector spaces which come from a constant Poisson tensor 
\cite{waldmann:2014a}. This theory was carried on in the master thesis of Matthias 
Sch\"otz \cite{schoetz:2014a}, who rephrased it using semi-inner product, which 
are a somehow more physical notion, since one can interpret spaces with such a 
topology as projective limits of pre-Hilbert spaces. This also allows a slightly 
coarser topology and hence a larger completion of the symmetric tensor algebra.

In this work, we will follow the second approach and apply it to 
another type of star product, the Gutt star product, which comes from a linear 
Poisson tensor on a vector space. Of course, this is the next logical step after 
constant Poisson tensors. However, these are also the first \emph{non-symplectic}
Poisson structures, which will be strictly quantized this way. Thus this 
master thesis really contributes something new to the theory of strict 
deformation quantization: a second example, in which Waldmann's locally convex 
topology on the tensor algebra leads to a continuous star product, when 
considered as a power series and not as an integral. Note that this also has a 
certain effect on Lie theory: the result can be seen as a functorial 
construction for a locally convex topology on the universal enveloping algebra 
$\algebra{U}(\lie{g})$ of a (possibly infitely-dimensional) Lie algebra 
$\lie{g}$. Therefore they may have applications to, for example, the 
representation theory of $\algebra{U}(\lie{g})$.

% Chapter 3
%

%
% Chapter 3 of my master thesis:
% The first real chapter
%

\chapter{Algebraic Preliminaries}

\section{Linear Poisson Structures in Infinite Dimensions}
\label{sec:chap3_LinearPoisson}

As we have seen before, there has already been done some work on how to 
strictly quantize Poisson structures on vector spaces. Star products of
exponential type on locally convex vector spaces were topologized by Stefan 
Waldmann in \cite{waldmann:2014a} and then investigated more closely by 
Matthias Schötz in \cite{schoetz:2014a}. The Poisson tensor which cooresponds to 
the Poisson bracket in these cases, is constant. Hence, as a the next step, we 
want to investigate linear Poisson structures on locally convex vector spaces. 
This will give a new big class of Poisson structures, which will be deformable 
in a strict way. Before we do so in the rest of this master thesis, we 
recall briefly some basics on linear Poisson structures.

In the following, we always consider a vector space $V$ and study Poisson 
structures on the coordinates which are elements of the dual space $V^*$. In order 
to cover most of the physically interesting examples by our reflections, we want 
to allow vector spaces of very general type and therefore assume that $V$ is a 
locally convex vector space. Every finite-dimensional vector space is normable and 
hence locally convex, so it fits in this framework. In this case, it is clear what 
$V^*$ should be and there is just one interesting topology on it. For 
infinite-dimensional spaces, the situation is more delicate: we have to think 
about what coordinates should be and how a Poisson structure on them could look 
like. A priori, it is not 
clear which dual we should consider: the algebraic dual $V^*$ of all linear 
forms on $V$, or the topological dual $V'$ which contains just the 
continuous linear forms? Here, one could argue that only $V'$ is of real 
interest, since otherwise we would encounter the very strange effect of 
having discontinuous polynomials, and the aim of constructing a continuous 
star product on them seems somehow pointless. 
But even if we stick to $V'$, the question of the topology still remains: 
do we want to consider the weak or the strong topology there and why one of 
them should be more interesting. In any case, we have to choose a topology 
on this space. Once this is done, we have to think about a good notion of 
Poisson structures in this context. However, we encounter quite a number of 
question, which have no trivial answer. For this reason, it is worth 
looking at some equivalent formulations of $\Pol^{\bullet}(V^*)$ in the 
finite-dimensional case, since they may allow better generalizations.

Let $V$ be a finite-dimensional vector space. Now, there is no question 
about the dual or its topology, since $V^* = V'$ is finite-dimensional, 
too, and we deal with polynomials on it. It is easy to see that a linear Poisson 
structure on $V^*$ is something very familiar: it is equivalent to a Lie algebra 
structure on $V$.
\begin{proposition}
	\label{Alg:Prop:LinPoissonIsLieAlg}
	Let $V$ be a vector-space of dimension $n \in \mathbb{N}$ and $\pi \in 
	\Secinfty(\Anti^2(TV^*))$. Then the two following things are 
	equivalent:
	\begin{propositionlist}
		\item
		$\pi$ is a linear Poisson tensor.
		
		\item
		$V$ has a uniquely determined Lie algebra structure.
	\end{propositionlist}
\end{proposition}
\begin{proof}
	We choose a basis $\{e_1, \ldots, e_n\} subset V$ and denote its dual basis 
	$\{e^1, \ldots, e^n\} \subset V^*$. Then we call the linear coordinates in 
	these bases $x_1, \ldots, x_n \in \Cinfty(V^*)$ and $\xi^1, \ldots, 
	\xi^n \in \Cinfty(V)$, such that for all $\xi \in V, x \in V^*$
	\begin{equation*}
		\xi
		=
		\xi^i(\xi) e_i
		\quad \text{ and } \quad
		x
		=
		x_i (x) e^i.
	\end{equation*}
	In these coordinates, the Poisson tensor reads
	\begin{equation*}
		\pi
		=
		\frac{1}{2}
		\pi_{ij}(x)
		\frac{\partial}{\partial x_i}
		\wedge
		\frac{\partial}{\partial x_j},
	\end{equation*}
	where $\pi$ is linear in the coordinates. Since $\pi_{ij}(x)$ is linear in the 
	coordinates, we can write it as
	\begin{equation*}
		\pi_{ij}(x)
		=
		c_{ij}^k x_k.
	\end{equation*}
	This gives for $f,g \in \Cinfty(V^*)$
	\begin{equation}
		\label{Alg:KksInCoordinates}
		\{f, g\} (x)
		=
		\pi(df, dg)(x)
		=
		x_k c_{ij}^k 
		\frac{\partial f}{\partial x_i}
		\frac{\partial g}{\partial x_j},
	\end{equation}
	using the identification $T^*V^* \cong V^{**} \cong V$.
	But now, the statement is obvious, since antisymmetry of $\pi$ means
	antisymmetry of the $c_{ij}^k$ in the indices $i$ and $j$ and the 
	Jacobi identity for the Poisson tensor gives
	\begin{equation}
		\label{Alg:JacobiInStructureConst}
		c_{ij}^\ell c_{\ell k}^m
		+
		c_{j k}^\ell c_{\ell i}^m
		+
		c_{ki}^\ell c_{\ell j}^m
		=
		0
	\end{equation}
	for all $i, j, k, m$, since it must be fulfilled for all smooth 
	functions. Vicely versa, \eqref{Alg:JacobiInStructureConst} ensures the 
	Jacobi identity of $\pi$ in \eqref{Alg:KksInCoordinates}. Hence the map
	\begin{equation}
		\label{Alg:LieBracketOfKks}
		[ \cdot, \cdot ]
		\colon
		V
		\times
		V
		\longrightarrow
		V
		\quad
		(e_i, e_j)
		\longmapsto
		c_{ij}^k e_k
	\end{equation}
	defines a Lie bracket, since the $c_{ij}^k$ are antisymmetric and 
	fulfil the Jacobi identity and are therefore structures constants. 
	Conversely, the structure constants of a Lie algebra on $V$ define a 
	Poisson tensor on $V^*$ via \eqref{Alg:KksInCoordinates}.
\end{proof}
A more detailed explanation of Poisson manifolds in general, their 
correspondence to Lie algebroids this special correspondence of linear Poisson 
structures and a vector space and Lie algebras can be found in the textbook of 
Waldmann \cite{waldmann:2007a} or in his lecture notes 
\cite{waldmann:2015a:script}. The fact that $V$ is a Lie algebra in the situation 
we want to consider motivates a change of notation: from now on, we will call the 
original vector space $\lie{g}$, which is more intuitive for a Lie algebra.
Since this kind of Poisson systems has a particular structure, there is a proper 
name for them.
\begin{definition}[Kirillov-Kostant-Souriau bracket]
	\label{Def:KKS}
	Let $\lie{g}$ be a finite-dimensional Lie algebra. We call the Poisson 
	bracket, which is given on $g^*$ by Equation~\eqref{Alg:KksInCoordinates}, 
	the Kirillov-Kostant-Souriau bracket and denote it by 
	$\{ \cdot , \cdot \}_{KKS}$.
\end{definition}

Proposition~\ref{Alg:Prop:LinPoissonIsLieAlg} gives us a hint how we could
think of infinite-dimensional vector spaces with linear Poisson 
tensor: We take $\lie{g}$ to be an infinite-dimensional Lie algebra, which 
gives something like a linear Poisson structure on $\lie{g}'$.
If we chose directly $\lie{g}'$ to have a linear Poisson tensor, we would get a 
Lie algebra structure on $\lie{g}''$. Of course, we could think of using this 
structure on $\lie{g}$, since it canonically injects into $\lie{g}''$.
The problem is that in general, this will \emph{not} be closed: taking the Lie 
bracket of $\xi, \eta \in \lie{g}$, we might drop out of $\lie{g}$ and have
$[\xi, \eta] \subseteq \lie{g}'' \backslash \lie{g}$. Usually, such
a behaviour will not be of interest, since the algebras of 
physical systems are closed objects and the double-dual is not what we are aiming 
for. This is why we will translate  the term ``linear Poisson structure on 
$\lie{g}^*$'' by ``$\lie{g}$ is a Lie algebra'' in infinite dimensions. Remark 
however that, from a mathematical point of view, this way of thinking about 
infinite-dimensional Poisson structures is a choice, not a 
logical necessity and other choices would have been possible.

The next task are the polynomials on $\lie{g}'$. As already mentioned, it 
is not easy to find a good generalization for them, since for a locally 
convex Lie algebra $\lie{g}$, even $\lie{g}'$ will be a huge vector space. 
Again, it is helpful to go back to the finite-dimensional case, where we 
have the following result:
\begin{proposition}
	\label{Alg:Prop:PolIsSym}
	Let $\lie{g}$ be a vector space of dimension $n \in \mathbb{N}$. Then 
	the algebras $\Sym^{\bullet}(\lie{g})$ and $\Pol^{\bullet}(\lie{g}^*)$ 
	are canonically isomorphic.
\end{proposition}
\begin{proof}
	Since this is a very well-known result, we just want to sketch the 
	proof briefly: take a basis $\{e_1, \ldots, e_n\}$ of $\lie{g}$ and its 
	linear coordinates $x_1, \ldots, x_n \in \Cinfty(\lie{g}^*)$ with 
	$x_i(x) = e_i(x)$ for $x \in \lie{g}^*$. On homogeneous symmetric 
	tensors this yields the map
	\begin{equation*}
		\mathcal{J}
		\colon
		\Sym^{\bullet}(\lie{g})
		\longrightarrow
		\Pol^{\bullet}(\lie{g}^*),
		\quad
		e_1^{\mu_1} \cdots e_n^{\mu_n}
		\longmapsto
		\xi_1^{\mu_1} \cdots \xi_n^{\mu_n}.
	\end{equation*}
	We see immediately that this is an isomorphism, but note 
	that we have used the identification $\lie{g}^{**} \cong \lie{g}$ via
	\begin{equation*}
		e_i(x)
		=
		\langle x, e_i \rangle.
	\end{equation*}
\end{proof}
In infinite dimensions, the last identification we used in the last step will not 
work in both directions any more, but just in one: we have a canonical injection
$\Sym^{\bullet}(\lie{g}) \subseteq \Pol^{\bullet}(\lie{g}')$, so every 
symmetric tensor still gives a polynomial. Anyway, this gives an idea 
how to avoid speaking about $\Pol^{\bullet}(\lie{g}^*)$ and its topology: we 
restrict from the beginning to $\Sym^{\bullet}(\lie{g})$.
For finite-dimensional systems, both points of view are equivalent,
but in infinite dimensions, this becomes a choice. However,
we have good reasons to think that this is enough: we get a closed 
and reasonably big subalgebra of the polynomials. Moreover, 
the symmetric tensor algebra is defined on infinite-dimensional spaces 
exactly in the same way as on finite-dimensional ones, and the 
construction is identical. 

So finally, we found a suitable way of speaking about our object of 
interest: we replace linear Poisson structures on $\Pol^{\bullet}
(\lie{g}^*)$ by $\Sym^{\bullet}(\lie{g})$.

\section{The Gutt Star Product}
\label{sec:chap3_GuttStar}

The aim of this chapter is to endow the symmetric algebra, and as a consequence 
the polynomial algebra, with a new, noncommutative product. This is possible in 
a very natural way, due to the Poincar\'e-Birkhoff-Witt theorem. It links 
the symmetric tensor algebra $\Sym^{\bullet}(\lie{g})$ of a Lie algebra 
$\lie{g}$ to its universal enveloping algebra $\algebra{U}(\lie{g})$.

\subsection{The Universal Enveloping Algebra}
\label{subsec:chap3_UniversalEnvelopingAlgebra}

If $\algebra{A}$ is an associative algebra, one can construct a Lie algebra 
out of it by using the commutator
\begin{equation*}
	[a,b]
	=
	a \cdot b - b \cdot a
	, \quad
	a,b \in \algebra{A}.
\end{equation*}
This construction is functorial, since it does not only map associative 
algebras to Lie algebras, but also morphisms of the former to those of the 
latter. While constructing a Lie algebra out of an associative algebra is 
easy, the reversed process is more complicated, but also possible. Every 
Lie algebra $\lie{g}$ can be embedded into a particular associative 
algebra, known as the universal enveloping algebra $\algebra{U}(\lie{g})$, 
which is uniquely determined (up to isomorphism) by the universal property: 
for every unital associative algebra $\algebra{A}$ and every Lie algebra 
homomorphism $\phi\colon \lie{g} \longrightarrow \algebra{A}$ using the 
commutator on $\algebra{A}$, one gets a unital homomorphism of associative 
algebras $\Phi \colon \algebra{U}(\lie{g}) \longrightarrow \algebra{A}$ 
such that the following diagram commutes:
\begin{center}
    \begin{tikzpicture}
        \matrix (m)[
        matrix of math nodes,
        row sep=2.5em,
        column sep=8em
        ]
        {
          \algebra{U}(\lie{g}) & \\
           & \algebra{A} \\
          \lie{g} &  \\
        };
        \draw
        [-stealth]
        (m-1-1) 	edge node 
        			[above] 
        			{$\Phi$}
        					(m-2-2)
        	(m-3-1)	edge node
        			[below]
        			{$\phi$}
        					(m-2-2)
        					
        			edge node
        			[left]
        			{$\iota$}
        					(m-1-1)
        	;
    \end{tikzpicture}
\end{center}
where $\iota$ denotes the embedding of $\lie{g}$ into $\algebra{U}(\lie{g})$.
The proof of existence and uniqueness of the universal enveloping algebra 
can be found in standard textbooks on Lie theory like 
\cite{hilgert.neeb:2012a} or \cite{varadarajan:1974a}, and we will not explain it 
here in detail. Just recall that existence is proven by an explicit 
construction: one takes the tensor algebra $\Tensor^{\bullet}(\lie{g})$ and 
considers the two-sided ideal
\begin{equation*}
	\mathfrak{I}
	=
	\langle \xi \tensor \eta - \eta \tensor \xi - [\xi, \eta] \rangle
	, \quad
	\text{ for }\xi, \eta \in \lie{g}
\end{equation*}
inside of it. Then one gets the universal enveloping algebra by the 
quotient
\begin{equation}
	\label{Alg:UnivEnvAlg}
	\algebra{U}(\lie{g})
	=
	\frac{\Tensor^{\bullet}(\lie{g})}{\mathfrak{I}}.
\end{equation}
To avoid confusion, we will always denote the multiplication in 
$\algebra{U}(\lie{g})$ by $\odot$, whereas the commutative product in 
$\Sym^{\bullet}(\lie{g})$ will be denoted without a sign.
It follows from this construction, that $\algebra{U}(\lie{g})$ is a 
filtered algebra
\begin{equation*}
	\algebra{U}(\lie{g})
	=
	\bigcup_{k \in \mathbb{N}}
	\algebra{U}^k(\lie{g})
	, \quad
	\algebra{U}^k(\lie{g})
	= 
	\Big\lbrace
		x 
		= 
		\sum_i
		\xi_1^i \odot \cdots \odot \xi_n^i
	\ \Big| \ 
		\xi_j^i \in \lie{g}
		, \
		1 \leq j \leq n,
		i \in I
	\Big\rbrace.
\end{equation*}
Generally, we just get a filtration, not a graded structure, since the 
ideal $\mathfrak{I}$ is not homogeneous in the symmetric degree. We will 
get a graded structure on $\algebra{U}(\lie{g})$, if and only if $\lie{g}$ 
was commutative. Then $\algebra{U}(\lie{g})$ is isomorphic to the symmetric
tensor algebra and thus also commutative. But $\algebra{U}(\lie{g})$
is much more than an associative algebra: it is also a Hopf algebra, 
since one can define a coassociative, cocommutative coproduct on it by
\begin{equation*}
	\Delta \colon
	\algebra{U}(\lie{g})
	\longrightarrow
	\algebra{U}(\lie{g})
	\tensor
	\algebra{U}(\lie{g})
	, \quad
	\xi
	\longmapsto
	\xi \tensor \Unit
	+
	\Unit \tensor \xi
	, \quad
	\text{ for } \xi \in \lie{g}
\end{equation*} 
which extends to $\algebra{U}(\lie{g})$ via algebra homomorphism, 
as well as an antipode
\begin{equation*}
	S \colon
	\algebra{U}(\lie{g})
	\longrightarrow
	\algebra{U}(\lie{g})
	, \quad
	\xi
	\longmapsto
	- \xi
	, \quad
	\text{ for } \xi \in \lie{g}
\end{equation*}
which extends to $\algebra{U}(\lie{g})$ via algebra antihomomorphism.
We will come back to those two maps and to the Hopf structure in Chapter 7,
when we talk about their continuity. More details on the algebraic
aspect of deformation theory using Hopf algebras can be found in 
\cite{chari.pressley:1994a} and \cite{majid:1995a}, for example.

\subsection{The Poincar\'e-Birkhoff-Witt Theorem}
\label{subsec:chap3_PoincareBirkhoffWitt}

The algebra $\algebra{U}(\lie{g})$ always admits a basis, which must be 
infinite. This result is due to the already mentioned theorem of 
Poincar\'e, Birkhoff and Witt:
\begin{theorem}[Poincar\'e-Birkhoff-Witt theorem]
	\label{Thm:Alg:PBW}
	Let $\lie{g}$ be a Lie algebra with a basis $\mathcal{B}_{\lie{g}} = \{ 
	\beta_i \}_{i \in I}$. Then the set
	\begin{equation*}
		\mathcal{B}_{\algebra{U}(\lie{g})}
		=
		\big\{
			\beta_{i_1}^{\mu_{i_1}}
			\odot \cdots \odot
			\beta_{i_n}^{\mu_{i_n}}
		\ \big| \
			n \in \mathbb{N}, \ 
			i_k \in I
			\text{ with } i_1 \earlier \cdots \earlier i_n 
			\text{ and } \beta_{i_k} \in \mathcal{B}_{\lie{g}}, \
			\mu_{i_1}, \ldots, \mu_{i_n} \in \mathbb{N}
		\big\}
	\end{equation*}
	defines a basis of $\algebra{U}(\lie{g})$.
\end{theorem}
There are different proofs for this statement. While a geometrical proof 
(like e.g. in \cite{waldmann:2015a:script}) is very convenient 
in the finite-dimensional case, a combinatorial argument must be used for 
infinite-dimensional Lie algebras. Most textbooks restrict to 
finite-dimensional Lie algebras and give a version of the latter one, except 
\cite{bourbaki:1989a}, which does it in full generality. The idea of most
of the combinatoric proof works with minor changes also for any Lie algebra,
since it  relies on ordered index sets which can be defined in any dimension.
The PBW theorem allows us to set up an isomorphism between $\Sym^{\bullet}
(\lie{g})$ and $\algebra{U}(\lie{g})$ immediately, since a basis of the 
former can be given by almost the same expression
\begin{equation*}
	\mathcal{B}_{\Sym^{\bullet}(\lie{g})}
	=
	\big\{
		\beta_{i_1}^{\mu_{i_1}}
		\cdots
		\beta_{i_n}^{\mu_{i_n}}
	\ \big| \
		n \in \mathbb{N}, \ 
		i_k \in I, 1 \leq k \leq n,
		i_1 \earlier \cdots \earlier i_n 
		\text{ and } \beta_{i_k} \in \mathcal{B}_{\lie{g}}, \
		\mu_{i_1}, \ldots, \mu_{i_n} \in \mathbb{N}
	\big\}
\end{equation*}
where we just have replaced the noncommutative product in $\algebra{U}
(\lie{g})$ by the symmetric tensor product in $\Sym^{\bullet}(\lie{g})$. This 
allows us to write down an isomorphism between the symmetric tensor algebra 
and the universal enveloping algebra, just by mapping the basis vectors to 
each other in a naive way. Of course, this can never be an isomorphism in 
the sense of algebras, but only of (filtered) vector spaces, because one of 
the algebras is commutative and the other isn't. Moreover, the symmetric 
algebra has a graded structure, which comes from the one on the tensor algebra, 
that the universal enveloping algebra does not have:
\begin{equation*}
	\Sym^{\bullet}(\lie{g})
	=
	\bigoplus\limits_{n = 0}^{\infty}
	\Sym^n(\lie{g})
	, \quad
	\Sym^n(\lie{g})
	=
	\underbrace{
		\lie{g} \vee \ldots \vee \lie{g}
	}_{
		n \text{ times}
	}.
\end{equation*}
We will denote by $\pi_n \colon \Sym^{\bullet}(\lie{g}) \longrightarrow
\Sym^n(\lie{g})$ the canonical projections of this grading.
This induces a filtration by $\Sym^{(k)}(\lie{g}) = \sum_{j=0}^k 
\Sym^j(\lie{g})$. Our simple isomorphism will respect the filtration, but 
not the grading. However, it is not the only isomorphism which one can write 
down. In \cite{berezin:1967a}, Berezin proposed another isomorphism which is 
more helpful to use:
\begin{equation}
	\label{Alg:BerezinQuantization}
	\mathfrak{q}_n
	\colon
	\Sym^n(\lie{g})
	\longrightarrow
	\algebra{U}^n(\lie{g})
	, \quad
	\beta_{i_1} \cdots \beta_{i_n}
	\longmapsto
	\frac{1}{n!}
	\sum\limits_{\sigma \in S_n}
	\beta_{i_{\sigma(1)}} 
	\odot \cdots \odot
	\beta_{i_{\sigma(n)}}
	, \quad
	\mathfrak{q}
	=
	\sum\limits_{n = 0}^{\infty}
	\mathfrak{q}_n.
\end{equation}
We will refer to it as the quantization map, for reasons that will soon 
become clear. It also respects the filtration and transfers the symmetric 
product to another symmetric expression. In this sense, we can now switch 
between both algebras and use the setting, which is more convenient in the 
current situation: the graded structure of $\Sym^{\bullet}(\lie{g})$, or the 
Hopf algebra structure of $\algebra{U}(\lie{g})$.

\subsection{The Gutt Star Product}
Since we know, that the universal enveloping and the symmetric tensor 
algebra are isomorphic as vector spaces, we have a good tool at hand to 
endow the symmetric tensor algebra, and hence the polynomials, with a 
noncommutative product by pulling back the product from $\algebra{U}
(\lie{g})$ to $\Sym^{\bullet}(\lie{g})$ via $\mathfrak{q}$. This is exactly 
what Gutt did in \cite{gutt:1983a}. She constructed a star product on 
$\Pol^{\bullet}(\lie{g}^*)$ from $\algebra{U}(\lie{g})$ while encoding the 
noncommutativity in a formal parameter $z \in \mathbb{C}$ in a convenient 
way. 
\begin{definition}[Gutt star product]
	\label{Def:GuttStar}
	Let $\lie{g}$ be a Lie algebra, $z \in \mathbb{C}$, and 
	$f, g \in \Sym^{\bullet}(\lie{g})$ of degree $k$ and $\ell$ 
	respectively. Then we define the Gutt star product by:
	\begin{equation}
		\label{Alg:GuttStar}
		\star_z
		\colon
		\Sym^{\bullet}(\lie{g})
		\times
		\Sym^{\bullet}(\lie{g})
		\longrightarrow
		\Sym^{\bullet}(\lie{g})
		, \quad
		(f,g)
		\longmapsto
		\sum\limits_{n = 0}^{k + \ell - 1}
		z^n
		\pi_{k + \ell - n}
		\left(
			\mathfrak{q}^{-1} \left(
				\mathfrak{q}(f) \odot \mathfrak{q}(g)
			\right)
		\right).
	\end{equation}
\end{definition}
This is the original way in which Gutt defined the star product in 
\cite{gutt:1983a}, but there are two more ways to do it. Define
\begin{equation*}
	\mathfrak{I}_z
	=
	\langle \xi \tensor \eta - \eta \tensor \xi - z [\xi, \eta] \rangle
\end{equation*}
for $z \in \mathbb{C}$. Then we set 
\begin{equation}
	\label{Alg:DeformedUnivEnvAlg}
	\algebra{U}(\lie{g}_z)
	=
	\frac{\Tensor^{\bullet}(\lie{g})}{\mathfrak{I}_z},
\end{equation}
and get the map
\begin{equation}
	\label{Alg:DeformedQuantization}
	\mathfrak{q}_{z,n}
	\colon
	\Sym^n(\lie{g})
	\longrightarrow
	\algebra{U}(\lie{g}_z)
	, \quad
	\beta_{i_1} \cdots \beta_{i_n}
	\longmapsto
	\frac{1}{n!}
	\sum\limits_{\sigma \in S_n}
	\beta_{i_{\sigma(1)}} 
	\odot \cdots \odot
	\beta_{i_{\sigma(n)}}
	, \quad
	\mathfrak{q}_z
	=
	\sum\limits_{n = 0}^{\infty}
	\mathfrak{q}_{z, n}.
\end{equation}
This way, we also get a star product:
\begin{equation}
	\label{Alg:DeformedStar}
	\widehat{\star}_z
	\colon
	\Sym^{\bullet}(\lie{g})
	\times
	\Sym^{\bullet}(\lie{g})
	\longrightarrow
	\Sym^{\bullet}(\lie{g})
	, \quad
	(f,g)
	\longmapsto
	\mathfrak{q}_z^{-1} 
	\left(
		\mathfrak{q}_z(f) \odot \mathfrak{q}_z(g)
	\right).
\end{equation}
In \cite{drinfeld:1983a}, Drinfel'd also constructed a star product using the 
Baker-Campbell-Hausdorff series: take $\xi, \eta \in \lie{g}$ and set
\begin{equation}
	\label{Alg:DrinfeldStar}
	\exp(\xi) \ast_z \exp(\eta)
	=
	\exp \left(
		\frac{1}{z}
		\bch{z \xi}{z \eta}
	\right),
\end{equation}
where the exponential series is understood a formal power series in $\xi$ 
and $\eta$. By formally differentiating, one gets the star product on all 
polynomials.

Of course, our aim is to show that these three maps are in fact identical 
and that they define a star product. Since this is a long way to go, we 
postpone the proof to the end of this chapter. It will be useful to learn 
something about the Baker-Campbell-Hausdorff series and the Bernoulli number 
first.

\section{The Baker-Campbell-Hausdorff Series}
\label{sec:chap3_BCH}

Since we have a formula for $\star_z$ which involves the
Baker-Campbell-Hausdorff series, we want to give a short overview
about it and introduce some results that will be helpful later on.
Note however, that there is not \emph{the BCH formula}, since one 
can always rearrange terms using antisymmetry or Jacobi identity, but for 
$\xi, \eta \in \lie{g}$, we can always write it as
\begin{equation}
	\label{Alg:NamesOfBCH}
	\bch{\xi}{\eta}
	=
	\sum\limits_{n = 1}^{\infty}
	\bchpart{n}{\xi}{\eta}
	=
	\sum\limits_{a,b = 0}^{\infty}
	\bchparts{a}{b}{\xi}{\eta},
\end{equation}
where $\bchpart{n}{\xi}{\eta}$ denotes all expressions having $n$ letters 
and $\bchparts{a}{b}{\xi}{\eta}$ denotes all expressions with $a$ times the letter 
$\xi$ and $b$ times the letter $\eta$. We have $\bchparts{0}{0}{\xi}{\eta} = 0$, 
$\bchparts{1}{0}{\xi}{\eta} = \xi$ and $\bchparts{0}{1}{\xi}{\eta} = \eta$. 
Clearly this gives
\begin{equation*}
	\bchpart{n}{\xi}{\eta}
	=
	\sum\limits_{a + b = n}
	\bchparts{a}{b}{\xi}{\eta}.
\end{equation*}
Of course, this only moves the problem of non-uniqueness to a later point 
when we will have to discuss the partial expressions. Yet, in the beginning, 
this will be helpful.

\subsection{Some General and Historical Remarks}
Assume $\lie{g}$ to be the Lie algebra of a finite-dimensional Lie group 
$G$. From the geometric point of view, the BCH formula is the infinitesimal 
counterpart of the multiplication law in $G$. Since the multiplication is 
smooth and the exponential function locally diffeomorphic around the 
unit element $e$, we would expect that there is a Lie algebraic analogon to 
the group multiplication, at least near the origin, which depends somehow 
``smoothly on the arguments''. Finding this expression is, however, a 
different task.

One approach to this would be the following: consider an algebra 
$\algebra{A}$ with the noncommuting elements $\xi, \eta$. We want to study 
the identity of formal power series 
\begin{equation*}
	\exp(\chi)
	=
	\exp(\xi)
	\exp(\eta).
\end{equation*}
There should be a $\chi \in \algebra{A}$, which fulfils this relation. 
We can rewrite the right hand side as
\begin{equation*}
	\exp(\xi)\exp(\eta) 
	= 
	\sum\limits_{n,m=0}^{\infty} 
	\frac{\xi^n \eta^m}{n! m!}
\end{equation*}
and use the formal power series for the logarithm
\begin{equation*}
	\log(\chi) 
	= 
	\sum\limits_{k = 1}^{\infty} 
	\frac{ (-1)^{k-1} }{ k } (\chi - 1)^k
\end{equation*}
in order to get an expression for $\chi$. This yields
\begin{equation}
	\label{Alg:BCHinRaw}
	\chi
	= 
	\log \left( \exp(\xi) \exp(\eta) \right) 
	= 
	\sum\limits_{k = 1}^{\infty} 
	\frac{ (-1)^{k - 1} }{k} 
	\sum\limits_{ \substack{ 
		i \in \{1, \ldots, k\} \\ 
		n_i, m_i \geq 0 \\ n_i + m_i \geq 1
	}}
	\frac{
		\xi^{n_1} \eta^{m_1} 
		\cdots 
		\xi^{n_k} \eta^{m_k}
	}
	{n_1! m_1! \cdots n_k! m_k!}.
\end{equation}
It is far from trivial, if and how this can be expressed using Lie brackets.
The first one who found a general way for this was Dynkin in the 1950's 
\cite{dynkin:1947a, dynkin:1950a}. Of course, the question of convergence 
still remains, although we would expect the expression to converge at least in 
a neighbourhood of $0$.

A different approach works via differential equations. We can consider flows 
on the Lie group. This gives also an expression of the group multiplication 
in logarithmic coordinates just using Lie brackets. One gets recursive 
relations for the $\bchpart{n}{\xi}{\eta}$ and the first formulas due to 
Baker \cite{baker:1905a}, Campbell \cite{campbell:1896a, campbell:1897a} and 
Hausdorff \cite{hausdorff:1906a} were of this kind. For the first terms, one 
finds
\begin{align}
	\nonumber
	\bch{\xi}{\eta}
	& =
	\log \left( \exp(\xi) \exp(\eta) \right)
	\\
	\nonumber
	& =
	\xi + \eta
	+ \frac{1}{2} [\xi, \eta]
	+ \frac{1}{12} ( [[\eta, \xi], \xi]  + [[\xi, \eta], \eta] )
	+ \frac{1}{24} [[[\eta, \xi], \xi], \eta]
	\\
	\nonumber
	& \quad
	+ \frac{1}{120}
	(
		[[[[\eta, \xi], \eta], \xi], \eta] + 
		[[[[\xi, \eta], \xi], \eta], \xi]
	)
	+ \frac{1}{360}
	(
		[[[[\eta, \xi], \xi], \xi], \eta] + 
		[[[[\xi, \eta], \eta], \eta], \xi]
	)
	\\
	\label{Alg:BCHSeriesLong}
	& \quad
	- \frac{1}{720}
	(
		[[[[\eta, \xi], \xi], \xi], \xi] + 
		[[[[\xi, \eta], \eta], \eta], \eta]
	)
	+ \cdots
\end{align}
which coincides of course with the result from Dynkin.

\subsection{Forms of the BCH Series}
As already mentioned, there are different forms of stating the BCH formula 
and depending on the problem one wants to solve, not every one is equally 
well suited. One can classify them roughly into four groups.
\begin{enumerate}
	\item
	There are recursive formulas, which calculate each term from the 
	previous one. The first expressions due to Baker, Campbell and Hausdorff 
	were of this kind. Though the idea is old, this approach is still much in
	use and allows powerful applications: Casas and Murua found an 
	efficient algorithm \cite{casas.murua:2009a} for calculating 	a form of 
	BCH series without redundancies based on a recursive formula, which was 
	given by Varadarajan in his textbook \cite{varadarajan:1974a}. 
	For such a non-redundant formula one needs a notion of basis of the free 
	Lie algebra. There are approaches to such (Hall or Hall-Viennot) basis, 
	which can e.g. be found in \cite{serre:2006a}.
	
	\item
	Most textbooks prove an integral form of the series, like 
	\cite{hall:2003a} and \cite{hilgert.neeb:2012a}. Since we will use it, 
	too, we want to introduce it briefly. Take the function
	\begin{equation}
		\label{Alg:DefinitionLogBernoullis}
		g \colon 
		\mathbb{C}
		\longrightarrow
		\mathbb{C},
		\quad
		z 
		\longmapsto
		\frac{ z \log(z) }{z - 1}
	\end{equation}
	and denote for $\xi \in \lie{g}$ by
	\begin{equation*}
		\ad_{\xi}
		\colon
		\lie{g}
		\longrightarrow
		\lie{g},
		\quad
		\eta
		\longmapsto
		[\xi, \eta]
	\end{equation*}
	the usual adjoint operator. Then one has for $\xi, \eta \in \lie{g}$
	\begin{equation}
		\label{Alg:BCHinIntegral}
		\bch{\xi}{\eta}
		=
		\xi + 
		\int\limits_0^1
		g \left( 
			\exp \left(
				\ad_{\xi}
			\right)
			\exp \left(
				t \ad_{\eta}
			\right)
		\right)
		(\eta) 
		dt.
	\end{equation}
	
	\item
	As already mentioned, Dynkin found a closed form for \eqref{Alg:BCHinRaw},
	which is the only one of this kind known so far. A proof can be found in 
	\cite{jacobson:1979a}, for example. It reads
	\begin{equation}
		\label{Alg:BCHinDynkin}
		\bch{\xi}{\eta}
		=
		\sum\limits_{k = 1}^{\infty} 
		\frac{ (-1)^{k - 1} }{k}
		\sum\limits_{\substack{
			i \in \{1, \ldots, k\} \\
			n_i, m_i \geq 0 \\ 
			n_i + m_i \geq 1
		}}
		\frac{1}{ \sum_{i = 1}^k(n_i + m_i) }
		\frac{ \left[
			\xi^{n_1} \eta^{m_1} \cdots \xi^{n_k} \eta^{m_k}
		\right]}
		{ n_1! m_1! \cdots n_k! m_k! },
	\end{equation}
	where the expression $[ \ldots ]$ denotes Lie brackets nested to the left:
	for instance, we have
	\begin{equation*}
		[\xi \eta \eta \xi]
		=
		[[[\xi, \eta], \eta], \xi].
	\end{equation*}
	Unfortunately, the combinatorics get extremely complicated for higher 
	degrees and increasingly many terms belong to the same Lie bracket 
	expression.
	
	\item
	Goldberg gave a form of the series which is based on words in two letters:
	\begin{equation}
		\bch{\xi}{\eta}
		=
		\sum\limits_{n = 1}^{\infty}
		\sum\limits_{|w| = n}
		g_w w.
	\end{equation}
	The $g_w$ are coefficients, which can be calculated using the recursively
	defined Goldberg polynomials (see \cite{goldberg:1956a}). It was put into 
	commutator form by Thompson in \cite{thompson:1982a}:
	\begin{equation}
		\label{Alg:BCHinGoldbergThompson}
		\bch{\xi}{\eta}
		=
		\sum\limits_{n = 1}^{\infty}
		\sum\limits_{|w| = n}
		\frac{g_w}{n} [w].
	\end{equation}
	Again, the $[w]$ are Lie brackets nested to the left. Of course, this 
	formula will also have redundancies, but its combinatorial aspect is much 
	easier than the one of Equation~\eqref{Alg:BCHinDynkin}. Since there are 
	estimates for the coefficients $g_w$, we will use this form for our Main 
	Theorem.
\end{enumerate}

\subsection{The Goldberg-Thompson Formula and some Results}

\subsubsection{Goldberg's Theorems}
We now introduce the results of Goldberg: he denoted a word in the letters $\xi$ 
and $\eta$ as
\begin{equation*}
	w_{\xi}(s_1, s_2, \ldots, s_m)
	=
	\xi^{s_1} \eta^{s_2} \cdots (\xi \vee \eta)^{s_m},
\end{equation*}
with $m \in \mathbb{N}$ and the last letter will be $\xi$ if $m$ is odd and 
$\eta$ if $m$ is even. The index $\xi$ of $w_{\xi}$ means that the word 
starts with a $\xi$. Now we can assign to each word $w_{\xi \vee \eta}(s_1, 
\ldots, s_m)$ a coefficient $c_{\xi \vee \eta}(s_1, \ldots, s_m)$. This is 
done by the following formula:
\begin{equation}
	\label{Alg:GoldbergCoeff1}
	c_{\xi}(s_1, \ldots, s_m) 
	= 
	\int\limits_0^1 
	t^{m'} (t - 1)^{m''} 
	G_{s_1}(t) \cdots G_{s_m}(t) dt,
\end{equation}
where we have $m' = \lfloor \frac{m}{2} \rfloor$, $m'' = \lfloor \frac{m-1}{2} 
\rfloor$ with $\lfloor \cdot \rfloor$ denoting the entire part of a real 
number and we have $n = \sum_{i = 1}^m s_i$. The $G_s$ are the recursively 
defined Goldberg polynomials
\begin{equation}
	\label{Alg:GoldbergPolynomials}
	G_s(t) 
	= 
	\frac{1}{s} 
	\frac{d}{dt} 
	t(t-1) G_{s-1}(t),
\end{equation}
for $s > 1$ and $G_1(t) = 1$. For $c_{\eta}$ we have
\begin{equation}
	\label{Alg:GoldbergCoeff2}
	c_{\eta}\left( s_1, \ldots, s_m \right)
	= 
	(-1)^{n-1} c_{\xi} \left( s_1, \ldots, s_m \right)
\end{equation}
and furthermore
\begin{equation*}
	c_{\eta}\left( s_1, \ldots, s_m \right)
	= 
	c_{\xi} \left( s_1, \ldots, s_m \right)
\end{equation*}
if $m$ is odd. This yields immediately
\begin{equation*}
	c_{\xi}\left( s_1, \ldots, s_m \right)
	= 
	c_{\eta}\left( s_1, \ldots, s_m \right) 
	= 
	0 
\end{equation*}
if $m$ is odd and $n$ is even. Of course, Goldberg found interesting identities 
which are fulfilled by the coefficients. A very remarkable one is that for all 
permutations $\sigma \in S_m$ one has
\begin{equation*}
	c_{\xi}\left( s_1, \ldots, s_m \right)
	=
	c_{\xi}\left( s_{\sigma(1)}, \ldots, s_{\sigma(m)} \right),
\end{equation*}
since \eqref{Alg:GoldbergCoeff1} obviously does not see the ordering of the 
$s_i$ and $m'$, $m''$ and $n$ are not affected by reordering. For words with 
$m = 2$, an easier formula can be found:
\begin{equation*}
	c_{\xi}(s_1, s_2) 
	= 
	\frac{ (-1)^{s_1} }{s_1! s_2!} 
	\sum\limits_{n = 1}^{s_2} 
	\binom{s_2}{n} B_{s_1 + s_2 - n},
\end{equation*}
where the $B_s$ denote the Bernoulli numbers, which will be explained more 
precisely in the next paragraph. First, we note that the only case which 
matters to us is of course $s_1 = 1$, since for $s_1, s_2 > 1$ we will find 
something like $[[\xi, \xi], \ldots ] = 0$. For simplicity, let's set 
$s_2 = 1$ and to permute $s_1 \leftrightarrow s_2$:
\begin{equation}
	\label{Alg:GoldbergCoeffBernoulli}
	c_{\xi}(1, s)
	=
	\frac{ (-1)^s }{s!}
	B_s.
\end{equation}

\subsubsection{Bernoulli numbers}
We have seen the Bernoulli numbers $B_n$ showing up and we will encounter them 
very often in the following. Hence it is useful to learn a few important 
things about them. They are defined by the series expansion of
\begin{equation}
	\label{Alg:DefBernoulli}
	g \colon
	\mathbb{C}
	\longrightarrow
	\mathbb{C}
	, \quad
	z \longmapsto
	\frac{z}{\E^z - 1}
	=
	\sum\limits_{n=0}^{\infty}
	\frac{B_n}{n!}z^n.
\end{equation}
Clearly, $g$ has poles at $z = 2 k \pi \I$, $k \in \mathbb{Z}\backslash \{0\}$.
Moreover, one can easily show that all odd Bernoulli numbers are zero, except 
for $B_1 = - \frac{1}{2}$ and since in some applications one wants $B_1$ to be 
positive, there is a different convention for naming them: one often 
encounters $B_n^* = (-1)^n B_n$ (which only differs for $n = 1$). The nonzero 
Bernoulli numbers alternate in sign. For their absolute value, one can show 
the asymptotic behaviour (see \cite{oeis:A027641, oeis:A027642})
\begin{equation*}
	|B_{2n}|
	\sim 
	(-1)^{n-1}
	\frac{2 (2n)!}{(2 \pi)^{2n}}.
\end{equation*}
This is not surprising, since we know that the generating function $g$ had 
poles at $\pm 2\pi \I$. The Bernoulli numbers can also be calculated by the 
recursion formula
\begin{equation}
	\label{Alg:BernoulliRecursive}
	B_n
	=
	- \frac{1}{n + 1}
	\sum\limits_{k = 0}^{n- 1}
	\binom{n + 1}{k}
	B_k,
\end{equation}
which is well-known in the literature (e.g. 
\cite{arakawa.ibukiyama.kaneko:2014a}). Since we will deal with them, 
we want to give the first numbers of this series here.
\begin{center}
	\begin{tabular}
	{c||c|c|c|c|c|c|c|c|c|c|c|c|c|c|c|c|c}
		$n$ & $0$ & $1$ & 
		$2$ & $3$ & $4$ & 
		$5$ & $6$ & $7$ & 
		$8$ & $9$ & $10$ & 
		$11$ & $12$ & $13$ & 
		$14$ & $15$ & $16$
		\\
		\hline 
		$B_n$ & $1$ & $-\frac{1}{2}$ & 
		$\frac{1}{6}$ & $0$ & $-\frac{1}{30}$ & 
		$0$ & $\frac{1}{42}$ & $0$ & 
		$-\frac{1}{30}$ & $0$ & $\frac{5}{66}$ & 
		$0$ & $\frac{691}{2730}$ & $0$ & 
		$\frac{7}{6}$ & $0$ & $\frac{3617}{510}$
	\end{tabular} 
\end{center}

\subsubsection{BCH up to first order}

\begin{proposition}
	\label{Alg:Prop:BCHFristOrder}
	Let $\lie{g}$ be a Lie algebra and the Bernoulli numbers as defined 
	before. Then we have for $\xi, \eta \in \lie{g}$
	\begin{align}
		\label{Alg:BCHFirstOrderXi}
		\bch{\xi}{\eta}
		& =
		\sum\limits_{n = 0}^{\infty}
		\frac{B_n^*}{n!}
		\left( \ad_{\xi} \right)^n (\eta)
		+ \mathcal{O} \left( \eta^2 \right)
		\\
		\label{Alg:BCHFirstOrderEta}
		& =
		\sum\limits_{n = 0}^{\infty}
		\frac{B_n}{n!}
		\left( \ad_{\eta} \right)^n (\xi)
		+ \mathcal{O} \left( \xi^2 \right)
	\end{align}
\end{proposition}
\begin{proof}
	We want to can calculate this using the Goldberg coefficients. Remind that 
	we put words to Lie brackets, and for computing the coefficients we 
	need the words $\eta \xi^n$ \emph{and} $\xi \eta \xi^{n-1}$ because 
	of antisymmetry and words of the form $\xi^k \eta \xi^{n-k}$ with $k > 1$ 
	give vanishing expressions. Now let $n \in \mathbb{N}$. We have
	\begin{equation*}
		c_{\eta}(1,n) 
		= 
		(-1)^n c_{\xi}(1,n) 
		= 
		(-1)^n 
		\frac{(-1)^n}{n!} 
		B_n 
		= 
		\frac{B_n}{n!}.
	\end{equation*}
	By $n$-fold skew-symmetry and \eqref{Alg:BCHinGoldbergThompson}, we get 
	the contribution
	\begin{equation*}
		\frac{(-1)^n}{(n+1)!} 
		B_n \left( \ad_{\xi} \right)^n (\eta)
		=
		\frac{1}{(n+1)!} 
		B_n^* \left( \ad_{\xi} \right)^n (\eta)		
	\end{equation*}
	Now we need $c_{\xi}(1, 1, n-1)$: let $n > 1$, then
	\begin{align*}
		c_{\xi}(1, 1, n - 1) 
		& = 
		\int\limits_0^1 
		t(t-1) G_{n-1}(t) 
		dt 
		\\
		& = 
		- \int\limits_0^1 
		t \frac{d}{dt} 
		\left( t(t-1) G_{n-1}(t) \right) 
		dt 
		\\
		& =
		- \int\limits_0^1 
		n t G_n(t) 
		dt 
		\\
		& = 
		-n 
		c_{\xi}(1,n) 
		\\
		& = 
		-n 
		\frac{(-1)^n}{n!} 
		B_n 
		\\
		& =  
		(-1)^{n + 1} 
		\frac{1}{(n-1)!} B_n,
	\end{align*}
	where we have done an integration by parts in the third step.
	So by using $n-1$ times the skew-symmetry of the Lie bracket, we get
	\begin{equation*}
		\frac{1}{n + 1} 
		\cdot (-1)^{n + 1} 
		\frac{1}{ (n - 1)! }
		B_n 
		[\ldots [\xi, \eta], \xi] \ldots ], \xi] 
		= 
		\frac{n}{(n+1)!} 
		B_n 
		\left( \ad_{\xi} \right)^n (\eta).
	\end{equation*}
	For $n > 1$ we add up those two and use the fact that $B_n = B_n^*$ and 
	find the result we want. For $n = 1$, there is just the first contribution 
	and $c_{\xi}(1,1) = - B_1$, which gives
	\begin{equation*}
		B_1^* \ad_{\xi}(\eta)
	\end{equation*}
	in total. For $n = 0$ , we get $c_{\xi}(1) = c_{\eta}(1) = 1$ and finally 
	get \eqref{Alg:BCHFirstOrderXi}. For \eqref{Alg:BCHFirstOrderEta}, note 
	that we need $c_{\xi}(1, n)$ and $c_{\eta}(1, 1, n-1)$. We have 
	$c_{\xi}(1,n) = (-1)^n c_{\eta}(1,n)$ and 
	$c_{\eta}(1, 1, n-1) = (-1)^n c_{\xi}(1,1,n-1$.
	This gives $(-1)^n$ and switches $B_n$ to $B_n^*$.
\end{proof}
\begin{remark}[Alternative Proof]
	Note that we could also have used the integral formula 
	\eqref{Alg:BCHinIntegral} to prove this. We want to sketch an alternative 
	proof here: if we write the second of the two exponential functions as a 
	series, we see that it can be cut after the constant term, since we are 
	looking for contributions which are linear in $\eta$. The function left to 
	integrate is then just $(g \circ \log)(z)$. Since we insert 
	$\exp(\ad_{\xi})$, we get
	\begin{equation*}
		\bch{\xi}{\eta}
		=
		\xi + 
		\int_0^1
		g \left( \ad_{\xi} \right)
		(\eta)
		dt
		+ \mathcal{O} \left( \eta^2 \right)
		=
		\xi +
		\sum\limits_{n = 1}^{\infty}
		\frac{B_n^*}{n!}
		\left( \ad_{\xi} \right)^n (\eta)
		+ \mathcal{O} \left( \eta^2 \right),
	\end{equation*}
	since there is no dependence on $t$ left and we get the same result.
\end{remark}

\subsubsection{Thompson's estimates on the coefficients}
We know that we can put the BCH series into the form of Equation 
\eqref{Alg:BCHinGoldbergThompson}:
\begin{equation}
	\bch{\xi}{\eta}
	=
	\sum\limits_{n = 1}^{\infty}
	\sum\limits_{|w| = n}
	\frac{g_w}{n} [w].
\end{equation}
Later, we will need estimates on the coefficients $g_w$ in order to show the 
continuity of the Gutt star product. The first simple estimate (together with a 
table of the first Goldberg coefficients and good explanation of Goldberg 
polynomials) was given by Newman and Thompson in \cite{newman.thompson:1987a}. 
The idea behind it was an analysis of the structure of the polynomials an its 
roots. This allows to put tight bounds on the values of $|g_w|$. We will need a 
bit more than that, but luckily, Thompson \cite{thompson:1989a} gave an estimate 
in exactly the form we will need.
\begin{proposition}
	\label{Alg:Prop:ThompsonsEstimate}
	Let $n \in \mathbb{N}$ and $g_w$ denotes the Goldberg coefficient of a 
	word in two letters. Then we have the estimate
	\begin{equation}
		\label{Alg:ThompsonsEstimate}
		\sum\limits_{|w| = n}
		\left| g_w \right|
		\leq
		2
	\end{equation}
\end{proposition}
\begin{proof}
	We want to sketch the proof here for convenience. One can see from the 
	recursion formula \eqref{Alg:GoldbergPolynomials}, that the $G_s(t)$ are 
	symmetric around $z = \frac{1}{2}$ (maybe up to a factor $t - \frac{1}
	{2}$): doing a shift $t \mapsto t \frac{1}{2}$, they are of the form
	\begin{equation*}
		G_s(t)
		=
		\frac{1}{s!}
		\frac{d}{dt}(t^2 - \frac{1}{4})
		\cdots
		\frac{d}{dt}(t^2 - \frac{1}{4}).
	\end{equation*}
	Their roots lie all in the interval $(0,1)$ and the polynomials are normed 
	(which can also be seen from the recursion formula). This means that they 
	can be written (in the shifted form) as
	\begin{equation*}
		G_s(t)
		=
		t^{s_0}
		(t + t_{s_1})(t - t_{s_1})
		\cdots
		(t + t_{s_r})(t - t_{s_r})
	\end{equation*}
	with $r = \lfloor \frac{s-1}{2} \rfloor$, $t_{s_i} \in (0, \frac{1}{2})$ 
	and $s_0 = 0$ if $s$ is odd and $s_0 = 1$ if $s$ is even. The symmetric, 
	quadratic terms are bounded in their absolute value by $\frac{1}{4}$ and 
	the linear term by $\frac{1}{2}$, since $t \in (-\frac{1}{2}, \frac{1}
	{2})$. Hence we have the estimate $|G_s(t)| \leq 2^{-s + 1}$ for the 
	integral domain. Now note that integration by parts yields
	\begin{equation*}
		\int\limits_{0}^1
		t^a (t-1)^b
		dt
		=
		\frac{a! b!}{(a + b + 1)!}.
	\end{equation*}
	We can put this together to get an estimate for $g_w = c_{\xi, \eta}(s_1, 
	\ldots, s_m)$. Note again by $n = s_1 + \cdots + s_m$. A slight 
	rearranging of the factor from integration by 
	parts using the fact that one of the numbers $m' = \frac{m}{2}$ and $m'' = 
	\frac{m - 1}{2}$ is not an integer and will therefore be rounded downwards 
	gives
	\begin{align*}
		|g_w|
		& \leq
		2^{-n + m}
		\frac{m' ! m'' !}{(m' + m'' + 1)!}
		\\
		& =
		2^{-n + m}
		\frac{1}{m}
		\binom{m-1}{m'}^{-1}.
	\end{align*}
	Now we just need to sum up all the expressions corresponding to 
	words of the length $n$. Note, that the words can start with $\xi$ or 
	$\eta$ and we therefore get a factor $2$ in front. The number of possible 
	arrangements $(s_1, \ldots, s_m)$ is due to a combinatorial argument 
	($n - m$ balls into $m$ buckets, since every bucket must contains at least 
	one ball) given by $\binom{n-1}{m-1}$ and we have to sum over all possible	
	$m$. We get
	\begin{align*}
		\sum\limits_{|w| = n}
		|g_w|
		& \leq
		\sum\limits_{m = 1}^n
		2 \binom{n-1}{m-1}
		2^{-n + m} \frac{1}{m}
		\binom{m-1}{m'}^{-1}
		\\
		& =
		2^{-n + 1}
		\sum\limits_{m = 1}^n
		\binom{n-1}{m-1}
		\underbrace{
			2^m
			\frac{1}{m}
			\binom{m-1}{m'}^{-1}
		}_{ \leq 2 \ (*)}
		\\
		& \leq
		2^{-n + 2}
		\sum\limits_{m = 1}^n
		\binom{n-1}{m-1}
		\\
		& =
		2.
	\end{align*}
	In $(*)$ we used the fact that $\binom{m-1}{m'}$ is the biggest term (or 
	one of the two) biggest terms in the binomial expansion of $(1 + 1)^{m-1}$ 
	and hence we have $m \binom{m-1}{m'} \geq 2^{m-1}$.
\end{proof}

\section{The Equality of the Star Products}
\label{sec:chap3_StarProductProof}

We want to prove the equality of the three star products. For a general and 
possibly infinite-dimensional Lie algebra, this is quite tedious. As a first 
step, it will be helpful to show their associativity.
\begin{remark}
	In the finite-dimensional case, there are different proofs for 
	Theorem~\ref{Alg:Thm:ThreeStarsAreOne}, the main theorem of this section, 
	which mostly rely on geometric arguments, like the one in 
	\cite{bordemann.neumaier.waldmann:1999a}. Unluckily, these techniques are 
	not at hand in infinite dimensions and one has to find an algebraic proof 
	instead. Since in the community of deformation quantization, this 
	statement is somehow folklore and believed for any Lie algebra, it 
	strongly seems like such a proof already exists. However, the we were 
	not able to trace it down in literature and therefore give an own proof.
\end{remark}
\begin{proposition}
	\label{Alg:Lemma:Associtivity}
	The three maps $\star_z$, $\widehat{\star}_z$ and $\ast_z$ from 
	\eqref{Alg:GuttStar}, \eqref{Alg:DeformedStar} and \eqref{Alg:DrinfeldStar}  
	respectively define	associative multiplications.
\end{proposition}
\begin{proof}
  All maps are defined as $\Sym^{\bullet}(\lie{g}) \times \Sym^{\bullet}
  (\lie{g}) \longrightarrow \Sym^{\bullet}(\lie{g})$, so we have to show 
  bilinearity and associativity.
  \begin{enumerate}
	\item 
	For $\widehat{\star}_z$, associativity and bilinearity are clear from the 
	construction, since we just pull-back the multiplication in $\algebra{U}
	(\lie{g}_z)$.
	
	\item
	For $\star_z$, bilinearity follows from the fact that all maps, that are 
	used in its definition, are (bi-)linear. For associativity, we have to 
	interchange sums and shift projections.
	Recall that $\pi_n(f \star_z g) = 0$, if $n > \deg(f) + \deg(g)$.
	Take homogeneous tensors $f, g, h \in \Sym^{\bullet}(\lie{g})$
	of degree $k, \ell, m \in \mathbb{N}$ respectively. Then we have
	\begin{align*}
		\big(
			f 
		&
			\star_z g
		\big)
		\star_z
		h
		\\
		& =
		\sum\limits_{j=0}^{k + \ell -1}
		\sum\limits_{i = 0}^{k + \ell + m -j - 1}
		z^i
		\left( 
			\pi_{k + \ell + m - j - i} 
			\circ 
			\mathfrak{q}^{-1}
		\right)
		\left(
			\mathfrak{q}
			\left(
				z^j
				\left(
					\pi_{k + \ell - j} 
					\circ 
					\mathfrak{q}^{-1}
				\right)
				\left(
					\mathfrak{q} (f) 
					\odot
					\mathfrak{q} (g)
				\right)
			\right)
			\odot
			\mathfrak{q}(h)
		\right)
		\\
		& =
		\sum\limits_{j=0}^{k + \ell -1}
		\sum\limits_{i = 0}^{k + \ell + m - 1}
		z^{i-j}
		\left( 
			\pi_{k + \ell + m - i} 
			\circ 
			\mathfrak{q}^{-1}
		\right)
		\left(
			\mathfrak{q}
			\left(
				z^j
				\left(
					\pi_{k + \ell - j} 
					\circ 
					\mathfrak{q}^{-1}
				\right)
				\left(
					\mathfrak{q} (f)
					\odot
					\mathfrak{q} (g)
				\right)
			\right)
			\odot
			\mathfrak{q}(h)
		\right)
		\\
		& =
		\sum\limits_{i = 0}^{k + \ell + m - 1}
		z^i
		\left( 
			\pi_{k + \ell + m - i} 
			\circ 
			\mathfrak{q}^{-1}
		\right)
		\left(
			\mathfrak{q}
			\left(
				\sum\limits_{j=0}^{k + \ell -1}
				z^{-j}
				z^j
				\left(
					\pi_{k + \ell - j} 
					\circ 
					\mathfrak{q}^{-1}
				\right)
				\left(
					\mathfrak{q} (f)
					\odot
					\mathfrak{q} (g)
				\right)
			\right)
			\odot
			\mathfrak{q}(h)
		\right)
		\\
		& =
		\sum\limits_{i = 0}^{k + \ell + m - 1}
		z^i
		\left( 
			\pi_{k + \ell + m - i} 
			\circ 
			\mathfrak{q}^{-1}
		\right)
		\left(
			\mathfrak{q}(f)
			\odot
			\mathfrak{q}(g)
			\odot
			\mathfrak{q}(h)
		\right),
	\end{align*}
	and we just need to do the reversed process on the right hand side to get 
	the wanted result.
	
	\item
	For $\ast_z$, we get associativity using the exponential function and the 
	logarithm. We have
	\begin{align*}
		\left( 
			\exp(\xi) \ast_z \exp(\eta) 
		\right) 
		\ast_z \exp(\chi)
		&=
		\exp
		\left(
			\frac{1}{z}
			\bch{
				\left(
					\frac{1}{z}
					\bch{z \xi}{z \eta}
				\right)
			}
			{z \chi}
		\right)
		\\
		&=
		\exp
		\left(
			\frac{1}{z}
			\bch{z \xi}
			{
				\left(
					\frac{1}{z}
					\bch{z \eta}{z \chi}
				\right)
			}
		\right)
		\\
		&=
		\exp(\xi)
		\ast_z
		\left( 
			\exp(\eta) \ast_z \exp(\chi) 
		\right), 
	\end{align*}
	since
	\begin{align*}
	\bch{
		\left(
			\frac{1}{z}
			\bch{z \xi}{z \eta}
		\right)
	}
	{z \chi}
	& =
	\log\left(
		\exp\left(
			\log\left(
				\frac{1}{z}
				\exp\left(
					z \xi
				\right)
				\exp\left(
					z \eta
				\right)
			\right)
		\right)
		\exp\left(
			z \chi
		\right)
	\right)
	\\
	& =
	\log\left(
		\frac{1}{z}
		\exp\left(
			z \xi
		\right)
		\exp\left(
			z \eta
		\right)
		\exp\left(
			z \chi
		\right)
	\right)
	\\
	& =
	\log\left(
		\exp\left(
			z \xi
		\right)
		\log\left(
			\left(
				\frac{1}{z}
				\exp\left(
					z \eta
				\right)
				\exp\left(
					z \chi
				\right)
			\right)
		\right)
	\right)
	\\
	& =
	\bch{z \xi}
	{
		\left(
			\frac{1}{z}
			\bch{z \eta}{z \chi}
		\right)
	}.
	\end{align*}
	Bilinearity follows from differentiating.
  \end{enumerate}	 
\end{proof}
Note that star products must fulfil the classical and the semi-classical limit. 
We will do this in Corollary~\ref{Formulas:Cor:LimitCases} and so just the 
equality is left to show. It is enough to prove the coincidence for terms of 
the form $\xi^k \star \eta$ with $\xi, \eta \in \lie{g}$ and $k \in 
\mathbb{N}$, because $\Sym^{\bullet}(\lie{g})$ is a commutative algebra and 
hence we get the coincidence on arbitrary monomials by polarization. The equality 
for the product of two monomials then follows by iteration, which is possible due 
to associativity. The next lemma presents a first big step.
\begin{lemma}
	\label{Alg:FuckingDeadlyLemma}
	Let $\xi, \eta \in \lie{g}$, then we have
	\begin{equation}
		\label{Alg:GuttStarFormula}
		\xi^k \widehat{\star}_z \eta
		=
	 	\sum\limits_{n=0}^k 
	 	z^n \binom{k}{n} B_n^*
	 	\xi^{k-n} 
	 	\left( \ad_{\xi} \right)^n (\eta).
	\end{equation}
\end{lemma}
\begin{proof}
	This proof is divided into the two following lemmata:
	\begin{lemma}
		\label{Alg:Lemma:BadSublemma1}
		Let $\xi, \eta \in \mathfrak{g}$ and $k\in \mathbbm{N}$. Then we have
		\begin{equation*}
			\lie{q}_z \left(\sum\limits_{n=0}^k z^n
			\binom{k}{n} B_n^* \xi^{k-n}
			\left(\ad_{\xi}\right)^n(\eta)\right)
			=
			\sum\limits_{s=0}^k\mathcal{K}(k,s)
			\xi^{k-s} \odot \eta \odot \xi^s
		\end{equation*}
		with
		\begin{equation*}
			 \mathcal{K}(k,s)
			 =
			 \frac{1}{k + 1} \sum\limits_{n=0}^k
			 \binom{k+1}{n} B_n^*
			 \sum\limits_{j=0}^n
			 (-1)^j \binom{n}{j}
			 \sum\limits_{\ell=0}^{k-n} \delta_{s, \ell+j}.
		\end{equation*}
	\end{lemma}
	\begin{subproof}
		Since the map $\mathfrak{q}_z$ is linear, we can pull out the
		constants and get
		\begin{equation*}
			\mathfrak{q}_z \left( \sum\limits_{n=0}^k
			z^n \binom{k}{n} B_n^* \xi^{k-n} 
			\left(
				\ad_{\xi}\right)^n(\eta)
			\right)
			=
			\sum\limits_{n=0}^k
			\binom{k}{n} B_n^*
			\mathfrak{q}_z
			\left(
				z^n \xi^{k-n}
				\left( \ad_{\xi} \right)^n
				(\eta)
			\right).
		\end{equation*}
		Now we need the two equalities
		\begin{align*}
			\mathfrak{q}_z
			\left( \xi^n \eta \right)
			& = 
			\frac{1}{n + 1}
			\sum\limits_{\ell = 0}^n 
			\xi^{k - \ell} \odot \eta \odot \xi^{\ell}
		\intertext{and}
			\mathfrak{q}_z
			\left(
				\left( z^n \ad_{\xi} \right)^n 
				(\eta) 
			\right)
			& = 
			\sum\limits_{j=0}^n
			(-1)^j \binom{n}{j}
			\xi^{n-j} \odot \eta \odot \xi^j
		\end{align*}
		which can easily be shown by induction. They give
		\begin{align*}
			\sum\limits_{n=0}^k
			\binom{k}{n} B_n^*
			&
			\mathfrak{q}_z
			\left(
				z^n \xi^{k-n}
				\left( \ad_{\xi} \right)^n
				(\eta)
			\right)
			\\
			& = 
			\sum\limits_{n=0}^k
			\binom{k}{n} \frac{B_n^*}{k-n+1}
			\sum\limits_{\ell=0}^{k-n}
			\xi^{k-n-\ell} \odot
			\left(
				\sum\limits_{j=0}^n
				(-1)^j \binom{n}{j} 
				\xi^{n-j} \odot \eta \odot \xi^j
			\right) 
			\odot \xi^{\ell}
			\\
			& =
			\sum\limits_{n=0}^k
			\binom{k}{n} \frac{B_n^*}{k-n+1}
			\sum\limits_{\ell=0}^{k-n}
			\sum\limits_{j=0}^n
			(-1)^j \binom{n}{j} 
			\xi^{k-\ell-j} \odot \eta \odot \xi^{\ell+j} 
			\\	
			& = 
			\frac{1}{k + 1} \sum\limits_{n=0}^k
			\binom{k+1}{n} B_n^*
			\sum\limits_{j=0}^n (-1)^j \binom{n}{j}
			\sum\limits_{\ell=0}^{k-n}
			\xi^{k-\ell-j} \odot \eta \odot \xi^{\ell+j}
		\end{align*}
		We just need to collect those terms for which we have $\ell + j = s$
		for all	$s = 0, \ldots, k$. If we do this with a Kronecker-delta, we 
		will get exactly the $\mathcal{K}(k, s)$.
	\end{subproof}
	For the second lemma, we need some statements on Bernoulli numbers and 
	binomial coefficients. Let $k,m,n \in \mathbbm{N}$. Then we have the 
	following identities:
	\begin{align}
		\label{Alg:Lemmastuff:1_B1}
		\sum\limits_{j=0}^k
		\binom{k+1}{n} B_j^* 
		& = k+1 \\
		\label{Alg:Lemmastuff:1_Carl}
		(-1)^k \sum\limits_{j=0}^k
		\binom{k}{j} B_{m+j} 
		& = 
		(-1)^m \sum\limits_{i=0}^m 
		\binom{m}{i} B_{k+i} \\
		\label{Alg:Lemmastuff:1_Bin1} 
		\sum\limits_{j=0}^m
		(-1)^j \binom{n}{j} 
		& = 
		(-1)^m \binom{n-1}{m} \\
		\label{Alg:Lemmastuff:1_Bin2}
		\binom{n}{m}\binom{m}{k}
		& =
		\binom{n}{k}\binom{n-k}{m-k}.
	\end{align}
	The first one can easily be proven using the recursive definition of the 
	Bernoulli numbers \eqref{Alg:BernoulliRecursive}.
	Equation~\eqref{Alg:Lemmastuff:1_Bin1} and \eqref{Alg:Lemmastuff:1_Bin2} 
	are standard identities in combinatorics and can be found in the textbook 
	of Aigner \cite{aigner:2000a}. Finally, 
	Equation~\eqref{Alg:Lemmastuff:1_Carl} is a theorem due to Carlitz 
	\cite{carlitz:1968a}. With them, we can show the next lemma which will 
	finish this proof.
	\begin{lemma}
		\label{Alg:Lemma:BadSublemma2}
		Let $\mathcal{K}(k,s)$ be defined as in Lemma 
		\ref{Alg:Lemma:BadSublemma1}, then we have for all $k\in \mathbbm{N}$
		\begin{equation*}
			\mathcal{K}(k, s)
			=
			\begin{cases}
				1 & s=0 \\
				0 & \text{else}.		
			\end{cases}
		\end{equation*}
	\end{lemma}
	\begin{subproof}
		This is divided into three parts. First, we show the statement for 
		$s=0$, then we show it for $s=1$ and then proceed by induction.
		\begin{enumerate}[(i)]
		  \item $s=0$:
			The Kronecker-delta will always be zero unless $l=j=0$. So we get
			\begin{equation*}
				\mathcal{K}(k,0)
				=
				\frac{1}{k+1} \sum\limits_{n=0}^k
				\binom{k+1}{n} B_n^*
				=
				\frac{k+1}{k+1}
				=
				1,
			\end{equation*}
			where we have used \eqref{Alg:Lemmastuff:1_B1}.

			\item $s=1$:
			To get a contribution from the $\delta$, we must have 
			$(j, \ell) = (1,0)$ or $(0,1)$. Except for $n=0$ and $n=k$, 
			both cases are possible. We split them off:
			\begin{align*}
				\mathcal{K}(k,1)
				& = 
				\underbrace{\frac{1}{k+1}}_{n=0}
				- \underbrace{k B_k^*}_{n=k}
				+ \frac{1}{k+1}
				\sum\limits_{n=1}^{k-1}
				\binom{k+1}{n}B_n^*
				\left(1 + (-1) \binom{n}{1} \right)\\
				& = 
				\frac{1}{k+1}
				+ \frac{1}{k+1} \sum\limits_{n=1}^{k-1}
				\binom{k+1}{n}B_n^* 
				- \frac{1}{k+1}	\sum\limits_{n=1}^{k-1}
				\binom{k+1}{n}n B_n^* 
				- k B_k^* \\
				& = 
				\underbrace{
					\frac{1}{k+1}
					\sum\limits_{n=0}^{k-1}
					\binom{k + 1}{n}
					B_n^*
				}_{ = 1 - B_k^*} 
				- \frac{1}{k+1} \sum\limits_{n=0}^k
				\binom{k+1}{n}n B_n^*
				\\
				& =
				1 - B_k^* 
				- \underbrace{\frac{k+1}{k+1}
				\sum\limits_{n=0}^k
				\binom{k+1}{n} B_n^*}_{ = k+1}
				+ \sum\limits_{n=0}^k 
				\underbrace{\frac{k+1-n}{k+1}
				\binom{k+1}{n}}_{\binom{k}{n}} B_n^*
				\\
				& = 
				1 - B_k^* - k - 1 
				+ \sum\limits_{n=0}^{k-1}
				\binom{k}{n} B_n^* + B_k^*
				\\
				& =
				- k + 
				\sum\limits_{n}^{k-1} 
				\binom{k}{n} B_n^*
				\\
				& = 
				0.
			\end{align*}
	
		  \item $s \mapsto s+1$:
			Due to the induction, it is sufficient to prove 
			$\mathcal{K}(k,s+1) - \mathcal{K}(k,s) = 0$.
			In order to do that, we must get rid of the $\delta$'s and 
			therefore rewrite $\mathcal{K}(k,s)$:
			\begin{align*}
				\mathcal{K}(k,s)
				& = 
				\frac{1}{k + 1}	\sum\limits_{n=0}^k
				\binom{k+1}{n} B_n^* \sum\limits_{j=0}^n
				(-1)^j \binom{n}{j} \sum\limits_{\ell=0}^{k-n}
				\delta_{s, \ell+j}\\
				& = 
				\frac{1}{k + 1} \sum\limits_{n=0}^s
				\binom{k+1}{n} B_n^* \sum\limits_{j=0}^n
				(-1)^j \binom{n}{j} \sum\limits_{\ell=0}^{k-n}
				\delta_{s, \ell+j} \\
				& \quad + \frac{1}{k + 1} \sum\limits_{n=s+1}^k
				\binom{k+1}{n} B_n^* \sum\limits_{j=0}^s
				(-1)^j \binom{n}{j} \sum\limits_{\ell=0}^{k-n}
				\delta_{s, \ell+j} \\
				& = \frac{1}{k + 1} \sum\limits_{n=0}^s
				\binom{k+1}{n} B_n^* 
				\sum\limits_{j=\max\{0, s+n-k\}}^n
				(-1)^j \binom{n}{j} \\
				& \quad + \frac{1}{k + 1}
				\sum\limits_{n=s+1}^k
				\binom{k+1}{n} B_n^*
				\sum\limits_{j=\max\{0, s+n-k\}}^s
				(-1)^j \binom{n}{j}.
			\end{align*}
			As long as $\max\{0, s+n-k\} = 0$, the first sum over $j$ will be 
			zero as it is just the binomial expansion of $(1-1)^n$, except for 
			$n = 0$. Hence we get a special case and a shorter first sum over 
			$n$. In the sums over $j$ we use again the binomial expansion of 
			$(1-1)^n$ and get 
			\begin{align*}
				\mathcal{K}(k,s)
				& = 
				\frac{1}{k+1} \left[ 
				1 + \sum\limits_{k+1-s}^s 
				\binom{k+1}{n} B_n^* \left(
				- \sum\limits_{j=0}^{s+n-k-1} 
				(-1)^j \binom{n}{j}\right) \right. \\
				& \quad + \left. 
				\sum\limits_{n=s+1}^k 
				\binom{k+1}{n} B_n^* \left( 
				- \sum\limits_{j=0}^{s+n-k-1} 
				(-1)^j\binom{n}{j} - 
				\sum\limits_{j=s+1}^n 
				(-1)^j \binom{n}{j} \right) \right].
			\end{align*}
			Now it is helpful to use \eqref{Alg:Lemmastuff:1_Bin1} and 
			$\binom{k}{n-k} = \binom{k}{n}$. We also get $(-1)^n$-terms which 
			we can put together with the $B_n^*$ to get $B_n$:
			\begin{align*}
				\mathcal{K}(k,s)
				& = 
				\frac{1}{k+1} \left[
				1+ \sum\limits_{n=k+1-s}^s
				\binom{k+1}{n} B_n (-1)^{k-s}
				\binom{n-1}{k-s} \right. \\
				& \quad + \left.
				\sum\limits_{n=s+1}^k
				\binom{k+1}{n}B_n \left(
				(-1)^{k-s}\binom{n-1}{k-s} + (-1)^{n+s}\binom{n-1}{s}
				\right) \right].
			\end{align*}
			We finally made the $\delta$ disappear. Hence we must compute
			$\mathcal{K}(k,s+1) - \mathcal{K}(k,s)$. Since we want to show 
			that it is $0$, we can multiply it with $k + 1$ in order to get 
			rid of the factor in front:
			\begin{align*}
				& \quad
				(k+1) \left( \mathcal{K}(k,s+1) - \mathcal{K}(k,s) \right)\\
				& =
				\sum\limits_{n=k-s}^{s+1}
				\binom{k+1}{n} B_n (-1)^{k-s-1}
				\binom{n-1}{k-s-1} - \sum\limits_{n=k+1-s}^s
				\binom{k+1}{n} B_n (-1)^{k-s} \binom{n-1}{k-s} \\
				& \quad
				+ \sum\limits_{n=s+2}^k
				\binom{k+1}{n} B_n \left(
				(-1)^{k-s-1}\binom{n-1}{k-s-1}
				+ (-1)^{n+s+1}\binom{n-1}{s+1}
				\right) \\
				& \quad
				- \sum\limits_{n=s+1}^k
				\binom{k+1}{n} B_n \left(
				(-1)^{k-s}\binom{n-1}{k-s} + (-1)^{n+s}\binom{n-1}{s}
				\right) \\
				& =
				-\sum\limits_{n=k-s}^k
				\binom{k+1}{n} B_n (-1)^{k-s}
				\binom{n-1}{k-s-1} -\sum\limits_{n=k-s+1}^k
				\binom{k+1}{n}B_n (-1)^{k-s} \binom{n-1}{k-s} \\
				& \quad -
				\sum\limits_{n=s+2}^k
				\binom{k+1}{n}B_n (-1)^{n+s}
				\binom{n-1}{s+1} - \sum\limits_{n=s+1}^k
				\binom{k+1}{n}B_n (-1)^{n+s} \binom{n-1}{s} \\
				& =
				-\sum\limits_{n=k-s}^k
				\binom{k+1}{n}B_n (-1)^{k-s} \left(
				\binom{n-1}{k-s-1} + \binom{n-1}{k-s}
				\right) \\
				& \quad -
				\sum\limits_{n=s+1}^k
				\binom{k+1}{n} B_n (-1)^{n+s}
				\left( \binom{n-1}{s+1} + \binom{n-1}{s} \right).
			\end{align*}
			We have rearranged the sums, added some zeros and shortened the 
			expression. Now we will use the recursion formula for the binomial 
			coefficients
			\begin{equation*}
				\binom{n-1}{k-1} + \binom{n-1}{k}
				=
				\binom{n}{k}
			\end{equation*}
			and our binomial multiplication equality 
			\eqref{Alg:Lemmastuff:1_Bin2}:
			\begin{align*}
				& =
				- \sum\limits_{n=k-s}^k
				\binom{k+1}{n} B_n (-1)^{k-s} \binom{n}{k-s}
				- \sum\limits_{n=s+1}^k
				\binom{k+1}{n} B_n (-1)^{n+s} \binom{n}{s+1} \\
				& = 
				- \sum\limits_{n=k-s}^k 
				\binom{k+1}{s+1} \binom{s+1}{n+s-k} B_n (-1)^{k-s}
				- \sum\limits_{n=s+1}^k
				\binom{k+1}{s+1} \binom{k-s}{n-s-1} B_n (-1)^{n+s}.
			\end{align*}
			Since we want to show that this is $0$, we can divide by
			$\binom{k+1}{s+1}$ which will never be zero because 
			$s \in \{0, 1, \ldots, k\}$. After doing so, can use $n > 1$ in 
			the second sum and thus only even $n$ will show up, because for 
			odd $n$ the Bernoulli numbers are zero. For this reason we have 
			$(-1)^n = 1$. Then we rewrite these sums by shifting the indices 
			and we add two zeros:
			\begin{align*}
				& \quad
				- \sum\limits_{n=k-s}^k
				\binom{s+1}{n+s-k} B_n (-1)^{k-s}
				+ \sum\limits_{n=s+1}^k
				\binom{k-s}{n-s-1} B_n (-1)^{s+1} \\
				& =
				(-1)^{s+1} \sum\limits_{\ell=0}^{k-s-1}
				\binom{k-s}{\ell} B_{\ell+s+1}
				- (-1)^{k-s}\sum\limits_{\ell=0}^s 
				\binom{s+1}{\ell} B_{\ell + k - s} \\
				& =
				(-1)^{s+1} \sum\limits_{\ell=0}^{k-s}
				\binom{k-s}{\ell} B_{\ell+s+1}
				- (-1)^{k-s} \sum\limits_{\ell=0}^{s+1}
				\binom{s+1}{\ell} B_{\ell + k - s}\\
				& \quad
				- (-1)^{s+1} \binom{k-s}{k-s} B_{k+1} 
				+ (-1)^{k-s} \binom{s+1}{s+1} B_{k+1}.
			\end{align*}
			The first two terms give the Carlitz-identity 
			\eqref{Alg:Lemmastuff:1_Carl} and vanish. So we are left with 
			the last two terms and get
			\begin{equation*}
				- (-1)^{s+1} B_{k+1} + (-1)^{k-s} B_{k+1}
				=
				(-1)^s B_{k+1} \left(1 + (-1)^k\right) = 0,
			\end{equation*}
			since the bracket will be zero if $k$ is odd and $B_{k+1}=0$ if 
			$k$ is even.
		\end{enumerate}
	\end{subproof}
	Finally, Lemma \ref{Alg:FuckingDeadlyLemma} is proven.
\end{proof}
In Lemma \ref{Formulas:Lemma:LinearMonomial1}, we will see that also	$\ast_z$ 
fulfils this identity. Hence $\ast_z = \widehat{\star}_z$. We only need to 
show $\widehat{\star}_z = \star_z$. For $z = 1$, the two 	maps are clearly 
identical and therefore we find
\begin{equation*}
	\xi^k \star_1 \eta
	=
	\sum\limits_{n=0}^k
	\binom{k}{n} B_n^* \xi^{k-n}
	\left( \ad_{\xi} \right)^k(\eta).
\end{equation*}
But now $\widehat{\star}_z = \star_z$ follows from the definition of $\star_z$: 
we just have to plug in powers of $z$ and find \eqref{Alg:GuttStarFormula}. 
So with the proofs in Chapter 4, we will have proven the following theorem:
\begin{theorem}
	\label{Alg:Thm:ThreeStarsAreOne}
	The three maps $\star_z$, $\widehat{\star}_z$ and $\ast_z$ 
	coincide on $\Sym^{\bullet}(\lie{g})$ and define star products.
\end{theorem}

% Chapter 4
%

%
% Chapter 4 of my master thesis:
% The formulas For the Gutt star product
%

\chapter{Formulas for the Gutt star product}

We have seen some results on the Baker-Campbell-Hausdorff series and an 
identity for the Gutt star product. The latter one, stated in Theorem 
\ref{Alg:Thm:ThreeStarsAreOne}, will be a very useful tool in the following, 
since we want to get explicit formulas for $\star_z$. There is still a part of 
the proof missing, but this will be caught up at the beginning of the first 
section of this chapter. From there, we will come to a first easy formula for 
$\star_z$. Afterwards, we will use the same procedure to find two more 
formulas for it: the first is a rather involved one for the $n$-fold star 
product of vectors. It will not be helpful for algebraic computations, 
but very useful for estimates. The second one is a more explicit formula for 
the product of two monomials.

From those formulas, we will be able to draw some easy, but nice 
conclusion in the second section and we will prove the classical and the 
semi-classical limit. Then, we will show how to calculate the Gutt star 
product explicitly by computing two easy examples.

\section{Formulas for the Gutt Star Product}
\label{sec:chap4_Formulas}

%
% An Iterative Approach from Linear Terms
%

\subsection{A Monomial with a Linear Term}

The easiest case for which we will develop a formula is surely the 
following one: for a given Lie algebra $\lie{g}$ and $\xi, \eta \in 
\lie{g}$ we would like to compute
\begin{equation*}
    \xi^k \star_z \eta
    =
    \sum\limits_{n=0}^k
    z^n C_n(\xi^k, \eta)
\end{equation*}
We have already done this for $\star_z$ and $\widehat{\star}_z$ in Lemma 
\ref{Alg:FuckingDeadlyLemma}, and want to do the same for $\ast_z$ now. This will finish the proof of the equality of 
the star products from Theorem~\ref{Alg:Thm:ThreeStarsAreOne}. We will use that
\begin{equation}\label{Formulas:MonomialDerivative}
    \xi^k
    =
    \frac{\partial^k}{\partial t^k}
    \At{t = 0} \exp(t \xi).
\end{equation}
Now we have all the ingredients to prove the following lemma:
\begin{lemma}
    \label{Formulas:Lemma:LinearMonomial1}
    Let $\lie{g}$ be a Lie algebra and $\xi, \eta \in \lie{g}$. We
    have the following identity for $\ast_z$:
    \begin{equation}
        \label{Formulas:LinearMonomial1}
        \xi^k \ast_z \eta
        =
        \sum\limits_{j=0}^k
        \binom{k}{j} z^j B_j^*
        \xi^{k-j}(\ad_{\xi})^j (\eta).
    \end{equation}
\end{lemma}
\begin{proof}
    We start from the simplified form for the Baker-Campbell-Hausdorff 
    series from Equation \eqref{Alg:BCHFirstOrderXi} in 
    Proposition~\ref{Alg:Prop:BCHFristOrder}:
    \begin{equation*}
		\bch{\xi}{\eta}
		=
		\xi 
		+ 
		\sum\limits_{n = 0}^{\infty}
		\frac{B_n^*}{n!}
		\left( \ad_{\xi} \right)^n (\eta)
		+
		\mathcal{O}(\eta^2).
    \end{equation*}
    If we insert this into the definition of the Drinfel'd star product 
    and use Equation~\eqref{Formulas:MonomialDerivative} we get
    \begin{align*}
        \xi^k \ast_z \eta
        & =
        \frac{\partial^k}{\partial t^k}
        \frac{\partial}{\partial s}
        \At{t=0, s=0}
        \exp \left(
            \frac{1}{z} \bch{z t \xi}{z s \eta}
        \right)
        \\
        & =
        \frac{\partial^k}{\partial t^k}
        \frac{\partial}{\partial s}
        \At{t=0, s=0}
        \exp \left(
            t \xi + \sum\limits_{j=0}^{\infty}
            z^j \frac{B_j^*}{j!}
            \left( \ad_{t \xi} \right)^j
            (s \eta)
            + \mathcal{O}(\eta^2)
        \right).
    \end{align*}
    We see that only terms which have exactly $k$ of the
    $\xi$'s in them and which are linear in $\eta$ will
    contribute. This means we can cut off the sum at $j = k$ and omit 
    higher orders in $\eta$. We now use the exponential series,
    cut it at $k$ for the same reason and get
    \begin{align*}
        \xi^k \ast_z \eta
        & =
        \frac{\partial^k}{\partial t^k}
        \frac{\partial}{\partial s}
        \At{t=0, s=0}
        \sum\limits_{n=0}^{k}
        \frac{1}{n!}
        \left(
            t \xi
            +
            \sum\limits_{j=0}^{k}
            (zt)^j \frac{B_j^*}{j!}
            \left(\ad_{\xi}\right)^j
            (s \eta)
        \right)^n
        \\
        & =
        \frac{\partial^k}{\partial t^k}
        \frac{\partial}{\partial s}
        \At{t=0, s=0}
        \sum\limits_{n=0}^{k}
        \frac{1}{n!}
        \sum\limits_{m = 0}^n
        \binom{n}{m}
        (t \xi)^{n - m}
        \left(
            \sum\limits_{j=0}^{k}
            (zt)^j \frac{B_j^*}{j!}
            \left(\ad_{\xi}\right)^j
            (s \eta)
        \right)^m
        \\
        & =
        \frac{\partial^k}{\partial t^k}
        \frac{\partial}{\partial s}
        \At{t=0, s=0}
        \left(
            \sum\limits_{n=0}^{k}
            \frac{1}{n!}
            (t \xi)^n
            +
            \sum\limits_{n=0}^{k}
            \sum\limits_{j=0}^k
            \frac{1}{(n - 1)!} t^{n + j - 1}
            z^j \frac{B_j^*}{j!}
            \xi^{n - 1}
            \left( \ad_{\xi} \right)^j
            (s \eta)
        \right).
    \end{align*}
    In the last step we set $m = 1$ since the other term have either too
    many or not enough $\eta$'s and will vanish because of the differentiation 
    with respect to $s$. We can finally differentiate to get the formula
    \begin{align*}
        \xi^k \ast_z \eta
        & =
        \sum\limits_{n=0}^k
        \sum\limits_{j=0}^k
        \delta_{k, n + j - 1}
        \frac{k!}{j! (n - 1)!}
        z^j B_j^*
        \xi^{n - 1}
        \left(\ad_{\xi}\right)^j
        (\eta)
        \\
        & =
        \sum\limits_{j=0}^k
        \binom{k}{j}
        z^j B_j^*
        \xi^{k - j}
        \left( \ad_{\xi} \right)^j
        (\eta),
    \end{align*}
    which is the wanted result.
\end{proof}
\begin{remark}
    We have proven the equality of the  star products $\widehat{\star}_z$
    $\ast_z$ by deriving an easy formula for both of them. From now on, we 
    will get all other formulas from $\ast_z$, since this is the 
    one which is easier to compute.
\end{remark}
Now it is actually easy to get the formula for monomials of the form $\xi_1 \ldots \xi_k$ with $\eta \in \lie{g}$:
\begin{proposition}
	\label{Formulas:Prop:LinearMonomial2}
    Let $\lie{g}$ be a Lie algebra and $\xi_1, \ldots, \xi_k, \eta \in 
    \lie{g}$. We have
    \begin{align}
    		\label{Formulas:LinearMonomial2}
	    	\xi_1 \cdots \xi_k \star_z \eta
    		& =
    		\sum\limits_{j=0}^k
    		\frac{1}{k!} \binom{k}{j}
    		z^j B_j^*
    		\sum\limits_{\sigma \in S_k}
    		[\xi_{\sigma(1)}, 
    			[ \ldots [\xi_{\sigma(j)}, \eta] \ldots ]
    		]
    		\xi_{\sigma(j+1)} \cdots \xi_{\sigma(k)} \text{ and}
    		\\
		\label{Formulas:LinearMonomial2T}
	    	\eta \star_z \xi_1 \cdots \xi_k
    		& =
    		\sum\limits_{j=0}^k
    		\frac{1}{k!} \binom{k}{j}
    		z^j B_j
    		\sum\limits_{\sigma \in S_k}
    		[\xi_{\sigma(1)}, 
    			[ \ldots [\xi_{\sigma(j)}, \eta] \ldots ]
    		]
    		\xi_{\sigma(j+1)} \cdots \xi_{\sigma(k)}.
    \end{align}
\end{proposition}
\begin{proof}
	We get the result by just polarizing the formula from Lemma 
	\ref{Formulas:Lemma:LinearMonomial1}. Let $\xi_1, \ldots, \xi_k \in 
	\lie{g}$, then we introduce the parameters $t_i$ for $i = 
	1, \ldots, k$ and set
	\begin{equation*}
		\Xi
		=		
		\Xi(t_1, \ldots, t_k)
		=
		\sum\limits_{i=1}^k t_i \xi^i.
	\end{equation*}
	Then we see that
	\begin{equation*}
		\xi_1 \cdots \xi_k
		=
		\frac{1}{k!}
		\frac{\partial^k}{\partial t_1 \cdots \partial t_k}
		\At{t_1, \cdots, t_k = 0}
		\Xi^k
	\end{equation*}
	since for every $i = 1, \ldots, k$ we have
	\begin{equation}\label{Formulas:LittleHelp1}
		\frac{\partial}{\partial t_i}
		\At{t_i = 0} \Xi
		=
		\xi_i.
	\end{equation}
	By writing out the $\Xi$'s and using multilinearity, we find
	\begin{align*}
		\xi_1 \cdots \xi_k \star_z \eta
		& =
		\frac{1}{k!}
		\frac{\partial^k}{\partial t_1 \cdots \partial t_k}
		\At{t_1, \cdots, t_k = 0}
        \sum\limits_{j=0}^k
        \binom{k}{j} z^j B_j^*
        \Xi^{k-j}(\ad_{\Xi})^j (\eta)
        \\
        & =
		\frac{1}{k!}
        \sum\limits_{j=0}^k
        \binom{k}{j} z^j B_j^*
		\sum\limits_{
			\{i_1, \ldots, i_k\}
			\in
			\{1, \ldots, k\}^k
		}
		\frac{\partial^k}{\partial t_1 \cdots \partial t_k}
		\At{t_1, \cdots, t_k = 0}
        t_{i_1} \cdots t_{i_k}
        \\
        & \qquad
        \cdot
        \xi_{i_1} \cdots \xi_{i_{k-j}}
        \ad_{\xi_{i_{k-j+1}}} 
        \circ \cdots \circ  
        \ad_{\xi_{i_k}}
        (\eta)
        \\
        & =
    		\sum\limits_{j=0}^k
    		\frac{1}{k!} \binom{k}{j}
    		z^j B_j^*
    		\sum\limits_{\sigma \in S_k}
    		[\xi_{\sigma(1)}, 
    			[ \ldots [\xi_{\sigma(j)}, \eta] \ldots ]
    		]
    		\xi_{\sigma(j+1)} \cdots \xi_{\sigma(k)}.
	\end{align*}
	In the last step, all expression which did not contain each $\xi_i$ 
	exactly once disappeared due to the differentiation. The proof of 
	Equation~\eqref{Formulas:LinearMonomial2T} works analogously.
\end{proof}
\begin{remark}
	\label{Formulas:Rem:EasierFormulaAlreadyKnown}
	This formula is actually not a new result: Gutt already gave it in her 
	paper \cite[Prop. 1]{gutt:1983a} and referred to Dixmier 
	\cite[part 2.8.12 (c)]{dixmier:1977a}, who already gave it in his 
	textbook. It can also be found in the diploma thesis of Neumaier 
	\cite[Rem. 5.2.8]{neumaier:1998a} and a work due to Kathotia 
	\cite[Eq. 2.23]{kathotia:1998a:pre}. Probably the first one to mention it 
	was Berezin in \cite[Eq. 30]{berezin:1967a}.
\end{remark}

% A first general Formula
%
\subsection{An Iterated Formula for the General Case}

Proposition \ref{Formulas:Prop:LinearMonomial2} allows theoretically 
to get a formula for the case of $\xi_1, \ldots, \xi_k \in \lie{g}$
\begin{equation*}
	\xi_1 \star_z \ldots \star_z \xi_k
	=
	\sum\limits_{j=0}^k
	C_{z,j} \left(
		\xi_1, \ldots, \xi_k
	\right)
\end{equation*}
which we will need to prove the functoriality of our later construction.
This could also be used to give an alternative proof for our main theorem.
Unluckily, this approach has a problem: iterating this formula, we get 
strangely nested Lie brackets, which would be very difficult to bring 
into a nice form with Jacobi identity. So this is not a 
good way to find a handy formula for the usual star product of two 
monomials. Nevertheless, we want to pursue it for a moment, since we 
will get an equality which will be, although rather involved looking, 
very useful in the following: for analytic observations, it will be 
enough to put (even rough) estimates on it and the exact nature of the 
combinatorics in the formula will not be important. Hence we rewrite 
Equation \eqref{Formulas:LinearMonomial2} in order to cook up such a 
formula.

Take $\xi_1, \ldots, \xi_k, \eta \in \lie{g}$, then we have
\begin{equation*}
	\xi_1 \cdots \xi_k \star_z \eta
	=
	\sum\limits_{n = 0}^k 
	C_n \left( \xi_1 \ldots \xi_k, \eta \right)
\end{equation*}
with the $C_n$ being as bilinear operators which are given explicitly 
on monomials by
\begin{align}
	\label{Formulas:MultipleStarDef1}
	C_n^k
	\colon
	\Sym^k(\lie{g})
	\times
	\lie{g}
	& 
	\longrightarrow
	\Sym^{k - n + 1}(\lie{g})
	\\
	\label{Formulas:MultipleStarDef2}
	(\xi_1 \cdots \xi_k, \eta)
	&
	\longmapsto
	\frac{1}{k!}
	\sum\limits_{\sigma \in S_k}
	\binom{k}{j} B_j^* z^j
	\xi_{\sigma(1)} \cdots \xi_{\sigma(k-j)}
	[\xi_{\sigma(k-j+1)}, [ \ldots, [\xi_{\sigma(k)}, \eta]]] 
\end{align}
with
\begin{equation*}
	C_n 
	= 
	\sum\limits_{k = 0}^{\infty}
	C_n^k.
\end{equation*}
This gives us a good way of writing the $n$-fold star product of vectors:
\begin{proposition}
	\label{Formulas:Prop:MultipleStars}
	Let $\lie{g}$, $2 \leq k \in \mathbb{N}$  and $\xi_1, \ldots, \xi_k 
	\in \lie{g}$. Then we have
	\begin{equation}
		\label{Formulas:MultipleStars}
		\xi_1 \star_z \ldots \star_z \xi_k		
		=
		\sum\limits_{\substack{
			1 \leq j \leq k-1 \\
			i_j \in \{0, \ldots, j\}
		}}
		z^{i_1 + \ldots + i_{k-1}}
		C_{i_{k-1}}
		\left(
			\ldots C_{i_2}
			\left(
				C_{i_1}
				\left( \xi_1, \xi_2 \right)
				, \xi_3	
			\right) 
			\ldots, \xi_{k}
		\right).
	\end{equation}
\end{proposition}
\begin{proof}
	This is an easy prof by induction over $k$. For $k = 2$ the statement is 
	clearly true. For the step $k \rightarrow k + 1$ we get
	\begin{align*}
		\xi_1 \star_z \ldots \star_z \xi_{k+1}
		& =
		\Bigg(
			\sum\limits_{\substack{
				1 \leq j \leq k-1 \\
				i_j \in \{0, \ldots, j\}
			}}
			z^{i_1 + \cdots + i_{k-1}}
			C_{i_{k-1}}
			\left(
				\ldots C_{i_2}
				\left(
					C_{i_1} 
					\left( \xi_1, \xi_2 \right)
					, \xi_3	
				\right) 
				\ldots, \xi_{k}
			\right)
		\Bigg)
		\star_z \xi_{k+1}
		\\
		& = 
		\sum\limits_{i_k = 0}^k
		z^{i_k}
		C_{i_k}
		\Bigg(
			\sum\limits_{\substack{
				1 \leq j \leq k-1 \\
				i_j \in \{0, \ldots, j\} \\
			}}
			z^{i_1 + \cdots + i_{k-1}}
			C_{i_{k-1}}
			\left(
				\ldots C_{i_2}
				\left(
					C_{i_1} 
					\left( \xi_1, \xi_2 \right)
					, \xi_3	
				\right) 
				\ldots, \xi_{k}
			\right)
			, \xi_{k+1}
		\Bigg)
		\\
		& = 
		\sum\limits_{\substack{
			1 \leq j \leq k \\
			i_j \in \{0, \ldots, j\} \\
		}}
		z^{i_1 + \cdots + i_k}
		\left(
			C_{i_{k-1}}
			\left(
				\ldots C_{i_2}
				\left(
					C_{i_1} 
					\left( \xi_1, \xi_2 \right)
					, \xi_3	
				\right) 
				\ldots, \xi_{k}
			\right)
			, \xi_{k+1}
		\right)
	\end{align*}
\end{proof}
\begin{remark}
	\mbox{}
	\begin{remarklist}
		\item
		Of course, Proposition~\ref{Formulas:Prop:MultipleStars}  is an easy 
		consequence from Proposition \ref{Formulas:Prop:LinearMonomial2}.
		It's value, however, is that we know how the $C_n$'s look like and 
		what the summation range in \eqref{Formulas:MultipleStars} is. This 
		will allow us to put estimates on things like iterated star products.
		
		\item
		As already mentioned, we would get an identity for the star product of 
		two monomials via
		\begin{equation}
			\label{Formulas:2MonomialsWeird}
			\xi_1 \cdots \xi_k \star_z \eta_1 \cdots \eta_{\ell}
			=
			\frac{1}{k! \ell!}
			\sum\limits_{\sigma \in S_k}
			\sum\limits_{\tau \in S_{\ell}}
			\xi_{\sigma(1)} \star_z \cdots \star_z \xi_{\sigma(k)}
			\star_z
			\eta_{\tau(1)} \star_z \cdots \star_z \eta_{\tau(\ell)}.
		\end{equation}
		This can be proven from the definition of the map $\mathfrak{q}_z$.
		Unfortunately, this would give a very clumsy formula to deal with.
	\end{remarklist}
\end{remark}

% A Formula for two Monomials
%

\subsection{A Formula for two Monomials}

If we want to get an identity for the star product of two monomials, we have 
to go back to Equation \eqref{Alg:DrinfeldStar}. This will not give a simple 
looking formula either, but we will at least be able to do some 
computations with concrete examples. As a first step, we must introduce 
a bit of notation:
\begin{definition}[G-Index]
	\label{Def:GuttIndex}
	Let $k, \ell, n \in \mathbb{N}$ and $r = k + \ell - n$. 
	Then we call an $r$-tuple $J$
	\begin{equation*}
		J = (J_1, \ldots, J_r) 
		= 
		((a_1, b_1), \ldots, (a_r, b_r)) 
	\end{equation*}
  	a G-index if it fulfils the following properties:
	\begin{enumerate}[(i)]
		\item
		$J_i \in \{0, 1, \ldots, k\} \times \{0, 1, \ldots, \ell\}$
		
  		\item 
		$|J_i| 
		= 
		a_i + b_i \geq 1 
		\quad \forall_{i = 1, \ldots, r}$
		
		\item 
		$\sum\limits_{i=1}^{r} a_i = k$ 
		and 
		$\sum\limits_{i=1}^{r} b_i = \ell$
		
		\item
		The tuple is ordered in the following sense:
		$i>j \Rightarrow |J_i| \geq |J_j| \quad \forall_{i,j = 1, 
		\ldots, r}$ and $|a_i| \geq |a_j|$ if $|J_i| = |J_j|$
		
		\item 
		If $a_i = 0$ [or $b_i = 0$] for some $i$, 
		then $b_i = 1$ [or $a_i = 1$].
	\end{enumerate}
	We call the set of all such $G$-indices	$\mathcal{G}_r(k,\ell)$.
\end{definition}
\begin{definition}[G-Factorial]
	\label{Def:GuttFactorial}
	Let $J = ((a_1, b_1), \ldots, (a_r, b_r)) \in \mathcal{G}_r(k,
	\ell)$ be a G-Index. We set for a given tuple $(a,b) \in \{0, 1, 
	\ldots, k\} \times \{0, 1, \ldots, \ell\}$
	\begin{equation*}
		\#_J (a,b)
		= 
		\textrm{ number of times that $(a,b)$ appears in } J.
	\end{equation*}
	Then we define the G-factorial of $J \in \{0, 1, \ldots, k\} \times 
	\{0, 1, \ldots, \ell\}$ as
	\begin{equation*}
		J!
		= 
		\prod\limits_{
			(a,b) \in 
			\{0, 1, \ldots, k\} 
			\times 
			\{0, 1, \ldots, \ell\}
		}
		\left( \#_J (a,b) \right)!
	\end{equation*}
\end{definition}
Each pair $(a,b)$ will later correspond to $\bchparts{a}{b}{\xi}{\eta}$.
Now we can state a good formula for the Gutt star product:
\begin{lemma}
	\label{Formulas:Lemma:2MonomialsFormula1}
	Let $\lie{g}$ be a Lie algebra, $\xi, \eta \in \lie{g}$ and $k, 
	\ell \in \mathbb{N}$. Then we have the following identity for the 
	Gutt star product:
	\begin{equation*}
    	\xi^k \star_z \eta^{\ell}
    	=
	    \sum\limits_{n=0}^{k + \ell -1}
    	z^n
    	C_n \left( \xi^k, \eta^{\ell} \right),
	\end{equation*}
	where the $C_n$ are given by
	\begin{align}
		\label{Formulas:2MonomialsExplicit}
        C_n \left( \xi^k, \eta^{\ell} \right)
        & =
        \frac{k! \ell!}{(k + \ell - n)!}
        \sum\limits_{\substack{a_1, b_1, \ldots, a_r, b_r \geq 0 \\
            a_i + b_i \geq 1 \\
            a_1 + \cdots + a_r = k \\
            b_1 + \cdots + b_r = \ell
        }}
        \bchparts{a_i}{b_i}{\xi}{\eta}
        \cdots
        \bchparts{a_r}{b_r}{\xi}{\eta}
        \\
        \label{Formulas:2MonomialsFormula1}
        & =
        \sum\limits_{J \in \mathcal{G}_{k + \ell - n}(k, \ell)}
        \frac{k! \ell!}{J!}
        \prod\limits_{i = 1}^{k + \ell - n}        
        \bchparts{a_i}{b_i}{\xi}{\eta}
	\end{align}
	and the product is taken in the symmetric tensor algebra.
\end{lemma}
\begin{proof}
	We want to see what the $C_n$ look like. Let's denote $r = k 
	+ \ell - n$ for brevity. Then we have
	\begin{equation*}
		C_n \left( \xi^k, \eta^{\ell} \right)
		\in \Sym^r(\lie{g}).
	\end{equation*}
	Of course, the only part of the series
	\begin{equation*}
		\exp\left(
			\frac{1}{z} \bch{z \xi}{z \eta}
		\right)
		=
		\sum\limits_{n = 0}^{k + \ell}
		\left(
			\frac{1}{z} \bch{z \xi}{z \eta}
		\right)^n
		+ \mathcal{O}(\xi^{k + 1}, \eta^{\ell + 1})
	\end{equation*} 
	which lies in $\Sym^r(\lie{g})$ is the summand for 
	$n = r $. We introduce the formal parameters $t$ and $s$.
	Since we differentiate with respect to them, we can omit 
	terms of higher orders in $\xi$ and $\eta$ than $k$ 
	and $\ell$ respectively.
    \begin{align*}
        z^n C_n \left( \xi^k, \eta^{\ell} \right)
        & =
        \frac{\partial^k}{\partial t^k}
        \frac{\partial^{\ell}}{\partial s^{\ell}}
        \At{t,s = 0}
        \frac{1}{z^r}
        \frac{\bch{z t \xi}{z s \eta}^r    }{r!}
        \\
        & =
        \frac{1}{z^r}
        \frac{1}{r!}
        \frac{\partial^k}{\partial t^k}
        \frac{\partial^{\ell}}{\partial s^{\ell}}
        \At{t,s = 0}
        \left(
            \sum\limits_{j = 1}^{k + \ell}
            \bchpart{j}{z t \xi}{z s \eta}
        \right)^r
        \\
        & =
        \frac{1}{z^r}
        \frac{1}{r!}
        \frac{\partial^k}{\partial t^k}
        \frac{\partial^{\ell}}{\partial s^{\ell}}
        \At{t,s = 0}
        \sum\limits_{\substack{
        	j_1, \ldots, j_r \geq 1 \\
            j_1 + \ldots + j_r = k + \ell
        }}
        \bchpart{j_1}{z t \xi}{z s \eta} 
        \cdots
        \bchpart{j_r}{z t \xi}{z s \eta}
        \\
        & =
        z^n
        \frac{k! \ell!}{r!}
        \sum\limits_{\substack{a_1, b_1, \ldots, a_r, b_r \geq 0 \\
            a_i + b_i \geq 1 \\
            a_1 + \cdots + a_r = k \\
            b_1 + \cdots + b_r = \ell
        }}
        \bchparts{a_i}{b_i}{\xi}{\eta}
        \cdots
        \bchparts{a_r}{b_r}{\xi}{\eta}
    \end{align*}
    We sum over all possible arrangements of the $(a_i, b_i)$. In order 
    to find an easier summation range, we put the ordering from Definition 
    \ref{Def:GuttIndex} on these multi-indices and avoid therefore double 
    counting. We loose the freedom of arranging the $(a_i, b_i)$ and 
    need to count the number of multi-indices $((a_1, b_1), \ldots, 
    (a_r, b_r))$ which belong to the same G-index $J$. This number will 
    be $\frac{r!}{J!}$, since we can not interchange the $(a_i, b_i)$ any 
    more (therefore $r!$), unless they are equal (therefore $J!^{-1}$). 
    Since the ranges of the $(a_i, b_i)$ in Equation 
    \eqref{Formulas:2MonomialsExplicit} and of  the elements in 
    $\mathcal{G}_r(k, \ell)$ are the same, we can change the summation there 
    to $J \in \mathcal{G}_r(k, \ell)$ and need to multiply by $\frac{r!}{J!}$. 
    This gives
    \begin{equation*}
    	z^n C_n \left( \xi^k, \eta^{\ell} \right)
    	=
    	z^n \frac{k! \ell!}{J!}
    	\sum\limits_{J \in \mathcal{G}_r(k, \ell)}
    	\bchparts{a_i}{b_i}{\xi}{\eta}
        \cdots
        \bchparts{a_r}{b_r}{\xi}{\eta}
    \end{equation*}
    which is equivalent to Equation \eqref{Formulas:2MonomialsFormula1}.
\end{proof}
Now we need to generalize this to factorizing tensors. To do so, we 
need a last definition:
\begin{definition}
	\label{Def:BCHTilde}
	Let $a,b \in \mathbb{N}$ and $\xi_1, \ldots, \xi_a, \eta_1, \ldots, 
	\eta_b \in \lie{g}$. Then we define by
	\begin{equation*}
		\widetilde{\mathrm{BCH}}_{a,b}
		\colon
		\lie{g}^{a + b}
		\longrightarrow
		\lie{g}
	\end{equation*}
	the map which we get when we replace in $\bchparts{a}{b}{\xi}{\eta}$ 
	the $i$-th $\xi$ by $\xi_i$ and the $j$-th $\eta$ by $\eta_j$
	for $i = 1, \ldots, a$ and $j = 1, \ldots, b$.
\end{definition}
\begin{proposition}
	\label{Formulas:Prop:2MonomialsFormula2}
	Let $\lie{g}$ be a Lie algebra, $k, \ell \in \mathbb{N}$ and $\xi_1, 
	\ldots, \xi_k, \eta_1, \ldots, \eta_{\ell} \in \lie{g}$. Then we 
	have the following identity for the Gutt star product:
	\begin{equation*}
    		\xi_1 \ldots \xi_k \star_z \eta_1 \ldots \eta_{\ell}
    		=
		\sum\limits_{n=0}^{k + \ell -1}
    		z^n C_n
    		\left( 
    			\xi_1 \cdots \xi_k, \eta_1 \cdots \eta_{\ell}
    		\right),
	\end{equation*}
	where the $C_n$ are given by
	\begin{align}
		\nonumber
        C_n
        \left( 
    			\xi_1 \cdots \xi_k, \eta_1 \cdots \eta_{\ell}
    		\right)
        & =
        \sum\limits_{J \in \mathcal{G}_{k + \ell - n}(k, \ell)}
        \frac{1}{J!}
        \sum\limits_{\sigma \in S_k}
        \sum\limits_{\tau \in S_{\ell}}
        \prod_{i=1}^{l + \ell - n}
        \widetilde{\mathrm{BCH}}_{a_i, b_i}
        \big( \xi_{\sigma(a_1 + \cdots + a_{i - 1} + 1)}, 
            \ldots, 
        \\
        \label{Formulas:2MonomialsFormula2}
        & \qquad            
            \ldots, \xi_{\sigma(a_1 + \cdots + a_i)} \big)
        \big( \eta_{\tau(b_1 + \cdots + b_{i - 1} + 1)}, 
            \ldots, \eta_{\tau(b_1 + \cdots + b_i)} \big).
        \\
        \nonumber
        & =
        \frac{1}{(k + \ell - n)!}
        \sum\limits_{\sigma \in S_k, \tau \in S_{\ell}}
        \sum\limits_{\substack{a_1, b_1, \ldots, a_r, b_r \geq 0 \\
            a_i + b_i \geq 1 \\
            a_1 + \cdots + a_r = k \\
            b_1 + \cdots + b_r = \ell
          }
        }
        \\
        \nonumber
        & \qquad
        \bchtilde{a_i}{b_i}
        {\xi_{\sigma(1)}, \ldots, \xi_{\sigma(a_1)}}
        {\eta_{\tau(1)}, \ldots, \eta_{\tau(b_1)}}
        \cdots
        \\
        \label{Formulas:2MonomialsFormula22}
        & \qquad
        \bchtilde{a_r}{b_r}
        {\xi_{\sigma(k - a_r + 1)}, \ldots, \xi_{\sigma(k)}}
        {\eta_{\tau(\ell - b_r + 1)}, \ldots, \eta_{\tau(\ell)}}.
	\end{align}
\end{proposition}
\begin{proof}
	The proof relies on polarization again and is completely analogous 
	to the one of Proposition \ref{Formulas:Prop:LinearMonomial2}. We set
	\begin{equation*}
		\Xi
		=
		\sum\limits_{i=1}^k t_i \xi^i
		\quad \text{ and } \quad
		\Eta
		=
		\sum\limits_{i=1}^{\ell} t_j \eta^j
		.
	\end{equation*}
	Then it is easy to see that we will get rid of the factorials in Equation 
	\eqref{Formulas:2MonomialsFormula1} since
	\begin{equation*}
		\xi_1 \cdots \xi_k \star_z \eta_1 \cdots \eta_{\ell}
		=
		\frac{1}{k! \ell!}
		\frac{\partial^{k + \ell}}
		{\partial_{t_1} \cdots \partial_{s_{\ell}}}
		\At{t_1, \ldots, s_{\ell} = 0}
		\Xi^k \star_z \Eta^{\ell}.
	\end{equation*}
	Instead of the factorials, we get symmetrizations over the 
	$\xi_i$ and the $\eta_j$ as we did in Proposition 
	\ref{Formulas:Prop:LinearMonomial2}, which gives the wanted result.
\end{proof}
\begin{remark}
	This formula was in some sense already found by Corti{\~n}as in 
	\cite{cortinas:2002a}. However, his formula is somewhat less transparent, 
	since he uses an explicit form of the Baker-Campbell-Hausdorff series, 
	namely the one due to Dynkin \eqref{Alg:BCHinDynkin}. As already 
	mentioned, this formula is very complicated from a combinatorial point of 
	view and hence also less adapted to put continuity estimates on it, as we 
	will have to do. This is why we derived another formula for it.
\end{remark}

\section{Consequences and Examples}
\label{sec:chap4_Consequences}

\subsection{Some Consequences}
Proposition \ref{Formulas:Prop:2MonomialsFormula2} allows us to get some 
algebraic results. For example, we would like to see that the Gutt star 
product fulfils the classical and the semi-classical limit from Definition 
\ref{Def:StarProduct}. We can prove this using Proposition 
\ref{Formulas:Prop:LinearMonomial2}. This will finish the proof of 
Theorem~\ref{Alg:Thm:ThreeStarsAreOne}.
\begin{corollary}
	\label{Formulas:Cor:LimitCases}
	Let $\lie{g}$ be a Lie algebra and $\Sym^{\bullet}(\lie{g})$ endowed with 
	the Gutt star product
	\begin{equation*}
		x \star_z y
		= 
		\sum\limits_{n = 0}^{\infty}
		z^n C_n(x,y).
	\end{equation*}	
	\begin{corollarylist}
		\item
		On factorizing tensors $\xi_1 \ldots \xi_k$ and $\eta_1 \ldots 
		\eta_{\ell}$, $C_0$ and $C_1$ give
		\begin{align}
			\label{Formulas:ClassicalLimit}
			C_0 \left(
				\xi_1 \cdots \xi_k, \eta_1 \cdots \eta_{\ell}
			\right)
			& = 
			\xi_1 \cdots \xi_k \eta_1 \cdots \eta_{\ell}
			\\
			\label{Formulas:SemiClassicalLimit}
			C_1 \left(
				\xi_1 \cdots \xi_k, \eta_1 \cdots \eta_{\ell}
			\right)
				& =
			\frac{1}{2}	
			\sum\limits_{i = 1}^k
			\sum\limits_{j = 1}^{\ell}
			\xi_1 \cdots \widehat{\xi_i} \cdots \xi_k
			\eta_1 \cdots \widehat{\eta_j} \cdots \eta_{\ell}
			[\xi_i \eta_j],
		\end{align}
		where the hat denotes elements which are left out.
		
		\item
		For $\lie{g}$ finite-dimensional and the canonical isomorphism 
		$\mathcal{J} \colon \Sym^{\bullet}(\lie{g}) \longrightarrow 
		\Pol^{\bullet}(\lie{g}^*)$ from Proposition~\ref{Alg:Prop:PolIsSym},
		we have for $f,g \in \Pol^{\bullet}(\lie{g}^*)$
		\begin{equation*}
			C_1 \left(
				\mathcal{J}^{-1} (f),
				\mathcal{J}^{-1} (g)
			\right)
			-
			C_1 \left(
				\mathcal{J}^{-1} (f),
				\mathcal{J}^{-1} (g)
			\right)
			= 
			\mathcal{J}^{-1} \left( 
				\{ f, g \}_{KKS}
			\right)
		\end{equation*}
		where $\{ \cdot, \cdot \}_{KKS}$ is the Kirillov-Kostant-Souriau 
		bracket.
		
		\item
		The map $\star_z$ fulfils the classical and the semi-classical limit 
		and is therefore a star product.
	\end{corollarylist}
\end{corollary}
\begin{proof}
	We take $\xi_1 \cdots \xi_k, \eta_1 \cdots \eta_{\ell} \in \Sym^{\bullet}
	(\lie{g})$ and consider the 	G-indices in $\mathcal{G}_{k + \ell}(k, \ell)$ 
	first. This is easy, since there is just one element inside:
	\begin{equation*}
		\mathcal{G}_{k + \ell}(k, \ell)
		=
		\Big\{
			( 
				\underbrace{(0,1), \ldots, (0,1)}_{
				\ell \text{ times}
				}
				,
				\underbrace{(1,0), \ldots, (1,0)}_{
				k \text{ times}
				}
			)
		\Big\}.
	\end{equation*}
	So we find
	\begin{align*}
		C_0 
		\left(
			\xi_1 \cdots \xi_k, \eta_1 \cdots \eta_{\ell}
		\right)
		& =
		\sum\limits_{\sigma \in S_k}
		\sum\limits_{\tau \in S_{\ell}}
		\frac{1}{J!}
		\bchparts{0}{1}{\varnothing}{\xi_{\sigma(1)}}
		\cdots
		\bchparts{0}{1}{\varnothing}{\xi_{\sigma(k)}}
		\\
		& \quad \cdot
		\bchparts{1}{0}{\eta_{\tau(1)}}{\varnothing}
		\cdots
		\bchparts{1}{0}{\eta_{\tau(\ell)}}{\varnothing}
		\\
		& =
		\sum\limits_{\sigma \in S_k}
		\sum\limits_{\tau \in S_{\ell}}
		\frac{1}{k! \ell!}
		\xi_{\sigma(1)} \cdots \xi_{\sigma(k)}
		\eta_{\tau(1)} \cdots \eta_{\tau(\ell)}
		\\
		& =
		\xi_1 \cdots \xi_k
		\eta_1 \cdots \eta_{\ell}
	\end{align*}
	where we used $J! = k! \ell!$  according to 
	Definition~\ref{Def:GuttFactorial}. We do the same for $C_1$. Also 
	here, we have just one element in $\mathcal{G}_{k + \ell - 1}(k, \ell)$:
		\begin{equation*}
		\mathcal{G}_{k + \ell}(k, \ell)
		=
		\Big\{
			( 
				\underbrace{(0,1), \ldots, (0,1)}_{
				\ell-1 \text{ times}
				}
				,
				\underbrace{(1,0), \ldots, (1,0)}_{
				k-1 \text{ times}
				}
				,
				(1,1)
			)
		\Big\}.
	\end{equation*}
	Using
	\begin{equation*}
		\bchparts{1}{1}{\xi}{\eta}
		=
		\frac{1}{2} [\xi, \eta]
	\end{equation*}
	and $J! = (k - 1)! (\ell - 1)!$, we find
	\begin{align*}
		C_1 
		\left(
			\xi_1 \cdots \xi_k, \eta_1 \cdots \eta_{\ell}
		\right)
		& =
		\frac{1}{2}
		\sum\limits_{\sigma \in S_k}
		\sum\limits_{\tau \in S_{\ell}}
		\frac{1}{(k-1)! (\ell-1)!}
		\xi_{\sigma(1)} \cdots \xi_{\sigma(k-1)}
		\eta_{\tau(1)} \cdots \eta_{\tau(\ell-1)}
		[\xi_{\sigma(k)}, \eta_{\tau(\ell)}]
		\\
		& =
		\frac{1}{2}
		\sum\limits_{i = 0}^k
		\sum\limits_{j = 0}^{\ell}
		\xi_1 \cdots \widehat{\xi_i} \cdots \xi_k
		\eta_1 \cdots \widehat{\eta_j} \cdots \eta_{\ell}
		[\xi_i, \eta_j].
	\end{align*}
	This finishes part one. From this, the anti-symmetry of the Lie bracket 
	yields
	\begin{equation*}
		C_1 
		\left(
			\xi_1 \cdots \xi_k, \eta_1 \cdots \eta_{\ell}
		\right)
		-
		C_1 
		\left(
			\eta_1 \cdots \eta_{\ell}, \xi_1 \cdots \xi_k
		\right)
		=
		\sum\limits_{i = 0}^k
		\sum\limits_{j = 0}^{\ell}
		\xi_1 \cdots \widehat{\xi_i} \cdots \xi_k
		\eta_1 \cdots \widehat{\eta_j} \cdots \eta_{\ell}
		[\xi_i, \eta_j].
	\end{equation*}
	We now need to compute the KKS brackets on polynomials. Because of the 
	linearity in both arguments, it is sufficient to check it on monomials of 
	coordinates. Let $e_1, \ldots, e_n$ be a basis of $\lie{g}$ with linear 
	coordinates 	$x_1, \ldots, x_n$ on $\lie{g}^*$. Now take $\mu_1, \ldots, 
	\mu_n, \nu_1, \ldots \nu_n \in \mathbb{N}$ and consider the monomials 
	$f = x_1^{\mu_1} \cdots x_n^{\mu_n}$ and $g = x_1^{\nu_1} \cdots 
	x_n^{\nu_n}$. We use the notation from 
	Proposition~\ref{Alg:Prop:LinPoissonIsLieAlg} and find for $x \in 
	\lie{g}^*$
	\begin{align*}
		\{ f, g \}_{KKS}(x)
		& =
		x_k c_{ij}^k
		\frac{\partial f}{\partial x_i}
		\frac{\partial g}{\partial x_j}
		\\
		& =
		\mu_i \nu_j c_{ij}^k
		x_k 
		x_1^{\mu_1} \cdots x_i^{\mu_i - 1} \cdots x_n^{\mu_n}
		x_1^{\nu_1} \cdots x_j^{\nu_j - 1} \cdots x_n^{\nu_n}.
	\end{align*}
	Applying $\mathcal{J}^{-1}$ to it gives
	\begin{equation}
		\label{Formulas:SemiClassicalIntermediate}
		\mathcal{J}^{-1}
		\left(
			\{ f, g \}_{KKS}
		\right)
		=
		\sum\limits_{i=0}^n
		\sum\limits_{j=0}^n
		\mu_i \nu_j
		e_1^{\mu_1} \cdots e_i^{\mu_i - 1} \cdots e_n^{\mu_n}
		e_1^{\nu_1} \cdots e_j^{\nu_j - 1} \cdots e_n^{\nu_n}
		[e_i, e_j].
	\end{equation}
	On the other hand, we have
	\begin{equation*}
		\mathcal{J}^{-1}(f)
		=
		e_1^{\mu_1} \cdots e_n^{\mu_n}
		\quad \text{ and } \quad
		\mathcal{J}^{-1}(g)
		=
		e_1^{\nu_1} \cdots e_n^{\nu_n}.
	\end{equation*}
	Together with \eqref{Formulas:SemiClassicalLimit} this gives
	\eqref{Formulas:SemiClassicalIntermediate} and proves part two.
	Due to the bilinearity of the $C_n$, the third part follows.
\end{proof}
It is clear, that the formulas from 
Proposition~\ref{Formulas:Prop:2MonomialsFormula2} and 
Proposition~\ref{Formulas:Prop:LinearMonomial2} should coincide. However, 
we want to check it, to have the evidence that everything works as we wanted.
\begin{corollary}
	\label{Formulas:Cor:FormulasCoincide}
	Given $\xi_1, \ldots, \xi_k, \eta \in \lie{g}$, the results of the 
	Equations \eqref{Formulas:2MonomialsFormula2} and 
	\eqref{Formulas:LinearMonomial2} are compatible.
\end{corollary}
\begin{proof}
	We have to compute sets of G-indices for $\xi_1, \ldots, \xi_k, 
	\eta_{\ell} \in \lie{g}$. Again, they only have one element:
	\begin{equation*}
		\mathcal{G}_{k + 1 - n}(k, 1)
		=
		\Big\{
			( 
				\underbrace{(1,0), \ldots, (1,0)}_{
				k - n \text{ times}
				}
				,
				(n,1)
			)
		\Big\}.
	\end{equation*}
	So we have with $J! = (k - n)!$ and $\bchparts{n}{1}{\xi}{\eta} = 
	\frac{B_n^*}{n!} \left( \ad_{\xi} \right)^n (\eta)$
	\begin{align*}
		z^n C_n
		\left(
			\xi_1 \ldots \xi_k, \eta_{\ell}
		\right)
		& =
		z^n
		\sum\limits_{\sigma \in S_k}
		\frac{1}{(k - n)!}
		\frac{B_n^*}{n!}
		\xi_{\sigma(1)} \ldots \xi_{\sigma(k-n)}
		[\xi_{\sigma(k-n+1)}, [
			\ldots, [\xi_{\sigma(k), \eta} ] \ldots 
		]]
		\\
		& =
		z^n
		\frac{1}{k!}
		\sum\limits_{\sigma \in S_k}
		\binom{k}{n} B_n^*
		\xi_{\sigma(1)} \ldots \xi_{\sigma(k-n)}
		[\xi_{\sigma(k-n+1)}, [
			\ldots, [\xi_{\sigma(k), \eta} ] \ldots 
		]]		
	\end{align*}
	Summing up over all $n$ gives Equation~\eqref{Formulas:LinearMonomial2}.
\end{proof}

\subsection{Two Examples}
Equation \eqref{Formulas:2MonomialsFormula2} is useful if one wants to do 
real computations with the star product, but it is maybe not intuitive to 
apply. This is why we will give two examples here. The easiest one which is 
not covered by the simpler formula~\eqref{Formulas:LinearMonomial2} will be 
the star product of two quadratic terms. The second one should be the a bit 
more complex case of a cubic term with a quadratic term.

\subsubsection*{Two quadratic terms}
Let $\xi_1, \xi_2, \eta_1, \eta_2 \in \mathfrak{g}$. We want to compute
\begin{equation*}
	\xi_1 \xi_2 \star_z \eta_1 \eta_2
	=
	C_0(\xi_1 \xi_2, \eta_1 \eta_2) 
	+ 
	z C_1(\xi_1 \xi_2, \eta_1 \eta_2) 
	+ 
	z^2 C_2(\xi_1 \xi_2, \eta_1 \eta_2) 
	+ 
	z^3 C_3(\xi_1 \xi_2, \eta_1 \eta_2).
\end{equation*}
The very first thing we have to do is computing the set of G-indices. Then we 
calculate the G-factorial and finally go through the permutations.
\begin{itemize}
	\item[$C_0$:]
	We already did this in Corollary \ref{Formulas:Cor:LimitCases}, and know 
	that the zeroth order in $z$ is just the symmetric product. Therefore we 
	have
	\begin{equation*}
		C_0(\xi_1 \xi_2, \eta_1 \eta_2)
		=
		\xi_1 \xi_2 \eta_1 \eta_2
	\end{equation*}
	
	\item[$C_1$:]
	We also did this one in Corollary \ref{Formulas:Cor:LimitCases}: There is 
	just one G-index and we finally get
	\begin{equation*}
		C_1(\xi_1 \xi_2, \eta_1 \eta_2)
		=
		\frac{1}{2} \left(
			\xi_2 \eta_2 [\xi_1, \eta_1] +
			\xi_2 \eta_1 [\xi_1, \eta_2] +
			\xi_1 \eta_2 [\xi_2, \eta_1] +
			\xi_1 \eta_1 [\xi_2, \eta_2]
		\right).
	\end{equation*}
	
	\item[$C_2$:]
	This is the first time, something interesting happens. We have three 
	G-indices:
	\begin{equation*}
		\mathcal{G}_2(2,2) 
		=
		\left\{ J_1, J_2, J_3 \right\}
		= 
		\big\{ 
			\big((0,1), (2,1)\big), 
			\big((1,0), (1,2)\big), 
			\big((1,1), (1,1)\big) 
		\big\}.
	\end{equation*}
	The G-factorials give $J_1! = J_2! = 1$ and $J_3! = 2$, since the index 
	$(1,1)$ appears twice in $J_3$. We take $\bchparts{a}{b}{\xi}{\eta}$ 
	from Equation \eqref{Alg:BCHSeriesLong} for $(a,b) \in \{(1,2), (2,1)\}$:
	\begin{equation*}
		\bchparts{1}{2}{\xi}{\eta}
		=
		\frac{1}{12}[[\xi, \eta], \eta]
		\quad \text{ and } \quad
		\bchparts{2}{1}{\xi}{\eta}
		=
		\frac{1}{12}[[\eta, \xi], \xi].
	\end{equation*}
	So we have to insert the $\xi_i$ and the $\eta_j$ into $\frac{1}{12} 
	\xi [[\xi, \eta], \eta]$ and $\frac{1}{12} \eta [[\eta, \xi], \xi]$ 
	respectively and then we go on with the last one, which is
	\begin{equation*}
		\frac{1}{2}
		\bchparts{1}{1}{\xi}{\eta}
		\bchparts{1}{1}{\xi}{\eta}
		=
		\frac{1}{8}
		[\xi, \eta][\xi, \eta].
	\end{equation*}
	We hence get
	\begin{align*}
		C_2(\xi_1, \xi_2, \eta_1, \eta_2) 
		& = 
		\frac{1}{12}
		\big( 
			\eta_1 [[\eta_2, \xi_1],\xi_2] + 
			\eta_1 [[\eta_2, \xi_2],\xi_1] + 
			\eta_2 [[\eta_1, \xi_1],\xi_2] + 
			\eta_2 [[\eta_1, \xi_2],\xi_1] +
		\\
		& \quad
			\xi_1 [[\xi_2, \eta_1],\eta_2] + 
			\xi_1 [[\xi_2, \eta_2],\eta_1] + 
			\xi_2 [[\xi_1, \eta_1],\eta_2] + 
			\xi_2 [[\xi_1, \eta_2],\eta_1] 
 		\big) +
 		\\
 		& \quad
		\frac{1}{4} 
		\big( 
			[\xi_1,\eta_1][\xi_2,\eta_2] + 
			[\xi_1,\eta_2][\xi_2,\eta_1] 
		\big)
	\end{align*}
	
	\item[$C_3$:]
	Here, we only have one G-index:
	\begin{equation*}
		\mathcal{G}_1(2,2) 
		=
		\big\{ 
			\big( (2,2) \big) 
		\big\}
	\end{equation*}
	The G-factorial is $1$. We take again Equation \eqref{Alg:BCHSeriesLong} 
	and see
	\begin{equation*}
		\bchparts{2}{2}{\xi}{\eta}
		=
		\frac{1}{24}
		[[[\eta, \xi], \xi], \eta].
	\end{equation*}
	This gives
	\begin{align*}
		C_3(\xi_1, \xi_2, \eta_1, \eta_2) 
		& = 
		\frac{1}{24}
		\big( 
			[[[\eta_1,\xi_1],\xi_2],\eta_2] + 
			[[[\eta_1,\xi_2],\xi_1],\eta_2] +
			[[[\eta_2,\xi_1],\xi_2],\eta_1] + 
			[[[\eta_2,\xi_2],\xi_1],\eta_1] 
		\big)
	\end{align*}
\end{itemize}
We just have to put all the four terms together and have the star product.

\subsubsection*{A cubic and a quadratic term}
Let $\xi_1, \xi_2, \xi_3, \eta_1, \eta_2 \in \mathfrak{g}$. We compute
\begin{equation*}
	\xi_1 \xi_2 \xi_3 \star_G \eta_1 \eta_2
	= 
	\sum\limits_{n = 0}^4
	z^n C_n(\xi_1 \xi_2 \xi_3, \eta_1 \eta_2)
\end{equation*}
\begin{itemize}
	\item[$C_0$:]
	The first part is again just the symmetric product:
	\begin{equation*}
		C_0(\xi_1 \xi_2 \xi_3, \eta_1 \eta_2)
		=
		\xi_1 \xi_2 \xi_3 \eta_1 \eta_2.
	\end{equation*}
	
	\item[$C_1$:]
	Here we have again the term from Corollary \ref{Formulas:Cor:LimitCases}:
	\begin{align*}
		C_1(\xi_1 \xi_2 \xi_3, \eta_1 \eta_2)
		& =
		\frac{1}{2} 
		\big( 
			\xi_2 \xi_3 \eta_2 [\xi_1, \eta_1] + 
			\xi_2 \xi_3 \eta_1 [\xi_1, \eta_2] + 
			\xi_1 \xi_3 \eta_2 [\xi_2, \eta_1] +
		\\
		& \quad 
			\xi_1 \xi_3 \eta_1 [\xi_2, \eta_2] + 
			\xi_1 \xi_2 \eta_2 [\xi_3, \eta_1] + 
			\xi_1 \xi_2 \eta_1 [\xi_3, \eta_2] 
 		\big)
	\end{align*}

	\item[$C_2$:]
	Here the calculation is very similar to the one of $C_2$ in the example 
	before. We have three G-indices:
	\begin{equation*}
		\mathcal{G}_3(3,2) 
		=
		\left\{
			J_1, J_2, J_3
		\right\}
		= 
		\big\{ 
			\big( (0,1), (1,0), (2,1) \big), 
			\big( (1,0), (1,0), (1,2) \big), 
			\big( (1,0), (1,1), (1,1) \big) 
		\big\}.
	\end{equation*}
	The G-factorials are now $J_1! = 1$ and $J_2! = J_3! = 2$. Again, we take 
	the BCH terms from Equation \eqref{Alg:BCHSeriesLong} and see, that we 
	must insert the $\xi_i$ and the $\eta_j$ into
	\begin{equation*}
		\frac{1}{12} \xi \eta [[\eta, \xi], \xi] +
		\frac{1}{24} \xi \xi [[\xi, \eta], \eta] +
		\frac{1}{8} \xi [\xi, \eta] [\xi, \eta].
	\end{equation*}
	Now we go through all the possible permutations and get
	\begin{align*}
		C_2(\xi_1 \xi_2 \xi_3, \eta_1 \eta_2) 
		& = 
		\frac{1}{12} 
		\big( 
			\xi_1 \xi_2 [[\xi_3, \eta_1], \eta_2] + 
			\xi_1 \xi_2 [[\xi_3, \eta_2], \eta_1] + 
			\xi_1 \xi_3 [[\xi_2, \eta_1], \eta_2] +
		\\
		& \quad 
			\xi_1 \xi_3 [[\xi_2, \eta_2], \eta_1] + 
			\xi_2 \xi_3 [[\xi_1, \eta_1], \eta_2] + 
			\xi_2 \xi_3 [[\xi_1, \eta_2], \eta_1] 
		\big) +
		\\ 
		& \quad
		\frac{1}{12} 
		\big( 
			\xi_1 \eta_1 [[\eta_2, \xi_2], \xi_3] + 
			\xi_1 \eta_2 [[\eta_1, \xi_2], \xi_3] + 
			\xi_1 \eta_1 [[\eta_2, \xi_3], \xi_2] +
		\\
		& \quad
			\xi_1 \eta_2 [[\eta_1, \xi_3], \xi_2] + 
			\xi_2 \eta_1 [[\eta_2, \xi_1], \xi_3] + 
			\xi_2 \eta_2 [[\eta_1, \xi_1], \xi_3] +
		\\
		& \quad
			\xi_2 \eta_1 [[\eta_2, \xi_3], \xi_1] + 
			\xi_2 \eta_2 [[\eta_1, \xi_3], \xi_1] + 
			\xi_3 \eta_1 [[\eta_2, \xi_2], \xi_1] +
		\\
		& \quad
			\xi_3 \eta_2 [[\eta_1, \xi_2], \xi_1] + 
			\xi_3 \eta_1 [[\eta_2, \xi_1], \xi_2] + 
			\xi_3 \eta_2 [[\eta_1, \xi_1], \xi_2] 
		\big) +
		\\
		& \quad
		\frac{1}{4} 
		\big( 
			\xi_1 [\xi_2, \eta_1] [\xi_3, \eta_2] + 
			\xi_1 [\xi_3, \eta_1] [\xi_2, \eta_2] + 
			\xi_2 [\xi_1, \eta_1] [\xi_3, \eta_2] + 
		\\
		& \quad
			\xi_2 [\xi_3, \eta_1] [\xi_1, \eta_2] + 
			\xi_3 [\xi_1, \eta_1] [\xi_2, \eta_2] + 
			\xi_3 [\xi_2, \eta_1] [\xi_1, \eta_2] 
 		\big).
	\end{align*}

	\item[$C_3$:]
	We first calculate the G-indices:
	\begin{equation*}
		\mathcal{G}_2(3,2) 
		= 
		\left\{
			J_1, J_2, J_3
		\right\}
		=
		\left\{ 
			\big( (0,1), (3,1) \big), 
			\big( (1,0), (2,2) \big), 
			\big( (1,1), (2,1) \big) 
		\right\}.
	\end{equation*}
	We can omit $J_1$, since $\bchparts{3}{1}{\xi}{\eta} = 0$. The 
	G-factorials for the other two indices are $1$. The BCH terms have been 
	computed before. So we have to fill in the expression
	\begin{equation*}
		\frac{1}{24} \xi [[[\eta, \xi], \xi], \eta] +
		\frac{1}{2 \cdot 12} [\xi, \eta] [[\eta, \xi], \xi].
	\end{equation*}
	Going through the permutations we get
	\begin{align*}
		C_3(\xi_1 \xi_2 \xi_3, \eta_1 \eta_2)
		& =
		\frac{1}{24}
		\big( 
			\xi_1[[[\eta_1, \xi_2], \xi_3], \eta_2] + 
			\xi_1[[[\eta_2, \xi_2], \xi_3], \eta_1] + 
			\xi_1[[[\eta_1, \xi_3], \xi_2], \eta_2] +
		\\
		& \quad
			\xi_1[[[\eta_2, \xi_3], \xi_2], \eta_1] + 
			\xi_2[[[\eta_1, \xi_1], \xi_3], \eta_2] + 
			\xi_2[[[\eta_2, \xi_1], \xi_3], \eta_1] +
		\\
		& \quad
			\xi_2[[[\eta_1, \xi_3], \xi_1], \eta_2] + 
			\xi_2[[[\eta_2, \xi_3], \xi_1], \eta_1] + 
			\xi_3[[[\eta_1, \xi_2], \xi_1], \eta_2] +
		\\
		& \quad
			\xi_3[[[\eta_2, \xi_2], \xi_1], \eta_1] + 
			\xi_3[[[\eta_1, \xi_1], \xi_2], \eta_2] + 
			\xi_3[[[\eta_2, \xi_1], \xi_2], \eta_1]
		\big) +
		\\
		& \quad
		\frac{1}{24}
		\big( 
			[\xi_1, \eta_1][[\eta_2, \xi_2], \xi_3] + 
			[\xi_1, \eta_2][[\eta_1, \xi_2], \xi_3] + 
			[\xi_1, \eta_1][[\eta_2, \xi_3], \xi_2] +
		\\
		& \quad
			[\xi_1, \eta_2][[\eta_1, \xi_3], \xi_2] + 
			[\xi_2, \eta_1][[\eta_2, \xi_1], \xi_3] + 
			[\xi_2, \eta_2][[\eta_1, \xi_1], \xi_3] +
		\\
		& \quad
			[\xi_2, \eta_1][[\eta_2, \xi_3], \xi_1] + 
			[\xi_2, \eta_2][[\eta_1, \xi_3], \xi_1] + 
			[\xi_3, \eta_1][[\eta_2, \xi_2], \xi_1] +
		\\
		& \quad
			[\xi_3, \eta_2][[\eta_1, \xi_2], \xi_1] + 
			[\xi_3, \eta_1][[\eta_2, \xi_1], \xi_2] + 
			[\xi_3, \eta_2][[\eta_1, \xi_1], \xi_2] 
		\big).
	\end{align*}

	\item[$C_4$:]
	Now there is only $C_4$ left. We have one G-index:
	\begin{equation*}
		\mathcal{G}_1(3,2) 
		= 
		\left\{ \big( (3,2) \big) \right\},
	\end{equation*}
	but there are more terms which belong to it. We have to go through
	\begin{equation*}
		\bchparts{3}{2}{\xi}{\eta}
		=
		\frac{1}{120} [[[[\xi, \eta], \xi], \eta], \xi] +
		\frac{1}{360} [[[[\eta, \xi], \xi], \xi], \eta].
	\end{equation*}
	So we permute and get
	\begin{align*}
		C_4(\xi_1 \xi_2 \xi_3, \eta_1 \eta_2)
		& =
		\frac{1}{120} 
		\big( 
			[[[[\xi_1, \eta_1], \xi_2], \eta_2], \xi_3] + 
			[[[[\xi_1, \eta_2], \xi_2], \eta_1], \xi_3] + 
			[[[[\xi_1, \eta_1], \xi_3], \eta_2], \xi_2] +
		\\
		& \quad
			[[[[\xi_1, \eta_2], \xi_3], \eta_1], \xi_2] + 
			[[[[\xi_2, \eta_1], \xi_1], \eta_2], \xi_3] + 
			[[[[\xi_2, \eta_2], \xi_1], \eta_1], \xi_3] +
		\\
		& \quad
			[[[[\xi_2, \eta_1], \xi_3], \eta_2], \xi_1] + 
			[[[[\xi_2, \eta_2], \xi_3], \eta_1], \xi_1] + 
			[[[[\xi_3, \eta_1], \xi_2], \eta_2], \xi_1] +
		\\
		& \quad
			[[[[\xi_3, \eta_2], \xi_2], \eta_1], \xi_1] + 
			[[[[\xi_3, \eta_1], \xi_1], \eta_2], \xi_2] + 
			[[[[\xi_3, \eta_2], \xi_1], \eta_1], \xi_2] 
		\big) +
		\\
		& \quad
		\frac{1}{360} 
		\big( 
			[[[[\eta_1, \xi_1], \xi_2], \xi_3], \eta_2] + 
			[[[[\eta_2, \xi_1], \xi_2], \xi_3], \eta_1] + 
			[[[[\eta_1, \xi_1], \xi_3], \xi_2], \eta_2] +
		\\
		& \quad
			[[[[\eta_2, \xi_1], \xi_3], \xi_2], \eta_1] + 
			[[[[\eta_1, \xi_2], \xi_1], \xi_3], \eta_2] + 
			[[[[\eta_2, \xi_2], \xi_1], \xi_3], \eta_1] +
		\\
		& \quad
			[[[[\eta_1, \xi_2], \xi_3], \xi_1], \eta_2] + 
			[[[[\eta_2, \xi_2], \xi_3], \xi_1], \eta_1] + 
			[[[[\eta_1, \xi_3], \xi_2], \xi_1], \eta_2] +
		\\
		& \quad
			[[[[\eta_2, \xi_3], \xi_2], \xi_1], \eta_1] + 
			[[[[\eta_1, \xi_3], \xi_1], \xi_2], \eta_2] + 
			[[[[\eta_2, \xi_3], \xi_1], \xi_2], \eta_1] 
		\big).
	\end{align*}
\end{itemize}
Now we only have to add up all those terms and we have finally computed the 
star product.

% Chapter 5
%

%
% Chapter 5 of my master thesis:
% The analysis part begins
%

\chapter{A Locally Convex Topology for the Gutt Star Product}

We have finished the algebraic part of this work, except for one little 
lemma concerning the Hopf theoretic chapter. Our next goal is setting up a 
locally convex topology on the symmetric tensor algebra, in which the Gutt 
star product will be continuous. At the beginning of this chapter, we will 
first give a motivation why the setting of locally convex algebras is 
convenient and necessary. In the second section, we will briefly recall the 
most important things on locally convex algebras and introduce the topology 
which we will work with. In the third section, the core of this chapter, the 
continuity of the star product and the analytic dependence on the formal 
parameter are proven. We also show that continuous representations of Lie 
algebras lift to  those of the deformed symmetric algebra and that our 
construction is in fact functorial. Part four treats the case when the formal 
parameter $z = 1$ and hence talks about a locally convex topology on universal 
enveloping algebras. We will show, that our topology is ``optimal'' in a 
specific sense.

\section{Why Locally Convex?}
\label{sec:chap5_Prelim}

There are very different types of topologies on vector spaces and algebras. On one 
hand, we want a topology on the symmetric tensor algebra -- which will 
after its completion be the algebra of observables -- that has as many 
convenient features as possible, of course. On the other hand, if we choose a too 
good topology, we will risk to have no representations of our algebra as 
(unbounded) operators on a Hilbert space, which means something like the $\hat 
q$ and $\hat p$ operators. Thus we have to think of the possible options first. 
Assume that $\lie{g}$ is a complex Lie algebra (if not, we can always complexify 
it), then the complex conjugation plays the role of a \emph{star involution} on 
the algebra of observables $\algebra{A}$, that means it is an automorphism
\begin{equation*}
	^* \colon
	\algebra{A}
	\longrightarrow
	\algebra{A}
	, \quad
	a 
	\longmapsto 
	a^*,
\end{equation*}
which fulfils for all $a \in \algebra{A}$ the identity $(a^*)^* = a$. We want to 
summarize briefly the options, which we would have in this case.
\begin{enumerate}
	\item
	We can try a $C^*$-algebra. This means that we have a norm $\norm{\cdot}$ 
	on our algebra, which fulfils the two properties
	\begin{equation}
		\label{LCAna:BanachAlg}
		\norm{a \cdot b}
		\leq
		\norm{a}
		\norm{b}
		, \quad
		\forall_{a,b \in \algebra{A}}
	\end{equation}
	and
	\begin{equation}
		\label{LCAna:CStarAlg}
		\norm{a^* \cdot a}
		=
		\norm{a}^2
		, \quad
		\forall_{a,b \in \algebra{A}}.
	\end{equation}
	The bounded operators on a Hilbert space form a $C^*$-algebra, for example.	
	In particular, we can define continuous functions of algebra elements in this
	case, see e.g. \cite{blackadar:2006a} or \cite{waldmann:2014a:script}.
		
	\item
	A Banach-$^*$-algebra is a bit less: we still have a norm, which fulfils
	\eqref{LCAna:BanachAlg}, but not \eqref{LCAna:CStarAlg} any more. This stil 
	allows holomorphic functions of algebra elements (like it can be seen in 
	\cite{rudin:1991a}).

	\item
	A weaker choice is a \emph{locally multiplicatively convex 
	$^*$-algebra}. We do not have a norm any more in this case, but a set of 
	submultiplicative \emph{seminorms} $\{p_i\}_{i \in I}$. A seminorm $p$ is 
	almost the same as a norm, except that we may have $p(x) = 0$ for some $x \in 
	\algebra{A}$ with $x \neq 0$. Submultiplicativity means that every seminorm
	$p$ fulfils
	\begin{equation}
		\label{LCAna:LMCAlg}
		p(a \cdot b)
		\leq
		p(a)
		p(b)
		, \quad
		\forall_{a,b \in \algebra{A}}.
	\end{equation}
	Examples for such algebras are the entire functions on $\mathbb{C}$ or the 
	continuous functions on $\mathbb{R}$. In this case, entire functions of 
	algebra elements can still be defined (see e.g. \cite{michael:1952a}).
	
	\item
	As a last possibility, we can have a \emph{locally convex $^*$-algebra}. 
	Here again, we have seminorms, but without \eqref{LCAna:LMCAlg}. Instead, we 
	have to use another seminorm $q$ to control $p$ of the product of two 
	elements. We will make this statement more precise later. Examples for such 
	structures are distributions or vector fields on vector space. Good books on 
	this kind of structure are \cite{hoermander:2003a}, \cite{rudin:1991a} or 
	\cite{jarchow:1981a}.
\end{enumerate}
So in principle, we have a lot of possibilities. Unfortunately, the canonical 
commutation relations, which should be fulfilled, do not allow 
submultiplicative structures: this is precisely what we show in the next 
proposition (and what can also be found in some textbooks, like 
\cite{reed.simon:1980a}, for example). So all the nice features like a continuous 
or an entire calculus will simply not be there, and therefore we have to aim at an 
honestly locally convex topology on the symmetric tensor algebra.
\begin{proposition}
	\label{Prop:LCAna:QMnotLMC}
	Let $\algebra{A}$ be a unital associative algebra which contains the 
	quantum mechanical observables $\hat{q}$ and $\hat{p}$ and in which 
	the canonical commutation relation
	\begin{equation*}
		[\hat{q}, \hat{p}]
		=
		i \hbar \Unit
	\end{equation*}
	is fulfilled. Then the only submultiplicative seminorm on it is 0.
\end{proposition}
\begin{proof}
	First, we need to show a little lemma:
	\begin{lemma}
		\label{Lemma:LCAna:NotLMCHelp}
		In the algebra given above, we have for $n \in \mathbb{N}$
		\begin{equation}
			\label{LCAna:NotLMCHelp}
			\left( \ad_{\hat{q}} \right)^n (\hat{p}^n)
			=
			(i \hbar)^n n! \Unit.
		\end{equation}
	\end{lemma}
	\begin{subproof}
		To show it, we use the fact that for $a \in \algebra{A}$ the 
		operator $\ad_a$ is a derivation. This is always true for the 
		commutator Lie algebra of an associative algebra, since for 
		$a, b, c \in \algebra{A}$ we have
		\begin{equation*}
			[a, bc]
			=
			a b c - b c a
			=
			a b c - b a c + b a c - b c a
			=
			[a, b] c + b [a, c].
		\end{equation*}
		For $n = 1$, Equation \eqref{LCAna:NotLMCHelp} is 
		certainly true. So let's look at the step $n \rightarrow n+1$.
		We make use of the derivation property and have
		\begin{align*}
			\left( \ad_{\hat{q}} \right)^{n+1}
			\left( \hat{p}^{n+1} \right)
			& =
			\left( \ad_{\hat{q}} \right)^{n}
			\left(
				i \hbar \hat{p}^n
				+
				\hat{p} 
				\ad_{\hat{q}} \left( \hat{p}^n \right)
			\right)
			\\
			& =
			(i \hbar)^{n + 1} n!
			+
			\left( \ad_{\hat{q}} \right)^{n}
			\left(
				\hat{p}
				\ad_{\hat{q}} \left( \hat{p}^n \right)
			\right)
			\\
			& =
			(i \hbar)^{n + 1} n!
			+
			\left( \ad_{\hat{q}} \right)^{n-1}
			\left(
				[\hat{q}, \hat{p}]
				\ad_{\hat{q}} \left( \hat{p}^n \right)
				+
				\hat{p}
				\left( \ad_{\hat{q}} \right)^2
				\left( \hat{p}^n \right)
			\right)
			\\
			& =
			(i \hbar)^{n + 1} n!
			+
			i \hbar
			\left( \ad_{\hat{q}} \right)^{n}
			\left( \hat{p}^n \right)
			+
			\left( \ad_{\hat{q}} \right)^{n-1}
			\left(
				\hat{p}
				\left( \ad_{\hat{q}} \right)^2
				\left( \hat{p}^n \right)
			\right)
			\\
			& =
			2 (i \hbar)^{n+1} n!
			+
			\left( \ad_{\hat{q}} \right)^{n-1}
			\left(
				\hat{p}
				\left( \ad_{\hat{q}} \right)^2
				\left( \hat{p}^n \right)
			\right)
			\\
			& \ot{($*$)}{=} 
			\quad \vdots
			\\
			& =
			n (i \hbar)^{n+1} n!
			+
			\ad_{\hat{q}}
			\left(
				\hat{p}
				\left( \ad_{\hat{q}} \right)^n
				\left( \hat{p}^n \right)
			\right)
			\\
			& =
			n (i \hbar)^{n+1} n!
			+
			i \hbar (i \hbar)^n n!
			\\
			& =
			(i \hbar)^{n + 1} (n + 1)!.
		\end{align*}
		At ($*$), we actually used another statement which is to be 
		proven by induction over $k$ and says
		\begin{equation*}
			\left( \ad_{\hat{q}} \right)^{n + 1}
			\left( \hat{p}^{n + 1} \right)
			=
			k (i \hbar)^{n + 1} n!
			+
			\left( \ad_{\hat{q}} \right)^{n + 1 - k}
			\left(
				\hat{p}
				\left( \ad_{\hat{q}} \right)^k
				\left( \hat{p}^n \right)
			\right).
		\end{equation*}
		Since this proof is analogous to the first lines of the 
		computation before, we omit it here and the lemma is proven.
	\end{subproof}	
	Now we can go on with the actual proof. Let $\norm{\cdot}$ be a 
	submultiplicative seminorm. Then we see from Equation 
	\eqref{LCAna:NotLMCHelp} that
	\begin{equation*}
		\norm{
			\left( \ad_{\hat{q}} \right)^n
			(\hat{p}^n)
		}
		=
		|\hbar|^n n! \norm{ \Unit }.
	\end{equation*}
	On the other hand, we have
	\begin{align*}
		\norm{
			\left( \ad_{\hat{q}} \right)^n
			(\hat{p}^n)
		}
		& =
		\norm{
			\hat{q}
			\left( \ad_{\hat{q}} \right)^{n-1}
			(\hat{p}^n)
			-
			\left( \ad_{\hat{q}} \right)^{n-1}
			(\hat{p}^n)
			\hat{q}
		}
		\\
		& \leq
		2 \norm{\hat{q}}
		\norm{
			\left( \ad_{\hat{q}} \right)^{n-1}
			(\hat{p}^n)
		}
		\\
		& \leq
		\quad \vdots
		\\
		& \leq
		2^n \norm{\hat{q}}^n
		\norm{ \hat{p}^n }
		\\
		& \leq
		2^n 
		\norm{\hat{q}}^n
		\norm{\hat{p}}^n
	\end{align*}
	So in the end we get
	\begin{equation*}
		|\hbar|^n n! \norm{ \Unit }
		\leq
		c^n
	\end{equation*}
	for $1 \norm{\hat q} \norm{\hat p} = c \in \mathbb{R}$. This cannot be 
	fulfilled for all $n \in \mathbb{N}$ unless $\norm{ \Unit } = 0$. But 
	then, by submultiplicativity, the seminorm itself must be equal to $0$.
\end{proof}
\begin{remark}
	The so called Weyl algebra, which fulfils the properties of the 
	foregoing proposition, can be constructed from a Poisson algebra with 
	constant a Poisson tensor. On one hand, it is a fair to ask the question, 
	why this restriction of not being locally m-convex should also be apply to
	linear Poisson systems. On the other hand, there is no reason to 
	expect that things become easier when we make the Poisson structure more 
	complex. Moreover, the Weyl algebra is actually nothing but a 
	quotient of the universal enveloping algebra of the so called 
	Heisenberg algebra, which is a particular Lie algebra. There is 
	no reason why the original algebra should have a ``better'' analytical 
	structure than its quotient, since the ideal, which is divided out by 
	this procedure, is a closed one.
\end{remark}
There is a second good reason why we should avoid our topology to be locally 
m-convex. The topology we set up on $\Sym^{\bullet}(\lie{g})$ for a Lie 
algebra $\lie{g}$ will also give a topology on $\algebra{U}(\lie{g})$.
In Proposition \ref{LCAna:Prop:NoBetterTopology}, we will show that, under 
weak (but for our purpose necessary) additional assumptions, there can be no 
topology on $\algebra{U}(\lie{g})$, which allows an entire holomorphic 
calculus. This underlines the results from Proposition 
\ref{Prop:LCAna:QMnotLMC}, since locally m-convex algebras always have such a 
calculus.

In this sense, we have good reasons to think that 
$\Sym^{\bullet}(\lie{g})$ will not allow a better setting than the
one of a locally convex algebra if we want the Gutt star product to 
be continuous. Before we go to this task, we have to recall some
technology from locally convex analysis.

\section{Locally Convex Algebras}
\label{sec:chap5_LCAlg}

\subsection{Locally Convex Spaces and Algebras}

Every locally convex algebra is also a locally convex space which 
is, of course, a topological vector space. To make clear what we talk about, 
we first give a definition which is taken from 
\cite[Definition 1.6]{rudin:1991a}.
\begin{definition}[Topological vector space]
	\label{Def:TVSpace}
	Let $V$ be a vector space endowed with a topology $\tau$. Then we call 
	$(V, \tau)$ (or just $V$, if there is no confusion possible) 
	a topological vector space, if the two following things  hold:
	\begin{definitionlist}
		\item
		for every point in $x \in V$ the set $\{x\} $ is closed and
		
		\item
		the vector space operations (addition, scalar multiplication) are 
		continuous.
	\end{definitionlist}
\end{definition}
Not all books require axiom $(i)$ for a topological vector space. It is, 
however, useful, since it assures that our topological vector 
space is Hausdorff -- a feature which we always want to have. The proof 
for this is not difficult, but since we do not want to go too much into details 
here, we refer to \cite{rudin:1991a} again, where it can be found in 
Theorem 1.10 and 1.12.

The most important class of topological vector spaces are, at least from a 
physical point of view, locally convex ones. Almost all 
interesting physical examples belong to this class: Finite-dimensional spaces, 
inner product (or pre-Hilbert) spaces, Banach spaces, Fr\'echet spaces, 
nuclear spaces and many more. There are at least two equivalent 
definitions of what is a locally convex space. While the first is more 
geometrical, the second is better suited for our analytic purpose.
\begin{theorem}
	\label{Thm:LCAna:LCSpace}
	For a topological vector space $V$, the following things are equivalent.
	\begin{theoremlist}
		\item
		$V$ has a local base $\algebra{B}$ of the topology whose members are 
		convex, i.e. for all $U in \algebra{B}$ we have for all $x,y \in U$
		and all $\lambda \in [0,1]$: $\lambda x + (1 - \lambda)y \in U$.
		\item
		The topology on $V$ is generated by a separating family of seminorms 
		$\mathcal{P}$.
	\end{theoremlist}
\end{theorem}
\begin{proof}
	This theorem is a very well-known result and can be found in standard 
	literature, such as \cite{rudin:1991a} again, where it is divided into 
	two 	Theorems (namely 1.36 and 1.37).
\end{proof}
\begin{definition}[Locally convex space]
	\label{Def:LCSpace}
	A locally convex space is a topological vector space in which one (and 
	thus all) of the properties from Theorem \ref{Thm:LCAna:LCSpace} are 
	fulfilled.
\end{definition}
The first property explains the term ``locally convex''. 
For our concrete computations, the second property is more helpful, 
since in this setting proving continuity just means putting estimates on 
seminorms. For this purpose, one often extends the set of seminorms 
$\mathcal{P}$ to the set of all continuous seminorms $\algebra{P}$ which 
contains all seminorms that are compatible with the topology:
\begin{equation*}
	p \in \algebra{P}
	\Longleftrightarrow
	\exists_{n \in \mathbb{N}}
	\exists_{c_1, \ldots, c_n > 0}
	\exists_{p_1, \ldots, p_n \in \mathcal{P}}
	\colon
	p 
	\leq 
	\sum\limits_{i = 1}^n
	c_i p_i.
\end{equation*}
In particular, one can take maxima of finitely many continuous seminorms to get 
again a continuous seminorm. From here, we can start looking at locally convex 
algebras.
\begin{definition}[Locally convex algebra]
	\label{Def:LCAlgebra}
	A locally convex algebra is a locally convex vector space with an 
	additional algebra structure which is continuous.
\end{definition}
More precisely, let $\algebra{A}$ be a locally convex algebra and 
$\algebra{P}$ the set of all continuous seminorms, then for all $p \in 
\algebra{P}$ there exists a $q \in \algebra{P}$ such that for all $x, y \in 
\algebra{A}$ one has
\begin{equation}
	\label{LCAna:ProductContinuity}
	p(a \cdot b)
	\leq
	q(a) q(b).
\end{equation}
Remind that we did not require our algebras to be associative. The product in 
this equation could also be a Lie bracket. If we talk about associative 
algebras, we will always say it explicitly.

\subsection{A Special Class of Locally Convex Algebras}

For our study of the Gutt star product, the usual continuity estimate 
\eqref{LCAna:ProductContinuity} will not be enough, since 
there will be an arbitrarily high number of nested brackets to control. 
We will need an estimate which does not depend on the number of Lie 
brackets involved. Since Lie algebras are just one type of algebras, 
we can define the property we need also for other locally convex algebras. 
\begin{definition}[Asymptotic estimate algebra]
    \label{Def:AE}
    Let $\algebra{A}$ be a locally convex algebra (not necessarily associative) 
    with $\cdot$ denoting the multiplication and $\algebra{P}$ the set 
    of all continuous seminorms. Let $p \in \algebra{P}$.
    \begin{definitionlist}
    \item \label{item:AsymptoticEstimate} 
        We call $q \in \algebra{P}$ an asymptotic estimate for $p$, if
        \begin{equation}
            \label{LCAna:AEProperty}
            p\left(
                w_n(x_1, \ldots, x_n)
            \right)
            \leq
            q(x_1) \cdots q(x_n)
        \end{equation}
        for all words $w_n(x_1, \ldots, x_n)$ made out of $n-1$
        products of the elements $x_1, \ldots, x_n \in \algebra{A}$
        with arbitrary position of placing brackets.
    \item \label{item:AEAlgebra} 
        A locally convex algebra is said to be an asymptotic estimate algebra 
        (AE algebra), if every $p \in \algebra{P}$ has an asymptotic estimate.
    \end{definitionlist}
\end{definition}
\begin{remark}[The notion ``asymptotic estimate'']
	\label{LCAna:Rem:AE1}
	\mbox{}
	\begin{remarklist}
		\item
		The term asymptotic estimate has, to the best of our knowledge, 
		first been used by Boseck, Czichowski and Rudolph in 
		\cite{boseck.czichowski.rudolph:1981a}. They gave a seemingly weaker
		definition of asymptotic estimates, which is in fact equivalent.
		They wanted that for every seminorm $p\in \algebra{P}$, there is a 
		$m \in \mathbb{N}$ and a $q \in \algebra{P}$ such that for all 
		$n \geq m$ one has
		\begin{equation}
			\label{LCAna:PseudiAE}
			p \left( x_1 \cdots x_n \right)
			\leq
			q \left(x_1\right) \cdots q \left(x_n\right)
			\quad
			\forall_{x_1, \ldots, x_n \in \algebra{A}}.
		\end{equation}
		But here we can set $m = 1$, since we just need to take the maximum of 
		a finite number of continuous seminorms. Let $q$ satisfy 
		\eqref{LCAna:PseudiAE}. By Continuity, we have for all 
		$i = 2, \ldots, m-1$
		\begin{equation*}
		 	p(x_1 \cdots x_i)
		 	\leq
		 	q^{(i)}(x_1) \cdots q^{(i)}(x_i)
		 	\quad
		 	\forall_{x_1, \ldots, x_i \in \algebra{A}}
		\end{equation*}
		Now we get an asymptotic estimate $q'$ for $p$ in our sense by setting
		\begin{equation*}
		 	q'
		 	=
		 	\max\{ 
		 		p, q^{(2)}, \ldots, q^{(m-1)}, q
		 	\}.
		\end{equation*}
		Also AE algebras were defined by Boseck et al. in 
		\cite{boseck.czichowski.rudolph:1981a}, but here the notion is really 
		different from ours: for them, in an AE algebra every continuous 
		seminorm admits a sequence of asymptotic estimates. This sequence must 
		fulfil two additional properties, which actually make the algebra 
		locally 	m-convex. Our definition is weaker,	since it does not imply, 
		a priori, the existence of an topologically equivalent set of 
		submultiplicative seminorms.
		
		\item
		In \cite{gloeckner.neeb:2012a}, Gl\"ockner and Neeb used a property to 
		which they referred as $(*)$ for associative algebras. It was then 
		used in 	\cite{bogfjellmo.dahmen.schmedig:2015a} by Bogfjellmo, Dahmen 
		and Schmedig, who called it the $GN$-property. It is rather easy to 
		see that it is equivalent to our AE condition. 
	\end{remarklist}
\end{remark}

There are, of course, a lot of examples of AE (Lie) algebras. 
All finite dimensional and Banach (Lie) algebras fulfil Inequality
\eqref{LCAna:AEProperty}, just as locally m-convex (Lie) algebras do. 
The same is true for nilpotent locally convex Lie algebras, 
since here again one just has to take the maximum of a finite 
number of semi-norms, analogously to the procedure in Remark 
\ref{LCAna:Rem:AE1}. On the other hand, it is not clear what the AE property 
implies exactly. Are there examples for associative algebras which are AE but 
not locally m-convex, for example? Are there Lie algebras which are truly AE 
and not locally m-convex or nilpotent? We don't have an answer to this 
questions, but we can make some simple observations, which allow us to give an 
answer for special cases.
\begin{proposition}[Entire calculus]
	Let $\algebra{A}$ be an associative AE algebra. Then it admits an
	entire calculus.
\end{proposition}
\begin{proof}
	The proof is the same as for locally m-convex algebras: let 
	$f \colon \mathbb{C} \longrightarrow \mathbb{C}$ be an entire 
	function with $f(z) = \sum_{n=0}^{\infty} a_n z^n$ and $p$
	a continuous seminorm with an asymptotic estimate $q$.
	Then one has $\forall_{x \in \mathcal{A}}$
	\begin{equation*}
		p(f(x))
		=
		p \left(
			\sum\limits_{n=0}^{\infty}
			a_n x^n
		\right)
		\leq
		\sum\limits_{n=0}^{\infty}
		|a_n| 
		p \left( x^n \right)
		\leq
		\sum\limits_{n=0}^{\infty}
		|a_n| q(x)^n
		<
		\infty.
	\end{equation*}
\end{proof}
\begin{remark}[Entire Calculus, AE and LMC algebras]
	The fact that associative AE algebras have an entire calculus makes them 
	very similar to locally m-convex ones.
	Now there is something we can say about associative algebras which have an 
	entire calculus: if such an algebra is additionally commutative and 
	Fr\'echet, then must be even locally m-convex. This statement was proven 
	in \cite{mitiagin.rolewicz.zelazko:1962a} by Mitiagin, Rolewicz and 
	{\.Z}elazko. Oudadess and El kinani extended this result to commutative, 
	associative algebras, in which the Baire category theorem holds 
	\cite{elkinani.oudadess:1997a}. For non-commutative algebras, the 
	situation is different. There are associative ''Baire algebras'' having an 
	entire calculus, which are not locally m-convex. {\.Z}elazko gave an 
	example for such an algebra in \cite{zelazko:1994a}. Unfortunately, his 
	example is also not AE. It seems to be is an interesting (and non-trivial) 
	question, if a non locally m-convex but AE algebra exists at all and if 
	yes, how an example could look like.
\end{remark}

\subsection{The Projective Tensor Product}

We want to set up a topology on $\Sym^{\bullet}(\lie{g})$. Therefore, 
we first construct a topology on the tensor algebra $\Tensor^{\bullet}
(\lie{g})$. Since the construction is possible starting from any locally 
convex vector space, we want to go back to this more general setting for a 
moment. $V$ will always denote a locally convex vector space and $\algebra{P}$ 
the set of its continuous seminorms. In this situation, we can use the 
projective tensor product $\tensor[\pi]$ to get a locally convex 
topology on each tensor power $V^{\tensor[\pi] n}$. The precise construction 
of the projective tensor product can be found in good textbooks on locally 
convex analysis like \cite[Chapter 15]{jarchow:1981a} or in the lecture notes 
\cite[Lemma 4.1.4]{waldmann:2014a:script}. 
Recall that for $p_1, \ldots, p_n \in \algebra{P}$ we have a continuous 
seminorm on $V^{\tensor[\pi] n}$ via
\begin{equation*}
	(p_1 \tensor[\pi] \cdots \tensor[\pi] p_n)(x)
	=
	\inf
	\left\{
	\left.
		\sum_i
		p_1 \left( x_i^{(1)} \right)
		\cdots
		p_n \left( x_i^{(n)} \right)
	\ \right| \
		x
		=
		\sum_i
		x_i^{(1)}
		\tensor \cdots \tensor
		x_i^{(n)}
	\right\}.
\end{equation*}
On factorizing tensors, we moreover have the property
\begin{equation}
	\label{LCAna:FactorTensor}
	(p_1 \tensor[\pi] \cdots \tensor[\pi] p_n)
	(x_1 \tensor[\pi] \cdots \tensor[\pi] x_n)
	=
	p_1(x_1) \cdots p_n(x_n)
\end{equation}
which will be extremely useful in the following and which can be proven by
the Hahn-Banach theorem. We also have 
\begin{equation*}
	(p_1 \tensor \cdots \tensor p_n) 
	\tensor 
	(q_1 \tensor \cdots \tensor q_m) 
	= 
	p_1 \tensor \cdots \tensor p_n \tensor q_1 \tensor \cdots \tensor q_m.
\end{equation*}
For a given $p \in \algebra{P}$ we will denote $p^n = p^{\tensor[\pi] n}$ and 
$p^0$ is just the absolute value on the field $\mathbb{K}$. The $\pi$-topology 
on $V^{\tensor[\pi] n}$ is set up by all the projective tensor products of 
continuous seminorms, or, equivalently, by all the $p^n$ for $p \in 
\algebra{P}$.

The projective tensor product has a very helpful feature: if we want to show a 
(continuity) estimate on the tensor algebra, it is enough to do so on 
factorizing tensors. We will use this very often and just refer to it as the 
``infimum argument''.
\begin{lemma}[Infimum argument for the projective tensor product]
	\label{LCAna:Lemma:InfimumArgument}
 	Let $V_1, \ldots, V_n$, \hfill  \\
 	$W$ be locally convex vector spaces and 
	\begin{equation*}
		\phi \colon
		V_1 \times \cdots \times V_n
		\longrightarrow
		W
	\end{equation*}
	a $n$-linear map, from which we get the linear map 
	$\Phi \colon V_1 \tensor[\pi] \cdots \tensor[\pi] V_n \longrightarrow W$.
	Then $\Phi$ is continuous if and only if this is true for $\phi$. 
	If the estimate
	\begin{equation}
		\label{LCAna:ContOnFactors}
		p \left(
			\Phi (x_1 \tensor \cdots \tensor x_n)
		\right)
		\leq
		q_i(x_1) \cdots q_i(x_n)
	\end{equation}
	is fulfilled for $p\in\algebra{P}_V$, $q_i \in \algebra{P}_{W_i}$ 
	and all $x_i \in V_i, i=1, \ldots, n$ then we have
	\begin{equation}
		\label{LCAna:ContGeneral}
		p \left(
			\Phi (x)
		\right)
		\leq
		\left( q_1 \tensor \cdots \tensor q_n \right)(x)
	\end{equation}
	for all $x \in V_1 \tensor \cdots \tensor V_n$.
\end{lemma}
\begin{proof}
	If $\Phi$ is continuous, the continuity of $\phi$ is clear. 
	The other implication is more interesting.
	Continuity for $\phi$ means, that for every continuous seminorm $p$ on 
	$W$ we have continuous seminorms $q_i$ on $V_i$ with $i = 1, \ldots, n$ 
	such that for all $x^{(i)} \in V_i$ the estimate
	\begin{equation}
		\label{LCAna:ContiFactorTensors}
		p \left(
			\phi\left( x^{(1)}, \ldots, x^{(n)} \right)
		\right)
		\leq
		q_1 \left( x^{(1)} \right) 
		\cdots 
		q_n \left( x^{(n)} \right)
	\end{equation}
	holds. Let $x \in V_1 \tensor[\pi] \ldots \tensor[\pi] V_n$, then it 
	has a representation in terms of factorizing tensors like
	\begin{equation*}
		x
		=
		\sum_j
		x_j^{(1)} 
		\tensor[\pi] \cdots \tensor[\pi] 
		x_j^{(n)}.
	\end{equation*}
	We thus have
	\begin{align*}
		p( \Phi(x) )
		& =
		p \left(
			\sum_j
			\Phi \left(
				x_j^{(1)} 
				\tensor[\pi] \cdots \tensor[\pi] 
				x_j^{(n)}
			\right)
		\right)
		\\
		& \leq
		\sum_j
		p \left(
			\phi \left(
				x_j^{(1)} 
				, \cdots ,
				x_j^{(n)}
			\right)
		\right)
		\\
		& \leq
		\sum_j
		q_1 \left( x_j^{(1)} \right) 
		\cdots 
		q_n \left( x_j^{(n)} \right).
	\end{align*}
	Now we take the infimum over all possibilities of writing $x$ as a sum of 
	factorizing tensors on both sides. While nothing happens on the left 
	hand side, on the right hand side we find $\left( q_1 \tensor[\pi] 
	\ldots \tensor[\pi] q_n \right)(x)$: 
	\begin{equation*}
		p( \Phi(x) )
		\leq
		\inf
		\left\{
		\left.
			\sum_j
			q_1 \left( x_j^{(1)} \right) 
			\cdots 
			q_n \left( x_j^{(n)} \right)
		\ \right| \
			x
			=
			\sum_j
			x_j^{(1)} \cdots x_j^{(n)}
		\right\}.
	\end{equation*}
	This is the estimate we wanted and the Lemma is proven.
\end{proof}

Most of the time, we will deal with the symmetric tensor algebra. Therefore, 
we want to recall some basic facts about $\Sym^n(V)$, when it inherits the 
$\pi$-topology from the $V^{\tensor[\pi] n}$. We will call it 
$\Sym_{\pi}^n(V)$ when we endow it with this topology.
\begin{lemma}
	\label{Lemma:LCAna:ProjTensSymm}
	Let $V$ be a locally convex vector space, $p$ a continuous seminorm 
	and $n,m \in \mathbb{N}$.
	\begin{lemmalist}
		\item
		The symmetrization map
		\begin{equation*}
			\Symmetrizer_n
			\colon
			V^{\tensor_{\pi} n}
			\longrightarrow
			V^{\tensor_{\pi} n}
			, \quad
			(x_1 \tensor \ldots \tensor x_n)
			\longmapsto
			\frac{1}{n!}
			\sum\limits_{\sigma \in S_n}
			x_{\sigma(1)}
			\tensor \ldots \tensor
			x_{\sigma(n)}
		\end{equation*}
		is continuous and we have for all $x \in V^{\tensor_{\pi} n}$ 
		the estimate
		\begin{equation}
			\label{LCAna:SymmProdCont}
			p^n(\Symmetrizer_n(x))
			\leq
			p^n(x).
		\end{equation}
		
		\item
		Each symmetric tensor power $\Sym_{\pi}^n(V) \subseteq 
		V^{\tensor_{\pi} n}$ is a closed subspace.
		
		\item
		For $x \in \Sym_{\pi}^n(V)$ and $y \in \Sym_{\pi}^m(V)$ we have
		\begin{equation*}
			p^{n + m}(xy)
			\leq
			p^n(x) p^m(y).
		\end{equation*}
	\end{lemmalist}
\end{lemma}
\begin{proof}
	The first part is very easy to see and uses most of the tools which are 
	typical for the projective tensor product. We have the estimate for 
	factorizing tensors $x_1 \tensor \ldots \tensor x_n$
	\begin{align*}
		p^n \left(
			\Symmetrizer \left(
				x_1 \tensor \ldots \tensor x_n
			\right)
		\right)
		& =
		p^n \left(
			\frac{1}{n!}
			\sum	\limits_{\sigma \in S_n}
			x_{\sigma(1)} 
			\tensor \ldots \tensor 
			x_{\sigma(n)}
		\right)
		\\
		& \leq
		\frac{1}{n!}
		\sum\limits_{\sigma \in S_n}
		p^n \left(
			x_{\sigma(1)} 
			\tensor \ldots \tensor 
			x_{\sigma(n)}
		\right)
		\\
		& =
		\frac{1}{n!}
		\sum\limits_{\sigma \in S_n}
		p \left( x_{\sigma(1)} \right)
		\ldots
		p \left( x_{\sigma(n)} \right)
		\\
		& =
		p \left( x_1 \right)
		\ldots
		p \left( x_n \right)
		\\
		& =
		p^n \left(
			x_1 \tensor \ldots \tensor x_n
		\right).
	\end{align*}
	Then we use the infimum argument from Lemma 
	\ref{LCAna:Lemma:InfimumArgument} and we are done.
	The second part is also easy since the kernel of a continuous map is
	always a closed subspace of the initial space and we have
	\begin{equation*}
		\Sym_{\pi}^n 
		= 
		\ker (\id - \Symmetrizer_n).
	\end{equation*}
	The third part is a consequence from the first.
\end{proof}
One could maybe think that the inequality in the first part of this lemma is 
just an artefact which is due to the infimum argument and should actually be 
an equality, if one looked at it more closely. It is very interesting to see, 
that this is \textit{not} the case, since it may happen that this 
inequality is strict. The following example illustrates this.
\begin{example}
	We take $V = \mathbb{R}^2$ with the standard basis $e_1,	e_2$ 
	and $V$ is endowed with the maximum norm. Hence we have $\norm{e_1} = 
	\norm{e_2} = 1$. Now look at $e_1 \tensor e_2$, which has the norm
	\begin{equation*}
		\norm{
			e_1 \tensor e_2
		}
		= 
		\norm{e_1}
		\tensor
		\norm{e_2}
		=
		1.
	\end{equation*}
	We now evaluate the symmetrization map on $V \tensor[\pi] V$:
	\begin{equation*}
		\Symmetrizer \left(
			e_1 \tensor e_2
		\right)
		=
		\frac{1}{2}
		\left(
			e_1 \tensor e_2
			+
			e_2 \tensor e_1
		\right).
	\end{equation*}
	Our aim is to show, that the projective tensor product of the norm of this 
	symmetrized vector is not $1$. Therefore we need to find another way of 
	writing it which has a norm of less than $1$. Observe that
	\begin{equation*}
		\frac{1}{2}
		(e_1 \tensor e_2 + e_2 \tensor e_1)
		=
		\frac{1}{4}
		(
			(e_1 + e_2) \tensor (e_1 + e_2)
			+
			(-e_1 + e_2) \tensor (e_1 - e_2)
		)
	\end{equation*}
	and we have
	\begin{align*}
		& 
		\frac{1}{4}
		\norm{
			(e_1 + e_2) \tensor (e_1 + e_2)
			+
			(-e_1 + e_2) \tensor (e_1 - e_2)
		}
		\\
		& \leq
		\frac{1}{4}
		(
			\norm{ (e_1 + e_2) \tensor (e_1 + e_2) }
			+
			\norm{ (-e_1 + e_2) \tensor (e_1 - e_2) }
		)
		\\
		& \leq
		\frac{1}{4}
		(
			\norm{e_1 + e_2}
			\cdot
			\norm{e_1 + e_2}
			+
			\norm{- e_1 + e_2}
			\cdot
			\norm{e_1 - e_2}
		)
		\\
		& =
		\frac{1}{4}
		(1 \cdot 1 + 1 \cdot 1)
		\\
		& =
		\frac{1}{2}.
	\end{align*}
	So we have $\norm{\Symmetrizer(e_1 \tensor e_2)} \leq \frac{1}{2} < 1$.
\end{example}

\subsection{A Topology for the Gutt Star Product}

The next step is to set up a topology on $\Tensor^{\bullet}(V)$ which 
has the $\pi$-topology on each component. A priori, there are a lot of such 
topologies and at least two natural ones: the direct sum topology which is 
very fine and has a very small closure, and the cartesian product topology 
which is very coarse and therefore has a very big closure. We need something 
in between, which we can adjust in a convenient way.
\begin{definition}[$\Tensor_R$-topology]
	Let $p$ be a continuous seminorm on a locally convex vector space 
	$V$ and $R \in \mathbb{R}$. We define by 
	\begin{equation*}
		p_R 
		= 
		\sum\limits_{n=0}^{\infty}
		n!^R p^n
	\end{equation*}
	a seminorm on the tensor algebra $\Tensor^{\bullet}(V)$. We write for 
	the tensor or the symmetric algebra endowed with all such seminorms  
	$\Tensor_R^{\bullet}(V)$ or $\Sym_R^{\bullet}(V)$ respectively.
\end{definition}
Now we want to collect the most important results on the locally convex 
algebras $\left( \Tensor_R^{\bullet}(V), \tensor \right)$ and $\left( 
\Sym_R^{\bullet}(V), \vee \right)$. 
\begin{lemma}[The $\Tensor_R$-topology]
    \label{Lemma:LCAna:Projections}%
    Let $R' \geq R \geq 0$ and $q, p$ are continuous semi-norms on $V$.
    \begin{lemmalist}
	  \item \label{Item:EstimateForSeminorms}
	    	If $q \geq p$ then $q_R \geq p_R$ and $p_{R'} \geq p_R$.
	  \item \label{Item:TensorProductContinuous}
	    	The tensor product is continuous and satisfies the following 
	    	inequality:
	    	 \begin{equation*}
	    		p_R(x \tensor y)
	    		\leq
	    		\left( 2^R p \right)_R(x)
	    		\left( 2^R p \right)_R(y)
	    	\end{equation*}
	  \item \label{Item:PitopologyOnComponents}
	    	For all $n \in \mathbb{N}$ the induced topology on 
	    	$\Tensor^n(V) \subset \Tensor_R^{\bullet}(V)$ and on 
	    	$\Sym^n(V) \subset \Sym_R^{\bullet}(V)$ is the $\pi$-
	    	topology.
	  \item \label{Item:ComponentProjectionsContinuous}
	    	For all $n \in \mathbb{N}$ the projection and the inclusion 
	    	maps
	        \begin{equation*}
	        	\begin{array}{ccccc}
		    	    \Tensor_R^{\bullet}(V)
		        	&
	    	   		\ot{$\pi_n$}{\longrightarrow}
	    	    		&
	    	    		V^{\tensor_{\pi} n}
	    	    		&
	    	    		\ot{$\iota_n$}{\longrightarrow}
	    		    &
	    		    \Tensor^{\bullet}(V)
	    		    \\
		        \Sym_R^{\bullet}(V)
		        &
	    	        \ot{$\pi_n$}{\longrightarrow}
	    	    		&
	    	    		\Sym_{\pi}^n(V)
	    	    		&
	    	    		\ot{$\iota_n$}{\longrightarrow}
	    		    &
	    		    \Sym_R^{\bullet}(V)
	        	\end{array}
	        \end{equation*}
	        are continuous.
	  \item \label{Item:CompletionExplicitly}
    		The completions $\widehat{\Tensor}_R^{\bullet}(V)$ of 
    		$\Tensor_R^{\bullet}(V)$ and $\widehat{\Sym}_R^{\bullet}(V)$ 
    		of $\Sym_R^{\bullet}(V)$ can be described explicitly as
    		\begin{equation*}
	    		\begin{array}{ccccc}
		    		\widehat{\Tensor}_R^{\bullet}(V)
		    		&
		    		=
		    		&
		    		\left\{
		    		\left.
		    			x
		    			=
		    			\sum\limits_{n=0}^{\infty}
		    			x_n
		    		\ \right| \ 
		    			p_R(x)
		    			<
		    			\infty
		    			, \text{ for all } p
		    		\right\}
		    		&
		    		\subseteq
		    		&
		    		\prod\limits_{n=0}^{\infty}
		    		V^{\hat{\tensor}_{\pi} n}
		    		\\
		    		&&&&
		    		\\
		    		\widehat{\Sym}_R^{\bullet}(V)
		    		&
		    		=
		    		&
		    		\left\{
		    		\left.
		    			x
		    			=
		    			\sum\limits_{n=0}^{\infty}
		    			x_n
		    		\ \right| \ 
		    			p_R(x)
		    			<
		    			\infty
		    			, \text{ for all } p
		    		\right\}
		    		&
		    		\subseteq
	    			&
	    			\prod\limits_{n=0}^{\infty}
	    			\Sym_{\hat{\tensor}_{\pi}}^n
	    		\end{array}
    		\end{equation*}
    		with $p$ running through all continuous semi-norms on $V$ and
    		the $p_R$ are extended to the Cartesian product allowing the 
	    	value $+ \infty$.
      \item \label{Item:StrictlyFinerForBiggerR}
    		If $R' > R$, then the topology on $\Tensor_{R'}^{\bullet}(V)$ 
    		is strictly finer than the one on $\Tensor_R^{\bullet}(V)$, 
    		the same holds for $\Sym_{R'}^{\bullet}(V)$ and 
    		$\Sym_R^{\bullet}(V)$. Therefore the completions get smaller 
    		for bigger $R$.
      \item \label{Item:ComponentInclusionsContinuous}
    		The inclusion maps $\widehat{\Tensor}_{R'}^{\bullet}(V) 
    		\longrightarrow \widehat{\Tensor}_R^{\bullet}(V)$ and 
    		$\widehat{\Sym}_{R'}^{\bullet}(\lie{g}) \longrightarrow 
    		\widehat{\Sym}_R^{\bullet}(\lie{g})$ are continuous.
	  \item \label{Item:LmcJustForZero}
    		The topology on $\Tensor_R^{\bullet}(V)$ with the tensor 
    		product and on $\Sym_R^{\bullet}(V)$ with the symmetric 
    		product is locally m-convex if and only if $R = 0$.
      \item \label{Item:FirstCountable}
    		The algebras $\Tensor_R^{\bullet}(V)$ and $\Sym_R^{\bullet}(V)$
    		are first countable if and only if this is true for $V$.
      \item \label{item:PointsArePoints} The evaluation functionals
        $\delta_{\varphi}\colon \Sym^\bullet_R(V) \longrightarrow
        \mathbb{K}$ for $\varphi \in \lie{g}'$ are continuous.
    \end{lemmalist}
\end{lemma}
\begin{proof}
	The first part is clear on factorizing tensors and extends to the whole 
	tensor algebra via the infimum argument. For part $(ii)$, take two 
	factorizing tensors
	\begin{equation*}
		x
		=
		x^{(1)} \tensor \cdots \tensor x^{(n)}
		\quad \text{ and }
		y
		=
		y^{(1)} \tensor \cdots \tensor y^{(m)}
	\end{equation*}
	and compute:
	\begin{align*}
		p_R \left(
			x \tensor y
		\right)
		& =
		(n + m)!^R
		p^{n + m} \left(
			x^{(1)} \tensor \cdots \tensor x^{(n)}
			\tensor
			y^{(1)} \tensor \cdots \tensor y^{(m)}
		\right)
		\\
		& =
		(n + m)!^R
		p^n \left(
			x^{(1)} \tensor \cdots \tensor x^{(n)}
		\right)
		p^m \left(
			y^{(1)} \tensor \cdots \tensor y^{(m)}
		\right)
		\\
		& =
		\binom{n + m}{n}^R
		n!^R m!^R
		p^n \left(
			x^{(1)} \tensor \cdots \tensor x^{(n)}
		\right)
		p^m \left(
			y^{(1)} \tensor \cdots \tensor y^{(m)}
		\right)
		\\
		& \leq
		2^{(n + m) R}
		p_R \left( x^{(1)} \tensor \cdots \tensor x^{(n)} \right)
		p_R \left( y^{(1)} \tensor \cdots \tensor y^{(m)} \right)
		\\
		& =
		\left( 2^R p\right)_R 
		\left( x^{(1)} \tensor \cdots \tensor x^{(n)} \right)
		\left( 2^R p \right)_R 
		\left( y^{(1)} \tensor \cdots \tensor y^{(m)} \right).
	\end{align*}
	The parts $(iii)$ and $(iv)$ are clear from the construction of the 
	$R$-	topology. In part $(v)$ we used the completion of the tensor product 
	$\hat{\tensor}$, the statement itself is clear and implies $(vi)$ 
	directly, since we have really more elements in the completion for $R < 
	R'$, like the series over $x^{n} \frac{1}{n!^t}$ for $t \in (R, R')$ and 
	$0 \neq x \in V$. Statement $(vii)$ follows from the first. For $(viii)$, 
	it is easy to see that $\Tensor_0^{\bullet}(V)$ and $\Sym_0^{\bullet}(V)$ 
	are locally m-convex. For every $R > 0$ we have
	\begin{equation*}
		p_R \left(x^n \right)
		=
		n!^R
		p(x)^n
	\end{equation*}
	for all $n \in \mathbb{N}$ and all $x \in V$. If we had a 
	submultiplicative seminorm $\varphi$ from an equivalent topology, 
	then we would have some $x \in V$, and a continuous seminorm $p$ with 
	$p(x) \neq 0$ such that $p_R \leq \varphi$, and hence
	\begin{equation*}
		n!^R p(x)^R
		\leq
		\varphi \left(x^n\right)
		\leq
		\varphi\left(x\right)^n.
	\end{equation*}
	Since this is valid for all $n \in \mathbb{N}$, we get a contradiction.
	For the ninth part, the tensor algebras can not be first countable if this 
	is not true for $V$ itself. On the other hand, if $V$ has a finite base of 
	the topology, then $\Tensor_R^{\bullet}(V)$ and $\Sym_R^{\bullet}(V)$ are 
	just a countable multiple of $V$ and stay therefore first countable.
	The last part finally assures, that every symmetric tensor really gives a 
	continuous function and comes from the continuity of $\varphi$.
\end{proof}
The projective tensor product obviously keeps a lot of important and 
strong properties of the original vector space $V$. But Lemma 
\ref{Lemma:LCAna:Projections} still leaves out some important things. 
We will not make use of them in the following, but it is worth naming 
them for completeness. To do this in full generality, we 
need one more definition, which will be also very important in Chapter 6.
\begin{definition}\label{ProjectiveLimit}
	For a locally convex vector space $V$ and $R \geq 0$ we set
	\begin{equation*}
		\Sym_{R^-}^{\bullet}(V)
		=
		\projlim\limits_{\epsilon \longrightarrow 0}
		\Sym_{1 - \epsilon}^{\bullet}(V)
	\end{equation*}
	and call its completion $\widehat{\Sym}_{R^-}^{\bullet}(V)$.
\end{definition}
The projective limit is the intersection of all the algebras 
$\Sym_R^{\bullet}(V)$ for $R < 1$. Its completion can be understood as all 
those series in the cartesian product, which converge for all $p_R$ for $R < 1$ 
and $p$ a continuous seminorm on $V$. We hence see that the completion 
$\widehat{\Sym}_{R^-}^{\bullet}(V)$ is bigger than $\widehat{\Sym}_R^{\bullet}
(V)$ without the projective limit. This will become important in Chapter 6.
Now we can state two more propositions. Since we won't use them, we omit the 
proofs here. They can be found in \cite{waldmann:2014a}.
\begin{proposition}[Schauder bases]
	\label{Prop:LCAna:Bases}
	Let $R \geq 0$ and $V$ a locally convex vector space.
	If $\{e_i\}_{i \in I}$ is an absolute Schauder basis of $V$ with 
	coefficient functionals $\{\varphi^i\}_{i \in I}$, i.e. for every 
	$x \in V$ we have
	\begin{equation*}
		x
		=
		\sum\limits_{i \in I}
		\varphi^i(x) e_i
	\end{equation*}
	such that for every $p \in \algebra{P}$ there is a 
	$q \in \algebra{P}$ such that
	\begin{equation}
		\label{LCAna:AbsoluteSchauder}
		\sum\limits_{i \in I}
		|\varphi^i(x)|
		p(e_i)
		\leq
		q(x),
	\end{equation}
	then the set 
	$\{e_{i_1} \tensor \ldots \tensor e_{i_n}\}_{i_1, \ldots, i_n \in I}$ 
	defines an absolute Schauder basis of $\Tensor_R^{\bullet}(V)$ 
	together with the linear functionals $\{\varphi^{i_1} \tensor \ldots 
	\tensor 	\varphi^{i_n}\}_{i_1, \ldots, i_n \in I}$ which satisfy
	\begin{equation*}
		\sum\limits_{n=0}^{\infty}
		\sum\limits_{i_1, \ldots, i_n \in I}
		\left| 
			\left(
				\varphi^{i_1} \tensor \ldots \tensor \varphi^{i_n}
			\right)
			(x)
		\right|
		p_R \left(
			e_{i_1} \tensor \ldots \tensor e_{i_n}
		\right)
		\leq
		q_R (x)
	\end{equation*}
	for every $x \in \Tensor_R^{\bullet}(V)$ whenever $p$ and $q$ satisfy 
	\eqref{LCAna:AbsoluteSchauder}. The same statement is true for 
	$\Sym_R^{\bullet}(V)$ and for $\Sym_{R^-}^{\bullet}(V)$ (for $R > 0$) when 
	we choose a maximal linearly independent subset out of the set $\{e_{i_1} 
	\ldots e_{i_n}\}_{i_1, \ldots, i_n \in I}$.
\end{proposition}
\begin{proposition}[Nuclearity]
	\label{LCAna:Nuclearity}
	Let $V$ be a locally convex space. For $R \geq 0$ the following statements 
	are equivalent:
	\begin{propositionlist}
	  \item
		$V$ is nuclear.
	  \item
	  	$\Tensor_R^{\bullet}(V)$ is nuclear.
	  \item
	  	$\Sym_R^{\bullet}(V)$ is nuclear.	  	
	\end{propositionlist}
	If moreover $R > 0$, then the following statements are equivalent:
	\begin{propositionlist}
	  \item
		$V$ is strongly nuclear.
	  \item
	  	$\Tensor_R^{\bullet}(V)$ is strongly nuclear.
	  \item
	  	$\Sym_R^{\bullet}(V)$ is strongly nuclear.
	\end{propositionlist}
\end{proposition}
Nuclearity is an important and natural property in physics. Actually every 
example of a locally convex vector space, which play a role in physics, is 
either a normed space or nuclear. Therefore we want this property to be 
conserved.
\begin{remark}[The $\Tensor_R$-topology in the literature]
	It seems that this topology is actually not a new construction, but it 
	was already found and refound in the past. A very similar construction is 
	due to Goodman \cite{goodman:1971a}, who looked at finite-dimensional 
	Lie algebras. He studied differential operators on the Lie group, that 
	means the universal enveloping algebra, which he topologised via the 
	tensor algebra using basically the same weights. He needed a basis to do so 
	and proved afterwards that the topology he found is independent of the 
	choice of the basis. He also showed the continuity of the multiplication, 
	but by going a different way: he proved that certain ideals are 
	closed in this topology and divided them out. In this sense his approach 
	is less explicit than ours and also does not apply to infinite-dimensional 
	Lie algebras, where the proofs are more involved. Goodman was more 
	interested in the representation theory which arises from his construction 
	and he also knew about its functoriality. However, he did not know 
	about the deformation aspect behind it, since the ideas of deformation 
	quantization did not exist at that time. It seems a bit like his work has 
	gone more or less unnoticed, we could not find a work which used this 
	topology later for studies in Lie theory. A very different way of 
	constructing a topology which coincides with ours for $R = 1$ was moreover 
	given before by Ra{\v{s}}evski{\u{i}} in \cite{rasevskii:1966a}, who also 
	restricted to the finite-dimensional case and used bases of his Lie algebra 
	to construct the topology. For this reason, also his ideas are bounded to 
	finite-dimensional Lie algebras.
\end{remark}

\section{Continuity Results for the Gutt Star Product}
\label{sec:chap5_TopologyStar}

From now on, we start with an AE Lie algebra $\lie{g}$ rather than with a 
general locally convex space $V$. We have most of the tools at hand to show the 
continuity of the Gutt star product, except for precise estimates on the 
$\widetilde{\mathrm{BCH}_{a, b}}$-terms. Therefore, we state the following lemma.
\begin{lemma}
	\label{LCAna:Lemma:BCHTermsEstiamte}
	Let $\lie{g}$ be a AE-Lie algebra, $\xi, \eta \in \lie{g}$, $p$ a 
	continuous seminorm, $q$ an	asymptotic estimate for it and $a, b, n \in 
	\mathbb{N}$ with $a + b = n$. Then, using the Goldberg-Thompson form of 
	the Baker-Campbell-Hausdorff series, we have the following estimates:
	\begin{lemmalist}
	  \item \label{Item:ThompsonEstimate}
		The coefficients $g_w$ from \eqref{Alg:BCHinGoldbergThompson} 
		fulfil the estimate
		\begin{equation}
			\label{LCAna:ThompsonEstimate}
			\sum\limits_{|w| = n}
			\left| \frac{g_w}{n} \right|
			\leq
			\frac{2}{n}.
		\end{equation}
		Recall that $|w|$ denotes the length of a word $w$ and $[w]$ is the
		word put in Lie brackets nested to  the left.
	\item \label{Item:AEonWords}
		For every word $w$, which consists of $a$ $\xi$ and $b$ $\eta$,
		we have
		\begin{equation}
		\label{LCAna:AEonWords}
			p([w])
			\leq
			q(\xi)^a q(\eta)^b.
		\end{equation}
	\item \label{Item:BCHEstimate}
		We have the estimate
		\begin{equation}
		\label{LCAna:BCHEstimate}
			p \left(
				\bchtilde{a}{b}
				{\xi_1, \ldots, \xi_a}
				{\eta_1, \ldots, \eta_b}
			\right)
			\leq
			\frac{2}{a + b}
			q\left( \xi_1 \right) \cdots q\left( \xi_a \right)
			q\left( \eta_1 \right) \cdots q\left( \eta_b \right)
		\end{equation}
	\end{lemmalist}
\end{lemma}
\begin{proof}
	We already showed part $(i)$ in Proposition 
	\ref{Alg:Prop:ThompsonsEstimate} and put in the factor $\frac{1}{n}$ 
	because these factors will appear later. The next estimate 
	\eqref{LCAna:AEonWords} is due to the AE property which does not see the 
	way how brackets are set but just counts the number of $\xi_i$ and 
	$\eta_j$ in the whole expression.
	Let us use the notation $|w|_{\xi}$ for the number of $\xi$'s appearing
	in a word $w$ and $|w|_{\eta}$ for the number of $\eta$'s. Clearly, $|w| =
	|w|_{\xi} + |w|_{\eta}$. With \eqref{LCAna:ThompsonEstimate} and the AE
	property of $\lie{g}$, we get
	\begin{align*}
		p \left(
			\bchtilde{a}{b}
			{\xi_1, \ldots, \xi_a}
			{\eta_1, \ldots, \eta_b}
		\right)
		& \leq
		\sum\limits_{\substack{
			|w|_{\xi} = a \\
			|w|_{\eta} = b
		}}
		p^{a + b} \left(
			\frac{g_w}{a + b}
			[w]
		\right)
		\\
		& \leq
		\sum\limits_{\substack{
			|w|_{\xi} = a \\
			|w|_{\eta} = b
		}}
		\frac{ |g_w| }{a + b}
		p([w])
		\\
		& \leq
		\sum\limits_{ |w| = a + b }
		\frac{ |g_w| }{a + b}
		q\left( \xi_1 \right) \cdots q\left( \xi_a \right)
		q\left( \eta_1 \right) \cdots q\left( \eta_b \right)
		\\
		& \leq
		\frac{2}{a + b}
		q\left( \xi_1 \right) \cdots q\left( \xi_a \right)
		q\left( \eta_1 \right) \cdots q\left( \eta_b \right).
	\end{align*}
\end{proof}
For showing the continuity of the star product we can either proceed via the 
big formula \eqref{Formulas:2MonomialsFormula2} for two monomials or via the 
smaller one \eqref{Formulas:LinearMonomial2} for a monomial with a vector and 
iterate it. While the results are a bit better for the first approach, the 
second makes the proof easier. Nevertheless, both approaches give strong 
results, and this is why we will to give both proofs here.

There will be a very general way how most of the proofs will work and which 
tools will be used in the following. If we want to show the continuity of a 
map $f \colon \Sym_R^{\bullet}(\lie{g}) \longrightarrow 
\Sym_R^{\bullet}(\lie{g})$, we will proceed most of the time like this:
\begin{enumerate}
	\item \label{Item:LCAna:Step1}
	First, we extend the map to the whole tensor algebra by symmetrizing 
	beforehand: $f = f \circ \Symmetrizer$. This is a real extension since the 
	symmetrization does not affect symmetric tensors. We do this here for the 
	Gutt star product
	\begin{equation*}
		\ostar_z
		\colon
		\Tensor^{\bullet}(\lie{g})
		\tensor
		\Tensor^{\bullet}(\lie{g})
		\longrightarrow
		\Tensor^{\bullet}(\lie{g})
		, \quad
		\ostar_z 
		= 
		\star_z \circ 
		\left( \Symmetrizer \tensor \right)
	\end{equation*}
	and for the $C_n$-operators analogously.

	\item \label{Item:LCAna:Step2}
	Then we start with an estimate which we prove only on 
	factorizing tensors in order to use the infimum argument 
	(Lemma \ref{LCAna:Lemma:InfimumArgument}).

	\item \label{Item:LCAna:Step3}
	During the estimation process, we find products of Lie 
	brackets. Those will be split up by the continuity of the symmetric 
	product \eqref{LCAna:SymmProdCont} from Lemma 
	\ref{Lemma:LCAna:ProjTensSymm} the AE property \eqref{LCAna:AEProperty}.
	
	\item \label{Item:LCAna:Step4}
	Finally, we rearrange the split up seminorms to the seminorm of a 
	factorizing tensor by \eqref{LCAna:FactorTensor}.
\end{enumerate}

\subsection{Continuity of the Product}
In the first proof, we want to give an estimate via the formula
\begin{equation}
    \label{eq:GstarOfXisAndEtas}
    \xi_1 \cdots \xi_k \star_z \eta_1 \cdots \eta_{\ell}
    =
    \sum\limits_{n=0}^{k + \ell -1}
    z^n
    C_n (\xi_1 \cdots \xi_k, \eta_1 \cdots \eta_{\ell}).
\end{equation}
To shorten the very long expression from Equation 
\eqref{Formulas:2MonomialsFormula22}, we abbreviate the summations by
\begin{equation}
	\label{LCAna:ShorterSumsInLongProduct}
    C_n\left(
        \xi_1 \cdots \xi_k,
        \eta_1 \cdots \eta_{\ell}
    \right)
    =
    \frac{1}{r!}
    \sum\limits_{\sigma, \tau}
    \sum\limits_{a_i, b_j}
    \bchtilde{a_1}{b_1}{\xi_{\sigma(i)}}{\eta_{\tau(j)}}
    \cdots
    \bchtilde{a_r}{b_r}{\xi_{\sigma(i)}}{\eta_{\tau(j)}},
\end{equation}
meaning the summations as given in
Proposition~\ref{Formulas:Prop:2MonomialsFormula2} and 
using $r = k + \ell - n$.
\begin{theorem}[Continuity of the star product]
    \label{Thm:LCAna:Continuity1}%
    Let $\lie{g}$ be an AE-Lie algebra, $R \geq 0$, $p$ a continuous seminorm
    with an asymptotic estimate $q$ and $z \in \mathbb{C}$.
    \begin{theoremlist}
    	\item \label{item:CnOperatorEstimate}
    	For $n \in \mathbb{N}$, the operator $C_n$ is continuous and for all
    	$x, y \in \Tensor_R^{\bullet}(\lie{g})$ we have the estimate:
    	\begin{equation}
    		\label{LCAna:CnOperators}
    		p_R \left( C_n(x,y) \right)
			\leq
        	\frac{n!^{1 - R}}{2 \cdot 8^n}
        	(16 q)_R (x)
        	(16 q)_R (y)
    	\end{equation}
    	\item \label{item:LCAna:Continuity1}
    	For $R \geq 1$, the Gutt star product is continuous and for all
    	$x, y \in \Tensor_R^{\bullet}(\lie{g})$ we have the estimate:
	    \begin{equation}
	   	    \label{LCAna:Continuity1}
	        p_R(x \ostar_z y)
	        \leq
	        (c q)_R(x) (c q)_R(y)
	    \end{equation}
	    with $c = 16(|z| + 1)$. Hence, the Gutt star product is continuous and
	    the estimate \eqref{LCAna:Continuity1} holds on
	    $\widehat{\Sym}_R^\bullet(\lie{g})$ for all $z \in \mathbb{C}$
	\end{theoremlist}
\end{theorem}
\begin{proof}
    Let us use $r = k + \ell - n$ as before and recall that the
    products are taken in the symmetric algebra. Then we can use
    Equation \eqref{LCAna:ShorterSumsInLongProduct} and put estimates on it.  
    Let $p$ be a continuous seminorm and let $q$ be an asymptotic estimate
    for it. Then we get
    \begin{align*}
        p_R \big(
            C_n \big(
                \xi_1 \tensor \cdots \tensor \xi_k, &
                \eta_1 \tensor \cdots \tensor \eta_{\ell}
            \big)
        \big)
        \\
        &=
        p_R \bigg(
         	\frac{1}{r!}
			\sum\limits_{\sigma, \tau}
			\sum\limits_{a_i, b_j}
			\bchtilde{a_1}{b_1}{\xi_{\sigma(i)}}{\eta_{\tau(j)}}
			\cdots
			\bchtilde{a_r}{b_r}{\xi_{\sigma(i)}}{\eta_{\tau(j)}}
        \bigg)
        \\
        & \ot{(a)}{\leq}
        \frac{1}{r!}
        r!^R
        \sum\limits_{\sigma, \tau}
		\sum\limits_{a_i, b_j}
        p \left(
            \bchtilde{a_1}{b_1}{\xi_{\sigma(i)}}{\eta_{\tau(j)}}
        \right)
        \cdots
        p \left(
            \bchtilde{a_r}{b_r}{\xi_{\sigma(i)}}{\eta_{\tau(j)}}
        \right)
        \\
        & \ot{(b)}{\leq}
        \frac{1}{r!^{1-R}}
        \sum\limits_{\sigma, \tau}
		\sum\limits_{a_i, b_j}
		\frac{2}{a_1 + b_1}
		\ldots
		\frac{2}{a_r + b_r}
		q(\xi_1) \cdots q(\xi_k)
		q(\eta_1) \cdots q(\eta_{\ell})
        \\
        & \ot{(c)}{\leq}
		q(\xi_1) \cdots q(\xi_k)
		q(\eta_1) \cdots q(\eta_{\ell})
        2^r
        \frac{k! \ell!}{r!^{1-R}}
        \sum\limits_{a_i, b_j}
		1,
    \end{align*}
    where we just used the continuity estimate for the symmetric
    tensor product in ($a$), Lemma~\ref{LCAna:Lemma:BCHTermsEstiamte},
    \refitem{Item:BCHEstimate}, in ($b$) and $\frac{2}{a_i + b_i} \leq
    2$ in ($c$). We estimate the number of terms in the sum and get
    \begin{equation*}
        \sum\limits_{\substack{a_1, b_1, \ldots, a_r, b_r \geq 0 \\
            a_i + b_i \geq 1 \\
            a_1 + \cdots + a_r = k \\
            b_1 + \cdots + b_r = \ell
          }}
		1
		\leq
        \sum\limits_{\substack{a_1, b_1, \ldots, a_r, b_r \geq 0 \\
            a_1 + b_1 + \cdots + a_r + b_r = k + \ell
          }}
		1
		=
		\binom{k + \ell + 2r - 1}{k + \ell}
		\leq
		2^{3(k + \ell) - 2n - 1}.
    \end{equation*}
    Using this estimate, we get
    \begin{align*}
        &p_R \big(
            C_n \big(
                \xi_1 \tensor \cdots \tensor \xi_k,
                \eta_1 \tensor \cdots \tensor \eta_{\ell}
            \big)
        \big) \\
        &\quad\leq
        q(\xi_1) \cdots q(\xi_k)
		q(\eta_1) \cdots q(\eta_{\ell})
        2^{k + \ell - n}
        \frac{k! \ell!}{(k + \ell - n)!^{1-R}}
        2^{3(k + \ell) - 2n - 1}
        \\
        &\quad=
        q_R \left(
            \xi_1 \tensor \cdots \tensor \xi_k
        \right)
        q_R \left(
            \eta_1 \tensor \cdots \tensor \eta_{\ell}
        \right)
        2^{4(k + \ell) - 3n - 1}
        \left(
        	\frac{k! \ell! n!}{(k + \ell - n)! n!}
        \right)^{1-R}
        \\
        &\quad\leq
        q_R \left(
            \xi_1 \tensor \cdots \tensor \xi_k
        \right)
        q_R \left(
            \eta_1 \tensor \cdots \tensor \eta_{\ell}
        \right)
        2^{4(k + \ell) - 3n - 1}
        2^{(1 - R)(k + \ell)}
        n!^{1 - R}
        \\
        &\quad=
        \frac{n!^{1 - R}}{2 \cdot 8^n}
        (16 q)_R \left(
            \xi_1 \tensor \cdots \tensor \xi_k
        \right)
        (16 q)_R \left(
            \eta_1 \tensor \cdots \tensor \eta_{\ell}
        \right).
    \end{align*}
    The estimate \eqref{LCAna:CnOperators} is now proven on factorizing
    tensors.
    For general tensors $x,y \in \Tensor_R^{\bullet}(\lie{g})$, we use the 
    infimum argument from Lemma~\ref{LCAna:Lemma:InfimumArgument} and the 
    first part is done. For the second part, let $x$ and $y$
    be tensors of degree at most $k$ and $\ell$ respectively. We use
    \eqref{LCAna:CnOperators} in (a), the fact that $R \geq 1$ in (b) and have
    \begin{align}
	    \nonumber
	    	p_R \left(
	    	x \star_z y
	    	\right)
	    	& =
	    	p_R \left(
	    		\sum\limits_{n=0}^{l + \ell - 1}
	    		z^n C_n(x, y)
	    	\right)
	    	\\
	    \nonumber
	    	& \leq
	    	\sum\limits_{n=0}^{k + \ell - 1}
	    	p_R \left(
	    		z^n C_n(x, y)
	    	\right)
	    	\\
	    \label{LCAna:2MonomialEstimate}
	    	& \ot{(a)}{\leq}
	    	\sum\limits_{n=0}^{k + \ell - 1}
	    	\frac{|z|^n}{2 \cdot 8^n}
	    	n!^{1 - R}
	    	(16 q)_R(x)
	    	(16 q)_R(y)
	    	\\
	    \nonumber
	    	& \ot{(b)}{\leq}
	    	\frac{(|z| + 1)^{k + \ell}}{2}
	    	(16 q)_R(x)
	    	(16 q)_R(y)
	    	\sum\limits_{n = 0}^{\infty}
	    	\frac{1}{8^n}
			\\
		\nonumber
			& \leq
	    	(16(|z| + 1) q)_R(x)
	    	(16(|z| + 1) q)_R(y).
    \end{align}
    Since estimates on $\Sym_R^{\bullet}(\lie{g})$ also hold for the
    completion, the theorem is proven.
\end{proof}
\begin{remark}[Uniform continuity]
	Part $(i)$ of the theorem makes clear, why exactly continuity will 
	only hold if $R \geq 1$: the estimate in \eqref{LCAna:CnOperators} shows, 
	that all the $C_n$ are indeed continuous for any $R \geq 0$, but only 
	for $R \geq 1$ there is something like a uniform continuity. When $R$ 
	decreases, the continuity of the $C_n$ ``gets worse'' and the uniform 
	continuity finally breaks down when the threshold $R = 1$ is trespassed. 
	But we need this uniform estimate, since we have to control the 
	operators up to an arbitrarily high order if we want to guarantee the 
	continuity of the whole star product. Continuity up to a previously fixed
	order $n$ does not suffice.
\end{remark}

Now, we want to give a second proof, which relies on 
\eqref{Formulas:LinearMonomial2}. Approaching like this, we do not account for 
the fact that we will encounter terms like $[\eta, \eta]$ which vanish, 
but we estimate more brutally. During this procedure, we also count the 
formal parameter $z$ more often than it is actually there. This is why we now 
have to make assumptions on $R$ and $z$ which are a bit stronger than before. 
Moreover, we split up tensor products and put them together again various 
times, which is the reason why an AE Lie algebra does not suffice any more: we 
need $\lie{g}$ to be locally m-convex. But if we make these assumptions, 
we get the following lemma, which finally simplifies the proof.
\begin{lemma}
    \label{LCAna:Lemma:PreContinuity2}%
    Let $\lie{g}$ be a locally m-convex Lie algebra and $R \geq1$. 
    Then if $|z| < 2 \pi$ or $R >1$ there exists, for $x \in
    \Tensor^{\bullet}(\lie{g})$ of degree at most $k$, $\eta \in \lie{g}$
    and each continuous submultiplicative seminorm $p$, a constant $c_{z,R}$ 
    only depending on $z$ and $R$ such that the following estimate holds:
    \begin{equation}
        \label{LCAna:PreContinuity2}
        p_R(x \ostar_z \eta)
        \leq
        c_{z,R} (k+1)^R p_R(x) q(\eta)
    \end{equation}
\end{lemma}
\begin{proof}
    We have for $\xi_1, \ldots, \xi_k, \eta \in \lie{g}$
    \begin{align*}
        p_R 
        &
        \left(
            \xi_1 \tensor
             \cdots \tensor \xi_k \ostar_z \eta
        \right)
        =
        p_R \left(
            \sum\limits_{n=0}^k
            \frac{B_n^* z^n}{n! (k-n)!}
            \sum\limits_{\sigma \in S_k}
            \xi_{\sigma(1)} \cdots \xi_{\sigma(k - n)}
            \left(
                \ad_{\xi_{\sigma(k - n + 1)}}
                \circ \cdots \circ
                \ad_{\xi_{\sigma(k)}}
            \right)
            (\eta)
        \right)
        \\
        & =
        \sum\limits_{n=0}^k
        \frac{|B_n^*| |z|^n}{n! (k - n)!}
        \sum\limits_{\sigma \in S_k}
        (k + 1 - n)!^R
        p^{k + 1 - n}
        \left(
            \xi_{\sigma(1)} \cdots \xi_{\sigma(k - n)}
            \left(
                \ad_{\xi_{\sigma(k - n + 1)}}
                \circ \cdots \circ
                \ad_{\xi_{\sigma(k)}}
            \right)
            (\eta)
        \right)
        \\
        & \leq
        (k + 1)^R
        \sum\limits_{n=0}^k
        \frac{|B_n^*| |z|^n}{n!}
        (k - n)!^{R - 1}
        k! p(\xi_1) \cdots p(\xi_k) p(\eta)
        \\
        & =
        (k + 1)^R
        \sum\limits_{n=0}^k
        \frac{|B_n^*| |z|^n}{n!^R}
        \left( \frac{(k - n)! n!}{k!} \right)^{R - 1}
        p_R \left(
         	\xi_1 \tensor \cdots \tensor \xi_k
        \right)
        p(\eta)
        \\
        & \leq
        (k + 1)^R
		p_R \left(
         	\xi_1 \tensor \cdots \tensor \xi_k
        \right)
        p(\eta)
        \sum\limits_{n=0}^k
        \frac{|B_n^*| |z|^n}{n!^R}.
    \end{align*}
    Now if $|z| < 2 \pi$ the sum can be estimated by extending it to a
    series which converges. So we get a constant $c_{z, R}$ depending
    on $R$ and on $z$ such that
    \begin{equation*}
        p_R \left(
            \xi_1
            \tensor \cdots \tensor \xi_k
            \ostar_z
            \eta
        \right)
        \leq
        (k + 1)^R c_{z, R}
		p_R \left(
         	\xi_1 \tensor \cdots \tensor \xi_k
        \right)
        p(\eta).
    \end{equation*}
    On the other hand, if $|z| \geq 2 \pi$ and $R > 1$ we can estimate
    \begin{equation*}
        p_R \left(
            \xi_1
            \tensor \cdots \tensor \xi_k
            \ostar_z
            \eta
        \right)
        \leq
        (k + 1)^R
		p_R \left(
         	\xi_1 \tensor \cdots \tensor \xi_k
        \right)
        p(\eta)
        \left(
            \sum\limits_{n=0}^k
            \frac{|B_n^*|}{n!}
        \right)
        \left(
            \sum\limits_{n=0}^k
            \frac{|z|^n}{n!^{R - 1}}
        \right).
    \end{equation*}
    Again, both series will converge and give constants depending only
    on $z$ and $R$. Hence, we have the estimate on factorizing tensors
    and can extend this to generic tensors of degree at most $k$ by using
    the infimum argument.
\end{proof}

In the following, we assume again that either $R > 1$ or $R \geq 1$
and $|z| < 2\pi$ in order the use
Lemma~\ref{LCAna:Lemma:PreContinuity2}. Now we can give a simpler proof of 
Theorem~\ref{Thm:LCAna:Continuity1} for the case of a locally m-convex 
Lie algebra:
\begin{proof}[Alternative Proof of Theorem~\ref{Thm:LCAna:Continuity1}]
    Assume that $\lie{g}$ is now even locally m-convex.  We want to
    replace $\eta$ in the foregoing lemma by an arbitrary tensor $y$
    of degree at most $\ell$. Again, we do that on factorizing tensors
    first and get
    \begin{align*}
        p_R \big(
         	\xi_1 \tensor \cdots \tensor \xi_k 
		&
         	\ostar_z 
         	\eta_1 \tensor \cdots \tensor \eta_\ell
        \big)
        =
        p_R \left(
	        \frac{1}{\ell !}
         	\sum\limits_{\tau \in S_\ell}
         	\xi_1 \tensor \cdots \tensor \xi_k 
         	\ostar_z
         	\eta_{\tau(1)} \tensor \cdots \tensor \eta_{\tau(\ell)}
        \right)
        \\
        & \leq
        \frac{1}{\ell !}
        \sum\limits_{\tau \in S_\ell}
        p_R \left(
         	\xi_1 \tensor \cdots \tensor \xi_k 
         	\ostar_z
         	\eta_{\tau(1)} \tensor \cdots \tensor \eta_{\tau(\ell)}
        \right)
        \\
        & \leq
        c_{z,R} (k + \ell)^R
        \frac{1}{\ell !}
        \sum\limits_{\tau \in S_\ell}
        p_R \left(
         	\xi_1 \tensor \cdots \tensor \xi_k 
         	\ostar_z 
         	\eta_{\tau(1)} \ostar \cdots \ostar_z \eta_{\tau(\ell - 1)}
        \right)
        p\left( \eta_{\tau(\ell)} \right)
        \\
        & \leq
        \quad \vdots
        \\
        & \leq
        c_{z,R}^{\ell} ((k + \ell) \cdots (k + 1))^R
        \frac{1}{\ell !}
        \sum\limits_{\tau \in S_\ell}
        p_R \left( \xi_1 \tensor \cdots \tensor \xi_k \right)
        p \left( \eta_{\tau(1)} \right)
        \cdots
        p \left( \eta_{\tau(\ell)} \right)
        \\
        & =
        c_{z,R}^{\ell} \left(
        \frac{(k + \ell)!}{k! \ell!}
        \right)^R
        p_R \left( \xi_1 \tensor \cdots \tensor \xi_k \right)
        p_R \left( \eta_1 \tensor \cdots \tensor \eta_\ell \right)
        \\
        & \leq
        (2^R p)_R 
        \left( \xi_1 \tensor \cdots \tensor \xi_k \right)
        (2^R c_{z,R}\ p)_R 
        \left( \eta_1 \tensor \cdots \tensor \eta_\ell \right).
    \end{align*}
    Once again, we extend the estimate via the infimum argument to the whole
    tensor algebra, since the estimate depends no longer on the degree
    of the tensors.
\end{proof}

Using this approach for continuity, it is easy to see that nilpotency
of the Lie algebra would change the estimate substantially: if we knew that
we had at most $N$ brackets because $N + 1$ brackets vanish,
then the sum in the proof of Lemma \ref{LCAna:Lemma:PreContinuity2}
would end at $N$ instead of $k$ and would therefore be independent of
the degree of $x$.

In both proofs, it is easy to see that we need at least $R \geq 1$ to
get rid of the factorials which come up because of the combinatorics
of the star product. It is nevertheless interesting to see that this
result is sharp, that means the Gutt star product really fails
continuity for $R < 1$:
\begin{example}[A counter-example]
    \label{LCAna:Ex:HeisenbergAlgebra}%
    Let $0 \leq R < 1$ and $\lie{g}$ be the Heisenberg algebra in three 
    dimensions, i.e. the Lie algebra generated by the elements 
    $P$, $Q$ and $E$ with the bracket $[P,Q] = E$ and all other brackets
    vanishing. This is a very simple example for a non-abelian Lie algebra
    and if continuity of the star product fails for this one, then we
    can not expect it to hold for more complex ones. We impose on
    $\lie{g}$ the $\ell^1$-topology with the norm $n$ and $n(P) = n(Q)
    = n(E) = 1$. This will be helpful, since here we really have the equality
    \begin{equation*}
    		n^{k + \ell} \left( X^k Y^\ell \right)
    		=
    		n^k \left( X^k \right)
    		n^\ell \left( Y^\ell \right)
    \end{equation*}
    for the symmetric product of vectors $X, Y$. Then we consider 
    the sequences
    \begin{equation*}
        a_k
        =
        \frac{P^k}{k!^{R + \epsilon}}
        \quad
        \textrm{and}
        \quad
        b_k
        =
        \frac{Q^k}{k!^{R + \epsilon}}
    \end{equation*}
    with $2 \epsilon < 1 - R$. It is easy to see that
    \begin{equation*}
        n_R(a_k)
        =
        n_R(b_k)
        =
        k!^{- \epsilon}
    \end{equation*}
    and hence we get the limit for any $c > 0$ by
    \begin{equation*}
     	\lim_{n \longrightarrow \infty}
     	(cn)_R(a_k)
     	=
     	\lim_{n \longrightarrow \infty}
     	(cn)_R(a_k)
     	=
     	0
    \end{equation*}
    We want to show that there is no $c > 0$ such that
    \begin{equation*}
        n_R(a_k \star_z b_k)
        \leq
        (c n)_R(a_k) (c n)_R(b_k).
    \end{equation*}
    In other words, $n_R(a_k \star_z b_k)$ grows faster than
    exponentially. But this is the case, since we can calculate the star
    product explicitly and see
    \begin{align*}
        n_R(a_k \star_z b_k)
        & =
        n_R \left(
        \sum\limits_{j=0}^k
        \binom{k}{j}
        \binom{k}{j}
        j! \frac{1}{k!^{2R + 2 \epsilon}}
        P^{k-j} Q^{k-j} E^j
        \right)
        \\
        & =
        \sum\limits_{j=0}^k
        \frac{k!^2 j! (2k - j)!^R}
        {(k-j)!^2 j!^2 k!^{2R + 2 \epsilon}}
        \underbrace{
        n^{2k-j}
        ( P^{k-j} Q^{k-j} E^j )
        }_{= 1}
        \\
        & =
        \sum\limits_{j=0}^k
        \underbrace{
        \binom{k}{j}^2 \binom{2k}{k} \binom{2k}{j}^{-1}
        }_{\geq 1}
        \frac{j!^{1-R}}{k!^{2 \epsilon}}
        \\
        & \geq
        \sum\limits_{j=0}^k
        \frac{j!^{1-R}}{k!^{2 \epsilon}}
        \\
        & \geq
        k!^{1 - R - 2\epsilon}.
    \end{align*}
    Finally, we get a contradiction to the continuity of the star product 
    since $a_k \star_z b_k$ is unbounded in the topology of $\Sym_R(\lie{g})$ 
    although the sequences themselves go to zero.
\end{example}

\subsection{Dependence on the Formal Parameter}

Now we want to analyse the completion $\widehat{\Sym}_R^{\bullet}(\lie{g})$ of
the symmetric algebra with $\star_z$. Concerning exponential functions, we 
get the following negative result:
\begin{proposition}
    \label{proposition:NoExponentialsSorry}%
    Let $\xi \in \lie{g}$ and $R \geq 1$, then $\exp(\xi) \not\in
    \widehat{\Sym}_R^{\bullet}(\lie{g})$, where $\exp(\xi) =
    \sum_{n=0}^{\infty} \frac{\xi^n}{n!}$.
\end{proposition}
\begin{proof}
    Take a seminorm $p$ such that $1 \geq p(\xi) \neq 0$. Then set 
    $c = p(\xi)^{-1} \geq 1$. For $\xi^n$ the powers in the sense of either 
    the usual tensor product, or the symmetric product or the star product
    are the same. So we have for $N \in \mathbb{N}$
    \begin{equation*}
        (cp)_R \left(
        \sum\limits_{n=0}^N
        \frac{1}{n!} \xi^n
        \right)
        =
        \sum\limits_{n=0}^N
        \frac{n!^R}{n!}
        c^n
        p(\xi)^{n}
        =
        \sum\limits_{n=0}^N
        n!^{R - 1}
        \geq
        N,
    \end{equation*}
    and hence $\exp(\xi)$ does not converge for the seminorm $(cp)_R$.
\end{proof}
When we go back to Theorem~\ref{Thm:LCAna:Continuity1}, we see that we have 
actually proven slightly more than we stated: we showed that the star product 
converges absolutely and locally uniform in $z$. This means that the star product 
is \emph{not only continuous}, but also that \emph{the formal series converges to 
the star product in the completion}. We can take a closer look at this proof in 
order to get a new result for the dependence on the formal parameter $z$:
\begin{proposition}[Dependence on $z$]
    \label{corollary:HolomorphicDependence}%
    Let $R \geq 1$, then for all $x, y \in \widehat{\Sym}_R^\bullet(\lie{g})$ 
    the map
    \begin{equation}
        \label{LCAna:Holomorphicity}
        \mathbb{K} \ni z
        \longmapsto
        x \star_z y \in
        \widehat{\Sym}_R^\bullet(\lie{g})
    \end{equation}
    is analytic with (absolutely convergent) Taylor expansion at $z = 0$ 
    given by Equation~\eqref{Formulas:2MonomialsFormula1}. For 
    $\mathbb{K} = \mathbb{C}$, the collection of algebras $\left\{ \left( 
    \widehat{\Sym}_R^\bullet(\lie{g}), \star_z \right) \right\}_{z \in 
    \mathbb{C}}$ is an entire holomorphic deformation of the completed 
    symmetric tensor algebra $\left( \widehat{\Sym}_R^\bullet(\lie{g}), \vee 
    \right)$.
\end{proposition}
\begin{proof}
	The crucial point is that for $x, y \in \widehat{\Sym}_R^\bullet
	(\lie{g})$ and every continuous seminorm $p$ we have an asymptotic
	estimate $q$ such that
	\begin{align*}
		p_R \left( x \star_z y \right)
		& =
		p_R
		\left(
			\sum\limits_{n=0}^{\infty}
			z^n C_n(x,y)
		\right)
		\\
		& =
		\sum\limits_{n=0}^{\infty}
		|z|^n
		p_R( C_n(x, y) )
		\\
		& \leq
		(16 q)_R (x)
        (16 q)_R (y)
		\sum\limits_{n=0}^{\infty}
        \frac{|z|^n n!^{1 - R}}{2 \cdot 8^n},
	\end{align*}
	where we used the fact that the estimate \eqref{LCAna:CnOperators}
	extends to the completion. For $R > 1$, this map is clearly analytic
	and absolutely convergent for all $z \in \mathbb{K}$.
	If $R = 1$, then for every $M \geq 1$ we go back to homogeneous,
	factorizing tensors $x^{(k)}$ and $y^{(\ell)}$ of degree $k$ and $\ell$
	respectively, and have
	\begin{align*}
		M^{n} p_R \left(
			C_n \left( x^{(k)}, y^{(\ell)} \right)
		\right)
		& \leq
		\frac{M^n}{2 \cdot 8^n}
		(16 q)_R \left( x^{(k)} \right)
		(16 q)_R \left( y^{(\ell)} \right)
		\\
		& \leq
		M^{k + \ell}
		(16 q)_R \left( x^{(k)} \right)
		(16 q)_R \left( y^{(\ell)} \right)
		\\
		& =
		(16M q)_R \left( x^{(k)} \right)
		(16M q)_R \left( y^{(\ell)} \right),
	\end{align*}
	where we used that $0 \leq n \leq k + \ell - 1$. The infimum argument 
	gives the estimate on all tensors $x,y \in \Tensor_R^{\bullet}(\lie{g})$ 
	and it extends to the completion such that
	\begin{equation*}
		p_R \left( z^n C_n(x, y) \right)
		\leq
		(16M q)_R(x) (16M q)_R(y)
		\frac{|z|^n}{2 \cdot (8M)^n}
	\end{equation*}
	and hence
	\begin{equation*}
		p_R(x \star_z y)
		\leq
		(16 M q)_R (x)
        (16 M q)_R (y)
		\sum\limits_{n=0}^{\infty}
        \frac{|z|^n}{2 \cdot (8M)^n}.
	\end{equation*}
	So the power series converges for all $z \in \mathbb{C}$ with
	$|z| < 8M$ and converges uniformly if $|z| \leq c M$ for $c < 8$.
	But then, the map \eqref{LCAna:Holomorphicity} converges on all open discs
	centered around $z = 0$, and it must therefore be entire.
\end{proof}
\begin{remark}[(Weakly) holomorphic maps with values in locally convex spaces]
	\label{LCAna:Rem:Holomorphicity} \hfill \\
	One could argue that the term ``holomorphic'' in a locally convex space 
	$V$ does not necessarily mean that a map has a absolutely convergent 
	Taylor expansion, but one should rather see the map 
	\eqref{LCAna:Holomorphicity} as a collection of paths $\mathbb{C} 
	\longrightarrow \widehat{\Sym}_R^{\bullet}(\lie{g})$ and ask for their 
	complex differentiability in the sense of a differential quotient. Of 
	course, this is also a valid formulation of the word ``holomorphic'' 
	and it is actually even the standard definition for holomorphic maps with 
	values in locally convex spaces. According to this definition one calls a map 
	``weakly holomorphic'', if every continuous linear form 
	$\lambda \colon V \longrightarrow \mathbb{C}$ applied to it gives a 
	holomorphic map $\mathbb{C} \longrightarrow \mathbb{C}$. Using this 
	terminology, we would just have proven the map \eqref{LCAna:Holomorphicity} to 
	be \emph{weakly holomorphic}. Yet, we have also
	proven ``strong'' holomorphicity, since in \cite[Theorem 3.31]{rudin:1991a} 
	Rudin proved that the two notions of holomorphicity coincide in locally convex 
	spaces.
\end{remark}

\subsection{Representations}

We now want to identify the deformed symmetric algebra $(\Sym_R^{\bullet}
(\lie{g}), \star_z)$ with the the universal enveloping algebra $\algebra{U}
(\lie{g}_z)$ with the scaled product.
\begin{definition}[The $\algebra{U}_R$-Topology]
	We denote by $\algebra{U}_R(\lie{g}_z)$ the locally convex algebra which we 
	get from $\algebra{U}(\lie{g}_z)$ by pulling back the $\Sym_R^{\bullet}
	(\lie{g})$-topology via the Poincar\'e-Birkholl-Witt isomophism
	$\mathfrak{q}$.
\end{definition}
We know this algebra has the universal property that Lie algebra 
homomorphisms into associative algebras can be lifted to unital homomorphisms 
of associative algebras. As a commutative diagram, this reads
\begin{center}
    \begin{tikzpicture}
        \matrix (m)[
        matrix of math nodes,
        row sep=6em,
        column sep=7em
        ]
        {
          \algebra{U}_R(\lie{g})
          & \Sym_R^{\bullet}(\lie{g}) \\
          \lie{g}
          & \algebra{A} \\
        };
        \draw
        [-stealth]
        (m-1-1) edge node [above] {$\lie{q}^{-1}$}
        (m-1-2) edge node [left] {$\Phi$}
        (m-2-2)
        (m-1-2) edge node [right] {$\widetilde{\Phi}$}
        (m-2-2)
        (m-2-1) edge node [left] {$\iota$}
        (m-1-1)
        (m-2-1) edge node [below] {$\phi$}
        (m-2-2);
    \end{tikzpicture}
\end{center}
where $\algebra{A}$ is the mentioned associative algebra. Since we endowed 
$\algebra{U}_R(\lie{g}_z)$ and $\Sym_R^{\bullet}(\lie{g})$ with a topology, we 
can now ask if the homomorphisms $\Phi$ and $\widetilde{\Phi}$ are continuous. 
This question is partly answered by the following result:
\begin{proposition}[Universal property]
    \label{LCAna:Prop:Semi-functoriality}%
    Let $\lie{g}$ be an AE-Lie algebra, $\algebra{A}$ an associative
    AE algebra and $\phi \colon \lie{g} \longrightarrow \algebra{A}$
    is a continuous Lie algebra homomorphism.  If $R \geq 1$, then the
    induced algebra homomorphisms $\Phi$ and $\widetilde{\Phi}$ are
    continuous.
\end{proposition}
\begin{proof}
    We define an extension of $\Phi$ on the whole tensor algebra
    again:
    \begin{equation*}
        \Psi \colon
        \Tensor_R^{\bullet}(\lie{g})
        \longrightarrow
        \algebra{A},
        \quad
        \Psi
        =
        \widetilde{\Phi} \circ \Symmetrizer
    \end{equation*}
    It is clear that if $\Psi$ is continuous on factorizing tensors,
    we will get the continuity of $\widetilde{\Phi}$ and of $\Phi$ via the
    infimum argument. So let $p$ be a continuous seminorm on
    $\algebra{A}$ with an asymptotic estimate $q$ and $\xi_1, \ldots,
    \xi_n \in \lie{g}$. Since $\phi$ is continuous, we find a
    continuous seminorm $r$ on $\lie{g}$ such that for all $\xi \in
    \lie{g}$ we have $q(\phi(\xi)) \leq r(\xi)$. Then we have
    \begin{align*}
        p \left(
        \Psi \left(
        \xi_1 \tensor \cdots \tensor \xi_n
        \right) \right)
        & =
        p \left(
        \widetilde{\Phi} \left(
        \xi_1 \star_z \cdots \star_z \xi_n
        \right) \right)
        \\
        & =
        p( \phi(\xi_1) \cdots \phi(\xi_n) )
        \\
        & \leq
        q( \phi(\xi_1) )
        \cdots
        q( \phi(\xi_n) )
        \\
        & \leq
        r(\xi_1) \cdots r(\xi_n)
        \\
        & \leq
        r_R(\xi_1 \tensor \cdots \tensor \xi_n),
    \end{align*}
    where the last inequality is true for all $R \geq 0$ and hence for all 
    $R \geq 1$.
\end{proof}
Although this is a nice result, our construction fails to be a universal object 
in the category of associative AE algebras, since the universal enveloping 
algebra endowed with our topology is \emph{not} AE in general. This is even 
very easy to see directly.
\begin{example}
    Take $\xi \in \lie{g}$, then we know that $\xi^{\tensor n} =
    \xi^{\ostar_z n} = \xi^n$ for $n \in \mathbb{N}$. 
    Let $R > 0$ and $p$ and $q$ be a continuous seminorms on $\lie{g}$ with 
    $q(\xi) \neq 0$, then we find
    \begin{equation}
        p_R(\xi^n)
        =
        n!^R p(\xi)^n
        =
        \frac{n!^R}{c^n} q(\xi)^n
    \end{equation}
    for $c = \frac{p(\xi)}{q(\xi)}$. But since $\frac{n!^R}{c^n}$ always
    diverges for $n \longrightarrow \infty$, we never get an asymptotic
    estimate for $p$.
\end{example}
Another argument would be that we know that we can not have exponential 
functions in $\widehat{\Sym}_R^{\bullet}(\lie{g})$. If 
$\widehat{\Sym}_R^{\bullet}(\lie{g})$ (and hence $\algebra{U}_R(\lie{g}_z)$) 
were AE, then we would have an entire calculus which clearly includes 
exponential functions. Although the construction is not universal, we can draw 
a strong conclusion from Proposition~\ref{LCAna:Prop:Semi-functoriality}:
\begin{corollary}[Continuous Representations]
    \label{LCAna:Coro:ContinuousRepresentations}%
    Let $R \geq 1$ and $\algebra{U}_R(\lie{g})$ the universal
    enveloping algebra of an AE-Lie algebra $\lie{g}$, then for every
    continuous representation $\phi$ of $\lie{g}$ into the bounded
    linear operators $\Bounded(V)$ on a Banach space $V$ the induced
    homomorphism of associative algebras $\Phi \colon
    \mathcal{U}(\lie{g}) \longrightarrow \Bounded(V)$ is continuous.
\end{corollary}
\begin{proof}
    This follows directly from
    Proposition~\ref{LCAna:Prop:Semi-functoriality} and $\Bounded(V)$
    being a Banach algebra.
\end{proof}
\begin{remark}
	From this, we get the special case that all finite-dimensional 
	representations of an AE-Lie algebra can be lifted to continuous 
	representations of $\algebra{U}_R(\lie{g}_z)$. 
	For infinite-dimensional Lie algebras, this statement will not be very 
	important, since there, one rather has strongly continuous 
	representations and no norm-continuous ones. Nevertheless, our 
	topology may help to think about continuous linear functionals on 
	$\algebra{U}_R(\lie{g}_z)$. We could then do GNS-representation 
	theory with it. This procedure would yield a representation of 
	$\algebra{U}_R(\lie{g}_z)$ on a (Pre-)Hilbert space and we could talk 
	about what strongly continuous representations of the universal 
	enveloping algebra should be. In this sense, the results of this 
	chapter may also open a door towards some new approaches in this field.
\end{remark}

\subsection{Functoriality}

Now let $\lie{g}, \lie{h}$ be two AE-Lie algebras. We know that a homomorphism 
of Lie algebras from $\lie{g}$ to $\algebra{U}(\lie{h}_z)$ would lift to a 
homomorphism $\algebra{U}(\lie{g}_z) \longrightarrow \algebra{U}(\lie{h}_z)$, 
if the latter one was AE, but this is not the case. Yet, we would like to have 
this result and prove that our construction is functorial, but it will be 
more complicated than that. We draw the commutative diagram to clarify the 
situation.
\begin{center}
    \begin{tikzpicture}
        \matrix (m)[
        matrix of math nodes,
        row sep=3em,
        column sep=7em
        ]
        {
          \Sym_R^{\bullet}(\lie{g})
          & \Sym_R^{\bullet}(\lie{h}) \\
          \algebra{U}_R(\lie{g}_z)
          & \algebra{U}_R(\lie{h}_z) \\
          \lie{g}
          & \lie{h} \\
        };
        \draw
        [-stealth]
        (m-1-1) edge node [above] {$\widetilde{\Phi}_z$}
        (m-1-2)
        (m-2-1) edge node [left] {$\mathfrak{q}_z^{-1}$}
        (m-1-1)
        (m-2-1) edge node [above] {$\Phi_z$}
        (m-2-2)
        (m-2-2) edge node [right] {$\mathfrak{q}_z^{-1}$}
        (m-1-2)
        (m-3-1)	edge node [below] {$\phi$}
        (m-3-2)
        (m-3-1)	edge node [left] {$\iota_z$}
        (m-2-1)
        (m-3-2)	edge node [right] {$\iota_z$}
        (m-2-2);
    \end{tikzpicture}
\end{center}
Assume that $\phi$ is a continuous Lie algebra homomorphism. We want to know if 
$\Phi_z$ and $\widetilde{\Phi}_z$ will be continuous, too. Luckily, the answer 
is yes and our construction is functorial. For the proof, we need the 
next two lemmas.
\begin{lemma}
    \label{LCAna:Lemma:NStarPrePreSubResult}%
    Let $\lie{g}$ be an AE-Lie algebra, $n \in \mathbb{N}$, $\xi_1, \ldots, \xi_n 
    \in \lie{g}$, $i_j \in \{0, \ldots, j\}$, $\forall_{j = 1, 
    \ldots, n-1}$ and denote $I = \sum_j i_j$. Then we have the formula
	\begin{align}
		\nonumber
		z^{i_{n-1}} 
		&
		C_{i_{n-1}}\left(
			\ldots
			z^{i_2} C_{i_2} \left(
				z^{i_1} c_{i_1} \left(
					\xi_1, \xi_2
				\right)
				,
				\xi_3
			\right)
			\ldots,
			\xi_n
		\right)
		=
		z^I B_{i_{n-1}}^* \cdots B_{i_1}^*
		\\
		& \quad
		\cdot
		\frac{\binom{1}{i_1}}{1!} 
		\frac{\binom{2 -i_1}{i_2}}
		{(2-i_1)!}
		\cdots 
		\frac{\binom{n-1 - i_1 - \cdots - i_{n-2}}{i_{n-1}}}
		{(n-1 - i_1 - \cdots - i_{n-2})!}
		\sum\limits_{\substack{
			\sigma_1 \in S_{2 - i_1} \\
			\ldots\\
			\sigma_{n-1} \in S_{n-1 - i_1 - \ldots - i_{n-2}}			
		}}
		[w_1] \cdots [w_{n-I}],
	\end{align}
	where the expressions $[w_i]$ denote nested Lie brackets in the $\xi_i$.
\end{lemma}
\begin{proof}
	The statement follows from the bilinearity of the $C_{i_j}$-operators. They 
	are always applied to homogeneous symmetric tensors. The proof itself is an 
	easy induction and follows very directly from
	Formula~\eqref{Formulas:LinearMonomial2}.
\end{proof}
\begin{lemma}
    \label{LCAna:Lemma:LemmaPreContinuityN}%
    Let $\lie{g}$ be an AE-Lie algebra, $R \geq 1$, $z \in \mathbb{C}$ and $p$ 
    a continuous seminorm with an asymptotic estimate $q$. Then for every $n 
    \in \mathbb{N}$ and all $\xi_1, \ldots, \xi_n \in \lie{g}$ the estimate
    \begin{equation}
        \label{LCAna:LemmaPreContinuityN}
        p_R \left(
            \xi_1 \star_z \cdots \star_z \xi_n
        \right)
        \leq
        c^n n!^R
        q^n(\xi_1 \tensor \cdots \tensor \xi_n)
    \end{equation}
    holds with $c = 8 \E (|z| + 1)$.
\end{lemma}
\begin{proof}
    We start with a continuous seminorm $p$:
    \begin{align}
        \nonumber
        p_R \big(
            \xi_1 \star_z 
        &
            \cdots \star_z \xi_n
        \big)
        =
        p_R \Bigg(
        \sum\limits_{\ell = 0}^{n-1}
        \sum\limits_{\substack{
			1 \leq j \leq n-1 \\
			i_j \in \{0, \ldots, j\} \\
			\sum_{j = 1}^{n - 1} i_j = \ell
		}}
		z^{i_{n-1}}
		C_{i_{n-1}}
		\left(
			\ldots z^{i_2} C_{i_2}
			\left(
				z^{i_1} C_{i_1}
				\left( \xi_1, \xi_2 \right)
				, \xi_3
			\right)
			\ldots, \xi_n
		\right)
        \Bigg)
        \\
        \nonumber
        & \leq
        \sum\limits_{\ell = 0}^{n-1}
        (n - \ell)!^R
        \sum\limits_{\substack{
			1 \leq j \leq n-1 \\
			i_j \in \{0, \ldots, j\} \\
			\sum_{j = 1}^{n - 1} i_j = \ell
		}}
        p^{n - \ell} \Big(
		z^{i_{n-1}}
		C_{i_{n-1}}
		\left(
			\ldots z^{i_2} C_{i_2}
			\left(
				z^{i_1} C_{i_1}
				\left( \xi_1, \xi_2 \right)
				, \xi_3
			\right)
			\ldots, \xi_n
		\right)
		\Big)
        \\
        \nonumber
        & \leq
        \sum\limits_{\ell = 0}^{n-1}
        (n - \ell)!^R
        \sum\limits_{\substack{
			1 \leq j \leq n-1 \\
			i_j \in \{0, \ldots, j\} \\
			\sum_{j = 1}^{n - 1} i_j = \ell
		}}
        |B_{i_{n-1}}^*| \cdots |B_{i_1}^*|
        |z|^{i_{n-1}} \cdots |z|^{i_1}
        \frac{\binom{1}{i_1}}{1!}
        \frac{\binom{2 - i_1}{i_2}}{(2-i_1)!}
        \cdots
        \\
        \label{LCAna:PreContinuityIntermediateN}
        & \quad \cdot
        \frac{\binom{n - 1 - i_1 - \cdots i_{n-2}}{i_{n-1}}}
        {(n - 1 - i_1 - \cdots - i_{n-2})!}
		\sum\limits_{\substack{
			\sigma_1 \in S_{2 - i_1} \\
			\ldots\\
			\sigma_{n-1} \in S_{n-1 - i_1 - \ldots - i_{n-2}}			
		}}       
        p^{n - \ell}
        \left( 
         	[w_1] \cdots [w_{n - \ell}]
        	\right).
    \end{align}
    In the last step, we used Lemma \ref{LCAna:Lemma:NStarPrePreSubResult}. Now
    we can apply the AE property to it and get $q(\xi_1) \cdots q(\xi_n)$. Thus 
    the sum over the $\sigma_j$ vanishes with the inverse factorials.
    By using the fact that $|B_m^*| \leq m!$ for all $m \in
    \mathbb{N}$ and grouping together the powers of $|z|$, we find
    \begin{align*}
        &p_R \left(
            \xi_1 \star_z \cdots \star_z \xi_n
        \right) \\
        & \quad \ot{(a)}{\leq}
        \sum\limits_{\ell = 0}^{n-1}
        (n - \ell)!^R
        \sum\limits_{\substack{
			1 \leq j \leq n-1 \\
			i_j \in \{0, \ldots, j\} \\
			\sum_{j = 1}^{n - 1} i_j = \ell
		}}
        |z|^{\ell}
        \frac{1!  (2 - i_1)! \cdots (n-1 - i_1 - \cdots - i_{n-2})!}
        {(1 - i_1)! \cdots (n-1 - i_1 - \cdots i_{n-1})!}
        q(\xi_1) \cdots q(\xi_n)
        \\
        &\quad\ot{(b)}{\leq}
        \sum\limits_{\ell = 0}^{n-1}
        (n - \ell)!^R
        \sum\limits_{\substack{
			1 \leq j \leq n-1 \\
			i_j \in \{0, \ldots, j\} \\
			\sum_{j = 1}^{n - 1} i_j = \ell
		}}
        |z|^{\ell}
        1^{i_1} 2^{i_2}
        \cdots (n-1)^{i_{n-1}}
        q(\xi_1) \cdots q(\xi_n)
        \\
        &\quad\ot{(c)}{\leq}
        \sum\limits_{\ell = 0}^{n-1}
        (n - \ell)!^R
        \sum\limits_{\substack{
			1 \leq j \leq n-1 \\
			i_j \in \{0, \ldots, j\} \\
			\sum_{j = 1}^{n - 1} i_j = \ell
		}}
        |z|^{\ell}
        n^{\ell}
        q(\xi_1) \cdots q(\xi_n)
        \\
        &\quad\ot{(d)}{\leq}
        \sum\limits_{\ell = 0}^{n-1}
        (n - \ell)!^R
        \sum\limits_{\substack{
			1 \leq j \leq n-1 \\
			i_j \in \{0, \ldots, j\} \\
			\sum_{j = 1}^{n - 1} i_j = \ell
		}}
        |z|^{\ell} (2 \E)^n \ell!
        q(\xi_1) \cdots q(\xi_n).
    \end{align*}
    In (a), the factorials coming from the Bernoulli numbers cancelled out with 
    the the $i_j$ from the binomial coefficients and in (b) we used
    \begin{equation*}
    		\frac{(k - i_1 - \cdots - i_{k-1})!}
    		{(k - i_1 - \cdots - i_k)!}
    		\leq
    		k^{i_k}.
    \end{equation*}
    We can estimate the product of all those expressions by $(n-1)^{i_1 + 
    \cdots + i_{n-1}}$ and this by $n^\ell$, since the sum over all $i_j$ 
    equals $\ell$. In the last step (d) we used $n^{\ell} \leq \E^n
    \frac{n!}{(n-\ell)!} = \E^n \binom{n} {\ell} \ell! \leq \E^n 2^n
    \ell!$. But now we can simply estimate $|z|^{\ell} \leq (|z| +
    1)^n$ and $(n - \ell)!^R \ell! \leq n!^R$ for $R \geq 1$. We just
    need to count the number of summands and get
    \begin{align*}
        p_R \left(
            \xi_1 \star_z \cdots \star_z \xi_n
        \right)
        & \leq
        n!^R (2 \E)^n (|z| + 1)^n
        q(\xi_1) \cdots q(\xi_n)
        \sum\limits_{\ell = 0}^{n-1}
        \sum\limits_{\substack{
			1 \leq j \leq n-1 \\
			i_j \in \{0, \ldots, j\} \\
			\sum_{j = 1}^{n - 1} i_j = \ell
		}}
		1
        \\
        &\ot{(a)}{\leq}
        n!^R (2 \E)^n (|z| + 1)^n
        q(\xi_1) \cdots q(\xi_n)
        \sum\limits_{\ell = 0}^{n-1}
        \binom{n - 1 + \ell - 1}{\ell - 1}
        \\
        &\ot{(b)}{\leq}
        n!^R (2 \E)^n (|z| + 1)^n
        q(\xi_1) \cdots q(\xi_n)
        2^{2n}
        \\
        &\leq
        c^n n!^R q^n(\xi_1 \tensor \cdots \tensor \xi_n),
    \end{align*}
    with $c = 8 \E (|z| + 1)$. In ($a$) the estimate for the big sum
    is the following: for every $j = 1, \ldots, n-1$ we have surely
    $i_j \in \{0, 1, \ldots, n-1\}$ and the sum of all the $i_j$ is
    $\ell$. If we forget about all other restrictions, we will just get
    more terms. But then the number of summands is same as there
    are ways to distribute $\ell$ items on $n-1$ places, which is
    given by $\binom{n - 1 + \ell - 1}{\ell - 1}$. Then in ($b$) we
    use
    \begin{equation*}
        \binom{n - 1 + \ell - 1}{\ell - 1} \leq \binom{2 n}{\ell - 1}
    \end{equation*}
    with the binomial coefficient being zero for $\ell = 0$. Then it
    is just the standard estimate for binomial coefficients via the
    sum over all $\ell$.
\end{proof}
\begin{proposition}[Functoriality]
	\label{LCAna:Prop:Functoriality}
	Let $R \geq 1$, $\lie{g}, \lie{h}$ be AE-Lie algebras and $\phi \colon
	\lie{g} \longrightarrow \lie{h}$ a continuous homomorphism between them.
	Then it lifts to a continuous unital homomorphism of locally convex
	algebras $\Phi_z \colon \algebra{U}_R(\lie{g}_z) \longrightarrow
	\algebra{U}_R(\lie{h}_z)$.
\end{proposition}
\begin{proof}
	First, if $\phi \colon \lie{g} \longrightarrow \lie{h}$ is continuous, 
	then for every continuous seminorm $q$ on $\lie{h}$, we have a continuous 
	seminorm $r$ on $\lie{g}$ such that for all $\xi \in \lie{g}$
	\begin{equation*}
		q\left( \phi(\xi) \right)
		\leq
		r(\xi).
	\end{equation*}
	Second, we define $\Psi_z$ on factorizing tensors via
	\begin{equation*}
		\Psi_z \colon
		\Tensor_R^{\bullet}(\lie{g})
		\longrightarrow
		\Sym_R^{\bullet}(\lie{h})
		, \quad
		\Psi_z
		=
		\widetilde{\Phi}_z \circ
		\Symmetrizer
	\end{equation*}
	and extend it linearly to $\Tensor_R^{\bullet}(\lie{g})$. Clearly, 
	$\Phi_z$ and $\widetilde{\Phi}_z$ will be continuous if $\Psi_z$ is 
	continuous. From this, we get for a seminorm $p$ on $\lie{h}$, an 
	asymptotic estimate $q$ and $\xi_1, \ldots, \xi_n$
	\begin{align*}
		p_R\left(
			\Psi_z \left(
				\xi_1 \tensor \cdots \tensor \xi_n
			\right)
		\right)
		& =
		p_R \left(
			\phi \left( \xi_1 \right)
			\star_z \cdots \star_z
			\phi \left( \xi_1 \right)
		\right)
		\\
		& \leq
		c^n n!^R
		q \left( \phi \left( \xi_1 \right) \right)
		\cdots
		q \left( \phi \left( \xi_n \right) \right)
		\\
		& \leq
		c^n n!^R
		r \left( \xi_1 \right)
		\cdots
		r \left( \xi_n \right)
		\\
		& =
		(c r)_R\left( 
			\xi_1 \tensor \ldots \tensor \xi_n
		\right).
	\end{align*}
	Again, we use the infimum argument and we have the estimate on all tensor 
	in $\Tensor_R^{\bullet}(\lie{g})$. It extends to the completion and the 
	statement is proven.
\end{proof}

\section{Alternative Topologies and an Optimal Result}
\label{sec:chap5_Optimality}

In this section we set the formal parameter $z = 1$ and make some 
observations. So far, we found a topology on $\Sym_R^{\bullet}(\lie{g})$ which 
gives a continuous star product and which has a reasonably large completion, 
but it is always fair to ask if we can do better than that: we have seen that 
our completed algebra does not contain exponential series, which would be a 
very good feature to have. So is it possible to put another locally convex 
topology on $\Sym_R^{\bullet}(\lie{g})$ which gives a completion with 
exponentials? The answer is no, at least under mild additional assumptions. 
\begin{proposition}[Optimality of the $\Tensor_R$-topology]
	\label{LCAna:Prop:NoBetterTopology}
	Let $\lie{g}$ be an AE Lie algebra in which one has elements $\xi, \eta$ 
	for which the Baker-Campbell-Hausdorff series does not converge. 
	Then there is no locally convex topology on $\Sym^{\bullet}(\lie{g})$ 
	such that all of the following things are fulfilled:
	\begin{propositionlist}
		\item
		The Gutt star product $\star_G$ is continuous.
		\item
		For every $\xi \in \lie{g}$ the series $\exp(\xi)$ converges 
		absolutely in the completion of $\Sym^{\bullet}(\lie{g})$.
		\item
		For all $n \in \mathbb{N}$ the projection and inclusion maps with 
		respect to the graded structure
		\begin{equation*}
			\Sym^{\bullet}(\lie{g})
	    		\ot{$\pi_n$}{\longrightarrow}
    	    		\Sym^n(\lie{g})
	    	    	\ot{$\iota_n$}{\longrightarrow}
	    		\Sym^{\bullet}(\lie{g})
		\end{equation*}
		are continuous.
	\end{propositionlist}
\end{proposition}
First of all, we should make clear what ``the Baker-Campbell-Hausdorff series 
does not converge'' actually means. This may be clear for a finite-dimensional 
Lie algebra, but in the locally convex setting, it is not that obvious. For 
simplicity, we assume our local convex space to be complete in the 
following, since we can always achieve this by completing it. First, we note 
here that a net or a sequence in a locally convex space is  convergent [or 
Cauchy], if and only if it is convergent [or Cauchy] with respect 
to all $p \in \algebra{P}$. Quite similar to a normed space, we can make the 
following definition.
\begin{definition}
	\label{Def:RadiusOfConvergenceLCS}
	Let $V$ be a locally convex vector space, $\algebra{P}$ the set of 
	continuous seminorms, $p \in \algebra{P}$ and 
	$\alpha = (\alpha_n)_{n \in \mathbb{N}}\subseteq V$ a sequence in $V$. 
	We set
	\begin{equation*}
		\varrho_p(\alpha)
		=
		\left(
			\limsup_{n \longrightarrow \infty}
			\sqrt[n]{p \left(\alpha_n \right)}
		\right)^{-1}
	\end{equation*}
	where $\varrho_p(\alpha) = \infty$ if $\limsup_{n \longrightarrow \infty} 
	\sqrt[n]{p \left(\alpha_n \right)} = 0$ as usual.
\end{definition}
From this, we immediately get the two following lemmas.
\begin{lemma}[Root test in locally convex spaces]
	\label{LCAna:Lemma:RootTest}
	Let $V$ be a complete locally convex vector space, $p \in \algebra{P}$ 
	and $\alpha = (\alpha_n)_{n \in \mathbb{N}} \subseteq V$ a sequence. 
	Then, if $\varrho_p(\alpha) > 1$, the series
	\begin{equation*}
		\mathcal{S}_n(\alpha)
		=
		\sum\limits_{j = 0}^n
		\alpha_j
	\end{equation*}
	converges absolutely with respect to $p$. Conversely, if this
	series converges with respect to $p$, then we have 
	$\varrho_p(\alpha) \geq 1$.	
\end{lemma}
\begin{proof}
	The proof is completely analogous to the one in finite dimensions.
\end{proof}
\begin{lemma}
	\label{LCAna:Lemma:PowerSeriesConvAbs}
	Let $V$ be a complete locally convex vector space, $p \in \algebra{P}$, 
	$\alpha = (\alpha_n)_{n \in \mathbb{N}} \subseteq V$ a sequence and
	$M > 0$. Then, if the power series
	\begin{equation*}
		\lim_{n \longrightarrow \infty}
		\sum\limits_{j = 0}^n
		\alpha_j z^j
	\end{equation*}
	converges with respect to $p$ for all $z \in \mathbb{C}$ with $|z| \leq M$, 
	it converges absolutely with respect to $p$ for all $z \in \mathbb{C}$ with 
	$|z| < M$.
\end{lemma}
\begin{proof}
	Like in the finite-dimensional setting, we use the root test: convergence 
	for $|z| \leq M$ means $\varrho_p(\alpha_z) \geq 1$, where we have set 
	$\alpha_z = (\alpha_n z^n)_{n \in \mathbb{N}}$. Hence, for every $z' < z$ 
	we get $\varrho_p(\alpha_{z'}) > 1$ and absolute convergence by Lemma 
	\ref{LCAna:Lemma:RootTest}.
\end{proof}
But this means that a convergent power series $\alpha_z$ in a locally convex space  
has $\varrho_p(\alpha_z) = \infty$ for all $p \in \algebra{P}$: if there was a $p 
\in \algebra{P}$ with a finite radius of convergence, we could take a multiple 
of $p$ in order to get a seminorm with arbitrarily small radius of convergence. 
Hence $\alpha_z$ could not converge for any $z \neq 0$. This helps us to make 
clear, how ``BCH does not converge'' can be interpreted, since we can read it 
as a power series in two variables:
\begin{equation*}
	\bch{t \xi}{s \eta}
	=
	\sum\limits_{a,b = 0}^{\infty}
	t^a s^b
	\bchparts{a}{b}{\xi}{\eta}.
\end{equation*}
So if it converges for $\xi$ and $\eta$, then it converges absolutely for all 
$t \xi$ and $s \eta$ with $t,s \in [0, 1)$. If it even converges for all $t 
\xi$ and $s \eta$ with $t, s \in \mathbb{K}$, then it converges absolutely for 
all $t, s$. If this is the case, we can reorder the sums as we like to, and the 
series
\begin{equation*}
	\sum\limits_{n = 0}^N
	\bchpart{n}{\xi}{\eta}
\end{equation*}
converges if and only if the series
\begin{equation*}
	\bch{t \xi}{s \eta}
	=
	\sum\limits_{a = 0}^{\infty}
	\sum\limits_{b = 0}^{\infty}
	\bchparts{a}{b}{\xi}{\eta}
\end{equation*}
converges. Now we can prove a Lemma, from which Proposition 
\ref{LCAna:Prop:NoBetterTopology} follows immediately.
\begin{lemma}
	Let $\lie{g}$ be an AE Lie algebra and $\Sym^{\bullet}(\lie{g})$ is 
	endowed with a locally convex topology, such that the conditions $(i) - 
	(iii)$ from Proposition~\ref{LCAna:Prop:NoBetterTopology} are fulfilled. 
	Then the Baker-Campbell-Hausdorff series converges absolutely for all 
	$\xi, \eta \in \lie{g}$.
\end{lemma}
\begin{proof}
	First, we complete the algebra to $\widehat{\Sym}^{\bullet}(\lie{g})$.
	We need the projection $\pi_1$ to the Lie algebra. Take 
	$\xi, \eta \in 	\lie{g}$. Now, since the the Gutt star product is 
	continuous and the exponential series is absolutely convergent for 
	all $\xi, \eta \in \lie{g}$, we get for $t, s \in \mathbb{K}$
	\begin{align*}
		\pi_1 \left( \exp(t \xi) \star \exp(s \eta) \right)
		& =
		\pi_1
		\left(
			\lim_{N \rightarrow \infty}
			\left(
				\sum\limits_{n=0}^N
				\frac{t^n \xi^n}{n!}
			\right)
			\star
			\lim_{M \rightarrow \infty}
			\left(
				\sum\limits_{m=0}^M
				\frac{s^m \eta^m}{m!}
			\right)
		\right)
		\\
		& \ot{(a)}{=}
		\pi_1
		\left(
			\lim_{N \rightarrow \infty}
			\lim_{M \rightarrow \infty}
			\left(
				\sum\limits_{n=0}^N
				\frac{t^n \xi^n}{n!}
			\right)
			\star_G
			\left(
				\sum\limits_{m=0}^M
				\frac{s^m \eta^m}{m!}
			\right)
		\right)
		\\
		& \ot{(b)}{=}
		\lim_{N \rightarrow \infty}
		\lim_{M \rightarrow \infty}
		\pi_1
		\left(	
			\left(
				\sum\limits_{n=0}^N
				\frac{t^n \xi^n}{n!}
			\right)
			\star_G
			\left(
				\sum\limits_{m=0}^M
				\frac{s^m \eta^m}{m!}
			\right)
		\right)
		\\
		& \ot{(c)}{=}
		\lim_{N \rightarrow \infty}
		\lim_{M \rightarrow \infty}
		\sum\limits_{n,m=0}^{N, M}
		\bchparts{t \xi}{s \eta}{n}{m}
	\end{align*}
	where we used the continuity of the star product in (a), the continuity of 
	the projection in (b) and evaluated the projection in (c). Since 
	$\exp(t \xi)$ and $\exp(s \eta)$ are elements in the completion, their 
	star product exists and hence the double series at the end of this 
	equation converges for any two elements $\xi, \eta \in \lie{g}$. But now 
	we can use the result that in this setting the BCH series converges 
	absolutely. We can rearrange the terms and get the convergence of
	\begin{equation*}
		\sum\limits_{n = 1}^N
		\bchpart{n}{t \xi}{s \eta}.
	\end{equation*}
	Therefore, the BCH series must converge globally.
\end{proof}
Obviously, this proves Proposition~\ref{LCAna:Prop:NoBetterTopology}.

Now we want to use a result due to Wojty\'nski \cite{wojtynski:1998a}, who 
showed that for Banach-Lie algebras, global convergence of the BCH series is 
equivalent to the fact that for any $\xi$ we have
\begin{equation*}
	\lim_{n \longrightarrow \infty}
	\norm{ 
		\left( \ad_{\xi} \right)^n 
	}^{\frac{1}{n}}
	=
	0.
\end{equation*}
A Banach-Lie algebra with this property is sometimes called quasi-nilpotent, 
radical or nil, see for example \cite{mueller:1994a} for various 
generalizations of nilpotency in the case of Banach algebras.
For finite-dimensional Lie algebras, quasi-nilpotency implies nilpotency. 
Hence for a finite-dimensional Lie algebra $\lie{g}$, BCH is globally 
convergent if and only if $\lie{g}$ is nilpotent.

From this, we see that at least for ``non-quasi-nilpotent'' Banach-Lie 
algebras, our result is in some sense optimal, at least if we want the grading 
structure to compatible with the topology.
\begin{remark}[Another topology in $\algebra{U}(\lie{g})$]
    \label{Rem:LCAnaBCHConvergence}
    In \cite{pflaum.schottenloher:1998a} Schottenloher and Pflaum mention an 
    alternative topology on the universal enveloping algebra for 
    finite-dimensional Lie algebras: they take the coarsest locally convex 
    topology, such that all finite-dimensional representations of $\lie{g}$ 
    extend to continuous algebra homomorphisms. This topology is in fact even 
    locally m-convex and has therefore an entire holomorphic calculus. In 
    particular, this completion contains exponential functions for all Lie 
    algebra elements. Therefore, as we have seen in 
    Proposition~\ref{LCAna:Prop:NoBetterTopology}, it can not respect the 
    grading structure, as our topology does. The $\Tensor_R$-topology must 
    hence be different from that. As we have seen in 
    Corollary~\ref{LCAna:Coro:ContinuousRepresentations}, it is finer for 
    $R \geq 1$, and since the topologies are different, it is even strictly 
    finer. One could argue that the $R$-topology is ``just'' locally convex, 
    but its advantage (for our purpose) is that the grading is necessary for 
    the holomorphic dependence on the formal parameter which is a feature that 
    we want.
\end{remark}

% Chapter 6
%

%
% Chapter 6 of my master thesis:
% The nilpotent case
%

\chapter{Nilpotent Lie Algebras}

At the end of the last chapter, we have seen that the Baker-Campbell-Hausdorff 
series and its convergence play an important role for a topology on the 
universal enveloping algebra. Thus it is natural to ask whether things will 
change, if we restrict our observations to Lie algebras with globally convergent 
BCH series. To make things not too complicated from the beginning, we focus on 
locally convex and truly nilpotent Lie algebras. Recall that a Lie algebra 
$\lie{g}$ is nilpotent, if there exists a $N \in \mathbb{N}$, such that for 
all $n > N$ and all $\xi_1, \ldots, \xi_n \in \lie{g}$ we have
\begin{equation}
	\label{Nilpot:StrongNilpotency}
	\ad_{\xi_1} \circ \cdots \circ \ad_{\xi_n}
	=
	0.
\end{equation}
In the infinite-dimensional case, this is a priori \emph{not} the same as 
\begin{equation}
	\label{Nilpot:WeakNilpotency}
	\left( \ad_{\xi} \right)^n
	=
	0
\end{equation}
for all $\xi \in \lie{g}$ and and $n > N$, but something
stronger. In the case of finite-dimensional Lie algebras, the notions 
\eqref{Nilpot:StrongNilpotency} and \eqref{Nilpot:WeakNilpotency} coincide due 
to the Engel theorem, which makes use of the existence of a 
finite descending series of nilpotent ideals in the Lie algebra. Such a 
terminating series does not need to exist in infinite dimensions. At least as 
soon as we are in the setting of Banach-Lie algebras, these two notions  
coincide again, but note that one can give different forms of 
\emph{quasi-nilpotency}, which weaken the statements from 
\eqref{Nilpot:StrongNilpotency} and \eqref{Nilpot:WeakNilpotency}, respectively. 
Those generalized notions do not coincide any more, so we have to be careful.

Before we look at this case more closely, we first want to make a list of things, 
which we expect to change or not when we go to this more particular setting.
\begin{enumerate}
	\item
	In Example~\ref{LCAna:Ex:HeisenbergAlgebra} we have seen that we can not 
	expect to get a continuous algebra structure for $R < 1$, even for very 
	simple nilpotent, but non-abelian Lie algebras. Therefore, we should not 
	expect to get much larger completions now.
	
	\item
	In \cite{waldmann:2014a}, Waldmann showed that for example the Weyl-Moyal 
	star product converges in the $\Tensor_R$-topology for $R \geq \frac{1}{2}$. 
	This so-called 	Weyl algebra is, however, nothing but a quotient of the 
	Heisenberg 	algebra. It would be interesting to understand this a bit better, 
	since we know that we need $R \geq 1$ for the Heisenberg algebra. The quotient 
	procedure must 
	therefore have some strong influence on this construction. Can we 
	reproduce the value $R \geq \frac{1}{2}$ somehow by dividing out an ideal?
	
	\item
	The argument we used in Proposition~\ref{LCAna:Prop:NoBetterTopology},
	namely the non-global convergence of BCH, is not given any more. Now there is 
	no longer a reason to expect that exponentials are not part 
	of the completion. In this sense, it would be at least nice to have 
	something more than ``just'' $R = 1$. Can we do that?
	
	\item
	As already mentioned, there are generalizations or weaker forms of 
	nilpotency in infinite-dimensions, especially for Banach-Lie algebras, 
	which are equivalent to the usual notion of nilpotency in finite 
	dimensions. If we get a stronger result for nilpotent Lie algebras, will it be 
	possible to extend it to some of these generalizations?
\end{enumerate}
The very fascinating and highly interesting answer to the three questions from 
$(ii) - (iv)$ is: yes, we can. The first section of this chapter will be devoted 
to the question from point $(iii)$: we get a bigger completion by using a 
projective limit. We will also see how to get again the good functorial 
properties we had before. The next section treats another phenomenon, which was 
not there in the generic case: we can observe bimodule-structures within 
$\Sym_R^{\bullet}(\lie{g})$. In third section, we will reproduce one of the 
results of Waldmann's, at least for the finite-dimensional case. The fourth part 
will take care of some generalizations of nilpotentcy for Banach-Lie algebras and 
will extend the result of the projective limit to a particular subcase there.

\section{The Projective Limit} 
\label{sec:chap6_ProjLim} 

\subsection{Continuity of the Product}

As already mentioned, it is possible to extend the continuity result. 
Therefore, we take a locally convex, nilpotent Lie algebra $\lie{g}$ and look at
\begin{equation*}
	\Sym_{1^-}^{\bullet}(\lie{g})
	=
	\projlim_{\epsilon \longrightarrow 0}
	\Sym_{1 - \epsilon}^{\bullet}(\lie{g}).
\end{equation*}
A tensor will be in the completed vector space $\widehat{\Sym}_{1^-}^{\bullet}
(\lie{g})$, if and only if it lies for every $\epsilon > 0$ in the completion of 
$\Sym_{1 - \epsilon}^{\bullet}(\lie{g})$. Otherwise stated: let $\algebra{P}$ 
be the set of all continuous seminorms of the Lie algebra $\lie{g}$, then 
\begin{equation*}
	f \in \widehat{\Sym}_{1^-}^{\bullet}(\lie{g})
	\quad \Longleftrightarrow \quad
	p_{1 - \epsilon} (f) 
	< 
	\infty
	\quad
	\forall_{p \in \algebra{P}}
	\forall_{\epsilon > 0}.
\end{equation*}
So, if we want to show, that the Gutt star product is continuous on 
$\widehat{\Sym}_{1^-}^{\bullet}(\lie{g})$, we need to show that for every 
$p \in \algebra{P}$ and  $R < 1$, there exists a $q \in \algebra{P}$ and a 
$R'$ with $R \leq R' < 1$, such that we have for all $x,y \in \Sym^{\bullet}
(\lie{g})$
\begin{equation*}
	p_R \left(
		x \star_z y
	\right)
	\leq
	q_{R'}(x)
	q_{R'}(y).
\end{equation*}
Before we prove the next theorem, we want to remind that locally convex, 
nilpotent Lie algebras are always AE Lie algebras. So the results we have 
found so far are valid in this case, too.
\begin{theorem}
    \label{Nilpot:Thm:ProjLimit}%
    Let $\lie{g}$ be a nilpotent locally convex Lie algebra with
    continuous Lie bracket and $N \in \mathbb{N}$ such that $N + 1$
    Lie brackets nested into each other vanish.
    \begin{theoremlist}
	    	\item \label{item:Nilpot:CnOperators}
	    	If $0 \leq R < 1$, the $C_n$-operators are continuous and fulfil the 
	    	estimate
	    	\begin{equation}
	    		\label{eq:Nilpot:CnOperators}
	    		p_R \left(
	    			C_n (x, y)
	    		\right)
	    		\leq
	    		\frac{1}{2 \cdot 8^n}
	    		(32 \E q)_{R + \epsilon}(x)
	    		(32 \E q)_{R + \epsilon}(y),
	    	\end{equation}
	    	for all $x, y \in \Sym_R^{\bullet}(\lie{g})$, where $p$ is a 
	    	continuous seminorm, $q$ an asymptotic estimate and 
	    	$\epsilon = \frac{N - 1}{N}(1 - R)$.
	    	
	    	\item \label{item_Nilpot:SEinsMinus}
	    	The Gutt star product $\star_z$ is continuous for the locally convex 
	    	projective limit $\Sym_{1^-}^\bullet(\lie{g})$ and we have
	    	\begin{equation}
	    		\label{eq:Nilpot:Continuity}
	    		p_R \left(
	    			x \star_z y
	    		\right)
	    		\leq
	    		(c q)_{R + \epsilon}(x)
	    		(c q)_{R + \epsilon}(y)
	    	\end{equation}
	    	with $c = 32 \E (|z| +1)$ and the $\epsilon$ from the first part.
	    	The Gutt star product extends continuously to
	    	$\widehat{\Sym}_R^{\bullet}(\lie{g})$, where it converges absolutely 
	    	and coincides with the formal series.
    \end{theoremlist}
\end{theorem}
\begin{proof}
    We use again $\star_z$ on the whole tensor algebra and compute the 
    estimate for $\xi_1 \tensor \cdots \tensor \xi_k$ and $\eta_1 \tensor \cdots 
    \tensor \eta_\ell$. The important point is that now, we get restrictions for 
    the values of $n$. Recall that $k + \ell - n$ is the symmetric degree of 
    $C_n \left( \xi_1 \tensor \cdots \tensor \xi_k, \eta_1 \tensor \cdots 
    \tensor \eta_\ell \right)$ and that we can have at most $N$ letters in one
    symmetric factor. This means
    \begin{equation*}
	    	(k + \ell - n) N
    		\geq
    		k + \ell
    		\quad 
    		\Longleftrightarrow 
    		\quad
    		n 
    		\leq
    		(k + \ell)
    		\frac{N - 1}{N}.
    \end{equation*}
    Hence we can estimate $n!^{1-R}$ in Equation \eqref{LCAna:CnOperators}: set 
    $\delta = \frac{N - 1}{N}$ and also denote a factorial where we have 
    non-integers, meaning the gamma function. We get
    \begin{align*}
        n!^{1-R}
        & \leq
        (\delta (k + \ell)!)^{1 - R}
        \\
        & \leq
        (\delta (k + \ell))^{(1 - R) \delta (k + \ell)}
        \\
        & \leq
        (k + \ell)^{(1 - R) \delta (k + \ell)}
        \\
        & =
        \left(
            (k + \ell)^{(k + \ell)}
        \right)^{(1 - R) \delta}
        \\
        & \leq
        \left(
            \E^{k + \ell} 2^{k + \ell} k! \ell!
        \right)^{(1-R) \delta}
        \\
        & =
        \left( (2 \E)^{\delta (1-R)} \right)^{k + \ell}
        k!^{\epsilon} \ell!^{\epsilon},
    \end{align*}
    using $\epsilon = \delta (1 - R)$. Hence
    \begin{align*}
        p_R \big(
        	C_n \big(
        		\xi_1 \tensor \cdots \tensor \xi_k, 
        &
        		\eta_1 \tensor \cdots \tensor \eta_\ell
        	\big)
        \big)
        \\
        & \leq
        \frac{
        	\left(
        		(2 \E)^{\delta (1 - R)}
        	\right)^{k + \ell}
        	k!^{\epsilon}
        	\ell!^{\epsilon}
        }{2 \cdot 8^n}
        (16 q)_R \left( \xi_1 \tensor \cdots \tensor \xi_k \right)
        (16 q)_R \left( \eta_1 \tensor \cdots \tensor \eta_\ell \right)
        \\
        & \leq
        \frac{1}{2 \cdot 8^n}
        (c q)_{R + \epsilon} 
        \left( \xi_1 \tensor \cdots \tensor \xi_k \right)
        (c q)_{R + \epsilon} 
        \left( \eta_1 \tensor \cdots \tensor \eta_\ell \right)
    \end{align*}
    with $c = 16 (2 \E)^{\delta (1 - R)} \leq 32 \E$.
    We then get the estimate on all tensors by the infimum argument and extend 
    it to the completion. Note, that for every $R < 1$ we also have 
    $R + \epsilon < 1$ with the $\epsilon = \delta(1-R)$ from above. Iterating 
    this continuity estimate, we get closer and closer to $1$ and it is not 
    possible to repeat this process an arbitrary number of times and stop at 
    some value strictly less than $1$. For the second part, we can conclude 
    analogously to the second part of Theorem~\ref{Thm:LCAna:Continuity1}.
\end{proof}
We have proven one of the four statements. This projective limit is interesting, 
because it has a bigger completion than just $R = 1$. For example, we get the 
following result.
\begin{corollary}
    \label{corollary:NilpotentCase}%
    Let $\lie{g}$ be a nilpotent, locally convex Lie algebra.
    \begin{corollarylist}
    \item \label{item:NilpotentHasExp} 
    	Let $\exp(\xi)$ be the
        exponential series for $\xi \in \lie{g}$, then we have $\exp(t
        \xi) \in \widehat{\Sym}_{1^-}^\bullet(\lie{g})$ for all $t
        \in \mathbb{K}$.
    \item \label{item:NilpotentExpGivesBCH} 
    	For $\xi, \eta \in
        \lie{g}$ and $z \neq 0$ we have $\exp(\xi) \ostar_z
        \exp(\eta) = \exp \left(\frac 1 z \bch{z \xi}{z \eta}
        \right)$.
    \item \label{item:NipotentOneParameterGroups}
    	For $s,t \in
        \mathbb{K}$ and $\xi \in \lie{g}$ we have $\exp(t \xi)
        \ostar_z \exp(s \xi) = \exp ((t + s) \xi)$.
    \end{corollarylist}
\end{corollary}
\begin{proof}
    For the first part, recall that the completion of the projective
    limit $1^-$ contains all those series $(a_n)_{n \in
      \mathbb{N}_0}$ such that
    \begin{equation*}
        \sum\limits_{n=0}^{\infty}
        a_n n!^{1 - \epsilon} c^n
        <
        \infty
    \end{equation*}
    for all $c > 0$.  This is the case for the exponential series of
    $t \xi$ for $t \in \mathbb{K}$ and $\xi \in \lie{g}$. The second
    part follows from the fact that all the projections $\pi_n$ onto
    the homogeneous subspaces $\Sym_{\pi}^n$ are continuous. The third
    part is then a direct consequence of the second.
\end{proof}
So the exponential series is what we want it to be, somehow. It generates a one 
parameter group. Recall that we can just exponentiate \emph{vectors}. If we wanted 
to exponentiate a quadratic tensor, this would yield something like a Gaussian, 
which is again \emph{not} part of the completion.

In the general case, we could find easier proof by assuming submultiplicativity of 
the seminorms. This is again the case for nilpotent Lie algebras. We get something 
like an alternative version of Lemma~\ref{LCAna:Lemma:PreContinuity2}.
\begin{lemma}
	\label{Nilpot:Lemma:PreContinuity2}
	Let $\lie{g}$ be a locally m-convex, nilpotent Lie algebra such that more than 
	$N$ Lie brackets nested into each other vanish. Let $p$ be a continuous 
	seminorm, $z \in \mathbb{K}$ and $R \geq 0$. Then, for every tensor 
	$x \in \Sym_R^{\bullet}(\lie{g})$ of degree at most $k \in \mathbb{N}$ and 
	$\eta \in \lie{g}$, we have the estimate
	\begin{equation}
		\label{Nilpot:PreContinuity}
		p_R \left( x \star_z \eta \right)
		\leq
		c (k + 1)^R k^{N (1-R)}
		p_R (x) p(\eta)
	\end{equation}
	with the constant $c = \sum_{n = 0}^N \frac{|B_n^*|}{n!} |z|^n$.
\end{lemma}
\begin{proof}
	We do the estimate on factorizing tensors and apply the infimum 
	argument later. So let $\xi_1, \ldots, \xi_k, \eta \in \lie{g}$, $k \in 
	\mathbb{N}$, $p$ a continuous seminorm on $\lie{g}$ and $z \in \mathbb{K}$. 
	Then, we have for $R \geq 0$
	\begin{align*}
		&
		p_R 
		\big(
			\xi_1 \tensor \cdots \tensor \xi_k \star_z \eta
		\big)
		\\
		& =
		\sum\limits_{n = 0}^k
		(k + 1 - n)!^R \binom{k}{n}
		|B_n^*| |z|^n
		p^{k + 1 - n} \left(
			\frac{1}{k!}
			\sum\limits_{\sigma \in S_k}
			\xi_{\sigma(1)} \cdots \xi_{\sigma_{k-n}}
			\left( 
				\ad_{\xi_{\sigma(k-n+1)}} 
				\circ \cdots \circ 
				\ad_{\xi_{\sigma(k)}} 
			\right) (\eta)
		\right)
		\\
		& \leq
		(k + 1)^R
		\sum\limits_{n = 0}^N
		\frac{ k! (k-n)!^R }{ (k-n)! n! }
		|B_n^*| |z|^n
		p\left( \xi_1 \right) \cdots p\left( \xi_k \right) 
		p(\eta)
		\\
		& =
		(k + 1)^R
		p_R \left(
			\xi_1 \tensor \cdots \tensor \xi_k
		\right)
		p(\eta)
		\sum\limits_{n = 0}^N
		\left( \frac{k!}{(k-n)!} \right)^{1-R}
		\frac{|B_n^*| |z|^n}{n!}
		\\
		& \leq
		(k + 1)^R
		k^{N (1-R)}
		p_R \left(
			\xi_1 \tensor \cdots \tensor \xi_k
		\right)
		p(\eta)
		\sum\limits_{n = 0}^N
		\frac{|B_n^*| |z|^n}{n!}.
	\end{align*}
\end{proof}
Now, we can iterate Lemma~\ref{Nilpot:Lemma:PreContinuity2} in the same way 
we did in Chapter 5:
\begin{proof}[Alternative Proof of Theorem~\ref{Nilpot:Thm:ProjLimit}]
	The calculation is done only on factorizing tensors. We need to 
	transform the $k^{N(1-R)}$ into a very small factorial somehow. This is 
	possible, since for given $N \in \mathbb{N}$ and $0 \leq R < 1$, the 
	sequence
	\begin{equation*}
		\left( \frac{k^N}{\sqrt{k!}} \right)^{1-R}
	\end{equation*}
	converges to $0$ for $k \longrightarrow \infty$ and is therefore bounded 
	by some $\kappa_N > 0$. Hence we get
	\begin{equation*}
		k^{ N (1-R) } 
		\leq
		\kappa_N \sqrt{k!}^{1-R},
	\end{equation*}
	and together with Lemma~\ref{Nilpot:Lemma:PreContinuity2} we find
	\begin{equation*}
		p_R \left( x \star_z \eta \right)
		\leq
		c \kappa_N
		(k + 1)^R k!^{\frac{1-R}{2}} 
		p_R (x) p(\eta)
	\end{equation*}
	for any tensor $x$ of degree at most $k$. Now, we can iterate this result
	for $\xi_1, \ldots, \xi_k, \eta_1, \ldots, \eta_\ell 
	\in \lie{g}$, $R \geq 0$:
	\begin{align*}
		&
		p_R \big(
			\xi_1 \tensor \cdots \tensor \xi_k 
			\star_z 
			\eta_1 \cdots \eta_\ell
		\big)
		\\
		& =
		p_R \left(
			\frac{1}{\ell!}
			\sum\limits_{\tau \in S_\ell}
			\xi_1 \tensor \cdots \tensor  \xi_k 
			\star_z
			\eta_{\tau(1)} \star_z \cdots \star_z \eta_{\tau(\ell)}
		\right)
		\\
		& \leq
		\frac{1}{\ell!}
		\sum\limits_{\tau \in S_\ell}
		c \kappa_N
		(k + \ell)^R
		(k + \ell - 1)!^{\frac{1-R}{2}}
		\\
		& \quad 
		\cdot
		p_R \left(
			\frac{1}{\ell!}
			\sum\limits_{\tau \in S_\ell}
			\xi_1 \tensor \cdots \tensor \xi_k 
			\star_z
			\eta_{\tau(1)} \star_z \cdots \star_z \eta_{\tau(\ell-1)}
		\right)
		p\left( \eta_{\tau(\ell)} \right)
		\\
		& \leq
		\quad \vdots
		\\
		& \leq
		(c \kappa_N)^{\ell}
		\left(
			\frac{(k + \ell)!}{k!}
		\right)^R
		(k + \ell - 1)!^{\frac{1-R}{2}}
		\ldots
		k!^{\frac{1-R}{2^N}}
		p_R \left( \xi_1 \tensor \cdots \tensor \xi_k \right)
		p\left( \eta_{\tau(1)} \right)
		\cdots
		p\left( \eta_{\tau(\ell)} \right)
		\\
		& \leq
		(c \kappa_N)^{\ell}
		\binom{k + \ell}{k}^R 
		\ell!^R
		(k + \ell)!^{\frac{(2^N - 1) (1-R)}{2^N} }
		p_R \left( \xi_1 \tensor \cdots \tensor \xi_k \right)
		p\left( \eta_{\tau(1)} \right)
		\cdots
		p\left( \eta_{\tau(\ell)} \right)
		\\
		& \leq
		(c \kappa_N)^{\ell}
		2^{(k + \ell) R}
		\ell!^R
		k!^{\frac{(2^N - 1) (1-R)}{2^N} }
		\ell!^{\frac{(2^N - 1) (1-R)}{2^N} }
		2^{ (k + \ell) \frac{(2^N - 1) (1-R)}{2^N} }
		\\
		& \quad
		\cdot
		p_R \left( \xi_1 \tensor \cdots \tensor \xi_k \right)
		p\left( \eta_{\tau(1)} \right)
		\cdots
		p\left( \eta_{\tau(\ell)} \right)
		\\
		& \leq
		(2 p)_{R + \epsilon} 
		\left( \xi_1 \tensor \cdots \tensor \xi_k  \right)
		(2 c \kappa_N p)_{R + \epsilon} 
		\left( \eta_1 \tensor \cdots \tensor \eta_\ell \right),
	\end{align*}
	where we have set $\epsilon = \frac{(2^N - 1)(1 - R)}{2^N}$. From this,
	we have $R + \epsilon < 1$ and get the wanted result for the 
	projective limit.
\end{proof}
Again, just like in the case of AE-Lie algebras for $R \geq 1$, we can show that 
the star product depends analytically on the formal parameter.
\begin{proposition}[Dependence on $z$]
    \label{Nilpot:corollary:HolomorphicDependence}%
    Let $\lie{g}$ be a nilpotent locally convex Lie algebra, $0 \leq R < 1$ and $z 
    \in \mathbb{K}$, then for all $x, y \in \widehat{\Sym}_{1^-}^\bullet(\lie{g})$ 
    the map
    \begin{equation}
        \label{Nilpot:Holomorphicity}
        \mathbb{K} \ni z
        \longmapsto
        x \star_z y \in
        \widehat{\Sym}_{1^-}^\bullet(\lie{g})
    \end{equation}
    is analytic with (absolutely convergent) Taylor expansion at $z = 0$ 
    given by Equation~\eqref{Formulas:2MonomialsFormula1}. For 
    $\mathbb{K} = \mathbb{C}$, the collection of algebras $\left\{ \left( 
    \widehat{\Sym}_{1^-}^\bullet(\lie{g}), \star_z \right) \right\}_{z \in 
    \mathbb{C}}$ is an entire holomorphic deformation of the completed 
    symmetric tensor algebra $\left( \widehat{\Sym}_{1^-}^\bullet(\lie{g}), \vee 
    \right)$.
\end{proposition}
\begin{proof}
	The proof is completely analogue to the case of AE Lie algebras when $R = 1$.
\end{proof}

%
% A bit Functioriality also in this case
%

\subsection{Representations and Functoriality}
\label{subsec:NilpotentFunctorialityRepresentations}

In the general AE case, we had some useful results concerning representations 
of Lie algebras and the functorialty of our construction. These results can be 
extended to the projective limit $\Sym_{1^-}^{\bullet}(\lie{g})$.
\begin{proposition}[Universal Property]
	\label{Nilpot:Prop:UnivProperty}
	Let $\lie{g}$ be a locally convex nilpotent Lie algebra, $\algebra{A}$ an 
	associative AE algebra and $\phi \colon \lie{g} \longrightarrow 
	\algebra{A}$ is a continuous homomorphism of Lie algebras. Then, the 
	lifted homomorphisms from $\Sym_{1^-}^{\bullet}(\lie{g})$ and 
	$\algebra{U}(\lie{g}_z)$ to $\algebra{A}$ are continuous.
\end{proposition}
\begin{proof}
	The proof is exactly the same as in the general AE case, since there, 
	we actually only needed $R \geq 0$.
\end{proof}
Again, this construction is not universal in the categorial sense, 
since $\Sym_{1^-}^{\bullet}(\lie{g})$ fails to be AE. But also here, we get 
the case of continuous representations into a Banach space (and in particular 
into a finite-dimensional space) as a corollary.
\begin{corollary}[Continuous Representations]
    \label{Nilpot:Coro:ContinuousRepresentations}%
    Let $\algebra{U}_R(\lie{g}_z)$ the universal enveloping algebra of locally 
    convex nilpotent Lie algebra $\lie{g}$ with bracket scaled by $z \in 
    \mathbb{C}$, then for every continuous 
    representation $\phi$ of $\lie{g}$ into the bounded linear operators 
    $\Bounded(V)$ on a Banach space $V$, the induced homomorphism of 
    associative algebras $\Phi \colon \algebra{U}_R(\lie{g}_z) \longrightarrow 
    \Bounded(V)$ is continuous.
\end{corollary}
We can also extend the functoriality statement to the projective limit, but we 
need to get another version of Lemma~\ref{LCAna:Lemma:LemmaPreContinuityN} for 
nilpotent Lie algebras, since this is the corner stone of the functoriality 
proof.
\begin{lemma}
    \label{Lemma:Nilpot:LemmaPreContinuityN}%
    Let $\lie{g}$ be locally convex nilpotent Lie algebra and $N \in 
    \mathbb{N}$ such that $N + 1$ Lie brackets vanish, $0 \leq R < 1$ and 
    $z \in \mathbb{C}$. Then for $p$ a continuous seminorm, $q$ an
    asymptotic estimate, $n \in \mathbb{N}$ and all $\xi_1, \ldots,
    \xi_n \in \lie{g}$ the estimate
    \begin{equation}
        \label{Nilpot:LemmaPreContinuityN}
        p_R \left(
            \xi_1 \star_z \cdots \star_z \xi_n
        \right)
        \leq
        c^n n!^{R + \epsilon}
        q^n(\xi_1 \tensor \cdots \tensor \xi_n)
    \end{equation}
    holds with $c = 16 \E^2 (|z| + 1)$ and $\epsilon = \frac{N-1}{N}(1 - R)$
    and the estimate is locally uniform in $z$.
\end{lemma}
\begin{proof}
    We take $R < 1$ and go directly into the proof of
    Lemma~\ref{LCAna:Lemma:LemmaPreContinuityN} at
    \eqref{LCAna:PreContinuityIntermediateN}. We know that, since we
    may have at most $N$ brackets, also the values for $\ell$ are
    restricted to
    \begin{equation*}
        \ell
        \leq
        \frac{N-1}{N} n
        =
        \delta n.
    \end{equation*}
    Using that in the proof of Lemma~\ref{LCAna:Lemma:LemmaPreContinuityN} 
    leads to
    \begin{align*}
        &p_R \left(
            \xi_1 \star_z \cdots \star_z \xi_n
        \right)
        \\
        &\quad\leq
        \sum\limits_{\ell = 0}^{\delta n}
        (n - \ell)!^R
        \sum\limits_{\substack{
			1 \leq j \leq n-1 \\
			i_j \in \{0, \ldots, j\} \\
			\sum_{j = 1}^{n - 1} i_j = \ell
		}}
        |z|^{\ell}
        \frac{1!  (2 - i_1)! \cdots (n-1 - i_1 - \cdots - i_{n-2})!}
        {(1 - i_1)! \cdots (n-1 - i_1 - \cdots - i_{n-1})!}
        q(\xi_1) \cdots q(\xi_n)
        \\
        &\quad\leq
        \sum\limits_{\ell = 0}^{\delta n}
        (n - \ell)!^R
        \sum\limits_{\substack{
			1 \leq j \leq n-1 \\
			i_j \in \{0, \ldots, j\} \\
			\sum_{j = 1}^{n - 1} i_j = \ell
		}}
        |z|^{\ell} (2 \E)^n \ell!
        q(\xi_1) \cdots q(\xi_n)
        \\
        &\quad\leq
        (2 \E)^n (|z| + 1)^n
        q(\xi_1) \cdots q(\xi_n)
        \sum\limits_{\ell = 0}^{\delta n}
        (n - \ell)!^R \ell!
        \binom{n + \ell - 2}{\ell - 1}.
    \end{align*}
    We have
    \begin{equation*}
        \ell!
        =
        \ell!^R
        \ell!^{1-R}
        \leq
        \ell!^R
        \left(
            (\delta n)^{\delta n}
        \right)^{1-R}
        \leq
        \ell!^R
        n^{\delta n (1-R)}
        \leq
        \ell!^R
        n!^{\delta (1 - R)}
        \E^{\delta n (1 - R)}.
    \end{equation*}
    Together with $\ell!^R (n - \ell)!^R \leq n!^R$ this gives
    \begin{align*}
        p_R \left(
            \xi_1 \star_z \cdots \star_z \xi_n
        \right)
        & \leq
        (2 \E)^n (|z| + 1)^n
        n!^R n!^{\delta (1 - R)}
        q(\xi_1) \cdots q(\xi_n)
        \sum\limits_{\ell = 0}^{\delta n}
        \binom{n + \ell - 2}{\ell - 1}
        e^{\delta n (1 - R)}
        \\
        & \leq
        (2 \E)^n (|z| + 1)^n
        \left(\E^{(1-R) \delta}\right)^n
        4^n n!^{R + \epsilon}
        q(\xi_1) \cdots q(\xi_n),
    \end{align*}
    with $\epsilon = \delta (1-R)$. It is clear that for all $R < 1$ we have
    $R + \epsilon < 1$. Set 
    \begin{equation*}
    	c 
    	= 
    	8 \E (|z|+1) \E^{(1-R)\delta}
    	\leq
    	16 \E^2 (|z|+1)
	\end{equation*}
	and note that the estimate is locally uniform in $z$.
\end{proof}
\begin{proposition}
	\label{Nilpot:Prop:Functoriality}
	Let $R \geq 1$, $\lie{g}, \lie{h}$ be two locally convex nilpotent Lie 
	algebras and $\phi \colon \lie{g} \longrightarrow \lie{h}$ a continuous 
	homomorphism between them. Then it lifts to a continuous unital 
	homomorphism of locally convex algebras $\Phi_z \colon 
	\algebra{U}_R(\lie{g}_z) \longrightarrow \algebra{U}_R(\lie{h}_z)$.
\end{proposition}
\begin{proof}
	The proof is analogous to the one of Proposition 
	\ref{LCAna:Prop:Functoriality}.
\end{proof}

\section{Module Structures}
\label{sec:chap6_Modules}

The projective limit $1^-$ is not the only additional structure we will get, if 
our Lie algebra $\lie{g}$ is nilpotent. For every $R \in \mathbb{R}$, the 
symmetric tensor algebra $\Sym_R^{\bullet}(\lie{g})$ is a locally convex vector 
space. For $R \geq 0$, the (symmetric) tensor product is continuous, which is 
very important for many estimates, and for $R \geq 1^-$, we have an algebra 
structure. In between however, we have more than ``only'' vector spaces: the 
spaces $\Sym_R^{\bullet}(\lie{g})$ form locally convex bimodules over the 
$\Sym_{R'}^{\bullet}(\lie{g})$ for certain values of $R'$. The next proposition 
makes this more exact.
\begin{proposition}[Bimodules in $\Sym_R^{\bullet}(\lie{g})$]
	\label{Nilpot:Prop:Bimodules}
	Let $\lie{g}$ be a nilpotent, locally convex Lie algebra, $N \in 
	\mathbb{N}$ such that $N + 1$ Lie brackets vanish, $z \in \mathbb{K}$ and 
	$0 \leq R < 1$. Then, for all $x, y \in \Sym^{\bullet}(\lie{g})$ and every 
	continuous seminorm $p$, we have a continuous seminorm $q$, such that the 
	estimates
	\begin{align}
		\label{Nilpot:BimoduleEstimate1}
		p_R \left(
			x \star_z y
		\right)
		& \leq
		\left(16 q\right)_R(x) 
		\left(16 c q\right)_{R + N(1-R)}(y)
		\\
	\intertext{and}
		\label{Nilpot:BimoduleEstimate2}
		p_R \left(
			x \star_z y
		\right)
		& \leq
		\left(16 c q\right)_{R + N(1-R)}(x)
		\left(16 q\right)_R(y) 
	\end{align}
	hold with $c = (N \E)^{N (1- R)}$.
	Hence, the vector space $\widehat{\Sym}_R^{\bullet}(\lie{g})$ forms a 
	bimodule over the algebra $\widehat{\Sym}_{R + N(1-R)}^{\bullet}(\lie{g})$. 
	In particular, if $\lie{g}$ is 2-step nilpotent, 
	the vector space $\widehat{\Sym}_0^{\bullet}(\lie{g})$ is a 
	$\widehat{\Sym}_1^{\bullet}(\lie{g})$-bimodule.
\end{proposition}
\begin{proof}
	Note that for every degree $n$ we loose, we get a bracket of $\xi$'s and 
	$\eta$'s. Since we can not have too highly nested brackets, we get the 
	following bounds:
	\begin{equation*}
		n
		\leq
		N k
		\quad \text{ and } \quad
		n 
		\leq 
		N \ell.
	\end{equation*}
	We prove the statement only on factorizing tensors again. We 
	want to show Estimate \eqref{Nilpot:BimoduleEstimate1} and take
	$x = \xi_1 \tensor \cdots \tensor \xi_k$ and $y = \eta_1 \tensor \cdots 
	\tensor \eta_\ell$. So
	\begin{align*}
		(\ell N)!^{1-R}
		& \leq
		(\ell N)^{(\ell N (1-R))}
		\\
		& \leq
		(N \E)^{\ell N (1- R)}
		\ell!^{N(1-R)}.
	\end{align*}
	This allows us again to go back to the proof of Theorem 
	\ref{Nilpot:Thm:ProjLimit} and we find
	\begin{align*}
		p_R \big(
			C_n \big(
				\xi_1 \tensor \cdots \tensor \xi_k,
		&
				\eta_1 \tensor \cdots \tensor \eta_\ell
			\big)
		\big)
		\\
		& \leq 
		\frac{ (N \E)^{\ell N (1- R)} \ell!^{N(1-R)} }
		{2 \cdot 8^n}
		(16 q)_R
		\left(
			\xi_1 \tensor \cdots \tensor \xi_k
		\right)
		(16 q)_R
		\left(
			\eta_1 \tensor \cdots \tensor \eta_\ell
		\right)
		\\
		& \leq
		\frac{1}{2 \cdot 8^n}
		(16 q)_R
		\left(
			\xi_1 \tensor \cdots \tensor \xi_k
		\right)
		(c q)_{R - N(1-R)}
		\left(
			\eta_1 \tensor \cdots \tensor \eta_\ell
		\right)
	\end{align*}
	with $c = 16 (N \E)^{N (1- R)}$. The rest of the proof is analogue to the 
	proofs of the Theorems \ref{Nilpot:Thm:ProjLimit} or
	 \ref{Thm:LCAna:Continuity1}.
\end{proof}
Once again, assuming submultiplicativity of the seminorms, it is possible to give 
an easier proof for Proposition~\ref{Nilpot:Prop:Bimodules} which relies on 
Lemma~\ref{Nilpot:Lemma:PreContinuity2}.
\begin{proof}[Alternative proof for Proposition~\ref{Nilpot:Prop:Bimodules}]
	We do the calculation on factorizing tensors: let 
	$\xi_1, \ldots, \xi_k, \eta_1, \ldots, \eta_\ell \in \lie{g}$, $R \geq 0$, $k, 
	\ell \in \mathbb{N}$. Using Lemma~\ref{Nilpot:Lemma:PreContinuity2}, we get
	\begin{align*}
		p_R \big(
			\xi_1 \tensor
		&			
			\cdots \tensor \xi_k
			\star_z
			\eta_1 \tensor \cdots \tensor \eta_\ell
		\big)
		\\
		&=
		p_R \left(
			\frac{1}{\ell!}
			\sum\limits_{\tau \in S_\ell}
			\xi_1 \tensor \cdots \tensor \xi_k
			\star
			\eta_{\tau(1)} \star_z \cdots \star_z \eta_{\tau(\ell)}
		\right)
		\\
		& \leq
		c (k + \ell)^R
		(k + \ell - 1)^{N (1-R)}
		\frac{1}{\ell!}
		\sum\limits_{\tau \in S_\ell}
		p_R \left(
			\xi_1 \tensor \cdots \tensor \xi_k
			\star
			\eta_{\tau(1)} \star_z \cdots \star_z \eta_{\tau(\ell-1)}
		\right)
		p\left( \eta_{\tau(\ell)} \right)
		\\
		& \leq
		\quad \vdots
		\\
		& \leq
		c^{\ell}
		\left(
			\frac{(k + \ell)!}{k!}
		\right)^R
		\left(
			\frac{(k + \ell - 1)!}{(k - 1)!}
		\right)^{N (1-R)}
		p_R \left(
			\xi_1 \tensor \cdots \tensor \xi_k
		\right)
		p\left( \eta_{\tau(1} \right)
		\cdots
		p\left( \eta_{\tau(\ell)} \right)
		\\
		& \leq
		c^{\ell}
		2^{k + l} 2^{N (k + \ell)}
		\ell!^{N (1-R)}
		p_R \left(
			\xi_1 \tensor \cdots \tensor \xi_k
		\right)
		p\left( \eta_{\tau(1} \right)
		\cdots
		p\left( \eta_{\tau(\ell)} \right)
		\\
		& =
		\left(2^{N + 1} p\right)_R 
		\left(
			\xi_1 \tensor \cdots \tensor \xi_k
		\right)
		\left(2^{N + 1} c p\right)_{R + N(1-R)}
		\left(
			\eta_1 \tensor \cdots \tensor \eta_\ell
		\right).
	\end{align*}
	The proof of the second estimate is analogous.
\end{proof}
\begin{remark}[Possible extensions]
	This result immediately raises new questions, like the one about possible 
	generalizations to ``weaker forms'' of nilpotency, for example. They may be 
	issues of some future work, but can not
	be addressed here, since we rather want to present something like  
	the ``big picture'', instead of getting 
	lost in its details too much. There are, without any doubt, questions that 
	are more significant than extending those estimates to very special cases 
	and finding sharp bounds there, although this is interesting and 
	important, too.
\end{remark}
Though it seems clear from the construction that these bimodules can not be 
there for general Lie algebras, we can give a concrete counter-example, which 
shows that there are Lie algebras, which don not allow them.
\begin{example}
	\label{Nilpot:Ex:NoModulesInGeneral}
	Choose $R < 1$ and take $\lie{g} = \mathbbm{R}^3$ with the basis $e_1, e_2, 
	e_3$ and the vector product as Lie bracket:
	\begin{equation*}
		[e_1, e_2] 
		= 
		e_3 
		\qquad 
		[e_2, e_3] 
		= 
		e_1 
		\qquad 
		[e_3, e_1] 
		= 
		e_2
	\end{equation*}
	Again, we take a $\ell^1$-norm $n$ such that $n(e_1) = n(e_2) = n(e_3) 
	= 1$. It has the nice property that for $k, \ell, m \in \mathbb{N}$ we get
	on the projective tensor product
	\begin{equation*}
		n^{k + \ell + m} \left(
			e_1^k e_2^{\ell} e_3^m
		\right)
		=
		1.
	\end{equation*}
	Now choose an $\epsilon > 0$ such that $R + \epsilon < 1$ and we define the 
	sequence $(a_k)_{k \in \mathbb{N}}$
	\begin{equation*}
		a_k 
		= 
		\frac{1}{k!^R} e_1^k,
	\end{equation*}
	for which we get $\lim_{k \longrightarrow \infty} n_R(a_k) = 0$. Now, we want 
	to show that $a_k \star_z e_2$ grows faster than exponentially:
	\begin{align*}
		n_R \left( a_k \star_z e_2 \right) 
		& = 
		n_R \left( 
			\sum\limits_{j = 0}^k 
			\binom{k}{j} B_j^* 
			\frac{1}{k!^{R + \epsilon}} 
			e_1^{n-j} 
			\left( 
				\operatorname{ad}_{e_1} 
			\right)^j(e_2) 
		\right)
		\\
		& =
		\sum\limits_{j = 0}^k 
		\binom{k}{j} 
		|B_j^*| 
		\frac{1}{k!^{R + \epsilon}} 
		(k-j+1)!^R 
		\underbrace{
			n^{k-j} \left( e_1 (e_2 \wedge e_3) \right)
		}_{ = 1}
		\\
		& = 
		\sum\limits_{j = 0}^k 
		(k-j+1)^R 
		\binom{k}{j}
		\frac{|B_j^*|}{j!} 
		\frac{(k-j)!^R j^R}{k!^{R + \epsilon}} 
		j!^{1-R}
		\\
		& = 
		\sum\limits_{j=0}^k 
		(k-j+1)^R 
		\binom{k}{j}^{1-R} 
		\frac{|B_j^*|}{j!} 
		\frac{j!^{1-R}}{k!^\epsilon}
		\\
		& \geq 
		\sum\limits_{j=0}^k 
		\frac{|B_j^*|}{j!} 
		\frac{j!^{1-R}}{k!^\epsilon}
		\\
		& \geq 
		\frac{|B_k^*|}{k!^{R + \epsilon}}.
	\end{align*}
	We know that for $R +\epsilon < 1$ and any $c > 0$
	\begin{equation*}
		\limsup_{n \longrightarrow \infty}
		\frac{|B_n^*|}{c^n n!^{R + \epsilon}}
		=
		\infty.
	\end{equation*}
	Hence the Limes superior of $n_R \left( a_k \star_z e_2 \right) $ grows 
	faster than any exponential function and can not be absorbed into the seminorm 
	of $e_2$. So the multiplication in the module can not be continuous.
\end{example}

\section{The Heisenberg and the Weyl Algebra}
\label{sec:chap6_HeisenbergWeyl}

Now we want to see how we get the link to the Weyl algebra from
\cite{waldmann:2014a}, since we have something like a discrepancy for the
parameter $R$ concerning the continuity of the product in the Weyl and the 
Heisenberg algebra. In the following, we will show that this gap actually makes 
a lot of sense. For simplicity, we consider the easiest case of
the Weyl/Heisenberg algebra with two generators $Q$ and $P$, but the 
calculation for the Heisenberg/Weyl algebra in $2n + 1$ [$2n$] dimensions 
is done exactly in the same way.
Recall that the Weyl algebra is a quotient of the enveloping algebra of the 
Heisenberg algebra $\lie{h}$ which one gets from dividing out its center. So 
let $h \in \mathbb{K}$ and we have a projection
\begin{equation}
    \label{Nilpot:WeylProjection}
    \pi \colon
    \Sym_R^\bullet(\lie{h})
    \longrightarrow
    \mathcal{W}_R(\lie{h})
    =
    \frac{\Sym_R^\bullet(\lie{h})}
    {\langle E - h \Unit \rangle}
\end{equation}
Of course we want to know if this projection is continuous.
\begin{proposition}
    \label{proposition:ProjectionWeylContinuous}%
    The projection $\pi$ is continuous for $R \geq 0$.
\end{proposition}
\begin{proof}
    We extend $\pi$ to the whole tensor algebra by symmetrizing beforehand. Let 
    then $p$ be a continuous seminorm on $\lie{h}$ and $k, \ell, m \in 
    \mathbb{N}_0$. We have
    \begin{align*}
        p_R(\pi (
        	Q^{\tensor k} \tensor
        	P^{\tensor \ell} \tensor
        	E^{\tensor m}
        ) )
        & =
        p_R( Q^k P^{\ell} h^m )
        \\
        & =
        |h|^m (k + \ell)!^R
        p^{k + \ell}(Q^k P^{\ell})
        \\
        & \leq
        (|h| + 1)^{k + \ell + m}
        (k + \ell + m)!^R
        p(Q)^k p(P)^{\ell} p(E)^m
        \\
        & =
        ((|h| + 1) p)_R
        (Q^{\tensor k} \tensor
        P^{\tensor \ell} \tensor
        E^{\tensor m}).
    \end{align*}
    Then we do the usual infimum argument and have the result on
    arbitrary tensors again.
\end{proof}

To establish the link to the continuity results of the Weyl algebra,
we need more: $\pi \circ \ostar_z$ should to be continuous for $R \geq \frac 
1 2$.
\begin{proposition}
    \label{proposition:ContinuousProductInWeyl}%
    Let $R \geq \frac{1}{2}$ and $\pi$ the projection from
    \eqref{Nilpot:WeylProjection}. Then the map $\pi \circ
    \ostar_z$ is continuous.
\end{proposition}
\begin{proof}
    Since we are in finite dimensions, we can choose a
    submultiplicative norm $p$ with $p(Q) = p(P) = p(E)$ without
    restrictions. Moreover, let $k, k', \ell, \ell', m, m' \in
    \mathbb{N}_0$. Then we have to get an estimate for $p_R \left(
    \pi\left( Q^k P^{\ell} E^m \ostar_z Q^{k'} P^{\ell'} E^{m'}
    \right) \right)$.  If we calculate the star product
    explicitly, we will see that we only get Lie brackets where we have
    $P$'s and $Q$'s. Let $r = k + \ell + m$ and $s = k' + \ell' + m'$,
    then we can actually simplify the calculations by
    \begin{align*}
        p_R \left(
         	\pi\left( 
         		Q^k P^{\ell} E^m
         		\right.
         	\right.
        &
         	\left.
         		\left.
         		\ostar_z
         		Q^{k'} P^{\ell'} E^{m'}
         	\right)
        \right)
        =
        (p_R \circ \pi) \left(
        \sum\limits_{n=0}^{r + s - 1}
        z^n C_n(Q^k P^{\ell} E^m,
        Q^{k'} P^{\ell'} E^{m'})
        \right)
        \\
        & \leq
        \sum\limits_{n=0}^{r + s - 1}
        |z|^n
        (p_R \circ \pi) \left(
        C_n(Q^k P^{\ell} E^m,
        Q^{k'} P^{\ell'} E^{m'})
        \right)
        \\
        & \leq
        \sum\limits_{n=0}^{r + s - 1}
        |z|^n
        (p_R \circ \pi) \left(
        C_n(Q^r, P^s)
        \right)
        \\
        & =
        \sum\limits_{n=0}^{r + s - 1}
        |z|^n
        \frac{r! s!}{(r-n)! (s-n)! n!}
        (p_R \circ \pi) \left(
        Q^{r-n} P^{s-n} E^n
        \right)
        \\
        & =
        \sum\limits_{n=0}^{r + s - 1}
        |z|^n |h|^n
        \frac{r! s!}{(r-n)! (s-n)! n!}
        p_R \left(
        Q^{r-n} P^{s-n}
        \right)
        \\
        & \leq
        \sum\limits_{n=0}^{r + s - 1}
        |z|^n |h|^n
        \frac{r! s!}{(r-n)! (s-n)! n!}
        \frac{(r + s - 2n)!^R}{r!^R s!^R}
        p_R \left(Q^{\tensor r} \right)
        p_R \left(P^{\tensor s} \right)
        \\
        & \leq
        \sum\limits_{n=0}^{r + s - 1}
        |z|^n |h|^n
        \binom{r}{n} \binom{s}{n}
        \frac{(r + s - 2n)!^R n!}{r!^R s!^R}
        p_R \left(Q^{\tensor r} \right)
        p_R \left(P^{\tensor s} \right)
        \\
        & \ot{(a)}{\leq}
        \sum\limits_{n=0}^{r + s - 1}
        |z|^n |h|^n
        \binom{r}{n} \binom{s}{n}
        \binom{r + s}{s}^R
        \binom{r + s}{2n}^{-R}
        p_R \left(Q^{\tensor r} \right)
        p_R \left(P^{\tensor s} \right)
        \\
        & \leq
        \sum\limits_{n=0}^{r + s - 1}
        (|z| + 1)^n (|h| + 1)^n
        4^{r + s}
        p_R \left(Q^{\tensor r} \right)
        p_R \left(P^{\tensor s} \right)
        \\
        & \ot{(b)}{ \leq }
        \underbrace{
        (8 (|z| + 1) (|h| + 1))^{r + s}
        }_{ = \tilde{c}^{r + s}}
        p_R \left(Q^{\tensor r} \right)
        p_R \left(P^{\tensor s} \right)
        \\
        & =
        (\tilde{c} p)_R \left(Q^{\tensor r} \right)
        (\tilde{c} p)_R \left(P^{\tensor s} \right)
        \\
        & \ot{(c)}{ = }
        (\tilde{c} p)_R \left(
        Q^{\tensor k} \tensor
        P^{\tensor \ell} \tensor
        E^{\tensor m} \right)
        (\tilde{c} p)_R \left(
        Q^{\tensor k'} \tensor
        P^{\tensor \ell'} \tensor
        E^{\tensor m'} \right),
    \end{align*}
    where we have set $c = 8 (|z| + 1) (|h| + 1)$. We rearranged the factorials in 
    ($a$) and used $R \geq \frac{1}{2}$. The estimates ($b$) are the standard 
    binomial coefficient estimates. In ($c$) we used $p(Q) = p(P) = p(E)$. Now
    we just use
    \begin{equation*}
        \left(
        	Q^{\tensor k} \tensor
        	P^{\tensor \ell} \tensor
        	E^{\tensor m}
        \right)
        \ostar_z
        \left(
        	Q^{\tensor k'} \tensor
        	P^{\tensor \ell'} \tensor
        	E^{\tensor m'}
        \right)
        =
        Q^k P^{\ell} E^m
        \ostar_z
        Q^{k'} P^{\ell'} E^{m'}
    \end{equation*}
    and the infimum argument to expand this estimate to all tensors.
    This concludes the proof.
\end{proof}
The previous proposition can be seen as something like the ``finite-dimensional 
version'' of Lemma 3.10 in \cite{waldmann:2014a}, just that we took a large 
detour for proving it. One could, most probably, redo some more results of this 
paper using finite-dimensional versions the Heisenberg algebra and the 
projection onto the Weyl algebra, but this would yield, also most probably, 
nothing new. It is good to know that this connections exists, but it is not 
something which is very helpful to pursue, since an evident generalization to 
infinite dimensions does not seem be obvious.

\section{Banach-Lie Algebras}
\label{sec:chap6_TheEProperty}

Now we want to focus a bit on weaker notions than true nilpotency. Since there 
are many of them, we want to restrict to the easier case of Banach-Lie algebras, 
where a somewhat developed theory already exists.

\subsection{Generalizations of Nilpotency}

In \cite{mueller:1994a}, M\"uller gives a list of weaker forms of nilpotency for 
associative Banach algebras. We can mostly copy the ideas and use them for 
Banach-Lie algebras, too
\begin{definition}
	Let $\lie{g}$ be a Banach-Lie algebra in which the Lie bracket fulfils the 
	estimate
	\begin{equation*}
		\norm{ [\xi, \eta] }
		\leq
		\norm{\xi}
		\norm{\eta}.
	\end{equation*}
	Denote by $\mathbb{B}_1(0)$ all elements $\xi \in \lie{g}$ with 
	$\norm{\xi} = 1$. We say that
	\begin{definitionlist}
		\item
		$\lie{g}$ is topologically nil (or radical, or quasi-nilpotent), if
		every $\xi \in \lie{g}$ is quasi-nilpotent, i.e.
		\begin{equation*}
			\lim_{n \longrightarrow \infty}
			\norm{\ad_{\xi}^n}^{\frac{1}{n}}
			=
			0.
		\end{equation*}
		
		\item
		$\lie{g}$ is uniformly topologically nil, if
		\begin{equation*}
			\lim_{n \longrightarrow \infty}
			\mathcal{N}_1(n)
			=
			0
		\end{equation*}
		for
		\begin{equation}
			\mathcal{N}_1(n)
			=
			\sup \left\{ 
			\left.
				\norm{ \ad_{\xi}^n}^{\frac{1}{n}} 
			\right|
				\xi \in \mathbb{B}_1(0)
			\right\}.
		\end{equation}
		
		\item
		$\lie{g}$ is topologically nilpotent, if for every sequence
		$(\xi_n)_{n \in \mathbb{N}} \subset \mathbb{B}_1(0)$ we have
		\begin{equation*}
			\lim_{n \longrightarrow \infty}
			\norm{ 
				\ad_{\xi_1} \circ \ldots \circ \ad_{\xi_n}
			}^{\frac{1}{n}}
			=
			0.
		\end{equation*}
		
		\item
		$\lie{g}$ is uniformly topologically nilpotent, if
		\begin{equation*}
			\lim_{n \longrightarrow \infty}
			\mathcal{N}(n)
			=
			0
		\end{equation*}
		for
		\begin{equation}
			\mathcal{N}(n)
			=
			\sup \left\{ 
			\left.
				\norm{ 
					\ad_{\xi_1} \circ \ldots \circ \ad_{\xi_n}
				}^{\frac{1}{n}} 
			\right|
				\xi_1, \ldots, \xi_n \in \mathbb{B}_1(0)
			\right\}.
		\end{equation}
	\end{definitionlist}
\end{definition}
It is clear that $(ii) \Rightarrow (i)$ and $(iv) \Rightarrow (iii)$. In the 
associative case, we have $(iii) \Leftrightarrow (iv)$ and hence $(iii) 
\Rightarrow (ii)$. Of course, it is a good question, if this remains true for 
Banach-Lie algebras. We have already encountered notion $(i)$: in 
\cite{wojtynski:1998a} Wojty\'nski gave a proof that it is equivalent to the 
global convergence of the BCH series. In the following, we will make use of notion 
$(iv)$: we will show, that it is possible to generalize the result of Theorem 
\ref{Nilpot:Thm:ProjLimit} to this case.

\subsection{An Adapted Topology for the Tensor Algebra}

The idea consists in changing the $\Tensor_R$-topology a bit: instead of taking 
$n!^R$ as weights with $0 \leq R < 1$, we take sequence $(\alpha_n)_{n \in 
\mathbb{N}}$ with a certain asymptotic 
behaviour for $n \longrightarrow \infty$ and use $\frac{n!}{\alpha_n}$ as 
weights. This will generalize the idea of $n!^R$ and will be the starting point 
for estimates.

First, we observe that every uniformly topologically nilpotent Banach-Lie 
algebra $\lie{g}$ comes with a characteristic, monotonously decreasing sequence
\begin{equation}
	\label{Nilpot:CharDownSequence}
	\omega_n
	=
	\sup_{m \geq n} \mathcal{N}(m).	
\end{equation}
If there exists a $N \in \mathbb{N}$, such that $\omega_n = 0$ for all $n \geq 
N$, then $\lie{g}$ is actually nilpotent and we can use the results of the 
first section in this chapter. We may hence restrict to those Banach-Lie 
algebras, where we have $\omega_n > 0$ for all $n \in \mathbb{N}$. This allows 
the next definition.
\begin{definition}[Rapidly increasing sequences]
	Let $\lie{g}$ be a uniformly topologically nilpotent Banach-Lie algebra and 
	$(\omega_n)_{n \in \mathbb{n}}$ the sequence defined in 
	\eqref{Nilpot:CharDownSequence}. Then we 
	define the characteristic sequence $(\chi_n^{\lie{g}})_{n \in \mathbb{N}}$ 
	of $\lie{g}$ by
	\begin{equation}
		\label{Nilpot:CharSequence}
		\chi_n^{\lie{g}}
		=
		\max
		\left\{ 
			\frac{1}{\omega_n}
			,
			2
		\right\}.
	\end{equation}
	Furthermore, we will say that a sequence $(\alpha_n)_{n \in \mathbb{N}}$ in 
	$(1, \infty)$ is $\lie{g}$-rapidly increasing, if it fulfils the following 
	properties:
	\begin{definitionlist}
		\item
		It grows fast than exponentially, i.e.
		\begin{equation*}
			\lim_{n \longrightarrow \infty}
			\frac{\log \left(\alpha_n \right)}{n}
			=
			\infty.
		\end{equation*}
		
		\item
		There exists a constant $c > 0$, such that
		\begin{equation*}
			\alpha_n 
			\leq 
			c^n  \chi_n^{\lie{g}}.
		\end{equation*}
	\end{definitionlist}
	We denote by $\mathfrak{I}_{\lie{g}}$ the set of all $\lie{g}$-rapidly 
	increasing sequences.
\end{definition}
Clearly, $(\chi_n^{\lie{g}})_{n \in \mathbb{N}}$ is a $\lie{g}$-rapidly 
increasing sequence itself. Note that with this definition we get for every 
sequence $(\xi_n)_{n \in \mathbb{N}} \subset \mathbb{B}_1(0)$
\begin{equation}
	\label{Nilpot:ChiEstimate}
	\norm{
		[ \ldots [[\xi_1, \xi_2], \xi_3], \ldots \xi_n]
	}
	\leq
	\frac{2^n}{\chi^{\lie{g}}_n}.
\end{equation}
\begin{remark}
	The number $2$ in \eqref{Nilpot:CharSequence} may look a bit confusing at 
	the first sight, since it is somehow arbitrary, but for technical reasons, 
	we will need $\chi^{\lie{g}}_n > 1$ later. 
	So actually every real number $c > 1$ 
	could have been used there. In this sense, the previous definition is rather 
	a technical tool than a ``general concept''.
\end{remark}
Each $(\alpha_n)_{n \in \mathbb{N}} \in \mathfrak{I}_{\lie{g}}$ gives a continuous 
seminorm.
\begin{definition}[Adapted seminorms]
	\label{Def:AdaptedBanachSeminorms}
	Let $(\alpha_n)_{n \in \mathbb{N}} \in \mathfrak{I}_{\lie{g}}$. Then
	\begin{equation*}
		p_{\alpha}
		=
		\sum\limits_{n = 0}^{\infty}
		\frac{n!}{\alpha_n}
		\norm{\cdot}^{\tensor[\pi] n}
	\end{equation*}
	defines a seminorm on the tensor algebra with the projective tensor product 
	$\Tensor_{\pi}^{\bullet}(\lie{g})$. We denote the set of all continuous 
	seminorms with respect to those	coming from such sequences by $\algebra{P}$.
\end{definition}
It will be important to see that every rapidly increasing sequence 
$(\alpha_n)_{n \in \mathbb{N}}$ yields a continuous function $f_{\alpha}$ by
\begin{equation}
	\label{Nilpot:IncreasingFunction}
	f_{\alpha}
	\colon
	\mathbb{R}_0^+
	\longrightarrow
	\mathbb{R}^+
	, \quad
	f_{\alpha}(0)
	=
	2
	,\
	f_{\alpha}(n)
	=
	\log \left( \alpha_n \right)
	,\
	\forall_{n \in \mathbb{N}}
\end{equation}
and linear interpolation between the values at the integers. The idea behind is 
that this will allow us to use a technical lemma, which we now introduce. This 
lemma shows that there are always ``many'' rapidly increasing functions in a 
certain sense. It is taken from a work \cite{mitiagin.rolewicz.zelazko:1962a} by 
Mitiagin, Rolewicz and \.{Z}elazko, where it is stated in Lemma~2.1 and Lemma~2.2.
\begin{lemma}
	\label{Nilpot:Lemma:MRZGrothLemma}
	Let $f\colon \mathbb{R}_0^+ \longrightarrow \mathbb{R}^+$ be a continuous 
	functions, such that
	\begin{equation}
		\label{Nilpot:GrowthProperty}
		\lim_{n \longrightarrow \infty}
		\frac{f(x)}{x}
		=
		\infty.
	\end{equation}
	Then there exists a convex, continuous function $g \colon \mathbb{R}_0^+ 
	\longrightarrow \mathbb{R}^+$, fulfilling \eqref{Nilpot:GrowthProperty} and
	\begin{equation}
		\label{Nilpot:SplittingProperty}
		g \left(
			t_1 + \cdots + t_n
		\right)
		\leq
		8 \left(
			g \left( t_1 \right)
			+ \cdots +
			g \left( t_n \right)
		\right)
		+
		f(n)
		, \quad
		\forall_{n \in \mathbb{N}}
		\text{ and all }
		t_i \in \mathbb{R}_0^+.
	\end{equation}
\end{lemma}
\begin{proof}
	We refer the reader to the paper \cite{mitiagin.rolewicz.zelazko:1962a},
	since we just want to use this lemma and do not want to go too much into
	details here.
\end{proof}
Note that if we have $(\alpha_n)_{n \in \mathbb{N}} \in \mathfrak{I}_{\lie{g}}$, 
then we can apply Lemma~\ref{Nilpot:Lemma:MRZGrothLemma} to 
the function $f_{\alpha}$, which is defined according to Equation 
\ref{Nilpot:IncreasingFunction}.

\subsection{A New Continuity Result}

Now, we have finally prepared our toolbox well enough to prove a new result.
\begin{proposition}
	\label{Nilpot:Prop:TopNilBanachLie}
	Let $\lie{g}$ be a uniformly topologically nilpotent Banach-Lie algebra,
	$(\alpha_n)_{n \in \mathbb{N}} \in \mathfrak{I}_{\lie{g}}$, $p_{\alpha}$ the 
	corresponding seminorm according to \eqref{Nilpot:IncreasingFunction} and 
	$z \in \mathbb{K}$. Then, there exists a series $(\beta_n)_{n \in \mathbb{N}} 
	\in \mathfrak{I}_{\lie{g}}$, such that for all $x,y \in \Tensor^{\bullet}
	(\lie{g})$ we have the estimate
	\begin{equation}
		\label{Nilpot:TopNilBanachLie}
		p_{\alpha} \left(
			x \star_z y
		\right)
		\leq
		(c p)_{\beta} (x)
		(c p)_{\beta} (y)
	\end{equation}
	with a $c > 0$, which only depends on $\alpha, z$ and the Lie algebra 
	$\lie{g}$.
\end{proposition}
\begin{proof}
	We compute the estimate on factorizing tensors and extend it with the 
	infimum argument later. Let hence $k, \ell \in \mathbb{N}$ and 
	$\xi_1, \ldots, \xi_k, \eta_1, \ldots, \eta_\ell \in \lie{g}$. 
	We need to estimate the $C_n$-operators for $n = 0, 1, \ldots, k + \ell - 1$. 
	Therefore we note $r = k + \ell - n$ and use the short-hand notation for the 
	sums, which appear in the Gutt star product, again.
	We take $(\alpha_n)_{n \in \mathbb{N}} \in \mathfrak{I}_{\lie{g}}$ and get 
	for $p_{\alpha}$
	\begin{align*}
		p_{\alpha} \big(
			C_n \big( 
			\xi_1 \tensor \cdots \tensor \xi_k,
		&
			\eta_1 \tensor \cdots \tensor \eta_\ell 
			\big)
		\big)
		\\
		& =
		p_{\alpha}
		\bigg(
         	\frac{1}{r!}
			\sum\limits_{\sigma, \tau}
			\sum\limits_{a_i, b_j}
			\bchtilde{a_1}{b_1}{\xi_{\sigma(i)}}{\eta_{\tau(j)}}
			\cdots
			\bchtilde{a_r}{b_r}{\xi_{\sigma(i)}}{\eta_{\tau(j)}}
        \bigg)
        \\
        & \leq 
    		\frac{1}{r!}
    		\frac{r!}{\alpha_r}
    		\sum\limits_{\sigma, \tau}
			\sum\limits_{a_i, b_j}
       	\norm{
    			\bchtilde{a_1}{b_1}{\xi_{\sigma(i)}}{\eta_{\tau(j)}}
    		}
    		\cdots
    		\norm{
    			\bchtilde{a_r}{b_r}{\xi_{\sigma(i)}}{\eta_{\tau(j)}}
    		}
        \\
        & \leq
    	\frac{k! \ell!}{\alpha_r}
    	\sum\limits_{a_i, b_j}
       	\frac{2}{\chi_{a_1 + b_1}^{\lie{g}}}
       	\ldots
       	\frac{2}{\chi_{a_r + b_r}^{\lie{g}}}
       	\norm{\xi_1} \cdots \norm{\xi_k}
       	\norm{\eta_1} \cdots \norm{\eta_\ell},
    \end{align*}
	where we have used the estimate from Lemma \ref{LCAna:Lemma:BCHTermsEstiamte}
	 $(\ref{Item:BCHEstimate})$ and Estimate \eqref{Nilpot:ChiEstimate}
	in the last step. Rearranging this, we have
	\begin{equation*}
		p_{\alpha} \left(
			C_n\left( \xi^{\tensor k}, \eta^{\tensor \ell} \right)
		\right)
		\leq
    	k! \ell! 2^r
       	\norm{\xi_1} \cdots \norm{\xi_k}
       	\norm{\eta_1} \cdots \norm{\eta_\ell}
    		\sum\limits_{\substack{a_1, b_1, \ldots, a_r, b_r \geq 0 \\
        		a_i + b_i \geq 1 \\
        		a_1 + \ldots + a_r = k \\
        		b_1 + \ldots + b_r = \ell
       	}}
       	\frac{1}
       	{
       		\alpha_r
       		\cdot
       		\chi^{\lie{g}}_{a_1 + b_1}
       		\cdots
       		\chi^{\lie{g}}_{a_r + b_r}
       	},
	\end{equation*}
	and we would like to find a $(\beta_n)_{n \in \mathbb{N}} \in 
	\mathfrak{I}_{\lie{g}}$ such that
	\begin{equation}
		\label{Nilpot:SavingSequence}
		\sup
		\left\{
		\left.
			\frac{\beta_k \cdot \beta_{\ell}}
	       	{
	       		\alpha_r
	       		\cdot
	       		\chi^{\lie{g}}_{a_1 + b_1}
	       		\ldots
	       		\chi^{\lie{g}}_{a_r + b_r}
	       	}
	    \right|
	    		k, \ell \in \mathbb{N},\
	    		a_i + b_i \geq 1,\
        		\sum_i a_i = k,\
        		\sum_j b_j = \ell
		\right\}
		\leq
		\kappa^{k + \ell}
	\end{equation}
	for some $\kappa > 0$, just depending on $(\alpha)$. Then we would have
	\begin{align*}
		p_{\alpha} \big(
			C_n \big( 
			\xi_1 \tensor \cdots \tensor \xi_k,
		&
			\eta_1 \tensor \cdots \tensor \eta_\ell 
			\big)
		\big)
		\\
		& \leq
    	\frac{k! \ell!}{\beta_k \beta_{\ell}} 
    	2^r
       	\norm{\xi_1} \cdots \norm{\xi_k}
       	\norm{\eta_1} \cdots \norm{\eta_\ell}
    		\sum\limits_{\substack{a_1, b_1, \ldots, a_r, b_r \geq 0 \\
        		a_i + b_i \geq 1 \\
        		a_1 + \ldots + a_r = k \\
        		b_1 + \ldots + b_r = \ell
       	}}
       	\kappa^{k + \ell}
       	\\
       	& =
       	2^{-n}
       	(2 \kappa)^{k + \ell}
       	p_{\beta}
       	\left( \xi_1 \tensor \cdots \tensor \xi_k \right)
       	p_{\beta}
       	\left( \eta_1 \tensor \cdots \tensor \eta_\ell \right)
       	\sum\limits_{\substack{a_1, b_1, \ldots, a_r, b_r \geq 0 \\
        		a_i + b_i \geq 1 \\
        		a_1 + \ldots + a_r = k \\
        		b_1 + \ldots + b_r = \ell
       	}}
       	1
       	\\
       	& \leq
       	\frac{1}{2 \cdot 8^n} 
       	(16 \kappa)^{k + \ell}
       	p_{\beta}
       	\left( \xi_1 \tensor \cdots \tensor \xi_k \right)
       	p_{\beta}
       	\left( \eta_1 \tensor \cdots \tensor \eta_\ell  \right)
       	\\ 
       	& =
       	\frac{1}{2 \cdot 8^n} 
       	(16 \kappa p)_{\beta}
       	\left( \xi_1 \tensor \cdots \tensor \xi_k \right)
       	(16 \kappa p)_{\beta}
       	\left( \eta_1 \tensor \cdots \tensor \eta_\ell  \right).
    \end{align*}
    From this we could conclude analogously to the procedure in the proof of 
    Theorem~\ref{Thm:LCAna:Continuity1} and the statement would follow.
    However, we need to show the existence of a $(\beta_n)_{n \in \mathbb{N}} 
    \in \mathfrak{I}_{\lie{g}}$ and a $\kappa > 0$, such that 
    \eqref{Nilpot:SavingSequence} holds.
    \begin{lemma}
    		For $f_{\alpha}$ defined as in \eqref{Nilpot:IncreasingFunction}, we 
    		take the function $g$ we get from 
    		Lemma~\ref{Nilpot:Lemma:MRZGrothLemma}. Then the sequence 
    		$((\beta_n)_{n \in \mathbb{N}}$ defined by
    		\begin{equation}
    			\beta_n
    			=
    			\exp \left(
    				\frac{g(n)}{8}
    			\right)
    		\end{equation}
    		and $\kappa = c \E^{8 g(1)}$ fulfil \eqref{Nilpot:SavingSequence}, 
    		where $c > 0$ is a constant such that $\alpha_n \leq c^n 
    		\chi^{\lie{g}}_n$.
    \end{lemma}
    \begin{subproof}
    		First, note that there is a fixed $c \geq 1$ such that
    		\begin{equation*}
    			\alpha_n
    			\leq
    			c^n \chi^{\lie{g}}_n
    			\quad \Longleftrightarrow \quad
    			\frac{1}{\chi(\lie{g})_n}
    			\leq
    			\frac{c'^n}{\alpha_n}
    		\end{equation*}
    		by the definition of a rapidly increasing sequence.
    		Denote $a_i + b_i = n_i$. Then we have
    		\begin{equation*}
    			\frac{\beta_k \cdot \beta_{\ell}}
	       	{
	       		\alpha_r
	       		\cdot
	       		\chi^{\lie{g}}_{n_1}
	       		\cdots
	       		\chi^{\lie{g}}_{n_r}
	       	}
	       	\leq
    			\frac{c^{k + \ell} \beta_k \cdot \beta_{\ell}}
	       	{
	       		\alpha_r
	       		\cdot
	       		\alpha_{n_1}
	       		\cdots
	       		\alpha_{n_r}
	       	}
    		\end{equation*}
    		and hence
    		\begin{align*}
    		\log 
    		&
    		\left(
    			\frac{\beta_k \cdot \beta_{\ell}}
		       	{
		       		\alpha_r
		       		\cdot
	  	     		\chi^{\lie{g}}_{n_1}
	 	      		\cdots
	  	     		\chi^{\lie{g}}_{n_r}
		       	}
	       	\right)
	       	\\
	       	& \leq
    		\log \left(
    			\frac{c^{k + \ell} \beta_k \cdot \beta_{\ell}}
		       	{
		       		\alpha_r
		       		\cdot
		       		\alpha_{n_1}
		       		\ldots
		       		\alpha_{n_r}
		       	}
	       	\right)
	       	\\
	       	& =
	       	(k + \ell) \log(c)
	       	+
	       	\frac{g(k)}{8}
	       	+
	       	\frac{g(\ell)}{8}
	       	-
	       	f_{\alpha}(r)
	       	-
	       	f_{\alpha}(n_1)
	       	- \cdots -
	       	f_{\alpha}(n_r)
	       	\\
	       	& \ot{(a)}{\leq}
	       	(k + \ell) \log(c)
	       	+
	       	\frac{g(k + \ell)}{8}
			-
	       	f_{\alpha}(r)
	       	-
	       	f_{\alpha}(n_1)
	       	- \cdots -
	       	f_{\alpha}(n_r)
	       	\\
	       	& \ot{(b)}{\leq}
			(k + \ell) \log(c)
	       	+
	       	g(n_1) + \cdots + g(n_r)
	       	+
	       	\frac{f_{\alpha}(r)}{8}
	       	\\
	       	& \qquad
	       	-
	       	f_{\alpha}(r)
	       	-
	       	f_{\alpha}(n_1)
	       	- \cdots -
	       	f_{\alpha}(n_r)
	       	\\
	       	& \ot{(c)}{\leq}
	       	(k + \ell) \log(c)
	       	+
	       	\left( g(n_1) - f_{\alpha}(n_1) \right)
	       	+ \cdots + 
	       	\left( g(n_r) - f_{\alpha}(n_r) \right)
	       	\\
	       	& \ot{(d)}{\leq}
	       	(k + \ell) \log(c')
	       	+
	       	8 g(1) n_1
	       	+ \cdots + 
	       	8 g(1) n_r
	       	\\
	       	& =
	       	(k + \ell)
	       	(\log(c') + 8 g(1)).
    		\end{align*}
    	In (a), we have used the convexity of $g$ and 
    	Lemma~\ref{Nilpot:Lemma:MRZGrothLemma} in (b). 
    	Then, (c) is just $\frac{f_\alpha(n)}{8} \leq f_\alpha(n)$.
    	In (d) we used Lemma~\ref{Nilpot:Lemma:MRZGrothLemma} again 
    	by estimating
    	\begin{equation*}
    		g(n)
    		=
    		g(1 + \ldots + 1)
    		\leq
    		8 n g(1) + f_{\alpha}(n)
    		\quad \Longleftrightarrow \quad
    		g(n) - f_\alpha(n)
    		\leq8 n g(1).
    	\end{equation*}
    	Now we need to exponentiate the inequality we just found and get
    	\begin{equation*}
   			\frac{\beta_k \cdot \beta_{\ell}}
	       	{
	       		\alpha_r
	       		\cdot
  	     		\chi^{\lie{g}}_{n_1}
 	      		\cdots
  	     		\chi^{\lie{g}}_{n_r}
	       	}
	       	\leq
	       	\E^{k + \ell)(\log(c') + 8 g(1))}
	       	=
	       	\left( c \E^{8 g(1)} \right)^{k + \ell}.
    	\end{equation*}
    	Since the $n_1, \ldots, n_r$ have been arbitrary, 
    	$(\beta_n)_{n \in \mathbb{N}}$ and $\kappa$ fulfil 
    	\eqref{Nilpot:SavingSequence} and the Lemma is proven.
    \end{subproof}
    Now, we only need to take the sum of the $C_n$. 
    We know this is possible from Chapter 5, and 
    we finally get a continuity estimate which is locally uniform in $z$.
\end{proof}
From the result of the previous proposition, we see that the next definition is 
useful:
\begin{definition}[The $\Sym_{1^{--}}$-Topology]
	Let $\lie{g}$ be a uniformly topologically nilpotent Banach-Lie algebra.
	Then we denote by $\Sym_{1^{--}}^{\bullet}(\lie{g})$ the symmetric tensor 
	algebra endowed with the topology, which is defined by the set of seminorms 
	$\algebra{P}$ from Definition \ref{Def:AdaptedBanachSeminorms} and which is 
	based $\lie{g}$-rapidly increasing sequences.
\end{definition}
Now we can rephrase the statement of Proposition~\ref{Nilpot:Prop:TopNilBanachLie} 
in the following way: $\Sym_{1^{--}}^{\bullet}(\lie{g})$ endowed with the Gutt 
star product is a locally convex algebra. We also see analogously to 
Corollary~\ref{corollary:NilpotentCase} that exponentials functions are in the 
completion $\widehat{\Sym}_{1^{--}}^{\bullet}(\lie{g})$.

At the end of Chapter 5, we showed that a finite-dimensional Lie algebra 
$\lie{g}$ is nilpotent if and only if its universal enveloping algebra 
$\algebra{U}(\lie{g})$ admitted a locally convex topology, such that the 
following three things are fulfilled.
\begin{enumerate}
	\item
	The product in $\algebra{U}(\lie{g})$ is continuous.

	\item
	For every $\xi \in \lie{g}$ the series $\exp(\xi)$ converges 
	absolutely in the completion of $\algebra{U}(\lie{g})$.

	\item
	Pulling back the topology to the symmetric tensor algebra, the projection 
	and inclusion maps with respect to the graded structure
	\begin{equation*}
		\Sym^{\bullet}(\lie{g})
    		\ot{$\pi_n$}{\longrightarrow}
   	    		\Sym^n(\lie{g})
    	    	\ot{$\iota_n$}{\longrightarrow}
    		\Sym^{\bullet}(\lie{g})
	\end{equation*}
	are continuous for all $n \in \mathbb{N}$.
\end{enumerate}
For Banach-Lie algebras, we came quite close to a similar statement: we know 
from Proposition~\ref{LCAna:Prop:NoBetterTopology} and the result of 
Wojty\`nski, that a Banach-Lie algebra must at least be \emph{topologically 
nil} to satisfy the three upper points, so being topologically nil is 
\emph{necessary}. We also know, that a uniformly topologically nilpotent 
Banach-Lie algebra $\lie{g}$ allows us to construct such a locally convex topology 
on $\algebra{U}(\lie{g})$ explicitly, hence uniform topological nilpotency is 
\emph{sufficient}. Maybe it is possible to find a notion of generalized 
nilpotency, which is equivalent to those three points.

% Chapter 7
%

%
% Chapter 7 of my master thesis:
% The Hopf algebra structure
%

\chapter{The Hopf Algebra Structure}

As already pointed out in Chapter 3, universal enveloping algebras of a Lie 
algebras are more than just associative, unital algebras: they constitute one of 
the most important types of \emph{Hopf algebras}, which are very common 
structures in mathematics. Hopf algebras are a particular \emph{bialgebras}, 
which in turn are on one hand associative, unital algebras, and on the other hand 
coassociative, counital coalgebras. Those two substructures of bialgebras must of 
course fulfil certain compatibility conditions. Note at this point, that 
coalgebras play a crucial role in the theory of formal deformation quantization, 
which is due to Kontsevich, but which has been a lot further developed since. 
Good references on this so-called formality theory are given by Esposito 
\cite{esposito:2015a} and Manetti \cite{manetti:2005a:script}.

In a Hopf algebra however, we have an 
additional map, called the \emph{antipode}, which again must fulfil some 
compatibility relations. So, in brief, a Hopf algebra over a field $\mathbb{K}$ 
is a tuple $(H, \cdot, \eta, \Delta, \varepsilon, S)$ of a vector space $H$ 
together with the following maps
\begin{equation*}
\begin{array}{rlrl}
	\cdot \colon
	&
	H \tensor H 
	\longrightarrow
	H, 
	& \quad &
	\text{ multiplication}
	\\
	\eta \colon
	&
	\mathbb{K}
	\longrightarrow
	H
	, 
	& \quad &
	\text{ unit}
	\\
	\Delta \colon
	&
	H
	\longrightarrow
	H \tensor H
	, & \quad &
	\text{ coproduct}
	\\
	\varepsilon \colon
	&
	H
	\longrightarrow
	\mathbb{K}
	, & \quad &
	\text{ counit}
	\\
	S \colon
	&
	H
	\longrightarrow
	H
	, & \quad &
	\text{ antipode},
\end{array}
\end{equation*}
such that we have a certain set of commuting diagrams. A very nice introduction 
to the theory of Hopf algebras with a lot of examples can e.g. be found in 
\cite{schweigert:2015a:script}.

In the previous chapters, we have mostly studied the properties of the 
multiplication in $\algebra{U}(\lie{g}_z)$, endowed with a particular topology. 
In this last chapter of this thesis, we will treat the remaining Hopf algebra 
structure maps. In the first section, we will see that the comultiplication and 
the antipode are not touched by our deformation procedure. Hence, it will be 
enough to show the continuity of the undeformed versions of those maps, which we 
will do in second section.

\section{An Undeformed Hopf Structure}

To get show the continuity of the remaining structure maps of 
$\algebra{U}_R(\lie{g}_z)$, we have to put estimates on seminorms again. For this 
purpose, we need explicit formulas for the antipode
\begin{equation}
    \label{eq:Antipode}
    S_z \colon
    \algebra{U}(\lie{g}_z)
    \longrightarrow
    \algebra{U}(\lie{g}_z)
\end{equation}
and the coproduct
\begin{equation}
    \label{eq:CoProduct}
    \coproduct_z \colon
    \algebra{U}(\lie{g}_z)
    \longrightarrow
    \algebra{U}(\lie{g}_z)
    \tensor
    \algebra{U}(\lie{g}_z)
\end{equation}
in the $\algebra{U}_R(\lie{g}_z)$ and in $\Sym^{\bullet}(\lie{g})$. We pull them 
back to the symmetric algebra and extend them to the whole tensor algebra by 
symmetrizing beforehand. We define
\begin{equation}
    \label{eq:AntipodeOnTensor}
    \widetilde{S}_z \colon
    \Tensor^\bullet(\lie{g})
    \longrightarrow
    \Sym^\bullet(\lie{g})
    , \quad
    \widetilde{S}_z
    =
    \mathfrak{q}_z^{-1}
    \circ
    S_z
    \circ
    \mathfrak{q}_z
    \circ
    \Symmetrizer
\end{equation}
and
\begin{equation}
    \label{eq:CoProductOnTensor}
    \widetilde{\ocoproduct}_z \colon
    \Tensor^\bullet(\lie{g})
    \longrightarrow
    \Sym^\bullet(\lie{g})
    \tensor
    \Sym^\bullet(\lie{g})
    , \quad
    \widetilde{\ocoproduct}_z
    =
    (\mathfrak{q}_z^{-1} \tensor \mathfrak{q}_z^{-1})
    \circ
    \coproduct_z
    \circ
    \mathfrak{q}_z
    \circ
    \Symmetrizer,
\end{equation}
to avoid that the maps on $\Sym^{\bullet}(\lie{g})$ and on $\algebra{U}
(\lie{g}_z)$ are denoted by the same symbols. The next lemma gives us the two 
explicit formulas we need.
\begin{lemma}
    \label{Thm:Hopf:Formulas}%
    For $\xi_1, \ldots, \xi_n \in \lie{g}$ we have the identities
    \begin{equation}
        \label{Hopf:AntipodeFormula}
        \widetilde{S}_z
        \left( \xi_1 \tensor \cdots \tensor \xi_n \right)
        =
        (-1)^n
        \xi_1 \cdots \xi_n
    \end{equation}
    and
    \begin{equation}
        \label{Hopf:CoproductFormula}
        \widetilde{\ocoproduct}_z
        \left(\xi_1 \tensor \cdots \tensor \xi_n \right)
        =
        \sum\limits_{
        	I \subseteq
        	\{1, \ldots, n\}
        }
        \xi_I
        \tensor
        \xi_1 \cdots
        \widehat{\xi_I}
        \cdots \xi_n
    \end{equation}
    where $\xi_I$ denotes the symmetric tensor product of all $\xi_i$ with 
    $i \in I$ and $\widehat{\xi_I}$ means that the $\xi_i$ with $i \in I$ 
    are left out.
\end{lemma}
\begin{proof}
	First, we derive Formula \ref{Hopf:AntipodeFormula}: the antipode gives
	$S_z(\xi) = - \xi$ for $\xi \in \lie{g}$ and extends to $\algebra{U}
	(\lie{g}_z)$ by algebra antihomomorphism, hence
	\begin{equation*}
		S_z \left( 
			\xi_1 \odot \cdots \odot \xi_n
		\right)
		=
		(-1)^n
		\xi_n \odot \cdots \odot \xi_1
	\end{equation*}
	in $\algebra{U}(\lie{g}_z)$ for $\xi_1, \ldots, \xi_n \in \lie{g}$.
	This means
	\begin{equation*}
		\widetilde{S}_z \left( 
			\xi_1 \star_z \cdots \star_z \xi_n
		\right)
		=
		(-1)^n
		\xi_n \star_z \cdots \star_z \xi_1
	\end{equation*}	
	in $\Sym^{\bullet}(\lie{g})$. But now, using the linearity of 
	$\widetilde{S}_z$ we get
	\begin{align*}
		\widetilde{S}_z \left( \xi_1 \cdots \xi_n \right)
		& =
		\widetilde{S}_z \left( 
		\frac{1}{n!}
			\sum\limits_{\sigma \in S_n}
			\xi_{\sigma(1)} 
			\star_z \cdots \star_z 
			\xi_{\sigma(n)}
		\right)
		\\
		& =
		\frac{1}{n!}
		\sum\limits_{\sigma \in S_n}
		\widetilde{S}_z \left( 
			\xi_{\sigma(1)} 
			\star_z \cdots \star_z 
			\xi_{\sigma(n)}
		\right)
		\\
		& =
		\frac{1}{n!}
		\sum\limits_{\sigma \in S_n}
		(-1)^n
		\xi_{\sigma(n)} 
		\star_z \cdots \star_z 
		\xi_{\sigma(1)}
		\\
		& =
		(-1)^n
		\xi_1 \cdots \xi_n.
	\end{align*}
	For the coproduct, we have well-known formula with shuffle permutations:
	\begin{equation*}
		\ocoproduct_z \left(
			\xi_1 \odot \cdots \odot \xi_n
		\right)
		=
		\sum\limits_{k=0}^n
		\sum\limits_{\substack{
			1 \leq i_1
			< \ldots <
			i_k \leq n
			\\
			I = \{i_1, \ldots, i_k\}
		}}
		\xi_{i_1} \odot \cdots \odot \xi_{i_k}
		\tensor
		\xi_1 
		\odot \cdots 
		\widehat{\xi_I}
		\cdots \odot
		\xi_n.
	\end{equation*}
	We can derive it from the fact that $\ocoproduct_z$ has the following form on 
	Lie algebra elements:
	\begin{equation*}
		\ocoproduct_z (\xi)
		=
		\xi \tensor \Unit
		+
		\Unit \tensor \xi.
	\end{equation*}
	The coproduct extends to $\algebra{U}(\lie{g}_z)$ by algebra 
	homomorphism. This yields
	\begin{align*}
		\ocoproduct_z \left(
			\xi_1 \odot \cdots \odot \xi_n
		\right)
		& =
		\ocoproduct_z \left( \xi_1 \right)
		\odot \cdots \odot
		\ocoproduct_z \left( \xi_n \right)
		\\
		&=
		\left( \xi_1 \tensor \Unit + \Unit \tensor \xi_1 \right)
		\odot \cdots \odot
		\left( \xi_n \tensor \Unit + \Unit \tensor \xi_n \right)
		\\
		&=
		\sum\limits_{k=0}^n
		\sum\limits_{\substack{
			1 \leq i_1
			< \ldots <
			i_k \leq n
			\\
			I = \{i_1, \ldots, i_k\}
		}}
		\xi_{i_1} \odot \cdots \odot \xi_{i_k}
		\tensor
		\xi_1 
		\odot \cdots 
		\widehat{\xi_I}
		\cdots \odot
		\xi_n.
	\end{align*}
	Now we can pull this back to $\Sym^{\bullet}(\lie{g}_z)$. We get a $\star_z$ 
	for every $\odot$. For symmetric tensors we have by linearity
	\begin{align*}
		&
		\widetilde{\ocoproduct}_z 
		\left( \xi_1 \cdots \xi_n \right)
		\\
		& =
		\widetilde{\ocoproduct}_z \left( 
			\frac{1}{n!}
			\sum\limits_{\sigma \in S_n}
			\xi_{\sigma(1)} 
			\star_z \cdots \star_z 
			\xi_{\sigma(n)}
		\right)
		\\
		& =
		\frac{1}{n!}
		\sum\limits_{\sigma \in S_n}
		\widetilde{\ocoproduct}_z \left( 
			\xi_{\sigma(1)} 
			\star_z \cdots \star_z 
			\xi_{\sigma(n)}
		\right)
		\\
		& =
		\frac{1}{n!}
		\sum\limits_{\sigma \in S_n}
		\sum\limits_{k=0}^n
		\sum\limits_{\substack{
			1 \leq i_1
			< \ldots <
			i_k \leq n
			\\
			I = \{i_1, \ldots, i_k\}
		}}
		\xi_{i_{\sigma(1)}} 
		\star_z \cdots \star_z
		\xi_{i_{\sigma(k)}}
		\tensor
		\xi_{\sigma(1)}
		\star_z \cdots
		\widehat{\xi_{\sigma(I)}}
		\cdots \star_z 
		\xi_{\sigma(n)}
		\\
		& =
		\frac{1}{n!}
		\sum\limits_{k=0}^n
		\sum\limits_{\sigma \in S_k}
		\sum\limits_{\tau \in S_{n - k}}
		\sum\limits_{
			\{i_1, \ldots, i_k\}
			\subseteq
			\{1, \ldots, n\}
		}
		\frac{n!}{k! (n - k)!}
		\cdot
		\xi_{i_{\sigma(1)}} 
		\star_z \cdots \star_z 
		\xi_{i_{\sigma(k)}}
		\tensor
		\xi_{\tau(1)}
		\star_z \cdots 
		\widehat{\xi_I}
		\cdots \star_z 
		\xi_{\tau(n)}
		\\
		& =
		\sum\limits_{k=0}^n
		\sum\limits_{
			\{i_1, \ldots, i_k\}
			\subseteq
			\{1, \ldots, n\}
		}
		\xi_{i_1} \cdots \xi_{i_k}
		\tensor
		\xi_1 \cdots \widehat{\xi_I} \cdots \xi_n.
	\end{align*}
\end{proof}
This lemma yields immediately the following result.
\begin{corollary}
	If we deform the symmetric tensor algebra of an AE-Lie algebra 
	$\Sym^{\bullet}(\lie{g})$ with the Gutt star product, the coproduct 
	and the antipode will remain undeformed.
\end{corollary}
\begin{remark}[Hopf structures on the tensor algebra]
	\mbox{}
	\begin{remarklist}
		\item
		On the first sight, this result looks astonishing, but actually it is 
		not. The Hopf algebra structure we found here is just the Hopf algebra 
		structure which comes from the tensor algebra of $\lie{g}$. To be more 
		precise: it is the Hopf algebra structure using the tensor product and 
		the shuffle coproduct. There is also a second Hopf algebra structure: the 
		one using the deconcatenation coproduct and the shuffle product. We find 
		these structures on every tensor algebra of a vector space, it does not 
		need to be a tensor algebra over a Lie algebra. Since the coalgebra 
		structure on $\algebra{U}(\lie{g}_z)$ is inherited from 
		$\Tensor^{\bullet}(\lie{g})$ and does not make use of the Lie bracket, 
		there is no reason why a rescaling of the Lie bracket should change it.
		
		\item
		What we have seen is hence just \emph{one} possibility to deform the 
		symmetric algebra. Another way to do so would be to deform the 
		costructure and leave the product untouched. Such a deformation of a Hopf 
		algebra would lead to the notion of \emph{quantum groups}.
	\end{remarklist}
\end{remark}

% AE-Lie Algebras
%

\section{Continuity of the Hopf Structure}
We need a topology on the tensor product in \eqref{Hopf:CoproductFormula},
since we want to prove the continuity of this map. As we have always used the 
projective tensor product for our construction so far, it seems just logic to do 
so again. The continuity of the two maps is then very easy to show.
\begin{proposition}
    \label{Prop:Hopf:CoproductContinuity}%
    Let $\lie{g}$ be an AE-Lie algebra and $R \geq 0$. For every continuous 
    seminorm $p$ and all $x \in \widehat{\Tensor}_R^\bullet(\lie{g})$
    the following estimates hold:
    \begin{equation}
        \label{Hopf:AntipodeContinuity}
        p_R \left( \widetilde{S}_z(x) \right)
        \leq
        p_R (x)
    \end{equation}
    and
    \begin{equation}
        \label{Hopf:CoproductContinuity}
        (p_R \tensor p_R)
        \left( \widetilde{\ocoproduct}_z(x) \right)
        \leq
        (2 p)_R (x).
    \end{equation}
\end{proposition}
\begin{proof}
	We just need to show both estimates on factorizing tensors and extend them 
	with the infimum argument. Inequality \eqref{Hopf:AntipodeContinuity} is 
	clear, since we only get a sign.
	To get Equation~\eqref{Hopf:CoproductContinuity}, we compute:
	\begin{align*}
		(p_R \tensor p_R)
        \left( \widetilde{\ocoproduct}_z
        	\left(
        		\xi_1 \tensor \cdots \tensor \xi_n
        	\right)
        \right)
        & =
		(p_R \tensor p_R)
        \left( 
			\sum\limits_{
        		I \subseteq
        		\{1, \ldots, n\}
        	}
        	\xi_I
        	\tensor
        	\xi_1 \cdots
        	\widehat{\xi_I}
        	\cdots \xi_n
        \right)
        \\
        & \leq
        \sum\limits_{
        	I \subseteq
        	\{1, \ldots, n\}
        }
        |I|!^R (n - |I|)!^R
        p^{|I|} \left( \xi_I \right)
        p^{n - |I|}
        \left( 
        	\xi_1 \cdots \widehat{\xi_I} \cdots \xi_n 
        \right)
        \\
        & \leq
        \sum\limits_{
        	I \subseteq
        	\{1, \ldots, n\}
        }
        |I|!^R (n - |I|)!^R
        p(\xi_1) \cdots p(\xi_n)
        \\
        & \leq
        \sum\limits_{
        	I \subseteq
        	\{1, \ldots, n\}
        }
        n!^R
        p(\xi_1) \cdots p(\xi_n)
        \\
        & =
        2^n n!^R
        p(\xi_1) \cdots p(\xi_n)
        \\
        & =
        (2p)_R \left(
        	\xi_1 \tensor \cdots \tensor \xi_n
        \right).
	\end{align*}
	This shows the statement.
\end{proof}
\begin{remark}
	We see that we get no dependence on the parameter $R$. Also this is not 
	surprising: the Hopf structure of the symmetric or the universal enveloping 
	algebra over a vector space is cocommutative. The symmetric tensor product is 
	commutative and its continuity estimate does not depend on the parameter $R$ 
	either. Note that for an abelian Lie algebra, also the product structure  
	remains undeformed and all Hopf algebra maps are continuous for $R \geq 0$. In 
	this sense, the independence of $R$ fits into the picture.
\end{remark}

The only maps which are left to consider are the unit and the counit. Since their 
continuity is clear by the definition of the $\Tensor_R$-topology, we get the 
a final result.
\begin{theorem}
    \label{Prop:Hopf:ContinuousHopf}%
    Let $\lie{g}$ be an AE-Lie algebra and $z \in \mathbb{K}$. Then, 
    if $R \geq 1$, $\widehat{\Sym}_R^\bullet (\lie{g})$ will be a locally convex 
    Hopf algebra. The same will hold for $\widehat{\Sym}_{1^-}^{\bullet}
    (\lie{g})$, if $\lie{g}$ is a nilpotent locally convex Lie algebra with 
    continuous Lie bracket and for $\widehat{\Sym}_{1^{--}}^{\bullet}(\lie{g})$, 
    if $\lie{g}$ is a uniformly topologically nilpotent Banach-Lie algebra.
\end{theorem}

% ============================================================================
% ////////////////////////////////////////////////////////////////////////////

% ============================================================================
% 			Bibliography Part
% ============================================================================

% Smaller than other text
%

{

  \footnotesize
  \renewcommand{\arraystretch}{0.5}

}

% ============================================================================
% ////////////////////////////////////////////////////////////////////////////

% Abschließende Formel
%

% ============================================================================
% ////////////////////////////////////////////////////////////////////////////

\end{document}